\documentclass{amsart}

\usepackage[utf8]{inputenc}
\usepackage[T1]{fontenc}
\usepackage[english]{babel}
\usepackage{graphicx}
\usepackage{amssymb}
\usepackage{amsmath}
\usepackage{tikz-cd}
\usepackage{adjustbox}
\usepackage{hyperref}
\usepackage{framed}
\usepackage{amsaddr}
\usepackage{mathtools}
\usepackage{mathscinet}
\usepackage{quiver}
\usepackage{MnSymbol}
\usepackage{bm}
\usepackage{placeins}
\usepackage[backend=bibtex,
style=alphabetic,
bibencoding=ascii
]{biblatex}
\usepackage{rotating}

\usepackage[all]{xy}

\usepackage{stmaryrd}
\SetSymbolFont{stmry}{bold}{U}{stmry}{m}{n}

\usepackage{etoolbox}
\makeatletter
\patchcmd{\part}{\normalfont}{\normalfont\Large}{}{}
\makeatother

\newtheorem{theorem}{Theorem}[subsection]
\newtheorem{lemma}[theorem]{Lemma}
\newtheorem{proposition}[theorem]{Proposition}
\newtheorem{corollary}[theorem]{Corollary}
\newtheorem{conjecture}[theorem]{Conjecture}

\theoremstyle{definition}
\newtheorem{definition}[theorem]{Definition}

\newtheorem{example}[theorem]{Example}

\theoremstyle{remark}
\newtheorem{remark}[theorem]{Remark}

\numberwithin{equation}{section}

\DeclareMathOperator{\Lie}{Lie}
\DeclareMathOperator{\Spec}{Spec}
\DeclareMathOperator{\Rep}{Rep}
\DeclareMathOperator{\Ind}{Ind}
\DeclareMathOperator{\rank}{rank}
\DeclareMathOperator{\ind}{\mathrm{ind}}

\DeclareMathOperator{\colim}{colim}
\DeclareMathOperator{\Spf}{Spf}

\newcommand{\U}{\mathbb{U}}
\newcommand{\ov}{\overline}

\newcommand{\Int}{\mathrm{Int}}

\newcommand{\der}{\mathrm{der}}
\newcommand{\mc}{\mathcal}

\newcommand{\Q}{\mathbb{Q}}

\newcommand{\mf}{\mathfrak}

\newcommand{\A}{\mathbb{A}}
\newcommand{\Sc}{\mathrm{sc}}
\newcommand{\ad}{\mathrm{ad}}

\newcommand{\bas}{\mathrm{bas}}

\newcommand{\bb}{\mathbb}
\newcommand{\C}{\mathbb{C}}
\newcommand{\R}{\mathbb{R}}
\newcommand{\Gal}{\mathrm{Gal}}
\newcommand{\Z}{\mathbb{Z}}

\newcommand{\SL}{\mathrm{SL}}

\newcommand{\RH}{\mathrm{RH}}

\newcommand{\Irr}{\mathrm{Irr}}

\newcommand{\Gm}{{\mathbb{G}_{\mathrm{m}}}}
\newcommand{\Hom}{\mathrm{Hom}}
\newcommand{\spec}{\mathrm{spec}}

\newcommand{\ab}{\mathrm{ab}}

\newcommand{\rel}{\mathrm{rel}}
\newcommand{\Res}{\mathrm{Res}}

\newcommand{\GL}{\mathrm{GL}}
\newcommand{\op}{\mathrm{op}}
\DeclareMathOperator{\End}{End}
\newcommand{\co}[1]{\prescript{#1}{}}
\newcommand{\wh}{\widehat}
\newcommand{\res}{\mathrm{res}}

\newcommand{\Bun}{\mathrm{Bun}}

\newcommand{\D}{\mathbb{D}}

\newcommand{\frob}{\mathrm{frob}}
\newcommand{\Par}{\mathrm{Par}}
\newcommand{\fss}{\mathrm{fss}}
\newcommand{\Fss}{\mathrm{Fss}}
\newcommand{\Sym}{\mathrm{Sym}}

\newcommand{\triv}{\mathrm{triv}}
\newcommand{\IndCoh}{\mathrm{IndCoh}}
\newcommand{\QCoh}{\mathrm{QCoh}}
\newcommand{\diag}{\mathrm{diag}}
\DeclareMathOperator{\Eis}{\mathrm{Eis}}
\newcommand{\RHom}{\mathrm{RHom}}

\newcommand{\Coh}{\mathrm{Coh}}
\newcommand{\Ad}{\mathrm{Ad}}
\DeclareMathOperator{\CT}{\mathrm{CT}}

\newcommand{\Sp}{\mathrm{Sp}}
\newcommand{\Spin}{\mathrm{Spin}}
\newcommand{\Cusp}{\mathrm{Cusp}}
\newcommand{\ur}{\mathrm{ur}}
\newcommand{\id}{\mathrm{id}}

\newcommand{\cusp}{\mathrm{cusp}}

\newcommand{\coh}{\mathrm{coh}}
\newcommand{\qc}{\mathrm{qc}}
\newcommand{\LCusp}{\mathrm{LCusp}}
\newcommand{\RCusp}{\mathrm{RCusp}}
\newcommand{\GS}{\mathrm{GS}}

\newcommand{\rc}{\mathrm{rc}}
\newcommand{\V}{\mathbb{V}}
\newcommand{\T}{\mathbb{T}}
\newcommand{\ProCoh}{\mathrm{ProCoh}}

\newcommand{\up}{\mathrm{up}}
\renewcommand{\L}{\mathbb{L}}
\newcommand{\WD}{\mathrm{WD}}
\newcommand{\Ga}{\mathbb{G}_{\mathrm{a}}}

\newcommand{\F}{\mathbb{F}}
\newcommand{\tor}{\mathrm{tor}}
\newcommand{\tf}{\mathrm{tf}}
\newcommand{\Alg}{\mathrm{Alg}}
\newcommand{\Ani}{\mathrm{Ani}}
\newcommand{\Stk}{\mathrm{Stk}}
\newcommand{\cl}{\mathrm{cl}}
\newcommand{\Mod}{\mathrm{Mod}}
\newcommand{\Perf}{\mathrm{Perf}}
\newcommand{\Map}{\mathrm{Map}}
\newcommand{\afp}{\mathrm{afp}}
\newcommand{\Pro}{\mathrm{Pro}}
\newcommand{\Aperf}{\mathrm{APerf}}
\newcommand{\LS}{\mathrm{LS}}
\newcommand{\Sect}{\mathrm{Sect}}
\newcommand{\adj}{\mathrm{adj}}
\newcommand{\red}{\mathrm{red}}

\newcommand{\RBas}{\mathrm{RBas}}
\newcommand{\LBas}{\mathrm{LBas}}
\newcommand{\Bas}{\mathrm{Bas}}

\DeclareMathOperator{\un}{\mathrm{un}}
\DeclareMathOperator{\gr}{gr}
\newcommand{\ev}{\mathrm{ev}}
\newcommand{\ush}{\lvert ! \rvert}

\newcommand{\wt}{\mathrm{wt}}

\newcommand{\rig}{\mathrm{rig}}
\newcommand{\SO}{\mathrm{SO}}

\newcommand{\magic}{\mathrm{magic}}
\newcommand{\MF}{\mathrm{MF}}
\newcommand{\mrL}{\mathrm{L}}
\newcommand{\M}{\mathrm{M}}
\newcommand{\cofib}{\mathrm{cofib}}
\newcommand{\Nm}{\mathrm{Nm}}

\newcommand{\shear}{{\mathbin{\mkern-6mu\fatslash}}}
\newcommand{\unshear}{\,{\mathbin{\mkern-6mu\fatbslash}}}

\title{Cuspidal sheaves on the stack of $L$-parameters}
\author{Alexander Bertoloni Meli, Teruhisa Koshikawa}

\begin{document}
\begin{abstract}
We introduce the notion of cuspidal coherent sheaves on the moduli stack of local $L$-parameters. We give a complete classification and relate them to Lusztig's cuspidal local systems (or character sheaves). Under the categorical local Langlands conjecture, these sheaves are expected to correspond to supercuspidal representations supported at basic points. We also construct semiorthogonal decompositions of the category of coherent sheaves on the moduli stack of $L$-parameters that conjecturally match the semiorthogonal decompositions of sheaves on $\Bun_G$ coming from the stratification by the Kottwitz set and Bernstein decomposition of representation categories.
\end{abstract}

\maketitle

\tableofcontents

\section{Introduction}
\subsection{Main results}
The local Langlands conjecture predicts a bijection between irreducible smooth representations of $p$-adic groups and so-called (enhanced) $L$-parameters consisting of a pair of a Weil--Deligne parameter $\phi$ and an irreducible representation of a suitable centralizer group of $\phi$ (see \cite{KalethaLLCnonQS} for instance). Several categorical versions of the local Langlands conjecture have been proposed, and in particular Fargues and Scholze \cite{fargues--scholze} conjecture that, in analogy with the geometric Langlands conjecture \cite{ArinkinGaitsgory}, 
\begin{equation}{\label{eqn: FSCLLC}}
    D(\Bun_G, \overline{\mathbb{Q}}_\ell) \simeq \IndCoh(\Par_G), 
\end{equation}
where $\Bun_G$ is the moduli v-stack of $G$-bundles on the Fargues--Fontaine curve for a quasi-split group $G$ over a $p$-adic field $F$, and $\Par_G$ is the moduli stack of Weil--Deligne parameters valued in the Langlands dual group $\wh G$ of $G$ over $\overline{\Q}_\ell$. The stack $\Par_G$ is a rather explicit classical (ind-)algebraic stack, and $\IndCoh$ denotes the $\infty$-category of ind-coherent sheaves \cite{GaitsgoryIndCohSheaves}, which is a ``renormalization'' of $\QCoh$ so that the coherent sheaves are compact. 

Some features of the left-hand side --- the geometric side --- are 
\begin{itemize}
    \item it admits a semiorthogonal decomposition indexed by the Kottwitz set $B(G)$ reflecting the geometry of $\Bun_G$ and the identification of the underlying topological space $\lvert \Bun_G\rvert \cong B(G)$,
    \item for each $b\in B(G)$, the piece at $b$ is equivalent to $D(G_b(F), \ov{\Q}_{\ell})$, the derived category of smooth representations of the reductive group $G_b(F)$ attached to $b$, and
    \item there are actually two semiorthogonal decompositions, corresponding to the two inclusion functors $i_{b!}$, $i_{b\sharp}$ (the left adjoint of $i_b^*$), and indexed by the descending partial order (i.e., basic points are maximal) respectively the ascending partial order (cf.~\cite[\S 9.0.2]{GHILZ}). If $b$ is basic, the embeddings agree. If $b$ is non-basic, the functors can be described using geometric Eisenstein series (\cite[Corollary 2.2.5]{HHS1}, \cite[Lemma 2.9]{Takayasecondadjointness}), e.g.\footnote{We will ignore the renormalization issue throughout.}
    \begin{equation*}
        i_{b\sharp} \simeq \Eis_{P_b} \circ i_{b_M!}\,\text{(up to a certain renormalization)}
    \end{equation*}
    for standard parabolic/Levi subgroups $P_b\supset M_b$ attached to $b$ so that $b$ lifts to a unique basic, strictly $P_b$-dominant element $b_M \in B(M_b)$.\footnote{The notion of dominance relative to a parabolic is central to the paper and explained precisely in \S \ref{s: notation}.}
    In particular, geometric Eisenstein series are fully faithful in this situation.  
\end{itemize}

By the categorical conjecture \eqref{eqn: FSCLLC}, one expects the ``same'' structures on the right-hand side --- the spectral side --- as well, and our first main result establishes them. 

Let us recall one definition first. For a parabolic subgroup $P \subset G$ with Levi quotient $M$, consider the diagram
\begin{equation*}
    \begin{tikzcd}
        \Par_M & \Par_P \arrow[l, swap,  "q"] \arrow[r, "p"] & \Par_G. 
    \end{tikzcd}
\end{equation*}
\begin{definition}
We define the \emph{spectral Eisenstein series} functor by
\begin{equation*}
    \Eis^{\spec}_P\colon \IndCoh(\Par_M) \rightarrow \IndCoh(\Par_G); \qquad \mc{F} \mapsto p_{\ast}q^{\ast, \IndCoh} \mc{F}\simeq p_* q^! \mc{F}.
\end{equation*}
\end{definition}
It is conjectured (\cite[Conjecture 1.6.2, Proposition 6.3.2]{HansenMann}) that geometric Eisenstein series and spectral Eisenstein series match under \eqref{eqn: FSCLLC}. Note also that $\IndCoh(\Par_M)$ admits a central weight decomposition indexed by $X^{\ast}(Z(\wh{M})^{\Gamma})$ (cf.~\S \ref{ss: weightdecomposition}).

Now, for every $G$, we define a certain full subcategory
\[
\Bas(\Par_G) \subset \IndCoh(\Par_G), 
\]
which corresponds to sheaves supported at basic points under the conjectural equivalence \eqref{eqn: FSCLLC}. 

\begin{definition}
   We define $\Bas(\Par_G) \subset \IndCoh(\Par_G)$ to be the full subcategory of objects $\mc{F}$ such that
    \begin{equation*}
       \Hom(\Eis^{\spec}_{P} \mc{F}', \mc{F})\simeq 0,
       \end{equation*}
       for each proper standard parabolic $P \subset G$ and each $\mc{F}' \in \IndCoh(\Par_M)$, with central weight $\chi$ for strictly $P$-dominant $\chi \in X^*(Z(\wh{M})^{\Gamma})$.

       We refer to sheaves in $\Bas(\Par_G)$ as \emph{basic}. The category of basic sheaves $\Bas(\Par_G)$ admits a central weight decomposition indexed by $\kappa\colon B(G)_{\bas} \cong X^* (Z(\wh G)^\Gamma) $ (it is simply induced by that of $\IndCoh(\Par_G)$). The summand for $b$ is denoted by $\Bas(\Par_G)^{\kappa(b)}$. 
\end{definition}
\begin{remark}
\leavevmode
\begin{itemize}
    \item The above is actually the definition of \emph{right}-basic sheaves. We also define left-basic sheaves, and show they agree (Corollary \ref{cor: BasParGgenbyCusp}).
    \item Assuming the existence of the equivalence \eqref{eqn: FSCLLC} satisfying its expected properties including compatibility with Eisenstein series, $\Bas(\Par_G)^{\kappa(b)}$, $b\in B(G)_{\bas}$ corresponds to the category of $i_{b!} \pi \simeq i_{b\sharp} \pi$ for $\pi \in D(G_b(F), \ov{\Q}_{\ell})$ (cf.~\cite[Proposition 6.7.3]{HansenMann}).
\end{itemize}
\end{remark}

As above, each element $b\in B(G)$ determines a unique standard Levi subgroup $M_b$ and standard parabolic subgroup $P_b$ with a basic lift $b_M \in B(M_b)$ such that the image under the Kottwitz map, $\kappa_{M_b}(b_M)$, is strictly $P_b$-dominant. 

\begin{theorem}\label{thm:IntroI}
The following statements hold:
\begin{itemize}
    \item the $\infty$-subcategories $\Eis^\spec_{P_b} \Bas(\Par_{M_b})^{\kappa_{M_b}(b_M)} \subset \IndCoh(\Par_G)$ for each $b \in B(G)$ form a semiorthogonal decomposition with the ascending partial order: $\Hom(\mc{F}_{b}, \mc{F}_{b'}) = 0$ unless $b \leq b'$ (Theorem \ref{thm: basSODparG}).
    \item a similar statement holds with $\Eis^\spec_{P_b^\op}$ in place of $\Eis^\spec_{P_b}$, with the opposite ordering (Theorem \ref{thm: basSODparG}).
    \item the functors $\Eis^\spec_{P_b}$, $\Eis^\spec_{P_b^\op}$ are fully faithful on $\Bas(\Par_{M_b})^{\kappa_{M_b}(b_M)}$ (Corollary \ref{cor: EisFFonBas}).
\end{itemize}
\end{theorem}

We will deduce Theorem \ref{thm:IntroI} from finer results.\footnote{After completing the first draft, we learned that Theorem \ref{thm:IntroI} itself can be derived essentially from \cite{BuPadurariuToda} by applying it to $\Par_G$.} To motivate this, let us return to consideration of the geometric side. There is a finer structure determined by the Bernstein decomposition of  $D(G_b(F), \ov{\Q}_{\ell})$:
\begin{itemize}
    \item the category $D(G_b(F), \ov{\Q}_{\ell})$ admits an orthogonal decomposition into Bernstein blocks determined up to inertial equivalence by a pair $[L, \sigma]$ of a Levi subgroup $L \subset G_b$ and supercuspidal representation $\sigma$ of $L(F)$,
    \item each Bernstein block has a natural compact projective generator constructed by Bernstein and given as the parabolic induction $\Ind^{G_b (F)}_{P(F)} \widetilde{\sigma}$, where $\widetilde{\sigma}$ is the universal unramified twist of $\sigma$ \cite{BernsteinRumelhart},
    \item if $G$ is semisimple and $b$ is basic, the category generated by supercuspidal Bernstein blocks (i.e., $L=G_b$) is semisimple,
    \item if $b_M \in B(M)$ is basic and $P$-dominant, with image $b\in B(G)$, then the extension of a parabolic induction can be realized as a geometric Eisenstein series (\cite[Corollary 2.2.14]{HHS1}, \cite[proof of Lemma 2.4]{Takayasecondadjointness}),
    \begin{equation*}
        \Eis_{P} \circ i_{b_M !} \simeq i_{b \sharp} \circ \Ind^{G_b(F)}_{Q_b(F)},
    \end{equation*}
    where $Q_b$ is the transfer of $P$ to $G_b$.
    In particular, geometric Eisenstein series are conservative in this setting.  
\end{itemize}

Our second main result establishes an analogous finer structure on the spectral side. Before stating it, we give the following definition which corresponds under the conjectural equivalence \eqref{eqn: FSCLLC} to the category generated by $i_{b!}\pi$, for $\pi$ a supercuspidal representation of $G_b(F)$, $b\in B(G)_{\bas}$.
\begin{definition}\label{def: IntroCusp}
       We define $\Cusp(\Par_G) \subset \IndCoh(\Par_G)$ to be the full subcategory of sheaves $\mc{F}$ such that
    \begin{equation*}
       \Hom(\Eis^{\spec}_{P} \mc{F}', \mc{F})\simeq 0,
       \end{equation*}
       for each proper standard parabolic $P \subset G$ and each $\mc{F}' \in \IndCoh(\Par_M)^{\chi}$, for each (not necessarily strictly!) $P$-dominant $\chi \in X^*(Z(\wh{M})^{\Gamma})$.\footnote{As before, this is actually the definition of \emph{right}-cuspidal sheaves. We also define left-cuspidal sheaves, and show they agree.}

       We refer to sheaves in $\Cusp(\Par_G)$ as \emph{cuspidal}. As before, $\Cusp(\Par_G)$ admits a central weight decomposition and we write $\Cusp(\Par_G)^{\kappa(b)}$ for the weight $b\in B(G)_\bas$ summand. 
\end{definition}

\begin{theorem}\label{thm:IntroII}
The following statements hold:
\begin{itemize}
    \item $\Bas(\Par_G)^{\kappa(b)}$ is generated under colimits by $\Cusp (\Par_G)^{\kappa(b)}$ and $\Eis^\spec_P \Cusp (\Par_M)^{\kappa_M(b_M)}$ with $b_M$ whose image in $B(G)$ is $b$, and they are orthogonal to each other (Corollary \ref{cor: BasParGgenbyCusp}).
    \item The $\infty$-subcategory generated under colimits by $\Eis^\spec_P \Cusp(\Par_M)^{\kappa_M(b_M)}$, $b_M$ $P$-dominant with a fixed image $b\in B(G)$ has natural mutually orthogonal compact generators whose endomorphism algebras are concentrated in degree 0, and,
    at least if $M=G$, the generators themselves are concentrated in degree 0  (Theorem \ref{thm: SODcuspsupp}, Corollary \ref{cor: cohspringerhoms}, Theorem \ref{thm: Generators}).
    \item If $G$ is semisimple, $\Cusp(\Par_G)$ is semisimple (Part (3) of Theorem \ref{thm: CVBclassification}).
    \item $\Eis^\spec_P$ is conservative on $\Cusp(\Par_M)^{\kappa_M(b_M)}$ with a fixed $P$-dominant $b_M$ (Corollary \ref{cor: ParGEisconservativity}). 
\end{itemize}
\end{theorem}

The natural compact ``projective'' generators above can be described explicitly. Similarly, it is also possible to classify simple objects in $\Cusp(\Par_G)$ explicitly, and this is our third main result. Under the conjectural equivalence \eqref{eqn: FSCLLC}, these are irreducible supercuspidal representations at basic points on the geometric side, and our classification is closely related to the cuspidal enhanced $L$-parameters of Aubert--Moussaoui--Solleveld, and to the fundamental work of Lusztig on cuspidal local systems and character sheaves. We call these simple objects (and their direct sums) \emph{cuspidal vector bundles}.

To describe our result more precisely, let $\lambda\colon W_F \rightarrow \wh G$ be a semisimple $L$-parameter and consider the stack $\V_{\lambda}$ defined by the classical fiber product
\begin{equation*}
    \begin{tikzcd}
        \V_{\lambda} \arrow[d, "\pi"] \arrow[r, hook, "i_\lambda"] & \Par_G \arrow[d] \\
       BH_{\lambda} \arrow[r, hook] &  (\text{stack of Weil parameters}).
    \end{tikzcd}
\end{equation*}
Then $\V_{\lambda}$ is the quotient of a prehomogeneous vector space $V_{\lambda}$ by the centralizer $H_{\lambda} \coloneq Z_{\wh{G}}(\lambda)$. We call $\V_{\lambda}$ a \emph{Vogan stack}. A Vogan stack has a unique open point corresponding to the parameter $\phi$ with ``maximal monodromy''. We say that the Vogan stack is (essentially) discrete if $\phi$ is. We denote by $\wh{T}_{H,\der}$ a fixed maximal torus of $H_{\lambda}\cap \wh G_\der$.

\begin{theorem}\label{thm:IntroIII}
The following statements hold:
\begin{itemize}
    \item any simple object in $\Cusp(\Par_G)$ is the pushforward along $i_\lambda$ of a vector bundle $\mc{V}$ on an essentially discrete $\V_{\lambda}$ (Theorem \ref{thm: Generators}).
    \item such a simple object is determined by the restriction to $BS_\phi$, and the representation of $S_\phi$ corresponding to the restriction is irreducible and cuspidal in the sense of Lusztig's cuspidal representations (more precisely, see Definition \ref{def: cuspgeneral}).
    \item every irreducible cuspidal representation of $S_\phi$ appears as the restriction of some simple object (Theorem \ref{thm: CVBclassification}).
    \item (in the graded setting) a vector bundle $\mc{V}$ on $\V_{\lambda}$ is the pullback along $\pi$ of a representation $\rho$ of $H_{\lambda}$ and $i_{\lambda\ast}\mc{V} \in \Cusp(\Par_G)$ exactly when the weights of $\rho$ lie in the interior of the polytope $\mc{P}_{\cusp} \subset X^{\ast}(\wh{T}_{H, \der})_{\R}$ given by
      \begin{equation*}
        \mc{P}_{\cusp} = \left(\bigoplus_{\alpha \in \wt(V_{\lambda})} [-\alpha, 0]\right) \bigominus \left(\bigoplus_{\beta \in \Phi(H_{\lambda})} [-\beta, 0]\right),
    \end{equation*}
where $\oplus$ and $\ominus$ denote Minkowski sum and difference and $[-\alpha, 0], [-\beta, 0] \subset X^{\ast}(\wh{T}_{H,\der})_{\R}$ are closed intervals (\S \ref{s: CVBs}).
\end{itemize} 
\end{theorem}

See \S \ref{s: statementsmain} and \S \ref{s: SODs} for further details of the main results. 

\subsection{Some perspectives}

\subsubsection{Classical local Langlands correspondence}
While the categorical local Langlands conjecture may be regarded as a categorification of the classical local Langlands conjecture, the relationship between the two conjectures is much more complex than one might naively expect. For instance, the conjectural coherent sheaf corresponding to an irreducible representation (at a basic point) is difficult to determine in general, and even if possible, there does not seem to be a direct relation to the (enhanced) $L$-parameter. Our Theorem \ref{thm:IntroIII} provides a clean answer for supercuspidal representations, and the relation in this case is as simple as possible in some sense. An interesting feature is that any irreducible cuspidal representation of $S_\phi$ admits a preferred (i.e., cuspidal in our sense) lift to an irreducible vector bundle over the closure, equivalently, an irreducible representation of $H_\lambda$. (For a discrete $L$-parameter associated with a map $\SL_2\to \wh G$, the existence of \emph{some} lift to $H_\lambda$ was already proven in \cite[Theorem 1.1]{Sommers}.)

Let us note that, after we essentially found the statement of Theorem \ref{thm:IntroIII}, we learned that Zhu also proves that an irreducible supercuspidal representation is sent to a vector bundle as above under his tame categorical equivalence \cite{ZhuTame}. 

\subsubsection{Generalized coherent Springer sheaves}
Beyond the supercuspidal case, our construction of compact ``projective'' generators (or rather a slightly refined fiberwise version) provides a further link to the classical local Langlands conjecture. We call them \emph{generalized coherent Springer sheaves} as they generalize those in \cite{BCHNspringer}, and, continuing the philosophy in \cite{BCHNconstructible}, we relate generalized coherent Springer sheaves and constructible Springer sheaves via linear Koszul duality and filtered $\mc D$-modules when they are expected to correspond to representations of pure inner forms. In fact, we use such a relation in the proof of our main results to access results of Lusztig. The generalized constructible Springer sheaves, or just simple perverse sheaves or $\mc D$-modules on Vogan stacks, have a clear relation with the classical local Langlands correspondence for pure inner forms via Vogan's conjecture \cite{VoganLLC}. In turn, generalized coherent Springer sheaves can be related to the classical local Langlands correspondence. A more detailed study will be worked out in a sequel. 

On a related note, generalized coherent Springer sheaves are expected to satisfy several properties, generalizing those in \cite{BCHNspringer, Propp, ZhuTame}, cf.~\cite{Beijingnotes}.\footnote{In these references, generalized coherent Springer sheaves may refer to a slightly different set of compact generators, but the difference is negligible in many discussions.} Theorems \ref{thm:IntroI} and \ref{thm:IntroII} contain results along this line, particularly the claim that the generalized coherent Springer sheaves have no nontrivial endomorphisms in nonzero degrees. 

\subsubsection{The perverse $t$-structure and generalized eigensheaves}
Our original motivation to study generalized coherent Springer sheaves comes from the spectral (sheared) ``eigensheaves'' we constructed for local $A$-parameters \cite{koshikawaoberwolfach2024}. In particular, assuming the categorical conjecture, we expect that they are perverse on the geometric side and their stalks at the neutral point are related to the $p$-adic Adams--Barbasch--Vogan packets \cite{CFMMX}. For instance, the perversity can be formulated using the semiorthogonal decomposition and compact projective generators, and hence it can be translated to the spectral side using Theorems \ref{thm:IntroI} and \ref{thm:IntroII}. We hope to establish these properties directly on the spectral side. 

\subsubsection{Proof of the categorical conjecture}
Not surprisingly, all our main results are useful and essential to \emph{proving} the categorical conjecture. David Hansen and Lucas Mann have formulated a program to prove the categorical local Langlands conjecture \cite{HansenMann}, making some assumptions including a certain compatibility of Eisenstein series on both sides. Theorems \ref{thm:IntroI}, \ref{thm:IntroII} (or their weaker forms) are used to establish an \emph{induction principle} \cite[Theorem 6.8.1]{HansenMann} (namely, induction on Levi subgroups). For classical groups, a further input is the ``generic categorical Langlands correspondence'' of \cite{HelmSolleveldXu}. The final step of their argument for a large class of classical groups is about cuspidal sheaves \cite[\S 8.3]{HansenMann}: they are able to show that the functor between cuspidal objects is fully faithful and invoke Theorem \ref{thm:IntroIII} (together with the known classification of irreducible supercuspidal representations) to conclude that it is an equivalence. 

As a technical note, \cite{HansenMann} also establishes several spectral results. But our paper is largely independent of theirs; we can reprove some of their spectral results. One common ground is the idea of the ``spectral geometric lemma'' that has already appeared in \cite{GLCIII}; however, we work it out in a different (i.e., Koszul dual) setting. 

\subsection{Idea of the proofs}
The proofs of Theorems \ref{thm:IntroI}, \ref{thm:IntroII} and \ref{thm:IntroIII} are long and complex. They involve several reductions to simpler situations or to related settings. 

Most of our main results and their variants roughly fall into the following three types:
\begin{description}
    \item[(1)] Construction of cuspidal/basic objects.
    \item[(2)] Hom computations such as the semiorthogonality relations and the full faithfulness of Eisenstein series. 
    \item[(3)] Generation of all compact objects in the ambient large category by cuspidal/basic objects and their Eisenstein series.
    \item[(3)'] Determination of compact objects of the categories of cuspidal/basic objects. 
\end{description}
(3) and (3)' are tightly related, and either one, together with (2), would imply the semiorthogonal decompositions. 

\subsubsection{Cotangent stacks}
The main body of the paper takes place in the setting of cotangent stacks of Vogan stacks. Fix a semisimple $L$-parameter $\lambda\colon W_F\to \wh G$ and let $V_\lambda$ denote the Vogan variety, the space of $N\in \Lie \wh G$ such that $(\lambda, N)$ is a Weil--Deligne parameter. Then the quotient $\V_\lambda$ of $V_\lambda$ by the natural action of $H_\lambda$ is the Vogan stack. 
Logically speaking, we first prove analogous results for the cotangent stacks $\T^*\V_\lambda$. In general, these are highly derived, though $\T^{\ast}\V_{\lambda}$ is classical if $G$ is semisimple and $\lambda$ underlies an essentially discrete $L$-parameter. However, there is always a concrete realization
\[
\T^*\V_\lambda\simeq\mu^{-1}(0)/H_\lambda, \qquad\mu\colon \T^*V_\lambda\to  (\Lie H_\lambda)^{\vee},
\]
where $\mu$ is the moment map and $\mu^{-1}(0)$ is the \emph{derived} $0$-fiber of $\mu$.

Assume now that $G$ is split, $\lambda$ is unramified, and $V_\lambda$ and $H_\lambda$ are determined by a cocharacter $\Gm\to \wh G$ in a certain way; we call such a situation the graded setting, and we reduce many things to the graded setting (\S \ref{s: cohreductiontograded}). It turns out that the semiorthogonal decomposition along the lines of Theorem \ref{thm:IntroI} in this setting is closely related to what algebraic geometers have studied (noncommutative resolution of singularities, Donaldson--Thomas theory, quasi-BPS categories), especially in the form of \cite{BuPadurariuToda},\footnote{It goes back to \cite[Proposition 7.4]{SpenkoVandenBerghSOD}.} in which context Eisenstein series are called (categorical) Hall products or induction; see \S \ref{ss: cotangenbybasic}. More precisely, if $G$ is semisimple, the category of basic coherent sheaves identifies with the category of graded matrix factorizations $\MF^{\gr}(\M((\T^*V_{\lambda}\times \Lie H_\lambda)/H_{\lambda}), f)$ of the magic\footnote{It essentially appears first in \cite{SpenkoVandenBergh}, and is called the magic window in \cite{HalpernLeistnerSam}; see \cite[footnote 14]{PTDolbeault}. The  window category has a physical origin \cite{HerbstHoriPage}.} category $\M((\T^*V_{\lambda}\times \Lie H_\lambda)/H_\lambda)$ under the Koszul equivalence (\S \ref{ss: limcat})
\[
\Coh(\T^*\V_\lambda) \simeq \MF^{\gr}(\Coh((\T^*V_{\lambda}\times \Lie H_\lambda)/H_{\lambda}), f). 
\]
This subcategory of $\Coh(\T^{\ast}\V_{\lambda})$ is called the limit\footnote{This is meant to be a classical limit of the category of $\mc D$-modules on $\V_\lambda$.} category $\mrL(\T^*\V_\lambda)$  by P\u{a}durariu--Toda \cite{PTDolbeault}, who generalize it to a more functorial notion beyond these quotient stacks. These results establish (2) and (3)/(3)' in the basic case. Interestingly, \cite[\S 4.4.3, \S 4.4.4, \S 4.4.5]{BuPadurariuToda} establish analogous decompositions for the spectral side of the Dolbeault/de Rham/Betti geometric Langlands conjecture respectively (see also \cite[Theorem 1.1]{PTHiggs}). 

To treat the cuspidal case, (1) is crucial. One essentially needs to describe the category of basic sheaves (or the limit category) using cuspidal sheaves and their Eisenstein series. There are at least three ingredients here:
\begin{enumerate}
    \item[(a)] For constructible sheaves or $\mc D$-modules, this is precisely the theory established by Lusztig \cite{LusztigGraded}. The generalized constructible Springer sheaves generate all sheaves and all simple cuspidal objects are classified by him \cite{LusztigIntersection}, and they are clean extensions of local systems $\mc L_\rho$ on the open orbit corresponding to certain irreducible representations $\rho$ of $S_\phi$, called cuspidal. 
    \item[(b)] It is possible to construct vector bundles $\mc V_\rho$ on $\T^*\V_\lambda$ that are cuspidal in a sense analogous to Definition \ref{def: IntroCusp} corresponding to cuspidal $\rho$ based on our earlier computation in \cite{BMKL} with Jonathan Leake, and they lift/deform to filtered $\mc D$-modules ``trivially'' in the spirit of \cite{BCHNconstructible}. 
    \item[(c)] These cuspidal vector bundles are clean extensions in the sense of P\u{a}durariu--Toda \cite{PTDolbeault}: they\footnote{We need a slight variant in our setting, but our argument is largely inspired by theirs.} construct two non-standard pushforward functors from the open orbit taking values in $\mrL (\T^*\V_\lambda)$, and $\mc V_\rho$ can be obtained as both types of pushforward of $\rho$ as a coherent sheaf on $BS_\phi \simeq \T^* BS_\phi$. 
\end{enumerate}
The items (a) and (b) solve (1) in the cuspidal case, but (c) is crucial for (1) in the basic case and (3). P\u{a}durariu and Toda rely on the work of Halpern-Leistner on derived $\Theta$-stratifications \cite{HalpernLeistner}. The upshot is that (c), or more generally \cite{HalpernLeistner}, enables us to run an argument reminiscent of those for constructible sheaves, and (b), via an inductive argument, allows us to compare generalized coherent Springer sheaves on $\T^*\V_\lambda$ and constructible Springer sheaves on $\V_\lambda$, generalizing \cite{BCHNconstructible}, so that we can transfer Lusztig's theory in (a) to the coherent side. We can combine them to establish (1) in the basic case and (3). 

We use a version of the spectral geometric lemma on cotangent stacks for (2); it is used for other properties of basic sheaves and crucial for the conservativity of Eisenstein series on cuspidal sheaves. 

\subsubsection{Moduli of $L$-parameters}
Once we establish the theory for cotangent stacks, we can deduce the results for moduli of $L$-parameters. The first key tool is the linear Koszul duality of Mirkovi\'{c}--Riche \cite{MRI, MRII, MRIwahoriMatsumoto}, up to shearing, 
\[
\Coh(\T^*\V_\lambda)^{\Gm}\simeq \wh\Coh{}^!(\wh{\T}\V_\lambda [-1])^\Gm
\]
in the form developed and used in \cite{Chen2023categorical} and \cite{BCHNconstructible}. Here, $\wh{\T}\V_\lambda [-1]$ is the completed odd tangent stack (completed along the zero section), and by \cite{BenZviNadler}, it is equivalent to the formal loop space
\[
\exp\colon \wh{\T}\V_\lambda [-1]\simeq \wh{\mc L}(\V_\lambda)\coloneq {\mc L}(\V_\lambda)^\wedge_{\V_\lambda}.
\]
The notation $\wh\Coh{}^!$ denotes the category of continuous ind-coherent sheaves, a good replacement for coherent sheaves in the formal setting. 
Note that $\V_\lambda$ is naturally a closed substack of $\Par_G$, sitting inside the fiber of $[\lambda]$ as a point of the coarse moduli space. 
It is essentially observed in \cite{BCHNconstructible} that
\[
\Par^\wedge_{G,\V_\lambda} \simeq \wh{\mc L}(\V_\lambda),
\]
while the completion along the fiber has the description
\[
\Par^\wedge_{G,\lambda} \simeq \mc L^u(\V_\lambda)\coloneq \Map(B\Ga, \V_\lambda), 
\]
$\mc L^u(\V_\lambda)$ being known as the unipotent loop space. 
The second key tool/fact is that the $!$-restriction
\[
\wh{\Coh}{}^!(\mc L^u(\V_\lambda)) \to \wh{\Coh}{}^!(\wh{\mc L}(\V_\lambda))
\]
is an equivalence, as a consequence of coherent completeness established in \cite{AlperHallRydhLuna, AlperHallRydhLocal, AlperHallLim}. These suggest that $\Coh(\T^*\V_\lambda)$ with $\lambda$ varying captures most of the information in $\Coh(\Par_G)$.  

Using these tools, (1), (2), and (3)' in the cuspidal setting are easy to transfer from the case of cotangent stacks; a fortunate fact is that, if $G$ is semisimple, every simple cuspidal sheaf on $\Par_G$ is supported on one Vogan stack $\V_\lambda$, and in fact, it is a vector bundle on $\V_\lambda$. This is the fundamental reason why our approach via cotangent stacks works. 

\begin{remark}
As outlined above, our argument ultimately reduces the classification of cuspidal sheaves to the setting of $\mc D$-modules or constructible sheaves via linear Koszul duality and ``deformation'' from coherent sheaves on the cotangent stack. We learned this fundamental picture from \cite{BCHNconstructible}. The idea that representations of pure inner forms are Koszul dual to $\mc D$-modules (or constructible sheaves) on the Langlands dual side goes back further to Soergel \cite{Soergel}, following the pioneering works of Adams--Barbasch--Vogan and Vogan \cite{ABV,VoganLLC}. 
\end{remark}
\subsection{Guide for the reader}
A reader looking to understand the main results of the paper should read the introduction, followed by \S \ref{s: cuspsheavesdef} and then \S \ref{s: statementsmain}, referring to the first two parts as necessary. 

Part \ref{p: preliminaries} gives an overview of the derived algebraic geometry and theories of coherent sheaves used in the paper. Part \ref{p: moduliofLparams} contains some new results on the moduli of Weil--Deligne parameters and Eisenstein series functors, but the reader may choose to skip it at first pass, referring back as necessary. Part \ref{p: cuspidalsheaves} is the heart of the paper, where our main results are proven. After reductions to the unramified graded setting and to cotangent stacks in \S \ref{s: cuspsheavesdef} and \S \ref{s: cohreductiontograded}, the key properties of the categories of basic and cuspidal sheaves are obtained in \S \ref{s: cotangentgeometry}, \S \ref{s: Springer}, and \S \ref{s: SODs} while \S \ref{s: CVBs} explicitly classifies cuspidal vector bundles.  A spectral geometric lemma for cotangent stacks is established in \S \ref{s: spectralgeometriclemma} and it is used to study spectral Eisenstein series and generalized coherent Springer sheaves. Finally, the construction of semiorthogonal decompositions is completed for cotangent stacks and $\Par_G$ in \S \ref{s: SODs}.

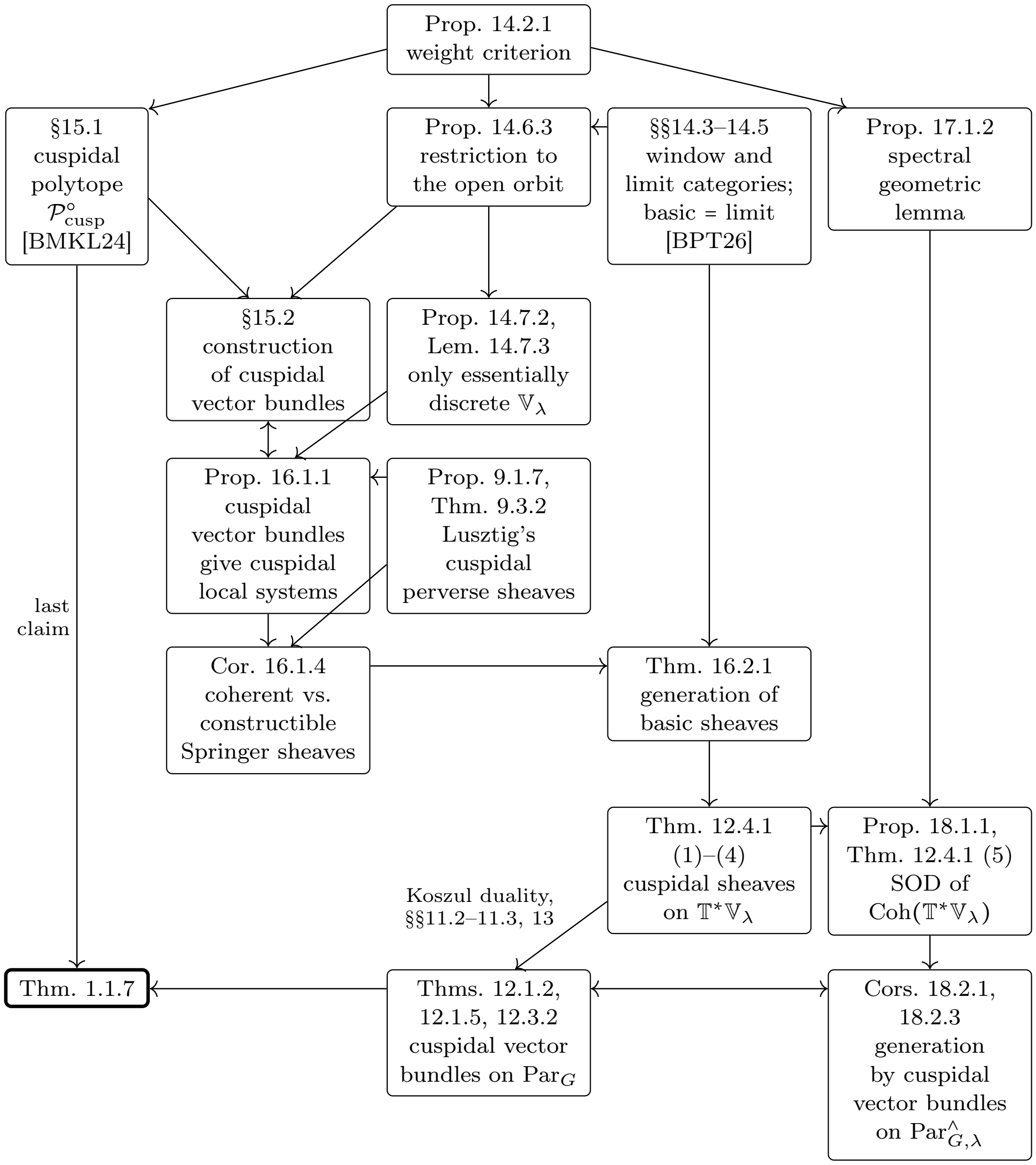
\begin{figure}[htbp]
\centering
\begin{tikzcd}[math mode=false, column sep=0.2cm, row sep=0.4cm,
  cells={nodes={draw, rounded corners=2pt, align=center, text width=2.2cm,
    font=\footnotesize, execute at begin node={\hyphenpenalty=10000\exhyphenpenalty=10000\relax}}}]
  & & Prop.~\ref{prop: basicweightcharacterization}\par weight criterion
      \arrow[dll] \arrow[d] \arrow[drr] & & \\
  |[text width=1.45cm]| \S\ref{ss: fundcuspineq}\par cuspidal polytope $\mc{P}^\circ_{\cusp}$ \cite{BMKL}
      \arrow[dr]
      \arrow[ddddd, swap, "{last\\ claim}"{font=\scriptsize, align=right}]
  & & Prop.~\ref{prop: cuspextopen}\par restriction to the open orbit
      \arrow[dl] \arrow[d]
  & \S\S\ref{ss: windowcat}--\ref{ss: cotangenbybasic}\par window and limit categories; basic $=$ limit \cite{BuPadurariuToda}
      \arrow[l] \arrow[ddd]
  & Prop.~\ref{prop: generaleisbasisbas}\par spectral geometric lemma
      \arrow[dddd] \\
  & \S\ref{ss: classificationunipotentcase}\par construction of cuspidal vector bundles
      \arrow[d, leftrightarrow]
  & Prop.~\ref{prop: cuspreductiontofinite}, Lem.~\ref{lem: zerosectiondirectsummand}\par only essentially discrete $\V_\lambda$
      \arrow[dl]
  & & \\
  & Prop.~\ref{prop: CVBimpliesCLS}\par cuspidal vector bundles give cuspidal local systems
      \arrow[d]
  & Prop.~\ref{prop: OrthogonalCuspidalConstructibleSheaves}, Thm.~\ref{thm: ClassificationCuspidalPerverseSheaves}\par Lusztig's cuspidal perverse sheaves
      \arrow[l] \arrow[dl]
  & & \\
  & Cor.~\ref{cor: coherentspringersummand}\par coherent vs.\ constructible Springer sheaves
      \arrow[rr]
  & & Thm.~\ref{thm: cotanbasgenbycuspeis}\par generation of basic sheaves
      \arrow[d]
  & \\
  & & & Thm.~\ref{thm: ClassificationCuspidalCotangent} (1)--(4)\par cuspidal sheaves on $\T^*\V_\lambda$
      \arrow[r]
      \arrow[dl, swap, "{Koszul duality,\\ \S\S\ref{ss: reductiontocompletedfiber}--\ref{ss: reductiontocotangent}, \ref{s: cohreductiontograded}}"{font=\scriptsize, align=right}]
  & Prop.~\ref{prop: refinedcotanSOD}, Thm.~\ref{thm: ClassificationCuspidalCotangent} (5)\par SOD of $\Coh(\T^*\V_\lambda)$
      \arrow[d] \\
  |[very thick, text width=1.45cm]| Thm.~\ref{thm:IntroIII}
  & & Thms.~\ref{thm: discretesupport}, \ref{thm: CVBclassification}, \ref{thm: Generators}\par cuspidal vector bundles on $\Par_G$
      \arrow[ll]
  & & Cors.~\ref{cor: SODhatCoh}, \ref{cor: CohCompletionSOD}\par generation by cuspidal vector bundles on $\Par^\wedge_{G,\lambda}$
      \arrow[ll, leftrightarrow]
\end{tikzcd}
\caption{Leitfaden for the proof of Theorem~\ref{thm:IntroIII}. An arrow $A\to B$ means that $A$ is used in the proof of $B$; not every such use is drawn. A two-headed arrow means that the two results are proven together: Proposition~\ref{prop: CVBimpliesCLS} and \S\ref{ss: classificationunipotentcase} by induction on the semisimple rank; Theorems~\ref{thm: discretesupport}, \ref{thm: CVBclassification}, \ref{thm: Generators} and Corollaries~\ref{cor: SODhatCoh}, \ref{cor: CohCompletionSOD} in alternation, since the cuspidal vector bundles of Theorem~\ref{thm: CVBclassification}\,(1) enter Corollaries~\ref{cor: SODhatCoh} and~\ref{cor: CohCompletionSOD}, which then give the remaining parts of Theorem~\ref{thm: CVBclassification} and Theorem~\ref{thm: Generators}.}
\label{fig: leitfaden}
\end{figure}

Figure~\ref{fig: leitfaden} summarizes the proof of Theorem~\ref{thm:IntroIII}.
 
\FloatBarrier

\subsection*{Acknowledgments}
The authors would like to thank David Hansen, whose ideas on $\ell$-adic categorical Langlands have been highly influential. We also thank him for spurring us to prove our results on semiorthogonal decompositions, for many mathematical conversations on this work, and for consistently encouraging us to complete our project. We would like to thank our collaborator on \cite{BMKL}, Jonathan Leake, whose computer computations with exceptional groups are used in \S \ref{sss: exceptionalgroups}. We thank Tasuki Kinjo who informed us of the work of P{\u{a}}durariu, Toda, and their collaborators, which turned out to be crucial for this project. 

We thank Noriyuki Abe, Anne-Marie Aubert, Ana Caraiani, Matthew Emerton, Linus Hamann, Naoki Imai, Adeel Khan, Lucas Mann, Masao Oi, Tudor P\u{a}durariu, David Schwein, Sug Woo Shin, Kari Vilonen, Ting Xue, Alex Youcis, Mingjia Zhang, and Xinwen Zhu for helpful conversations related to this project.

The first named author was partially supported by NSF grant DMS-2502131 during the completion of this work. The second named author was supported by JSPS KAKENHI Grant Numbers 24K16895 and 22H00093. 

Work on this project began in 2023, and no AI was used until September, 2026 when the project was in an advanced state. Since that time, the authors had free access to Anthropic models Opus 5 and Opus 5.5 via the Claude team plan for scientists and these were extensively used to proofread the manuscript. In several instances, such as Lemma \ref{lem: ParGdomswitch}, Proposition \ref{prop: refinedcotanSOD}, Lemma \ref{lem: zerosectiondirectsummand} an LLM identified flawed arguments and gaps, helping the authors fully spell out and correct the write-up. However, the entire document except for the leitfaden in Figure~\ref{fig: leitfaden} was written manually. We also used OpenAI model GPT-6 Astra via ChatGPT for Academic Researchers for further proofreading. The authors take full responsibility for the correctness of the results in this article.

\part{Preliminaries}{\label{p: preliminaries}}

\section{Notation}{\label{s: notation}}

Fix a prime number $p$, a $p$-adic local field $F$ and a choice of algebraic closure $\ov{F}$, and let $\Gamma$ denote the Galois group $\Gal(\ov{F} / F)$. We let $W_F$ denote the Weil group of $F$ and fix a lift, $\frob$, of \emph{arithmetic} Frobenius. We use the convention that $\lvert \frob \rvert =q$, where $q$ is the size of the residue field of $F$. 

We assume throughout that our coefficient ring $\Lambda$ is an algebraically closed field of characteristic $0$. We adopt the notation that all functors are derived unless otherwise specified. In particular, $\Hom$ denotes $\RHom$; in the $\Lambda$-linear setting, it also denotes the enriched mapping complex, which is not necessarily connective. Affine algebraic groups are assumed to be classical (of finite type), and they are smooth since we will work only in characteristic 0.

\subsection{Recollections on \texorpdfstring{$L$}{L}-groups}
Unless otherwise indicated, $G$ denotes a quasi-split connected reductive group over $F$. Since we work entirely on the dual side and inner twists induce isomorphisms on dual groups, we could just as well assume that $G$ is an arbitrary connected reductive group. However, we keep the quasi-split assumption for simplicity.  We implicitly fix a pinning $(B, T, \{x_{\alpha}\}_{\alpha \in \Delta})$. We let $\wh{G}$ denote the dual group of $G$ with respect to our pinning and $\co{L}{G} \coloneq \wh{G} \rtimes W_F$ the Weil form of the $L$-group. It comes with a $W_F$-stable Borel pair $(\wh{B}, \wh{T})$ and a pinning. 

We quickly recall the parabolic and Levi subgroups of $L$-groups from \cite[\S 3.3, 3.4]{Bor2}, \cite[\S 2.5]{SilbergerZink}. A subgroup of $\co{L}{G}$ is parabolic if it surjects onto $W_F$ and its intersection with $\wh{G}$ is a parabolic subgroup. A parabolic subgroup is called standard (resp.~semistandard) if it contains $\wh{B}\rtimes W_F$ (resp.~$\wh{T}\rtimes W_F$). A standard (resp.~semistandard) parabolic subgroup has the form $\wh{P}\rtimes W_F$, where $\wh{P}$ is a $W_F$-stable standard (resp.~semistandard) parabolic subgroup of $\wh{G}$. Every parabolic subgroup is $\wh{G}$-conjugate to a unique standard parabolic subgroup.  

A Levi subgroup of ${}^L G$ is a Levi complement of the unipotent radical of some parabolic subgroup of $\co{L}{G}$. A Levi subgroup is standard (resp.~semistandard) if it is the unique Levi complement of the unipotent radical of some standard (resp.~semistandard) parabolic subgroup of $\co{L}{G}$ containing $\wh{T}\rtimes W_F$. Any Levi subgroup is $\wh{G}$-conjugate to at least one standard Levi subgroup. By \cite[Lemma 3.5]{Bor2}, for any Levi subgroup $\wh M\rtimes W_F$, we have $\wh M\rtimes W_F =Z_{\co{L}{G}}(Z(\wh M)^{\Gamma,\circ})$.  

Given a semistandard parabolic (resp.~Levi) subgroup $P$ (resp.~$M$) of $G$ over $F$, we let $\wh{P}$ (resp.~$\wh{M}$) denote the corresponding parabolic (resp.~Levi) subgroup of $\wh{G}$ such that $\wh{P}\rtimes W_F$ (resp.~$\wh{M}\rtimes W_F$) is a semistandard parabolic (resp.~Levi) subgroup \cite[Proposition 2.5.4 (i)]{SilbergerZink}. This correspondence preserves the standard ones. 

Given a maximal torus $T \subset G$, we denote the Weyl group by $W(T,G) \coloneq N_G(T)/T$.  When the choice of maximal torus is clear, we denote $W(T,G)$ by $W_G$. We use analogous notation for $\wh{G}$ and we have a $\Gamma$-equivariant isomorphism $W_{G} \cong W_{\wh{G}}$ induced by our fixed pinnings. Let $A_T \subset T$ be the maximal split subtorus and let $A_{\wh{T}} \subset \wh{T}$ be the maximal $\Gamma$-fixed subtorus. Then we define $W_G(A_T) \coloneq N_G(A_T)/Z_G(A_T)$ and analogously for $\wh{G}$ and have identifications
\begin{equation*}
    W_G(A_T) \cong W^{\Gamma}_G \cong W^{\Gamma}_{\wh{G}} \cong W_{\wh{G}}(A_{\wh{T}}).
\end{equation*}
The above bijection between the sets of semistandard parabolic subgroups of $G$ and $\wh{G}$ is equivariant with respect to these identifications of Weyl groups \cite[Proposition 2.5.4 (ii)]{SilbergerZink}. 

We denote by $B(G)$ the Kottwitz set defined as $G(\breve{F})$ modulo $\frob$-conjugacy. The set $B(G)$ admits maps
\begin{equation*}
    \kappa\coloneq\kappa_G\colon B(G) \rightarrow X^{\ast}(Z(\wh{G})^{\Gamma}), \qquad \nu_G: B(G) \rightarrow \ov{C}_{\Q} \subset X_{\ast}(A_T)_{\Q},
\end{equation*}
where $\ov{C}_{\Q}$ denotes the closed Weyl chamber of $X_{\ast}(A_T)_{\Q}$. We say that $b \in B(G)$ is \emph{basic} if $\nu_G(b)$ is central in $G$ and recall that $\kappa_G$ induces an isomorphism $B(G)_{\bas} \cong X^{\ast}(Z(\wh{G})^{\Gamma})$ (\cite[Proposition 5.6]{kot8}), where $B(G)_{\bas}$ denotes the subset of basic elements. For $M \subset P \subset G$ a semistandard Levi and parabolic subgroup, we denote by $X^{\ast}(Z(\wh{M})^{\Gamma})^{+} \subset X^{\ast}(Z(\wh{M})^{\Gamma})$ (resp.~$X^{\ast}(Z(\wh{M})^{\Gamma})^{\circ +}$)   the set of (resp.~strictly) $P$-dominant elements. Namely, those that pair non-negatively (resp.~positively) with the roots of $A_M$ in $\Lie U_P$ under the map
\begin{equation*}
    X^{\ast}(Z(\wh{M})^{\Gamma}) \to X^{\ast}(Z(\wh{M})^{\Gamma})_{\Q} \cong X_{\ast}(A_M)_{\Q}.
\end{equation*}
We let $B(M)^+ \subset B(M)$ denote the subset of $b \in B(M)$ such that $\kappa_M(b) \in X^{\ast}(Z(\wh{M})^{\Gamma})^+$ and analogously for $B(M)^{\circ +}$. For $b, b' \in B(G)$, we say $b \leq b'$ if $\kappa_G(b) = \kappa_G(b')$ and $\nu_G(b') - \nu_G(b)$ is a non-negative $\Q$-combination of positive coroots. We call this the ascending order and it satisfies that basic elements are minimal.
 
\section{Derived algebraic geometry}{\label{s: DAG}}
We will give an informal overview of derived algebraic geometry and omit fully precise formulations. Let $\Lambda$ be an algebraically closed field of characteristic $0$. We sometimes also consider a (classical) noetherian ring $\Lambda'$ over $\Lambda$ as a base ring; we will omit the obvious compatibility between different bases. 

\subsection{Basic local definitions}{\label{ss: DAGbasiclocaldefns}}
We recall relevant notions in derived algebraic geometry over $\Lambda$. In this setup, one may work with connective commutative differential graded $\Lambda$-algebras; introductory texts are \cite{EugsterPridham, ToenDAG}. This is the approach taken in the usual geometric Langlands program such as \cite{GaitsgoryIndCohSheaves, GRInd, GRI}, while these treatments also rely on foundational work of Lurie on higher category theory \cite{LurieHTT, LurieHA}. Over more general rings, the theory is developed using simplicial commutative $\Lambda$-algebras, or animated (commutative) $\Lambda$-algebras \cite{CesnaviciusScholze} in modern terms; these ideas originated in Lurie's thesis \cite{LurieThesis} and works of To\"en--Vezzosi \cite{TVI, TVII}, and \cite{HansenMann} is a useful reference too. We freely use the equivalence between the underlying $\infty$-categories of connective commutative differential graded $\Lambda$-algebras and animated $\Lambda$-algebras via the monoidal Dold--Kan correspondence. 

In the rest of the paper, a ring or $\Lambda$-algebra means an animated ring or $\Lambda$-algebra. We call the usual rings and $\Lambda$-algebras classical (these are called static in the condensed mathematics of Clausen--Scholze). A ring $A$ is noetherian if $\pi_0(A)$ is noetherian and $\pi_i(A)$ is a finitely generated $\pi_0(A)$-module for all $i$. Let $\Alg_{\Lambda}$ denote the $\infty$-category of $\Lambda$-algebras. 

Let $f\colon A\to B$ be a map of $\Lambda$-algebras. The relative cotangent complex is well-defined as an animated $B$-module. By \cite[Corollary 3.2.17]{LurieThesis}, $f$ is an equivalence if $\pi_0(f)$ is an isomorphism and $\bb L_{B/A}=0$. 

The map $f$ is almost of finite presentation if $\tau_{\leq n} B$ is compact in the $\infty$-category of $n$-truncated $A$-algebras for all $n$; if $A$ is noetherian, it is equivalent to say that $\pi_n(B)$ is a finitely generated $\pi_0(B)$-module for all $n$ and $\pi_0(B)$ is a finitely generated $\pi_0(A)$-algebra in the classical sense \cite[Proposition 3.1.5]{LurieThesis}. More generally, it is equivalent to \cite[Definition 2.1.6]{HansenMann} by Remark 2.1.8 therein. The map $f$ is flat (resp.~faithfully flat) if $\pi_0(f)$ is flat (resp.~faithfully flat) and $\pi_i(A)\otimes_{\pi_0(A)}\pi_0(B)\cong \pi_i(B)$ for all $i$. The map $f$ is fppf if it is faithfully flat and almost of finite presentation.\footnote{Or fpppf might be better but fppf is commonly used.} We will consider the fppf topology as in e.g.~\cite{DrinfeldGaitsgory}. The default topology in \cite{GRI} is the \'etale topology, but everything works also for the fppf topology \cite[\S 2.2.3.2]{GRI}. 

Recall also that the map $f$ is \'etale if $\bb L_{B/A}=0$ and it is almost of finite presentation \cite[Definition 3.4.7]{LurieThesis}. Equivalently, $f$ is flat and $\pi_0(f)$ is \'etale in the classical sense \cite[Proposition 3.4.10]{LurieThesis}. In fact, the functor $B\to \pi_0(B)$ induces an equivalence of $\infty$-categories between \'etale $A$-algebras and \'etale $\pi_0(A)$-algebras \cite[Theorem 3.4.13]{LurieThesis} (or \cite[Proposition 5.2.4]{CesnaviciusScholze}). An \'etale map $f$ is an open immersion if $\pi_0(f)$ is an open immersion. Similarly, there is a notion of smooth morphisms, and $f$ is smooth if and only if $f$ is flat and $\pi_0(f)$ is smooth in the classical sense \cite[Proposition 3.4.9]{LurieThesis}, in which case $\bb L_{B/A}$ is a finite projective $B$-module. Weakening the last condition, we say that the map $f$ is quasi-smooth if it is almost of finite presentation and $\bb L_{B/A}$ has Tor amplitude in $[-1, 0]$ (these conditions imply $\bb L_{B/A}$ is a perfect complex). 

The property of a map being quasi-smooth (resp.~smooth, resp.~\'etale, resp.~an open immersion, resp.~almost of finite presentation) satisfies flat descent \cite[Proposition 3.5.3]{LurieThesis}, and it is \'etale-local on the source except in the case of open immersions, cf.~\cite[Lemma 2.4.15]{HansenMann}.

\subsection{Algebraic stacks}
A (derived) prestack is an (accessible\footnote{We will ignore the set-theoretic issues.}) functor $\Alg_{\Lambda}\to \Ani$. A (derived) stack is a prestack that is a sheaf with respect to the fppf topology. Let $\Stk$ denote the $\infty$-category of stacks. The fppf topology is subcanonical \cite[Example 4.2.4]{LurieThesis}, and the Yoneda functor induces a fully faithful functor $\Alg_{\Lambda}^\op\to \Stk$. An object in the essential image is called a (derived) affine scheme. For a $\Lambda$-algebra $A$, its image is denoted by $\Spec A$. Its underlying classical affine scheme is $\Spec \pi_0 (A)$ with a canonical map $\Spec \pi_0 (A)\to \Spec A$. The notion of open immersions makes sense for maps of stacks, and the notion of open covers makes sense as well. A stack is a (derived) scheme if it admits an open cover by affine schemes. For any derived scheme $X$, we have the underlying classical scheme $X_{\cl}$ with a canonical map $X_{\cl}\to X$. Similarly, the notion of \'etale (resp.~smooth, resp.~quasi-smooth) maps makes sense for maps representable by schemes, and a stack $X$ is an algebraic space if it admits an \'etale cover by (affine) schemes and its diagonal is representable by a scheme, and we have the underlying classical algebraic space $X_{\cl}$. The above classes of maps extend to maps representable by algebraic spaces. Define the underlying topological space $\lvert X\rvert \coloneq \lvert X_{\cl}\rvert$. A scheme or algebraic space $X$ is locally noetherian if it admits an open or \'etale cover by noetherian affine schemes, and it is noetherian if moreover $\lvert X\rvert$ is quasi-compact. 

A stack $X$ is (1-)Artin if there exists a smooth surjective map $U\to X$, called an atlas, from an algebraic space $U$ and the diagonal $X\to X\times_{\Lambda}X$ is representable by an algebraic space. Such an $X$ is quasi-compact (resp.~smooth) if $U$ can be chosen to be quasi-compact (resp.~smooth). An Artin stack is (0-)geometric if the diagonal is affine. A map of Artin stacks $X\to Y$ is smooth if $U\to Y$ is smooth for one, equivalently any, atlas $U\to X$ of $X$.  

\begin{remark}
In \cite{GRI}, the default topology is the \'etale topology. So, the notion of Artin stack is a priori different. As in \cite[Remark 2.4.1.4]{GRI}, their definition of Artin stack satisfies fppf descent by a result of To\"en \cite[Theorem 2.1]{Toen}. 
\end{remark}

\subsection{Almost of finite type stacks}
The theory of ind-coherent sheaves is developed with finiteness conditions. Let $\Lambda'$ be a classical noetherian ring over $\Lambda$. A scheme or algebraic space $X$ over $\Lambda'$ is locally almost of finite type (over $\Lambda'$) if it admits an open or \'etale covering by affine schemes $\Spec A_i$ such that $A_i$ is a $\Lambda'$-algebra almost of finite presentation. Equivalently, it is a scheme or algebraic space that is, on $n$-truncated $\Lambda'$-algebras for any fixed $n$, left Kan-extended from a functor $\Alg_{\Lambda'}^{\afp, \leq n}\to \Ani$, where $\Alg_{\Lambda'}^{\afp, \leq n}$ denotes the $\infty$-category of $n$-truncated $\Lambda'$-algebras almost of finite presentation. 

Recall that a prestack $X$ is nilcomplete if the natural map
\[
X(A) \to \lim_n X(\tau_{\leq n}A ) 
\]
is an equivalence for any $\Lambda$-algebra $A$. (This property is called convergent in \cite{GRI, DrinfeldGaitsgory}.) Algebraic spaces are nilcomplete and more generally, an Artin stack is nilcomplete \cite[Proposition 2.4.4.9]{GRI} and so is its formal completion defined below.

A nilcomplete stack over $\Lambda'$ is called locally almost of finite type (over $\Lambda'$) if its restriction to $n$-truncated $\Lambda'$-algebras is left Kan-extended from a functor $\Alg_{\Lambda'}^{\afp, \leq n}\to \Ani$ for any $n$ \cite[\S 2.1.7.2]{GRI}, and almost of finite type if it is in addition quasi-compact. We may regard prestacks locally almost of finite type as functors from (eventually coconnective) $\Lambda'$-algebras almost of finite type \cite[Proposition 2.1.7.6]{GRI}. An Artin stack is almost of finite type if and only if there is an atlas from an algebraic space almost of finite type \cite[Proposition 2.4.5.2]{GRI}. (Note that \cite{GRI} works only over $\Lambda$, but the argument works for $\Lambda'$.)

Recall also the notion of QCA stacks from \cite[Definition 1.1.8]{DrinfeldGaitsgory}.\footnote{We warn the reader that \cite[Definition 2.6.1]{HansenMann} uses a different definition.} For us, it is enough to know that any scheme (or algebraic space) almost of finite type over $\Lambda'$ and the quotient of a quasi-affine scheme almost of finite type over $\Lambda'$ by an affine algebraic group over $\Lambda$ are QCA. More generally, any geometric Artin stack almost of finite type over $\Lambda'$ is QCA. 

\subsection{Some examples}
Many stacks we will face arise from the following types of stacks and taking fiber products. 

\begin{example}
\begin{enumerate}
    \item (Classifying stacks) Let $G$ be an affine algebraic group over $\Lambda$, which is automatically smooth over $\Lambda$. The classifying stack $BG$ is a geometric Artin stack.
    \item (Quotient stacks) More generally, if $X=\Spec A$ for a $\Lambda$-algebra $A$ almost of finite presentation with an action of $G$, the quotient stack $X/G$ is a geometric Artin stack. If $G^{\circ}$ is reductive, we define $X\sslash G$ to be $\Spec A^G$, which is a good moduli space in the sense of \cite{GoodModuli}. If $X$ is classical, it identifies with the GIT quotient.
    \item (Mapping stacks) For any stacks $X$, $Y$, we can consider the stack $\Map(X, Y)$ of maps $X\to Y$, with the evaluation map $\ev\colon X\times \Map(X,Y)\to Y$. Such stacks are studied in detail in \cite{HalpernLeistnerPreygel}. 
\end{enumerate}    
\end{example}

Here are somewhat more concrete examples that are relevant for us. 

\begin{example}\label{ex: ExamplesDerivedStacks}
Let $X$ be a stack almost of finite type. 

\begin{enumerate}
    \item The loop space of $X$ is defined as 
    \[
    \mc L(X)\coloneq \Map (S^1, X)\simeq X\times_{X\times X}X
    \]
    with a projection $\mc L(X)\to X$ corresponding to restriction to a fixed base point of $S^1$ \cite[Definition 2.1.3]{Chenequivariant}. If $X=Y/G$,  the loop space has the following presentation as a fiber product (\cite[Proposition 2.1.8]{Chenequivariant})
    \begin{equation}{\label{eqn: chenfiberprod}}
  \mc{L}(Y/G) \simeq \left( (Y \times G)\times_{Y \times Y} Y \right) / G,
    \end{equation}
    where the map $Y \rightarrow Y \times Y$ is the diagonal and $Y \times G \rightarrow Y \times Y$ maps $(y,g) \mapsto (y, gy)$. This leads to the informal description:
    \[
    \mc L(Y/G)\simeq \{ (y, g)\in Y\times G \mid g\cdot y=y\} /G \to Y/G,
    \]
    with a canonical map to $\mc L(BG)\simeq G/^\adj G$, where the quotient is taken for the adjoint action. 
    \item The unipotent loop space \cite[Definition 2.1.3]{Chenequivariant} is 
    \[
    \mc L^u (X)\coloneq\Map (B\Ga, X)
    \]
    with a canonical map to $\mc L(X)$ via the affinization map $S^1 = B\Z \to B\Ga$ induced by the map $\Z \rightarrow \Ga$ on groups. There is a natural action of $B\Ga \rtimes \Gm$. This is \emph{not} an Artin stack in general even if $X$ is so, and we will give a description in \S \ref{ss: loopspaces}. 
    \item Let $G, H$ be classical group schemes. The stack $\Map (BG, BH)$ is the quotient of the space $\Hom(G, H)$ of group homomorphisms $G\to H$ by the adjoint action of $H$, cf.~\cite[Proposition 2.1.22]{Chenequivariant} (in loc.~cit., the affineness of $G$ and $H$ is not used). If $G^{\circ}$ is reductive, $\pi_0(G)$ is finite, and $H$ is affine, it is a geometric Artin stack locally almost of finite type by \cite[Example 1.2.4]{HalpernLeistnerPreygel}. 
    
    In the case where $G$ is discrete and $H$ is affine, $\Hom(G, H)$ is shown to be affine in \cite[\S 3.2]{ToenDAG}; see also \cite[\S5]{GalatiusVenkatesh}, \cite[Definition 2.8, Remark 2.9, Proposition 2.10]{Zhu} and \cite[Proposition 3.1.5]{HansenMann}. 
    \item The stack $\mathrm{Grad}(X)\coloneq\Map(B{\Gm}, X)$ is sometimes called the stack of graded points. If $X$ is a geometric Artin stack, then $\mathrm{Grad}(X)$ is a geometric Artin stack locally almost of finite type by \cite[Proposition 4.3.4, Theorem 5.1.1]{HalpernLeistnerPreygel}. 
    \item Let $\Gm$ act on $\bb A^1$ by weight $1$ geometrically. The stack $\mathrm{Filt}(X)\coloneq\Map(\bb A^1 /\Gm, X)$ is sometimes called the stack of filtered points. If $X$ is a geometric Artin stack, then $\mathrm{Filt}(X)$ is a geometric Artin stack locally almost of finite type by \cite[Example 1.2.2]{HalpernLeistnerPreygel}.  
\end{enumerate}
\end{example}

\begin{example}[Fixed points]
Suppose $G$ is an affine algebraic group whose identity component $G^{\circ}$ is a torus, and $G$ acts on a stack $X$. We will consider the stack $X^G$ of ``fixed points''. 
\begin{enumerate}
    \item If $X$ is classical, we have the classical notion of fixed points $X^{G, \cl}\to X$. If $X=\Spec A$ is affine and $G=\Gm$, then $X^{G,\cl}=\Spec A/I$, where $I$ is the ideal generated by the homogeneous elements with nonzero weights. If $A$ is either non-positively or non-negatively graded, then $A/I$ identifies with the homogeneous elements of weight $0$. 
    
    Even when $X$ is not classical, $G$ acts on $X_{\cl}$ and $X_{\red}$, and the notations $X_{\cl}^{G, \cl}$, $X_{\red}^{G, \cl}$ make sense. 
    \item The stack of homotopy fixed points $X^G$ is defined by the following cartesian square
    \[
    \begin{tikzcd}
    X^G \arrow[r] \arrow[d] & \Map(BG, X/G) \arrow[d] \\
    \{\id_{BG}\} \arrow[r] & \Map (BG, BG), 
    \end{tikzcd}
    \]
    see \cite[Appendix A]{AranhaKhanLatyntsevParkRavi}, \cite[\S A.2.2]{BCHNconstructible}. If $X$ is a geometric Artin stack locally almost of finite type, then so is $X^G$ \cite[Example 4.3.5, Theorem 5.1.1]{HalpernLeistnerPreygel}. If $X$ is smooth (resp.~quasi-smooth), so is $X^G$ \cite[Corollary A.30]{AranhaKhanLatyntsevParkRavi}. In particular, if $X$ is a smooth algebraic space, then $X^G$ agrees with the classical fixed points $X^{G, \cl}$ of $G$ and $X^G\to X$ is a closed immersion, cf.~\cite[Proposition A.21]{AranhaKhanLatyntsevParkRavi} and its proof. 
\end{enumerate}
\end{example}

\subsection{Quasi-coherent sheaves}\label{ss: QCoh}
We discuss quasi-coherent sheaves first. Some general references are \cite{BenZviFrancisNadler}, \cite[Ch.~3]{GRI}, \cite{DrinfeldGaitsgory}, \cite[Appendix A]{HalpernLeistnerPreygel}, \cite[\S 9.2]{ZhuTame}, \cite[\S 2.3]{HansenMann}. 

\begin{remark}
We only follow a classical treatment. A modern sophisticated approach uses solid quasi-coherent sheaves of Clausen--Scholze to construct a 6-functor formalism. This is done in \cite[\S 2.7]{HansenMann}, for instance. 
\end{remark}

On affine schemes $\Spec A$, we set $\QCoh(\Spec A)\coloneq\Mod_A$, the symmetric monoidal presentable $\Lambda$-linear stable $\infty$-category of $A$-modules. Inside $\Mod_A$, we have the idempotent-complete $\infty$-category $\Perf_A$ of perfect complexes. Both $\Mod_A$ and $\Perf_A$ satisfy fppf descent \cite[Corollary 1.3.5, Corollary 3.6.5]{GRI}. By right Kan extension, we obtain $\QCoh(X)$ and its subcategory $\Perf(X)$ for any prestack $X$. The category $\QCoh(X)$ has a symmetric monoidal structure and $\Perf(X)$ agrees with the dualizable objects, equivalently, locally compact \cite[Lemma 3.3.6.7]{GRI}.  

For any map $f\colon X\to Y$, we have the symmetric monoidal $*$-pullback functor $f^*\colon \QCoh(Y)\to \QCoh(X)$, $\Perf(Y)\to \Perf(X)$. By Lurie's adjoint functor theorem, the former admits a right adjoint $f_*\colon \QCoh(X)\to \QCoh(Y)$. For any map $f$ whose base change to an affine scheme is QCA, $f_*$ is continuous and it satisfies both the projection formula and base change \cite[Corollary 1.4.5]{DrinfeldGaitsgory}. 

A perfect complex may not be compact in general, while it is always locally compact. By \cite[Theorem 1.3.2]{DrinfeldGaitsgory}, any perfect complex is compact if $X$ is a QCA stack. For such $X$, $\QCoh(X)$ is compactly generated if and only if $\Ind(\Perf(X))\simeq \QCoh(X)$ canonically by the proof of \cite[Proposition 3.9]{BenZviFrancisNadler}. If $X$ has an affine diagonal and $\Ind(\Perf(X))\simeq \QCoh(X)$ holds, $X$ is called perfect \cite[Definition 3.2]{BenZviFrancisNadler}. 

If $X$ is the quotient of an affine scheme almost of finite type over $\Lambda'$ by an affine group $G$, then $X$ is perfect by \cite[Proposition 3.21]{BenZviFrancisNadler}. Let us remark that any quasi-compact algebraic space with affine diagonal is also known to be perfect. 

\begin{lemma}
Let $X$ be a perfect stack. 
\begin{enumerate}
    \item The fiber products of perfect stacks are perfect. In particular, the loop space $\mc L(X)$ is perfect. 
    \item For any affine map $Y\to X$, $Y$ is perfect. In particular, the total space of any derived vector bundle over $X$ is perfect.  
\end{enumerate}
\end{lemma}

\begin{proof}
(1) follows from \cite[Proposition 3.24]{BenZviFrancisNadler}. (2) follows from \cite[Proposition 3.21]{BenZviFrancisNadler}. 
\end{proof}

\begin{proposition}[{\cite[Theorem 1.2]{BenZviFrancisNadler}}]\label{prop:QCoh_satisfies_Kunneth}
For any maps $X\to Y\leftarrow X'$ of perfect stacks, there is a natural equivalence
\[
\QCoh(X\times_Y X') \simeq \QCoh(X)\otimes_{\QCoh(Y)}\QCoh(X'), 
\]
where $\otimes$ is the Lurie tensor product of symmetric monoidal presentable $\Lambda$-linear stable $\infty$-categories. 
\end{proposition}

Recall that there is a canonical $t$-structure on $\QCoh(X)$, at least for Artin stacks. If $X=\Spec A$ is affine and almost of finite type over $\Lambda'$, we have the bounded above category $D_{\coh}^{-}(X)\subset \QCoh(X)$ with coherent cohomology; it identifies with the category of almost perfect complexes. It in turn defines the category $\Aperf^*(X)\coloneq D_{\coh}^{-}(X)\subset \QCoh(X)$ for an Artin stack $X$ almost of finite type over $\Lambda'$. One can also define the $\infty$-category $\Coh(X)\coloneq D_{\coh}^{b}(X)$ of bounded complexes with coherent cohomology. 

\begin{remark}
Let $X$ be a classical Artin stack. By \cite[Proposition 1.3]{HallRydh}, $\QCoh(X)$ is equivalent to the unbounded derived category $D_{\qc}(X_{\textnormal{lis-\'et}}, \mc O_X)$ of sheaves with quasi-coherent cohomology on the lisse-\'etale site. Moreover, if $X$ is QCA and perfect, then it is known to be further canonically equivalent to the unbounded derived category of the abelian category of quasi-coherent sheaves \cite[Theorem 1.2]{HallNeemanRydh}. If $X$ has the form $Y/G$, this Grothendieck abelian category identifies with the abelian category of $G$-equivariant quasi-coherent sheaves on $Y$. 
\end{remark}

Let $f\colon X\to Y$ be a map of schemes almost of finite type. We shall say that $f$ is of finite Tor amplitude (or eventually coconnective) if $f^*\colon \QCoh(Y)\to \QCoh(X)$ has bounded cohomological amplitude, equivalently, $f^*$ sends $\Coh(Y)$ to $\Coh(X)$ \cite[Lemma 3.6.3]{GaitsgoryIndCohSheaves}. A map $f\colon X\to Y$ of Artin stacks is locally of finite Tor amplitude if it is of finite Tor amplitude fppf locally on $X$. Any representable quasi-smooth map is locally of finite Tor amplitude. 

\begin{remark}
In the 6-functor formalism of solid quasi-coherent sheaves, maps locally of finite Tor amplitude are exactly suave maps \cite[Remark 2.7.13]{HansenMann}. Moreover, duality for suave objects recovers Grothendieck--Serre duality \cite[Corollary 2.8.4]{HansenMann}.   
\end{remark}

\subsection{Cotangent complexes}
The theory of cotangent complexes (relative to $\Lambda$) extends to (pre)stacks \cite[\S 3.2]{LurieThesis} using quasi-coherent sheaves. If $X$ is an Artin stack, the cotangent complex $\bb L_X$ exists and $\bb L_X[1]$ is connective \cite[Proposition 5.1.5]{LurieThesis}, \cite[Proposition 2.4.8]{HansenMann}. If $X$ is moreover smooth, $\bb L_X$ is perfect and has Tor amplitude in $[0,1]$.

\begin{example}[{\cite[Example 2.1.10]{Chenequivariant}, \cite[Proposition 2.4.13]{HansenMann}}]{\label{ex: cotangentcomplexes}}
If $X=BG$, $\QCoh(X)$ is the derived category of $\Rep G$, and $\bb L_X\simeq (\Lie G)^{\vee}[-1]$ with the coadjoint action. More generally, if $X=(\Spec A) /G$ with $\Spec A$ being smooth, then $\bb L_{X}$ is represented by a complex $\Omega^1_A\to A\otimes (\Lie G)^{\vee}$ as $G$-equivariant sheaves on $\Spec A$. 
\end{example}

A map $f\colon X\to Y$ of Artin stacks is quasi-smooth if $\bb L_{X/Y}$ is a perfect complex of Tor amplitude in $[-1, 1]$. If $X$, $Y$ are almost of finite type and $f$ is representable by algebraic spaces, then this agrees with the previous definition. 

The following general computation will be freely used: suppose that $f\colon X \to Y$ is the limit of some diagram of maps $X_i \rightarrow Y_i$ for $i$ in some index category $I$ in the $\infty$-category of stacks. If all $\L_{X_i/Y_i}$ exist and $\varinjlim\limits_i g_i^*\L_{X_i/Y_i}$ for $g_i\colon X\to X_i$ is bounded above, then $\L_{X/Y}$ exists and 
$\L_{X/Y} \simeq \varinjlim\limits_i g_i^*\L_{X_i/Y_i}$;
see \cite[Proposition 3.2.9]{LurieThesis}, \cite[Proposition 2.4.8]{HansenMann}.

\begin{lemma}[{\cite[Proposition 5.1.10]{HalpernLeistnerPreygel}}]\label{lem: MappingStackCotangentComplex}
Let $\pi\colon X\to \Spec \Lambda$ be a stack of finite Tor amplitude. Assume that $\pi^*$ admits a left adjoint $\pi_{\natural}$ after base change to any affine scheme $\Spec A$. For any stack $Y$ that admits a cotangent complex, the mapping stack $\Map(X, Y)$ admits a cotangent complex. Let $\eta\colon \Spec A \to \Map(X,Y)$ be a point corresponding to $f\colon X_A\to Y$. Then
\[
\eta^*\bb L_{\Map(X, Y)}\simeq \pi_\natural f^* \bb L_{Y}.
\]
In particular, if the left adjoint $\pi_\natural$ preserves perfect complexes and $Y=BG$, then it identifies with
\[
\pi_\natural f^* (\Lie G)^\vee[-1] \simeq \Gamma (X_A, f^* \Lie G)^{\vee}[-1]. 
\]
\end{lemma}

\begin{remark}\label{rem: LeftAdjointOfPullBack}
If the left adjoint $\pi_\natural$ preserves perfect complexes, then it sends a perfect complex $F$ to $(\pi_* F^\vee)^\vee$, where $(-)^\vee =\underline{\Hom}(-, \mc O_X)$ because $\Perf(X)$ is self-dual with respect to $(-)^{\vee}$ and $\pi^*$ commutes with $(-)^{\vee}$. See also the proof of \cite[Proposition 8.12]{ScholzeSixFunctor}.  
\end{remark}

\begin{example}\label{ex: GroupHomologyExample}
\begin{enumerate}
\item If $X$ is a proper algebraic space of finite Tor amplitude, then $\pi_\natural$ exists and preserves perfect complexes \cite[Theorem 8.10]{ScholzeSixFunctor}, \cite[Proposition 9.7]{ZhuTame}. 
\item Let $X=BH$ for a group scheme $H$ locally of finite type with $H^\circ$ affine. First, if $H$ is discrete, the left adjoint exists and is given by $-\otimes_{A[H]}A$. Assume $\pi_0(H)$ is an extension of $\Z$ by a finite group. It follows easily from the Hochschild--Serre type spectral sequence that $\pi_\natural$ exists and preserves perfect complexes. Therefore, 
\[
\bb L_{\Map (BH, BG)} \simeq \Gamma (H, \ev^* \Lie G)^{\vee}[-1].
\] 
\item If $X=B\Gm$, or more generally $BH$ with $H^{\circ}$ reductive and $\pi_0(H)$ finite, then the left adjoint is given by $\pi_*$. Therefore, $\bb L_{\Map(X, Y)}\simeq \Gamma (H, \ev^*\bb L_Y)$. 
\item A formula for the left adjoint for $X=\bb A^1/\Gm$ is given in \cite[Lemma 1.3.1]{HalpernLeistner}, and $\bb L_{\mathrm{Filt}(Y)}\simeq \Gamma (\bb A^1/\Gm, \ev^*\bb L_{Y}\otimes \Lambda [t, t^{-1}]/t\Lambda [t])$, where $t$ is the coordinate of $\bb A^1$. 
\end{enumerate}
\end{example}

\subsection{Ind-coherent sheaves}
We now move to ind-coherent sheaves. Some references are \cite{GaitsgoryIndCohSheaves, GRI, DrinfeldGaitsgory, ArinkinGaitsgory, Preygel, ZhuTame,HansenMann}. 

\begin{remark}
Starting with the 6-functor formalism of solid quasi-coherent sheaves, Hansen and Mann \cite[\S 2.9]{HansenMann} work out the 6-functor formalism of ind-coherent sheaves over a general regular noetherian ring $\Lambda$, essentially $\IndCoh^!$ or $\mathrm{QC}^!$ for Artin stacks below. However, it seems the ind-completion of coherent sheaves (or its subcategory) is more appropriate for the categorical local Langlands conjecture with general coefficients, and Zhu constructs such a theory in \cite[\S 9.3]{ZhuTame} and uses it for his fully faithful unipotent categorical local Langlands functor in the torsion case. Beyond the case where $\Lambda$ is regular, Hansen and Mann remark that the category of coherent sheaves has to be replaced \cite[Remark 2.8.5]{HansenMann}.
\end{remark}

\subsubsection{The case of schemes}\label{ss: indcoh for schemes}
Let $X$ be a scheme almost of finite type over $\Lambda'$. Recall $\Coh(X)\coloneq D_{\coh}^b (X)$ and set 
\[
\IndCoh(X)\coloneq \Ind(\Coh (X)),  
\]
with a canonical $t$-structure compatible with filtered colimits. It is $\QCoh (X)$-linear. 
There is a canonical $t$-exact $\QCoh (X)$-linear functor
\[
\Psi_X\colon \IndCoh(X)\to \QCoh (X),  
\]
which induces an equivalence on bounded below objects \cite[Proposition 1.2.4]{GaitsgoryIndCohSheaves}. It follows that the category of compact objects of $\IndCoh(X)$ is $\Coh(X)$ \cite[Corollary 1.2.6]{GaitsgoryIndCohSheaves}. 

The functor $\Psi_X$ exhibits $\QCoh(X)$ as the left completion of $\IndCoh(X)$ \cite[Proposition 1.3.4]{GaitsgoryIndCohSheaves}. If $X=\Spec A$, it is an equivalence if and only if $A$ is classical and regular \cite[Lemma 1.1.6, Proposition 1.6.4]{GaitsgoryIndCohSheaves}. 

For any quasi-compact almost of finite type schemes $X$ and $Y$, there is a natural equivalence
\[
\IndCoh(X)\otimes\IndCoh(Y) \simeq \IndCoh(X\times Y)
\]
that is compatible with $\QCoh(X)\otimes\QCoh(Y) \simeq \QCoh(X\times Y)$ under $\Psi$. 

For any map $f\colon X\to Y$ of schemes almost of finite type, $f_*\colon \Coh(X)\to \QCoh (Y)$ induces a functor
\[
f_*\colon \Coh(X)\to D_{\qc}^+(Y)\to \IndCoh (Y),  
\]
and its ind-extension $f_*^{\IndCoh}\colon \IndCoh (X)\to \IndCoh(Y)$, which is left $t$-exact and $\QCoh (Y)$-linear; by its construction, $f_*^{\IndCoh}$ and $f_*$ are compatible under $\Psi_X, \Psi_Y$ \cite[Proposition 3.1.1, Proposition 3.1.3]{GaitsgoryIndCohSheaves}. This construction naturally lifts to $\infty$-categories \cite[Proposition 3.2.4]{GaitsgoryIndCohSheaves}.

The next step is one of the most difficult parts of the theory \cite[Theorem 5.2.2]{GaitsgoryIndCohSheaves}, \cite{GRI}: construct homotopy-coherently $\QCoh(Y)$-linear $f^!\colon \IndCoh(Y)\to \IndCoh(X)$ and (a priori) additional data for any map $f\colon X\to Y$ of quasi-compact schemes almost of finite type satisfying the following:
\begin{itemize}
    \item If $f$ is proper (i.e., $\pi_0(f)$ is proper), $f^!$ identifies with the $\QCoh(Y)$-linear right adjoint of $f_*^{\IndCoh}$ canonically. (The right adjoint exists by the adjoint functor theorem, and it is moreover continuous since $f_*^{\IndCoh}$ preserves compact objects \cite[Corollary 3.3.6]{GaitsgoryIndCohSheaves}, and the right adjoint is $\QCoh(Y)$-linear \cite[Corollary 3.3.11]{GaitsgoryIndCohSheaves}.)  
    \item If $f$ is an open immersion, $f^!$ identifies with the left adjoint of $f_*^{\IndCoh}$ canonically. (The functor $f^*\colon \QCoh(Y)\to \QCoh(X)$ sends $\Coh(Y)$ to $\Coh(X)$, and it induces a $\QCoh(Y)$-linear functor $f^{*, \IndCoh}\colon \IndCoh(Y)\to \IndCoh(X)$ that is left adjoint to $f_*^{\IndCoh}$.)
    \item Additional data of base change isomorphisms such that if $f$ is proper, it arises from the adjunction and proper base change \cite[Proposition 3.4.2]{GaitsgoryIndCohSheaves}, and, if $f$ is an open immersion, it arises from the adjunction and base change \cite[Lemma 3.6.9]{GaitsgoryIndCohSheaves}.
    \item And compatibilities with the symmetric monoidal structure $\overset{!}{\otimes}$ defined using $\Delta^!$ for the diagonal $\Delta\colon X\to X\times X$. 
\end{itemize}

It is also an fppf sheaf:

\begin{proposition}[{\cite[Theorem 8.3.2]{GaitsgoryIndCohSheaves}}]
In the fppf topology, $\IndCoh$ satisfies $!$-descent.      
\end{proposition}

If $Y$ is eventually coconnective, the (relative) dualizing complex $\omega_{X/Y}$ is defined to be $f^! \mc O_Y$.  

\begin{remark}
After a recent axiomatization of 3-functor formalisms \cite[\S A.5]{MannPadic}, Scholze pointed out in \cite[Lecture VIII]{ScholzeSixFunctor} that
\[
\ProCoh(X)\coloneq \Pro(\Coh (X))\simeq \Pro(\Coh(X)^{\op}) \simeq \IndCoh(X)^{\op},  
\]
where $\Coh (X)\simeq \Coh(X)^{\op}$ is induced by the Grothendieck--Serre duality functor $\bb D_{\GS}$ recalled below, fits more naturally in the 3-functor formalism \cite[Definition 3.1.1]{HeyerMann}. More concretely, this means we instead consider the pro-versions of $\otimes$, $f^*$, $f_*$ for proper maps, and the left adjoint of $f^*$ for open immersions. One disadvantage is that $\ProCoh (X)$ is not presentable, but once we pass to $\IndCoh(X)$, we can extend it further to a 6-functor formalism \cite[Definition 3.2.1]{HeyerMann}. A recent result of \cite{DauserKuijper} implies that this construction with $\ProCoh(X)$ is essentially unique: one can show that proper maps (resp.~open immersions) are cohomologically proper (resp.~\'etale), cf.~\cite[Theorem 2.9.12]{HansenMann}.  
\end{remark}

\begin{remark}
The construction of \cite{GaitsgoryIndCohSheaves, GRI}\footnote{A part is also constructed in a different way in \cite[\S 5.4.3]{GRI}, but base change is missing, for instance.} relied on unproven statements about $(\infty, 2)$-categories \cite[\S 10.0.4.2]{GRI}, but proofs of these statements have been announced (\cite{LoubatonRuit} and references therein).   

The above type of construction can also be done rigorously in different ways:  
\begin{enumerate}
    \item \cite[Lecture VIII, Option 5]{ScholzeSixFunctor} gives such a construction using $\ProCoh(X)$ as above. (As written, only separated morphisms are allowed there, but it can be relaxed.) The construction is $(\infty, 1)$-categorical and relies on \cite{LiuZheng}. 
    \item Similarly, \cite[(9.13)]{ZhuTame} gives such a construction (extended to almost of finite type algebraic spaces), again relying on \cite{LiuZheng}. See also \cite[Remark 9.41(2)]{ZhuTame}. 
    \item Most conceptually, the universality of the $(\infty, 2)$-category of spans established in \cite{CnossenLenzLinskens} provides a 2-categorically enhanced construction. 
\end{enumerate}
However, the inputs for these constructions and the precise formulations of formalisms are somewhat different; see also \S \ref{ss: GSduality} below. 
\end{remark}

\subsubsection{More properties}{\label{sss: moreproperties}}
As in the case of open immersions, we can construct a left adjoint $f^{*, \IndCoh}\colon \IndCoh(Y)\to \IndCoh(X)$ of $f_*^{\IndCoh}$ for any map of schemes almost of finite type $f\colon X\to Y$ of finite Tor amplitude, and it is compatible with $f^*$ under $\Psi_X$, $\Psi_Y$. Both $f^{*, \IndCoh}$ and $f^!$ preserve compact objects in this setting \cite[Lemma 7.1.2]{GaitsgoryIndCohSheaves}. In particular, $\omega_{X/Y}\in \Coh (X)$ if $Y$ is eventually coconnective. The pullback $f^{*, \IndCoh}$ (resp.~$f^!$) satisfies base change with pushforward for any map \cite[Proposition 4.3.2.2]{GRI} (resp.~\cite[Proposition 7.1.6]{GaitsgoryIndCohSheaves}). In the smooth topology, $\IndCoh$ satisfies $*$-descent \cite[Proposition 8.3.8]{GaitsgoryIndCohSheaves}. Moreover, in the cartesian diagram
\[
\begin{tikzcd}
X' \arrow[r, "f'"] \arrow[d, "g'"] & Y' \arrow[d, "g"] \\
X \arrow[r, "f"] & Y
\end{tikzcd}
\]
with any $g$, we have the identification $(f')^{*, \IndCoh}\circ g^!\simeq (g')^{!}\circ f^{*, \IndCoh}$ \cite[Proposition 7.1.6]{GaitsgoryIndCohSheaves}; in fact, the case where $f$ is an open immersion was implicitly used in the construction of the formalism. 

Still assuming $f$ is of finite Tor amplitude, we have a continuous $\QCoh(Y)$-linear functor $f^{!,\QCoh}\colon \QCoh(Y)\to \QCoh(X)$ that is compatible with $f^!$ under $\Psi_X$, $\Psi_Y$ \cite[Proposition 7.2.2]{GaitsgoryIndCohSheaves}. If $f$ is proper, $f^{!, \QCoh}$ is the right adjoint of $f_*$.  

\begin{remark}    
Assume $X$ is eventually coconnective. Then $\Perf(X)\subset \Coh(X)$, and we obtain a right $t$-exact $\QCoh(X)$-linear fully faithful functor $\Xi_X\colon \QCoh(X)\to \IndCoh(X)$, and $\Xi_X$ is left adjoint to $\Psi_X$ \cite[\S 1.5]{GaitsgoryIndCohSheaves}. If $X$, $Y$ are both eventually coconnective, the two pullbacks $f^{*, \IndCoh}$, $f^{*}$ are compatible under $\Xi_X$, $\Xi_Y$ \cite[Lemma 3.5.11]{GaitsgoryIndCohSheaves}, and the two pushforwards $f_*^{ \IndCoh}$, $f_*$ are compatible under $\Xi_X$, $\Xi_Y$ \cite[Proposition 3.6.7]{GaitsgoryIndCohSheaves}. 
\end{remark}

\subsubsection{Grothendieck--Serre duality}\label{ss: GSduality}
In \cite[\S 9]{GaitsgoryIndCohSheaves}, \cite[\S 5.4.2]{GRI}, a duality is deduced from their formalism itself, and then it is shown to recover Grothendieck--Serre duality \cite[Proposition 9.5.7]{GaitsgoryIndCohSheaves}; i.e., the internal Hom functor $\underline{\Hom}(-, \omega_X)$ on $\Coh(X)$ for the dualizing complex $\omega_X$. As in \cite[Remark 4.4.5]{DrinfeldGaitsgory}, this characterizes the duality by the formula
\begin{equation}\label{eqn: GSdualityformula}
\Hom_{\IndCoh(X)} (\bb D_{\GS}\mc F, \mc G) \simeq f_*^{\IndCoh}(\mc F \overset{!}{\otimes} \mc G)    
\end{equation}
for $f\colon X\to \Spec \Lambda$, $\mc F, \mc G \in \Coh(X)$. In particular, the proper $*$-pushforward commutes with Grothendieck--Serre duality on $\Coh \subset \IndCoh$ by the projection formula \cite[(3.4)]{DrinfeldGaitsgory}, which follows from base change. Similarly, for maps of finite Tor amplitude, the duality exchanges the $*$-pullback and $!$-pullback on $\Coh$. 

\begin{remark}
In \cite[\S 9.3.4]{ZhuTame}, Grothendieck--Serre duality is deduced from the case of classical algebraic spaces, once the theory is constructed in \cite[Theorem 9.39]{ZhuTame}. Then, the author deduces the formula \eqref{eqn: GSdualityformula} above.      
\end{remark}

\begin{remark}
The construction of \cite[Lecture VIII, Option 5]{ScholzeSixFunctor} presumes Grothendieck--Serre duality, while, as mentioned earlier, Grothendieck--Serre duality can be realized first in the formalism of solid quasi-coherent sheaves, for instance. Then, the above formula \eqref{eqn: GSdualityformula} can be deduced exactly as in the proof of \cite[Theorem 9.42]{ZhuTame}. 
\end{remark}

We say that the map $f$ of finite Tor amplitude is Gorenstein if $\omega_{X/Y}\coloneq f^{!, \QCoh}\mc O_Y$ is a graded line bundle. This is equivalent to saying that $f^{!}(-)\simeq f^{*,\IndCoh} (-)\otimes \omega_{X/Y}$ \cite[Proposition 7.3.8, Proposition 7.4.5]{GaitsgoryIndCohSheaves}.\footnote{If $Y$ is eventually coconnective, then $\mc O_Y\in \Coh(Y)$ and $\omega_{X/Y}\simeq f^! \mc O_Y$.} Any quasi-smooth map of schemes almost of finite type is Gorenstein \cite[Corollary 2.2.7]{ArinkinGaitsgory}. Any map of finite Tor amplitude between quasi-smooth schemes almost of finite type is Gorenstein. 

\begin{proposition}\label{prop: QuasismoothDualizingScheme}
Let $f\colon X\to Y$ be a map of schemes almost of finite type. Let $\det(-)[-]$ denote the graded determinant of a perfect complex \cite{Heleodoro}, \cite[Appendix B]{HalpernLeistner}.  
\begin{enumerate}
    \item If $f\colon X\to Y$ is quasi-smooth, $\omega_{X/Y}\simeq \det (\bb L_{X/Y})[\rank(\bb L_{X/Y})]$. 
    \item If $X$, $Y$ are quasi-smooth, $\omega_{X/Y} \simeq \det (\bb L_{X/Y})[\rank(\bb L_{X/Y})]$.  
\end{enumerate}
\end{proposition}

\begin{proof}
(1) implies (2) by \cite[(46), (48)]{HalpernLeistner}. If $f$ is a quasi-smooth closed immersion, (1) is \cite[Corollary 9.7.3.2]{GRII}. A clearer, but still similar, argument is given in the proof of \cite[Lemma 2.4.5]{HalpernLeistner} using a version of the deformation to the normal cone/bundle \cite{KhanRydh}. The separated smooth case of (1), which is Serre duality for derived schemes, is \cite[Proposition 9.7.3.4]{GRII} (see also \cite[Lecture VIII, Option 1]{ScholzeSixFunctor}). These imply (1) for compositions of quasi-smooth closed immersions and separated smooth morphisms. Now (2) follows if the structure maps $X\to \Spec \Lambda$, $Y\to \Spec \Lambda$ are such compositions.

To treat the general derived case of (1), we assume the extension to stacks recalled below. Then, the argument of \cite[Lemma 2.4.5]{HalpernLeistner} with the deformation to the normal bundle stack \cite[Theorem A]{HekkingKhanRydh}, which is a 1-Artin stack in the current setting, would still work (see also the comments before \cite[Proposition 3.1.5]{HalpernLeistner}): we end up first with the dualizing complex of the zero section of the derived vector bundle stack $\bb L_{X/Y}^{\vee}[1]$ over $X$. Arguing as in \cite[Proposition 8.14]{ScholzeSixFunctor}, it identifies with the dual of the dualizing complex of the diagonal of the zero section. This diagonal identifies with the zero section of $\bb L_{X/Y}^{\vee}$ over $X$, and it further reduces to the dualizing complex of $\bb L_{X/Y}^{\vee}\to X$. Since $\bb L_{X/Y}^{\vee}$ is a scheme that has a quasi-smooth closed immersion into a smooth scheme over $X$, it follows from the established case of (1). 
\end{proof}

\begin{remark}
The functor $f^{!, \QCoh}$ and its properties would be more natural from the 6-functor formalism of solid quasi-coherent sheaves, since maps locally of finite Tor amplitude (resp.~Gorenstein) are suave (resp.~cohomologically smooth) in that setup. In the 6-functor formalism of $\IndCoh$, Gorenstein maps of eventually coconnective schemes are prim and the codualizing complex is given by the dual of $f^{!, \QCoh}\mc O_Y$ under $\Xi_X$. 
\end{remark}

\subsubsection{Extension to stacks}
All prestacks below are assumed to be almost of finite type.\footnote{We will treat the non-quasi-compact case only via ind-stacks.} There are a few approaches to ``ind-coherent'' sheaves in general. First, we consider a $\Lambda$-linear category $\Ind(\Coh(X))$ for an Artin stack $X$ almost of finite type with a natural functor $\Psi_X\colon \Ind(\Coh(X))\to \QCoh(X)$. 

As in the case of schemes, we have the following:
\begin{itemize}
    \item A natural $t$-structure such that $\Psi_X$ is $t$-exact. 
    \item For any map $f\colon X\to Y$ of Artin stacks almost of finite type, a left $t$-exact continuous functor $f_*^{\IndCoh}\colon \Ind(\Coh(X))\to \Ind(\Coh(Y))$ that is compatible with $f_*$ under $\Psi_X$, $\Psi_Y$. 
    \item If $f\colon X\to Y$ is of finite Tor amplitude, $f_*^{\IndCoh}$ admits a left adjoint $f^{*, \IndCoh}$, and it is compatible with $f^*$ under $\Psi_X$, $\Psi_Y$. 
    \item Base change holds for $*$-pushforward and $*$-pullback \cite[Lemma 9.26]{ZhuTame}. 
    \item If $f$ is schematic and proper, then $f_*^{\IndCoh}$ preserves compact objects and admits a continuous right adjoint $f^!\colon \Ind(\Coh(Y))\to \Ind(\Coh(X))$. Base change holds for $*$-pushforward and proper $!$-pullback, and the $*$-pullback commutes with the $!$-pullback in the cartesian diagram when these are defined \cite[Lemma 9.27]{ZhuTame}. 
    \item Compatibilities with $\Ind(\Coh(X))\otimes \Ind(\Coh(Y))\to \Ind(\Coh (X\times Y))$ \cite[Proposition 9.31]{ZhuTame}. 
\end{itemize}

One precise formalism is given in \cite[Theorem 9.32]{ZhuTame}. Note that $\Ind(\Coh(-))$ may not satisfy the schematic smooth descent with respect to the $*$-pullbacks, while it does so for QCA Artin stacks. 

\begin{remark}
It would not be difficult to construct the above 3-functor formalism from the 6-functor formalism of (solid) quasi-coherent sheaves using $f_*$ for all maps and $f^*$ for maps of finite Tor amplitude. 
\end{remark}

Next, we consider $\IndCoh^*(X)$ (sometimes called $\mathrm{QC}^*(X)$) which is extended from affine schemes almost of finite type along the $*$-pullback defined for maps of finite Tor amplitude. It is $\QCoh(X)$-linear. By \cite[Proposition 11.7.5]{GaitsgoryIndCohSheaves}, there is a canonical $t$-structure compatible with filtered colimits, and a $t$-exact functor
\[
\Psi_X\colon \IndCoh^*(X)\to \QCoh(X)
\]
inducing an identification on bounded below objects and exhibiting $\QCoh(X)$ as the left completion of $\IndCoh^* (X)$. In particular, the subcategory of bounded objects contains $\Coh(X)$, and we obtain a $t$-exact functor
\begin{equation}\label{eqn: IndCohComparison}
\Ind(\Coh(X))\to \IndCoh^*(X).     
\end{equation}
For any map of finite Tor amplitude $f\colon X\to Y$ of Artin stacks almost of finite type, we have a right $t$-exact $\QCoh(Y)$-linear functor $f^{*, \IndCoh}\colon \IndCoh^*(Y)\to \IndCoh^*(X)$ compatible with $f^{*,\IndCoh}\colon \Ind(\Coh(Y))\to \Ind(\Coh(X))$. 

The functor \eqref{eqn: IndCohComparison} may not be an equivalence. It is an equivalence if $X$ is QCA \cite[Theorem 3.3.5]{DrinfeldGaitsgory}, which includes the case of algebraic spaces almost of finite type. In particular, $f_*^{\IndCoh}\colon \IndCoh^* (X)\to \IndCoh^*(Y)$ is defined if $X$, $Y$ are QCA, and it agrees with the one in \cite[\S 3.6]{DrinfeldGaitsgory}. The functor $f_*^{\IndCoh}\colon \IndCoh^* (X)\to \IndCoh^*(Y)$ is also defined if $f$ is schematic (and quasi-compact, which is automatic here) \cite[\S 10.6]{GaitsgoryIndCohSheaves}. If $f$ is of finite Tor amplitude, the fact that $f_*^{\IndCoh}$ is the right adjoint of $f^{*,\IndCoh}$ is also stated in \cite[Corollary 3.7.13]{DrinfeldGaitsgory}. Let us also remark that if either $X$ or $Y$ is QCA, 
\[
\IndCoh^*(X)\otimes\IndCoh^* (Y)\simeq \IndCoh^* (X\times Y)
\]
by \cite[Corollary 4.2.3]{DrinfeldGaitsgory}.  If $ X\to Y$, $Y'\to Y$ are maps of perfect Artin stacks almost of finite type and $f\colon Y'\to Y$ is smooth and its base change to any affine scheme is QCA, then a natural map
\[
f^*\colon \QCoh(Y')\otimes_{\QCoh(Y)}\IndCoh^*(X) \to \IndCoh^*(Y'\times_Y X)
\]
is an equivalence.\footnote{Indeed, using that $\QCoh(Y)$, $\QCoh(Y')$ are rigid dualizable, we can reduce it to the case where $X$ is an affine scheme. And then it further reduces to the case $X\simeq Y$. As $\IndCoh^*(X)\simeq \Ind(\Coh(X))$ is dualizable, we can further reduce it to the case $Y'$ is an affine scheme. This case follows from \cite[Proposition 4.5.3]{GaitsgoryIndCohSheaves}.}

The category $\IndCoh^*(X)$ (sometimes called $\mathrm{QC}^*(X)$) has an alternative description, and it is extended from (affine) schemes almost of finite type along the $!$-pullback for all maps by smooth descent and the commutativity of the $!$- and $*$-pullbacks in cartesian diagrams. In fact, this makes sense for any prestack, and is the main approach of \cite{GaitsgoryIndCohSheaves, GRI}.\footnote{In these references, the functors are extended first from eventually coconnective affine schemes of finite type, but they are equivalent \cite[Corollary 10.2.2]{GaitsgoryIndCohSheaves}.} It is denoted by $\IndCoh^!(X)$ (or $\mathrm{QC}^!(X)$), and the comparison
\[
\IndCoh^* (X)\simeq \IndCoh^! (X)
\]
is stated in \cite[Corollary 11.5.3]{GaitsgoryIndCohSheaves} (see also \cite[Theorem 2.9.12(i), Lemma 2.10.8(i)]{HansenMann}). For any map $f\colon X\to Y$, a functor $f^!\colon \IndCoh^! (Y)\to \IndCoh^!(X)$ is then induced; it satisfies several properties and compatibilities:
\begin{itemize}
    \item It commutes with the $*$-pullbacks in cartesian diagrams, cf.~\cite[Corollary 11.5.3]{GaitsgoryIndCohSheaves}. 
    \item If $f$ is schematic, base change holds with $f_*^{\IndCoh}$. If $f$ is moreover proper, it is the right adjoint of $f_*^{\IndCoh}$ and the base change map is induced from the adjunction \cite[\S 5.3.4.4]{GRI}. 
    \item Base change also holds if $X$, $Y$ are QCA \cite[Corollary 3.7.14]{DrinfeldGaitsgory}. 
    \item If $f$ is Gorenstein smooth locally on $X$, then $f^!\simeq f^{*, \IndCoh}\otimes \omega_{X/Y}$, where $\omega_{X/Y}$ is a graded line bundle, regarded as an object of $\QCoh(X)$, called the relative dualizing complex (or sheaf) \cite[(3.12)]{DrinfeldGaitsgory}. 
    \item $\IndCoh^!(X)$ admits a symmetric monoidal structure $\overset{!}{\otimes}$, the $!$-pullbacks are symmetric monoidal, and the projection formula holds for the $*$-pushforward and $!$-pullback \cite[\S 5.4.1]{GRI}. 
\end{itemize}

A modern detailed exposition of these constructions (in slightly different form and generality) appears in \cite{HansenMann}. In fact, the construction itself is immediate from the case of schemes by applying \cite[Theorem 3.4.11]{HeyerMann}. 

\begin{remark}
Let $\omega_X$ denote the dualizing complex $f^! \Lambda$ for $f\colon X\to \Spec \Lambda$. It is the unit for the symmetric monoidal structure. Then, the following functor
\[
\Upsilon_X\colon \QCoh(X) \to \IndCoh^!(X);\quad \mc F \mapsto \mc F \otimes \omega_X
\]
is $\QCoh(X)$-linear and functorial with respect to the $*$-pullbacks and $!$-pullbacks \cite[\S 6.3.3]{GRI}. It in turn provides the projection formula involving $\QCoh$.     
\end{remark}

\begin{remark}[Grothendieck--Serre duality for stacks]
Grothendieck--Serre duality for $\Coh(X)$ for Artin stacks almost of finite type follows from the case of schemes by using the $*$-pullback and $!$-pullback for smooth maps \cite[(4.4)]{DrinfeldGaitsgory}. Phrased differently, the extension of a 6-functor formalism automatically encodes duality for stacks as well. Duality is again characterized by the formula \eqref{eqn: GSdualityformula} \cite[Proposition 4.4.4]{DrinfeldGaitsgory}, but with $f^{\IndCoh}_*$ replaced by $f_*\circ \Psi_X$ if $X$ is not QCA. 
\end{remark}

\begin{remark}
A version of the formalism for the $!$-pullback is constructed in \cite[Theorem 9.39]{ZhuTame} with $\Ind(\Coh(X))$ in place of $\mathrm{QC}^!(X)$ above. It would also work for our purpose. Note that the formalisms should agree on QCA stacks whenever relevant functors are defined. 
\end{remark}

\begin{proposition}\label{prop: QuasismoothDualizingStack}
Let $f\colon X\to Y$ be a map of Artin stacks almost of finite type. Let $\det(-)[-]$ denote the graded determinant of a perfect complex \cite{Heleodoro}, \cite[Appendix B]{HalpernLeistner} as in Proposition \ref{prop: QuasismoothDualizingScheme}.
\begin{enumerate}
    \item If $f\colon X\to Y$ is quasi-smooth, $\omega_{X/Y}\simeq \det (\bb L_{X/Y})[\rank(\bb L_{X/Y})]$. 
    \item If $X$, $Y$ are quasi-smooth, $\omega_{X/Y} \simeq \det (\bb L_{X/Y})[\rank(\bb L_{X/Y})]$.  
\end{enumerate}
\end{proposition}

\begin{remark}
In particular, $f^!$ is a left adjoint of $f_*^{\IndCoh}$ if $f$ is \'etale. The case of open immersions is \cite[\S 3.5.5]{GRII}. 
\end{remark}

\begin{proof}
As in the proof of Proposition \ref{prop: QuasismoothDualizingScheme}, (2) follows from (1). If $f$ is a quasi-smooth closed immersion (resp.~separated and smooth), (1) is \cite[Corollary 9.7.3.2]{GRII} (resp.~\cite[Proposition 9.7.3.4]{GRII}). The case where $Y$ is classical is also proven in \cite[Proposition 3.1.5]{HalpernLeistner}, and it holds for base changes of such maps. Therefore, (1) is true for compositions of such maps.\footnote{We only need these cases later.} 

The argument along the lines of the proof of Proposition \ref{prop: QuasismoothDualizingScheme} would work in the general case; a main issue is that we have to work with 2-Artin stacks. The details will appear in a forthcoming paper of A.~Jiang and A.~A.~Khan as announced in \cite[\S 0.4.3]{HekkingKhanRydh}.
\end{proof}

For QCA stacks $X$, we will simply write $\IndCoh(X)$ for
\[
\Ind(\Coh(X)) \simeq \IndCoh^*(X)\simeq \IndCoh^! (X). 
\]
If a map of QCA stacks $f\colon X\to Y$ is either of finite Tor amplitude or proper and schematic and $\mc F\in \Coh(Y)$, then we sometimes write $f^{!,\QCoh}\mc F$ for $\Psi_X (f^! \mc F)$. In the proper schematic case, it can also be described using the right adjoint of $f_*\colon \QCoh(X)^+ \to \QCoh(Y)^+$ on bounded below objects. 

\section{Formal derived algebraic geometry}
We will work with several ind-algebraic stacks (including $\Par_G$ itself) such as formal completions and certain loop spaces. We recall the foundations and discuss results on continuous coherent sheaves. Our references include \cite{chendhillon, DottoEmertonGee, ElmantoKubrakSosnilo, GRInd, HalpernLeistnerPreygel, ZhuTame}.  

\subsection{Ind-algebraic stacks}

Let us emphasize that we will consider finiteness conditions only relative to $\Lambda$. 

\begin{definition}
A prestack $X$ is an \emph{ind-algebraic stack} (resp.~ind-(QCA algebraic) stack) over $\Lambda$ if it is nilcomplete and is the filtered colimit of (resp.~QCA) Artin stacks $X_i, i\in I$ almost of finite type over $\Lambda$ such that the transition maps are closed immersions.  
\end{definition}

\begin{remark}
Following \cite[Definition 9.4]{ZhuTame}, we assumed $X_i$ are quasi-compact.  
\end{remark}

\begin{lemma}
If $X$ is an ind-algebraic stack, then $X$ is an fppf sheaf. 
\end{lemma}

\begin{proof}
As in the proof of \cite[Proposition 1.5.2]{GRInd}, this follows from the nilcompleteness and \cite[Lemma 1.3.3]{GRInd}. 
\end{proof}

\begin{example}
If an Artin stack is a disjoint union of open and closed substacks almost of finite type over $\Lambda$, then it is an ind-algebraic stack. 
\end{example}

The $*$-pushforward, $!$- and $*$-pullbacks in the theory $\IndCoh(-)$ extend to ind-algebraic stacks \cite{GRInd} (technically, for ind-schemes), \cite[Theorem 9.39]{ZhuTame}, or their arguments. More precisely, this is literally true with the approach in \cite{ZhuTame} as stated there, while one should consider ind-(QCA algebraic) stacks for $\IndCoh^!(-)$, and the two approaches are equivalent in that case. In our main setup, there is no ambiguity.  

\subsection{Formal completions}
The most typical example of ind-algebraic stacks is formal completions. Recall that, for a prestack $X$, the (thin) de Rham prestack $X_{\mathrm{dR}}$ of $X$ is defined by
\[
X_{\mathrm{dR}}(R)=X(\pi_0(R)_{\red}),  
\]
cf.~\cite[Definition 7.1.1]{ElmantoKubrakSosnilo}. 

\begin{definition}
Let $Z\to X$ be a map of prestacks. The (thin) formal completion $X^\wedge_Z$ of $X$ along $Z$ is defined to be $X\times_{X_\mathrm{dR}}Z_{\mathrm{dR}}$. 
\end{definition}

\begin{remark}
We usually consider Artin stacks almost of finite type over $\Lambda$ and their closed substacks. 
In this setting, the thin de Rham stack and the absolute de Rham stack are canonically equivalent \cite[Lemma 7.1.15]{ElmantoKubrakSosnilo}, and hence the corresponding notions of formal completions agree. 
\end{remark}

\begin{remark}
The formal completion depends only on $X$ and the underlying reduced classical stack of $Z$ (with the embedding into $X$). If $X$ is classical, the formal completion agrees with the classical one \cite[Corollary 2.1.5]{HalpernLeistnerPreygel}. 
\end{remark}

\begin{remark}{\label{rem: formalcompletionpullback}}
The following compatibility with pullback is clear from the definition: 
If $Z$ is the pullback of $Z'\subset X'$ along a map $X\to X'$, then 
\[
X^\wedge_Z \simeq X'^{\wedge}_{Z'}\times_{X'}X. 
\]
\end{remark}

\begin{remark}
Let $Z_2\to Z_1 \to X$ be maps of stacks. Then, $(X^\wedge_{Z_1})^\wedge_{Z_2}\simeq X^{\wedge}_{Z_2}$ holds. 
\end{remark}

\begin{proposition}
If $X$ is a scheme almost of finite type over $\Lambda'$ and $Z$ is a closed subscheme, then $X^\wedge_Z$ is the filtered colimit of closed subschemes $Z_\alpha$ almost of finite type over $\Lambda'$ in $X$ set-theoretically supported on $Z$ with a factorization $Z\to Z_{\alpha}\to X$. 
\end{proposition}

\begin{proof}
This follows from \cite[Proposition 6.5.5]{GRInd}, which works for $\Lambda'$ as well. 
\end{proof}

\begin{example}
Let $A$ be a classical $\Lambda$-algebra of finite type, and $I$ an ideal. The completion of $\Spec A$ along $\Spec A/I$ canonically identifies with $\Spf A^\wedge_I$, and it is an ind-algebraic stack in the above sense.   
\end{example}

\begin{corollary}
Let $X\to S$ be a map from an Artin stack almost of finite type over $\Lambda'$ to a scheme almost of finite type over $\Lambda'$. Given a closed subscheme $T\to S$, let $Z$ denote the pullback of $T$. The formal completion $X^\wedge_Z$ is the filtered colimit of $T_\alpha\times_{S} X$, where $T_\alpha \to S$ is a closed subscheme of $S$ almost of finite type over $\Lambda'$ set-theoretically supported on $T$ with a factorization $T\to T_\alpha\to S$.  
\end{corollary}

\begin{example}
Typically, we consider a map $f\colon X\to S$ from an Artin stack to a classical scheme, both almost of finite type over $\Lambda$. Pick a closed point $s\in S$. Then we have the base change
\[
X\times_S \Spec \mathcal{O}_{S,s}^\wedge \to \Spec \mathcal{O}_{S,s}^\wedge,  
\]
where the source is almost of finite type over $\mathcal{O}_{S,s}^\wedge$, but the target need not be almost of finite type over $\Lambda$. 
The completion of the source along $f^{-1}(s)$ identifies with $X^\wedge_{f^{-1}(s)}$, which is an ind-algebraic stack over $\Lambda$. 
\end{example}

\begin{lemma}
Let $X$ be an Artin stack almost of finite type over $\Lambda'$ with a closed substack $Z$. The $\infty$-category of factorizations $Z\to Y_i\to X$, where $Y_i$ is a closed substack almost of finite type over $\Lambda'$ set-theoretically supported on $Z$, is filtered, and $X^\wedge_Z\times_X U\simeq \colim_i (Y_i\times_X U)$ for any atlas $U\to X$ that is quasi-compact (as a map). 
\end{lemma}


\begin{proof}
We first claim that the $\infty$-category in the statement admits finite colimits, and hence it is filtered. For schemes, this is \cite[Proposition 3.1.3(i)]{GRInd}. These finite colimits commute with flat localizations on $X$ (for Zariski localizations, this is \cite[Proposition 3.1.3(ii)]{GRInd}). The case of stacks follows by taking an atlas $U\to X$ via the descent of closed substacks almost of finite type. 

The second part for schemes is \cite[Proposition 6.5.5]{GRInd}. The general case follows as in the proof of \cite[Proposition 2.2.21]{chendhillon} by passing to \v{C}ech nerves and taking scheme-theoretic images; we assume $U\to X$ is quasi-compact to be safe. 
\end{proof}

It is very likely that any formal completion is an ind-algebraic stack in the above sense under reasonable assumptions. This is easy to check when $X$ is a classical stack and $Z$ is a closed substack of $X$. More generally, at least the formal completion $X^\wedge_Z$ of a quotient stack has the expected form:

\begin{example}[{\cite[Proposition 7.1.16]{ElmantoKubrakSosnilo}}]
Let $G$ be an affine algebraic group over $\Lambda$, and let $Z\hookrightarrow X$ be a $G$-equivariant closed immersion of affine schemes almost of finite type over $\Lambda$. Then, 
\[
(X/G)^\wedge_{Z/G} \simeq X^\wedge_Z/G.  
\]
In particular, for any $G$-equivariant presentation of $X^\wedge_Z=\colim_i X_i$, $X^\wedge_Z/G\simeq \colim_i X_i/G$, and it is an ind-(QCA algebraic) stack. 
\end{example}

\begin{remark}
Note however that when $X=\Spec A$ and $Z=\Spec A/I$, there may not be an action of $G$ on $\Spec A^\wedge_I$ compatible with those on $\Spec A/I^n$, $n\geq 1$,~\cite[footnote 6]{chendhillon}. 
\end{remark}

Note that for any ind-algebraic stack $X$, 
\[
\Aperf^* (X) \subset \QCoh(X)
\]
makes sense and, for any presentation $X\simeq \colim_i X_i$, 
\[
\QCoh(X)\simeq \lim_i \QCoh(X_i), \quad \Aperf^*(X)\simeq \lim_i\Aperf^*(X_i), 
\]
and $\QCoh(X)$ is presentable. There are compatible $t$-structures on these stable $\infty$-categories (see also \cite[Proposition 2.1.11]{HalpernLeistnerPreygel}), and we consider 
\[
\wh{\Coh}{}^*(X)\subset \Aperf^*(X),
\]
called continuous pro-coherent sheaves in \cite{chendhillon}, consisting of bounded objects. 

\begin{example}
If $X$ is an affine scheme $\Spec A$ and $Z=\Spec \pi_0(A)/I$ with $I\subset \pi_0(A)$, then $\Aperf^*(X^\wedge_Z)$ is $t$-exactly equivalent to $\Aperf^*(\Spec A^\wedge_I)$ \cite[Lemma 2.1.10]{HalpernLeistnerPreygel}, where $A^\wedge_I$ is the derived $I$-adic completion. 

Any $A^\wedge_I$-module $M$ in $\Aperf^*(\Spec A^\wedge_I)$ is $I$-adically derived complete and $M\otimes_A \pi_0(A)/I$ is in $\Aperf^*(\Spec \pi_0(A)/I)$. If $A$ is classical, then $A^\wedge_I$ is the classical $I$-adic completion and the heart $\Aperf^*(X^\wedge_Z)^\heartsuit$ is the abelian category of finitely generated $A^\wedge_I$-modules.
\end{example}

\begin{lemma}\label{lem: completion_is_texact}
Suppose $X$ is an Artin stack almost of finite type over $\Lambda'$, and $i\colon Z\to X$ is a closed immersion. Then, $\Aperf^*(X^\wedge_Z)$ is left $t$-complete, right $t$-bounded, and the functor
\[
\wh{i}^*\colon \Aperf^* (X)\to \Aperf^*(X^\wedge_Z)
\]
is $t$-exact.  
\end{lemma}

For a similar $t$-exact claim for classical $X$ with affine diagonal, see also \cite[Appendix B]{DottoEmertonGee}, e.g.~Proposition B.3.21 and (B.3.28). 

\begin{proof}
See \cite[Proposition 2.1.11]{HalpernLeistnerPreygel} for the left $t$-completeness and the right $t$-boundedness. Let $U=\Spec A\to X$ be a smooth atlas, and suppose the reduced pullback of $Z$ is defined by $I\subset \pi_0(A)$. By the definition of $t$-structure, all the claims reduce to the case $U=X$ by flat descent. In the affine case, the $t$-exactness of $\wh{i}^*$ follows from \cite[Proposition 4.3.8, Corollary 4.3.9]{LurieDAGXII} adapted to animated $\Lambda'$-algebras. 
\end{proof}

\begin{remark}
The functor $\wh{i}^*\colon \QCoh(X)\to \QCoh(X^\wedge_Z)$ is left $t$-bounded and right $t$-exact, cf.~\cite[Lemma 2.2.1]{HalpernLeistnerPreygel}. 
\end{remark}

\subsection{Unipotent and formal loop spaces}{\label{ss: loopspaces}}
Formal completions also arise from certain loop spaces. Let $X$ be an Artin stack almost of finite type. Recall the unipotent loop space $\mc L^u(X)=\Map(B\Ga ,X)$. 

For an affine algebraic group $G$ over $\Lambda$, the unipotent cone $\mc U$ as a scheme  \cite[\S4.2]{Cotner} is defined either as the Zariski closure of unipotent elements of $G(\Lambda)$ regarded as a reduced scheme, or as the fiber of the composite $G\to \overline{G} \to \overline{G} \sslash \overline{G}$ over the image of $e\in G$, where $\overline{G}$ is the reductive quotient of $G$.  

\begin{proposition}[{\cite[Proposition 2.1.25]{Chenequivariant}}]
Let $G$ be an affine algebraic group. Then, there is an identification
\[
\mc L^u (BG) \simeq \wh{\mc U}/^\adj G, 
\]
where $\wh{\mc U}$ is the formal completion of $G$ along the unipotent cone $\mc U$. If $G$ acts on a scheme $X$ almost of finite type, then
\[
\mc L^u (X/G) \simeq \wh{\mc U}/^\adj G\times_{G/^\adj G}\mc L(X/G). 
\]
In particular, $\mc L^u (X/G)$ is an ind-algebraic stack. 
\end{proposition}

We can complete the loop space further. Suppose $X$ is geometric. The map $X\to \mc L(X)$ of constant loops is a closed immersion \cite[Lemma 2.1.13]{Chenequivariant}.

\begin{definition}
The \emph{formal loop space} $\wh{\mc L}(X)$ is the formal completion of $\mc L(X)$ along the constant loops $X\to \mc L(X)$. 
\end{definition}

\begin{example}[{\cite[Example 2.1.21]{Chenequivariant}}]
Let $G$ be an affine algebraic group over $\Lambda$. Then, there is an identification
\[
\wh{\mc L} (BG) \simeq G^\wedge_e/^\adj G.  
\]
If $G$ acts on a scheme $X$ almost of finite type, then
\[
\wh{\mc L} (X/G) \simeq G^\wedge_e/^\adj G\times_{G/^\adj G}\mc L(X/G). 
\]
In particular, $\wh{\mc L} (X/G)$ is an ind-algebraic stack.  
\end{example}

Now suppose $G$ is an affine algebraic group such that $G^\circ$ is reductive, and it acts on a scheme $X$ almost of finite type. Let $z\in G(\Lambda)$ be a semisimple element with $[z]\in G\sslash G$. Following \cite{Chenequivariant, BCHNconstructible}, we define certain variants of loop spaces: 

\begin{definition}
We define the $z$-formal loop space $\wh{\mc L}_z (X/G)$ to be the formal completion of $\mc L(X/G)\to \mc L(BG)$ along the fiber of $BZ_G (z)\hookrightarrow\mc L(BG)$, and the $z$-unipotent loop space $\mc L^u_z(X/G)$ to be the formal completion of 
\[
\mc L (X/G)\to \mc L(BG)\to G\sslash G
\]
along the fiber of $[z]\in G\sslash G$. 
\end{definition}

\begin{remark}
For $z=e$, $\wh{\mc L}_z (X/G)\simeq \wh{\mc L}(X/G)$ and $\mc{L}^u_z (X/G)\simeq \mc L^u (X/G)$. 
\end{remark}

The following seems implicit in \cite{BCHNconstructible}:

\begin{lemma}
Suppose $X$ is a smooth scheme. Let $A=\langle z\rangle \subset G$ denote the closed subgroup generated by $z$. Let $X^z$ denote the homotopy fixed points of the action map $z\colon X\to X$. Then, there is a natural equivalence $X^A\cong X^z_{\cl}$. 
\end{lemma}

\begin{proof}
The scheme $X^z_{\cl}$ is just the classical fixed points of $z$, and hence it agrees with the classical fixed points of $A$. 
\end{proof}

Set $X(z)\coloneq X^A$ and $G(z)\coloneq Z_G (z)$. If $X$ is smooth, $X(z)\cong X^z_{\cl}$ as above.

\begin{theorem}[{\cite[Theorem 3.1.12]{Chenequivariant}, \cite[Theorem 2.9]{BCHNconstructible}}]\label{thm: LoopSpaceLocalization}
The natural localization maps
\[
l^u_z \colon \mc L^u_z (X(z)/G(z)) \to \mc L^u_z (X/G), \quad \wh{l}_z\colon \wh{\mc L}_z (X(z)/G(z))\to \wh{\mc L}_z(X/G)
\]
are equivalences. 
\end{theorem}

Suppose $G$ is an affine algebraic group and $z\in G(\Lambda)$ is semisimple. We define $z$-formal loops $\wh{\mc L}_z(X/G)$ as in the reductive case. To define $z$-unipotent loops, we let $\overline{G}$ denote the reductive quotient of $G$ and let $\overline{z}$ denote the image of $z$. There are maps
\[
\mc L(X/G) \to \mc L(BG) \to \overline{G}\sslash \overline{G}, 
\]
and we define $\mc L^u_z (X/G)$ to be the formal completion of $\mc L(X/G)$ along the fiber of $[\overline{z}]\in \overline{G}\sslash \overline{G}$. For $z=e$, we have again $\mc{L}^u_z (X/G)\simeq \mc L^u (X/G)$. 

\begin{lemma}\label{lemma: NonreductiveLocalization}
The maps 
\[
l^u_z\colon \mc L^u_z (BG(z))\to \mc L^u_z (BG), \quad \wh{l}_z\colon \wh{\mc L}_z (BG(z))\to \wh{\mc L}_z (BG)
\]
are equivalences. 
\end{lemma}

\begin{proof}
It suffices to treat the $z$-unipotent loops. The sequence
\[
1 \to Z_{U_G}(z) \to G(z) \to Z_{\overline{G}}(\overline{z})\to 1,
\]
where $U_G$ is the unipotent radical, is exact by the semisimplicity of $z$: take $u \in G$ with $\overline{u}\in Z_{\overline{G}}(\overline{z})$. Then $uzu^{-1}$ is semisimple and maps to $\overline{z}$. By Lemma \ref{lemma: semisimpleconjugation} below, there exists $u'\in U_G$ such that $u' uzu^{-1}u'^{-1}=z$, and then $u'u$ is another lift of $\overline{u}$ and it belongs to $G(z)$. Therefore, $Z_G (\overline{z})=Z_G(z)U_G$. It follows that $\mc L^u_z (BG)$ is the formal completion of $G/^\adj G$ along $z\mc U_{G(z)}U_G /^\adj G(z)U_G$ and $\mc L^u_z (BG(z))$ is the formal completion of $G(z)/^\adj G(z)$ along $ z \mc U_{G(z)}/^\adj G(z)$. Again by Lemma \ref{lemma: semisimpleconjugation}, $ z \mc U_{G(z)}/^\adj G(z)\simeq z\mc U_{G(z)}U_G /^\adj G(z)U_G$. We deduce that $\mc L^u_z (BG(z))\simeq \mc L^u_z (BG)$, e.g.~by looking at the cotangent complexes.  
\end{proof}

\begin{lemma}\label{lemma: semisimpleconjugation}
Let $G$ be an affine algebraic group with reductive quotient $\overline{G}$. If $s$, $s'$ are semisimple elements of $G(\Lambda)$ and $\overline{s}=\overline{s'}$ in $\overline{G}$, then $s$ and $s'$ are conjugate under $U_G$.   
\end{lemma}

\begin{proof}
This follows from \cite[Theorem VIII.4.3]{hochschild1981basic}, \cite{mostow1956fully}; these results imply that any semisimple element is contained in a Levi subgroup of the (possibly disconnected) group $G$, and Levi subgroups are conjugate under $U_{G}$. Finally, note that the projection $G \rightarrow \overline{G}$ is an isomorphism when the source is restricted to a Levi subgroup of $G$.
\end{proof}

The following claim is clear:

\begin{lemma}\label{lem: ShiftingMap}
Suppose $z$ is central in $G$. The shifting map
\[
sh_z\colon \mc L^u (BG) \simeq \wh{\mc U}_G /^\adj G \to z\wh{\mc U}_G /^\adj G \simeq \mc L^u_z (BG)
\]
is an equivalence, and $\wh{\mc L}(BG)\simeq \wh{\mc L}_z(BG)$ similarly. 
\end{lemma}

\subsection{Derived coherent completeness and equivariant localization}
We recall results on coherent completeness, especially in the graded unipotent case discussed in \cite{AlperHallLim}, which is related to equivariant localization with respect to a torus action.\footnote{A related idea of equivariant localization appears, for instance, in \cite{Chenequivariant}. We originally had an argument inspired by \cite{Chenequivariant}, and only later became aware of the literature on coherent completeness.}

\begin{definition}
An Artin stack $X$ almost of finite type over $\Lambda'$ is complete along a closed subset $Z\subset \lvert X\rvert$ if the functor
\[
\wh{i}^* \colon \Aperf^*(X) \to \Aperf^*(X^\wedge_Z)
\]
is an equivalence for $\wh{i}\colon X^\wedge_Z\to X$.
\end{definition}

\begin{remark}
This notion (in the derived setting) is studied in \cite{HalpernLeistnerPreygel}. In the non-derived stacky setting, related notions (see \cite[Definition 3.2, Remark 3.4]{AlperHallLim}) are studied in \cite{GeraschenkoZureickBrown, AlperHallRydhLuna, AlperHallRydhLocal, AlperHallLim}. The statement closest to our interest is probably \cite[Theorem 1.3]{AlperHallLim} as in the discussion below. 
\end{remark}

\begin{lemma}\label{lem: completeness_reduction_to_classical}
If the underlying classical stack of $X$ is complete along $Z\subset \lvert X\rvert$, then $X$ is also complete along $Z$. 
\end{lemma}

\begin{proof}
This follows from \cite[Lemma 2.3.4]{HalpernLeistnerPreygel}. 
\end{proof}

\begin{lemma}\label{lem: reduction_to_formalfunctions}
Suppose $X$ is classical. Then $X$ is complete along $Z$ if and only if the following two conditions are both satisfied:
\begin{itemize}
    \item $(X, Z)$ satisfies \emph{formal functions}, i.e., for any coherent sheaf $\mc F\in \Coh(X)^\heartsuit$, the natural map
    \[
    H^i (X, \mc F) \to H^i (X^\wedge_Z, \wh{\mc F})
    \]
    is an isomorphism for all $i\geq 0$. 
    \item $(X, Z)$ is \emph{coherently complete}, i.e., the natural functor
    \[
    \Coh (X)^\heartsuit \to \varprojlim_n \Coh(X_n)^\heartsuit
    \]
    is an equivalence, where $X_n$ is the $n$-th infinitesimal thickening (with the 1st being the reduced closed substack).  
\end{itemize}
Moreover, $\Coh(X)\simeq \wh{\Coh}{}^*(X^\wedge_Z)$ in this case.
\end{lemma}

\begin{proof}
\cite[Proposition 3.5]{AlperHallLim} shows that the two conditions in the statement above are together equivalent to $(X,Z)$ being derived coherently complete in the sense of \cite[Definition 3.2]{AlperHallLim}. 

If $(X, Z)$ is derived coherently complete, then $X$ is complete along $Z$ as $\wh{i}^*$ is $t$-exact.
Conversely, if $X$ is complete along $Z$, then, by Lemma~\ref{lem: completion_is_texact}, it follows that $\Coh(X)\simeq \wh{\Coh}{}^*(X^\wedge_Z)$. Moreover, by arguing as in the proof of \cite[Proposition 3.5]{AlperHallLim}, it follows that $(X, Z)$ is derived coherently complete.
\end{proof}

Let us recall results on quotients by actions of reductive groups. 

\begin{theorem}[{\cite{HalpernLeistnerPreygel}}]\label{thm: completeness_along_fiber}
Suppose $A$ is a classical $\Lambda'$-algebra of finite type with $X=\Spec A$ and $G$ is reductive. If $G$ acts on $\Spec A$ and $A^G$ is a complete local ring with respect to the closed point $0\in \Spec A^G$, then $X/ G$ is complete along $f^{-1}(0)/G$, where $f\colon X\to \Spec A^G$ denotes the GIT quotient map. 
\end{theorem}

\begin{proof}
This follows from \cite[Theorem 4.2.1, Proposition 4.2.3]{HalpernLeistnerPreygel} (see also \cite[Theorem 2.3.16]{chendhillon}). 

Alternatively, this follows from the previous lemma together with \cite[Theorem 1.1]{GeraschenkoZureickBrown} and \cite[Corollary 4.2]{AlperHallRydhLocal} (or just \cite[(2.1)]{AlperHallRydhLuna}). 
\end{proof}

In fact, the following stronger result holds; see also \cite[Theorem B.4.17]{DottoEmertonGee}. Recall that in the above setup, there is a unique closed point $*\in X/G$. 

\begin{theorem}[{\cite[Theorem 4.5]{AlperHallRydhLocal}}]\label{thm: completeness_along_closedpoint}
In the above setting, $(X/ G, *)$ satisfies formal functions and is coherently complete. In particular, $X/G$ is complete along $*$. The same holds with $*$ replaced by any nonempty closed substack of $f^{-1}(0)/G$. 
\end{theorem}

\begin{proof}
This follows from \cite[Theorem 4.5]{AlperHallRydhLocal} (see also \cite[Theorem 1.3]{AlperHallRydhLuna}) and the previous lemma, together with \cite[Lemma 5.6]{AlperHallLim}. 
\end{proof}

We turn our attention to the non-reductive case, or more precisely, the graded unipotent case. 

Assume $P$ is a parabolic subgroup of a reductive group $G$ with a Levi subgroup $M$. Note that the following holds:
\begin{quotation}
There exists a cocharacter $\nu \colon \Gm\to M\subset P$ such that $Z_P(\nu)=M$ and the unipotent radical $U=U_P$ of $P$, or its Lie algebra $\Lie U$, has positive $\nu$-weights.  
\end{quotation}

\begin{theorem}\label{thm: completeness_along_fiber_parabolic}
Assume $\Lambda'$ is a complete noetherian local $\Lambda$-algebra. Suppose $P$ acts on an affine scheme $X=\Spec A$ almost of finite type over $\Lambda'$, and the GIT quotient (with respect to $M$) $f\colon X\to X\sslash M$ canonically identifies with $\Spec \Lambda'$. If the $\nu$-weights of $\pi_0 (A)$ are $\leq 0$, then $X/P$ has a unique closed point $*$, and $X/P$ is complete along any nonempty closed substack of $f^{-1}(0)/P$ for the closed point $0\in \Spec \Lambda'$. 
\end{theorem}

\begin{remark}
This is a slight variant of \cite[Theorem 1.3]{AlperHallLim}. 
\end{remark}

\begin{proof}
By Lemma \ref{lem: completeness_reduction_to_classical}, we may assume $X$ is classical. By \cite[Lemma 5.6]{AlperHallLim}, it suffices to show that $(X/P, *)$ satisfies formal functions and is coherently complete. This follows from \cite[Theorem 1.3]{AlperHallLim}: Conjecture 1.1 in \cite{AlperHallLim} holds by Theorem \ref{thm: completeness_along_closedpoint} (and Lemma \ref{lem: reduction_to_formalfunctions}). 
\end{proof}

\subsection{Ind-coherent sheaves and continuous ind-coherent sheaves}
For any ind-algebraic stack $X=\colim_i X_i$, the presentable stable $\infty$-category of ind-coherent sheaves (with a $t$-structure) can be defined as
\[
\IndCoh^!(X)=\lim_i \IndCoh^!(X_i);
\]
recall that the limit is taken with respect to $!$-pullbacks. If all $X_i$ are QCA, then
\[
\lim_i \IndCoh(X_i) \simeq \colim_i \IndCoh(X_i), 
\]
where the colimit is taken with respect to $*$-pushforwards (=$!$-pushforwards), and $\IndCoh(X)\coloneq\IndCoh^!(X)$ is compactly generated by $\Coh(X)=\colim_i \Coh(X_i)$ in this case. 

There is a dual version of $\Aperf^*(X)$. The full subcategory $\Aperf^! (X)\subset \IndCoh^! (X)$ consists of objects $\mc F$ whose $!$-restriction to each $X_i$ is left $t$-bounded, i.e., in 
\[
\IndCoh(X_i)^{+} \simeq \QCoh(X_i)^{+}
\]
and has coherent truncations, cf.~\cite[Proposition 2.2.14]{chendhillon}. The $t$-structure on $\IndCoh^!(X)$ induces a $t$-structure on $\Aperf^!(X)$, since the truncations preserve the left $t$-boundedness. 

Next, following \cite{chendhillon, Chen2023categorical, BCHNconstructible}, we introduce a dual version of the continuous pro-coherent sheaves $\wh{\Coh}{}^*$, i.e., continuous ind-coherent sheaves $\wh{\Coh}{}^!$.

\begin{definition}
Suppose $X$ is an ind-algebraic stack. We define $\wh{\Coh}{}^!(X)\subset \Aperf^!(X)$ to be the full subcategory consisting of $t$-bounded objects.  
\end{definition}

\begin{lemma}\label{lem:right_adjoint_from_hatCoh}
Suppose $X$ is a QCA Artin stack over $\Lambda$, and $Z$ is a closed substack with $\wh{i}\colon X^\wedge_Z\to X$. Then
the image of $\Coh(X)\subset \IndCoh(X)$ under $\wh i^!$ is contained in $\wh{\Coh}{}^!(X^\wedge_Z)$. 
The resulting functor
\[
\wh{i}^!\colon \IndCoh^!(X) \to \Ind\wh{\Coh}{}^! (X^\wedge_Z)
\]
admits a right adjoint $\wh{i}_?$. If $X^\wedge_Z\simeq \colim_i X_i$ with $X_i$ being QCA, and  $\mc F_j \simeq i_{ij*}\mc F_i$ for some $\mc F_i\in \Coh(X_i)$ and all $j\geq i$, $i_{ij}\colon X_i\to X_j$, then $\wh{i}_? \mc F\simeq i_{i*} \mc F_i$, where $i_i\colon X_i\to X$.  
\end{lemma}

\begin{proof}
The first claim is \cite[Lemma 2.3.7]{chendhillon} (whose proof works under the current assumptions). The second claim follows from the adjoint functor theorem (cf.~\cite[\S A.2.1]{BCHNconstructible}). For the last claim, we have equivalences
\[
\RHom (\mc G, i_{i*}\mc F_i)\simeq
\RHom (\wh{i}_* \wh{i}^!\mc G, i_{i*}\mc F_i)\simeq
\RHom (\wh{i}^!\mc G, \wh{i}^! i_{i*}\mc F_i)\simeq 
\RHom (\wh{i}^!\mc G, \mc F). 
\]
\end{proof}

One can relate $\wh{\Coh}{}^!$ and $\wh{\Coh}{}^*$ via Grothendieck--Serre duality following \cite{chendhillon}, but this is slightly non-canonical when working over $\Lambda'$. Suppose $\Lambda'$ is a complete noetherian local ring and fix a dualizing complex, which exists by \cite[\href{https://stacks.math.columbia.edu/tag/0BFR}{Tag 0BFR}]{StacksProject}.

Suppose $X$ is a geometric Artin stack almost of finite type over $\Lambda'$ and $Z$ is a closed substack. 
As in \cite[Corollary 2.3.5]{chendhillon}, the above choice determines compatible contravariant equivalences
\[
\Aperf^* (X^\wedge_Z) \simeq \Aperf^! (X^\wedge_Z), \quad \wh{\Coh}{}^*(X^\wedge_Z)\simeq \wh{\Coh}{}^!(X^\wedge_Z),
\]
exchanging $\wh{i}^*$ and $\wh{i}^!$. 

\begin{corollary}{\label{cor: cohhatcompletion}}
In the above setup, if $X$ is complete along $Z$, then
\[
\Aperf^! (X)\xrightarrow{\wh{i}^!} \Aperf^! (X^\wedge_Z), \quad
\Coh (X) \xrightarrow{\wh{i}^!} \wh{\Coh}{}^! (X^\wedge_Z)
\]
are equivalences. 
\end{corollary}

\begin{remark}
Note that the above statements are independent of the choice of the dualizing complex. 
\end{remark}

\section{Derived vector bundles and linear Koszul duality}

We recall derived vector bundles and linear Koszul duality in our relevant setting, following Chen \cite{Chen2023categorical}. 

\subsection{Derived vector bundles}
Let $X$ be a prestack over $\Lambda$. The category of derived vector bundles on $X$ is the category of perfect complexes on $X$. If $\mathcal{E}$ is a coconnective perfect complex on $X$, then we can consider the relative spectrum
\[
\mathbb{E}_X\coloneq \Spec_X \mathrm{Sym}_{\mathcal{O}_X} \mathcal{E}^\vee,
\]
a prestack that is affine over $X$, and we call $\mathbb{E}_X$ a derived vector bundle as well. It is equipped with the zero section $z\colon X\to \mathbb{E}_X$, and we define the completed vector bundle $\widehat{\mathbb{E}}_X$ as the (thin) formal completion of $\mathbb{E}_X$ along $z\colon X\to \mathbb{E}_X$. If $n\geq 1$ is an integer, we define $\mathbb{E}^{\leq n}_X$ to be the $(n+1)$-st infinitesimal neighborhood, i.e., it is the relative spectrum of the $(n+1)$-st truncated algebra of symmetric powers of $\mathcal{E}^\vee$. 

\begin{remark}
The above definition of completed vector bundles agrees with \cite[Definition 7.1.24]{ElmantoKubrakSosnilo} if $X$ is an Artin stack almost of finite type over $\Lambda$, and this assumption holds in our main setup.  
\end{remark}

\begin{lemma}\label{lem:colimit=convergent colimit}
Assume $X$ is a classical Artin stack almost of finite type over $\Lambda$ and $\mathcal{E}$ is concentrated in degrees $[0,1]$. Then, 
\[
\widehat{\mathbb{E}}_X \simeq \colim_n \, \mathbb{E}_X^{\leq n}. 
\]
\end{lemma}

\begin{proof}
For any $\mathcal{E}$, \cite[Proposition A.1.1]{ElmantoKubrakSosnilo} shows that 
\[
\widehat{\mathbb{E}}_X \simeq \colim_n^{\mathrm{conv}} \, \mathbb{E}_X^{\leq n},  
\]
where $\colim_n^{\mathrm{conv}}$ denotes the colimit as convergent prestacks. Under the assumption that $\mathcal{E}$ is concentrated in degrees $[0, 1]$, there exists an integer $m$ such that the $n$-th truncated algebra of symmetric powers of $\mathcal{E}^{\vee}$ is concentrated in degrees $[-m,0]$ for all $n$, and it follows that the colimit in prestacks agrees with the colimit in convergent prestacks. 
\end{proof}

\subsection{Completed odd tangent stacks and formal loop spaces}

Let $X$ be a quasi-smooth Artin stack over $\Lambda$ so that the cotangent complex $\L_X$ is perfect. We write $\T X[-1]$ for the derived vector bundle corresponding to $(\L_X [1])^{\vee}=\L_X^{\vee}[-1]$, and $\widehat{\T}X[-1]$ for its completion along the zero section. 

The following theorem is \cite[Theorem A.3.4]{ElmantoKubrakSosnilo}, which clarifies the earlier work of Ben-Zvi--Nadler \cite[Theorem 6.9]{BenZviNadler}. We will make use of it only in the case where $X$ is smooth, in which the situation is clearer \cite[Remark 1.3.2]{ElmantoKubrakSosnilo} (and by Lemma \ref{lem:colimit=convergent colimit}). 

\begin{theorem}\label{thm:exp}
If $X$ is geometric, there is a functorial equivalence of prestacks 
\[
\exp_X\colon \widehat{\T} X[-1] \xrightarrow{\simeq} \widehat{\mathcal{L}}(X). 
\]
\end{theorem}

\subsection{Cotangent stacks}
Suppose $X$ is an Artin stack smooth over $\Lambda$ so that $\T^*X$ is a derived vector bundle concentrated in $[0,1]$, i.e., $\mathcal{E}=\L_X$ in the notation above. 

We further assume that $X=Y/G$ for a connected smooth affine scheme $Y$ of finite type with an action of an affine algebraic group $G$. Recall the moment map
\[
\mu\colon \T^* Y \to (\Lie G)^{\vee}. 
\]

The following results are well-known. 

\begin{lemma}
The cotangent stack $\T^* (Y/G)$ identifies with the quotient of the (derived) zero fiber $\mu^{-1}(0)$ by $G$. 
\end{lemma}

\begin{proof}
This is shown in \cite[\S4.3]{BeniniSafronovSchenkel}, \cite[\S3.1]{BeniniPridhamSchenkel}.    
\end{proof}

\begin{lemma}
The zero fiber $\mu^{-1}(0)$ is topologically the union of conormal bundles $\T^*_{\mc O}Y$ for all $G$-orbits $\mc O\subset Y$. If there are only finitely many $G$-orbits in $Y$, then the underlying classical scheme of $\mu^{-1}(0)$ is equidimensional of dimension $\dim Y$, and its irreducible components are given by the closures of conormal bundles of $G^\circ$-orbits. 
\end{lemma}

\begin{proof}
See the beginning of Chapter 19 of \cite{ABV}, especially Lemma 19.2.
\end{proof}

\subsection{Linear Koszul duality}{\label{ss: koszuloddtangent}}
We will use the linear Koszul duality of Mirkovi\'c--Riche \cite{MRI, MRII, MRIwahoriMatsumoto}, formulated as in \cite{Chen2023categorical}. We will specialize to the completed odd tangent stack and the cotangent stack.  

Let $X$ be a smooth QCA stack with the trivial action of $\Gm$. We will regard $\L_X$ as a $\Gm$-equivariant sheaf by giving it weight $-2$\footnote{Beware that \cite{Chen2023categorical} uses the weight $-1$ instead.}, and its dual $\L_X^\vee$ has the weight $2$. The prestacks $\wh{\T}X[-1]$, $\T^\ast X$ admit the resulting $\Gm$-actions, and we will speak about the graded sheaves on these spaces. 

We define $\QCoh(\T^\ast X[2])$ (resp.~$\QCoh(\T^\ast X[2]/\Gm)$) to be the category of (resp.~graded) $\Sym_{\mathcal{O}_X} (\L_X^\vee[-2])$-modules in $\QCoh(X)$; this agrees with the definition in \cite{Chen2023categorical} by \cite[Proposition 2.1.11]{Chen2023categorical}. 

Recall also that $\IndCoh(\wh{\T}(-)[-1])$ (resp.~$\QCoh(\T^\ast (-)[2])$) has $!$-pullback for all maps $X\to Y$ of smooth QCA stacks (resp.~maps of finite Tor amplitude), and $*$-pushforward for all maps $X\to Y$, defined via correspondences as in \cite[Definition 2.3.1]{Chen2023categorical}. For $\QCoh (\bb T^*(-)[2]/ \Gm)$, the $!$-pullback is the composite of the $!$-pullback for $(\bb T^* Y[2]\times_Y X) / \Gm\to \bb T^* Y[2]/ \Gm$ (defined by the assumption) and the $*$-pushforward along $(\bb T^* Y[2]\times_Y X)/ \Gm \to \bb T^* X[2]/ \Gm$. Similarly, for $\QCoh (\bb T^*(-)[2]/ \Gm)$, the $*$-pushforward is the composite of the $*$-pullback for $(\bb T^* Y[2]\times_Y X) / \Gm\to \bb T^* X[2]/ \Gm$ and the $*$-pushforward along $(\bb T^* Y[2]\times_Y X) / \Gm\to \bb T^* Y[2]/ \Gm$. 

The first form of the linear Koszul duality is the following:

\begin{theorem}[{\cite[Theorems 2.1.19, 2.3.6]{Chen2023categorical}}]{\label{thm: gradedkoszulduality}}
Recall that $z\colon X\to \wh{\T}X[-1]$ is the zero section. The functor $\kappa\coloneq z^!$ induces an equivalence of categories
\begin{equation*}
        \IndCoh(\wh{\T}X[-1] / \Gm) \simeq  \QCoh(\T^{\ast}X [2]/ \Gm). 
\end{equation*}
Moreover, this duality is compatible with $!$-pullback for all maps and $*$-pushforward for proper schematic maps. 
\end{theorem}

\begin{proof}
The last part follows from \cite[Theorem 2.1.33]{Chen2023categorical} and the definitions of $!$-pullback and $*$-pushforward.
\end{proof}

\begin{remark}
As in the proof in \cite{Chen2023categorical}, the equivalence is a simple application of the Barr--Beck--Lurie theorem. 
\end{remark}

It is more convenient to replace $\T^\ast X [2]$ by $\T^\ast X$, and this is done via shearing \cite{Chen2023categorical} (see also \cite[\S 6]{BZSV}). Our normalization is that the shearing $\shear=[\![-1]\!]$, and it sends a sheaf in degree 0 to the line in weight-degree $(k,k)$. Its inverse is $\unshear=[\![1]\!]$.\footnote{In \cite[Definition 2.1.7]{Chen2023categorical}, $\shear=[\![-2]\!]$, $\unshear=[\![ 2]\!]$, and he uses $(k, 2k)$ in place of $(k,k)$.}
The $!$-pullback and proper schematic $*$-pushforward are also defined on $\IndCoh(\T^\ast (-))$. 

The second form of the linear Koszul duality is the following:

\begin{theorem}[{\cite[Theorems 2.1.19, 2.3.6]{Chen2023categorical}}]{\label{thm: koszulduality}}
The functor $ \unshear \circ \kappa$ induces an equivalence of categories
\begin{equation*}
\IndCoh(\wh{\T}X[-1] / \Gm) \simeq  \QCoh(\T^{\ast}X / \Gm). 
\end{equation*}
Moreover, this duality is compatible with $!$-pullback for all maps and $*$-pushforward for proper schematic maps. 
\end{theorem}

There is also a renormalized version, and in fact it is the most convenient form for us. We simply call it (linear) Koszul duality. 

\begin{theorem}[{\cite[Theorems 2.2.10, 2.3.6]{Chen2023categorical}}]{\label{thm: renormalizedcontinuousKoszul}}
The functor $\unshear \circ \kappa$ induces an equivalence of categories
\begin{equation*}
\wh{\Coh}{}^!(\wh{\T}X[-1] / \Gm) \simeq  \Coh(\T^{\ast}X / \Gm),
\end{equation*}
and it is compatible with $!$-pullback for smooth maps and $*$-pushforward for proper schematic maps. The same holds for the $\Ind$-completions. 
\end{theorem}

\begin{example}[{\cite[Theorem 2.3.6.5]{Chen2023categorical}}]
Let $p\colon \T^\ast X\to X$ denote the projection. The functor $\unshear\circ \kappa$ exchanges the functors
\begin{align*}
z_* &\colon \Coh(X/\Gm)=\Perf(X/\Gm) \to \Coh(\wh{\T}X[-1]/\Gm)\subset \wh{\Coh}{}^!(\wh{\T}X[-1]/\Gm), \\
p^* &\colon \Coh(X/\Gm)=\Perf(X/\Gm) \to \Perf(\T^\ast X/\Gm) \subset \Coh(\T^\ast X/\Gm). 
\end{align*}
\end{example}

\begin{example}
The image of $\omega_{\wh{\T}X[-1]}$ is $z_*\omega_X\in \Coh(\T^\ast X)$ supported on the zero section of $\T^\ast X$ (ignoring the grading). 
\end{example}

\begin{example}{\label{ex: derivedVB}}
Suppose $X=V/H$ for a (possibly disconnected) reductive group $H$ and a representation $V$ of $H$ such that the action of $H$ factors through $H'\coloneq H/Z$ for a central torus $Z$. Let us consider the following form of ``truncation'' of the odd cotangent complex 
\[
[V^\vee \otimes \mc O_X\to (\Lie H')^{\vee}\otimes \mc O_X]\to
[V^\vee \otimes \mc O_X\to (\Lie H)^{\vee}\otimes \mc O_X]\simeq\L_X[1]
\]
and dually, up to a shift, 
\[
\L_X^\vee \simeq [\Lie H\otimes \mc O_X\to V\otimes \mc O_X] \to
[\Lie H'\otimes \mc O_X\to V\otimes \mc O_X].
\]
Note that the truncation is also a derived vector bundle. Geometrically, let us write them as
\[
\delta\colon \wh{\T}X[-1] \to \wh{\V}[-1], \quad 
\check{\delta}\colon \V^\vee\to \T^* X. 
\]
Let us consider a vector bundle $\mc V$ on $X$ and regard it as a sheaf on the zero section of $\wh{\V}[-1]$. Its Koszul dual is the $*$-pullback $\mc V|_{\V^\vee}$ of $\mc V$ to $\V^\vee$. The pushforward $\check{\delta}_* (\mc V|_{\V^\vee})$ is Koszul-dual to $\delta^! \mc V$ (with its natural grading). The relative dualizing complex of $\T X[-1]\to \V[-1]$ is trivial in degree $-\dim Z$. Therefore, $\delta^! \mc V$ agrees with the $!$-pullback of $\mc V[\dim Z]$ on $\Lie Z \times X\subset \T X[-1]$ to $\wh{\T}X[-1]$. In particular, $\delta^! \mc V[-\dim Z]$ is in degree $0$ with respect to the $t$-structure.   
\end{example}

\section{Filtered \texorpdfstring{$\mc{D}$}{D}-modules}
For a formalism of filtered $\mc D$-modules for smooth Artin stacks, we follow \cite[\S 3.3]{Chen2023categorical}, but see also \cite{MixedHodgeProject} with more details in the non-stacky case. Let us emphasize that we will work throughout with \emph{right} $\mc D$-modules, with increasing filtrations.\footnote{The side-changing functor in \cite[\S 3.3.5]{Chen2023categorical} seems to be normalized incorrectly, and there should be no shift coming from $\widetilde{\omega}_X$.}

Let $X$ be a smooth scheme. Let $F\mc D (X)$ denote the $\QCoh(\bb A^1/\Gm)$-linear (stable) $\infty$-category of filtered $\mc D$-modules defined in \cite[Definition 3.3.1]{Chen2023categorical}. It is functorial with respect to $!$-pullbacks \cite[Proposition 3.3.10]{Chen2023categorical} and forms a sheaf for the smooth topology with values in $\QCoh(\bb A^1/\Gm)$-linear (stable) $\infty$-categories \cite[Proposition 3.3.18]{Chen2023categorical}. The definition extends to smooth Artin stacks via descent \cite[Definition 3.3.20]{Chen2023categorical}. It comes with the following functors
\[
\mc D(X) \xleftarrow{\un} F\mc{D} (X) \xrightarrow{\gr} \QCoh(\bb T^* X/\Gm)=\QCoh(\T^*X)^\Gm, 
\]
where $\mc D(X)$ denotes the category of right $\mc D$-modules, the left functor forgets filtrations, and the right functor takes the associated graded sheaf.
\footnote{Let us recall that our normalization of $\Gm$-action is different from \cite{Chen2023categorical}.}

There is also a functor $\res\colon F\mc D (X) \to \QCoh(X)^{\textnormal{Fil}}$ taking the underlying filtered quasi-coherent sheaf on $X$, called restriction in \cite{Chen2023categorical}. The restriction admits a left adjoint $\ind\colon \QCoh(X)^{\textnormal{Fil}} \to F\mc{D}(X)$, called induction. The functors $\un$ and $\gr$ on $\QCoh(X)^{\textnormal{Fil}}$ are compatible with induction and restriction in the sense that we have a commutative diagram 
\begin{equation*}
    \begin{tikzcd}
       \QCoh(X) \arrow[d, "\ind", shift left]   & \QCoh(X)^{\textnormal{Fil}} \arrow[d, "\ind", shift left] \arrow[r, "\gr"] \arrow[l, "\un", swap] & \QCoh(X)^{\gr} \arrow[d, "\ind", shift left]  \\
        \mc{D}(X) \arrow[u, "\res"] & F\mc{D}(X) \arrow[u, "\res"] \arrow[r, "\gr", swap] \arrow[l, "\un"] & \QCoh(\T^{\ast}X / \Gm) \arrow[u, "\res"].
    \end{tikzcd}
\end{equation*}

A filtered $\mc D$-module $\mc M$ is called coherent if $\gr (\mc M)$ is coherent and the (increasing) filtration is separated.\footnote{The first sentence of \cite[Remark 3.3.4]{Chen2023categorical} does not seem completely correct.} For quasi-compact Artin stacks, the subcategory $F\mc D_c(X)$ of coherent filtered $\mc D$-modules can be strictly larger than the subcategory of compact objects. On the other hand, the objects induced from $\Coh(X)^{\textnormal{Fil}}$ are compact and generate $F\mc D (X)$ as in the corresponding results for $\mc D$-modules in \cite{DrinfeldGaitsgory}. Therefore, we sometimes consider
\[
\mc D_c(X) \xleftarrow{\un} F\mc{D}_c (X) \xrightarrow{\gr} \Coh(\bb T^* X/\Gm) 
\]
or their ind-completions instead.

\begin{remark}{\label{rem: Dmodulefunctoriality}}
The above diagram is functorial with respect to $!$-pullbacks for smooth maps and $*$-pushforward along proper schematic maps \cite[Theorem 3.4.9]{Chen2023categorical} (the former is claimed only for schematic maps, but being schematic is not used). 
\end{remark}

\subsection{The Riemann--Hilbert correspondence}
Assume $\Lambda=\C$. For smooth schemes $X$, the Riemann--Hilbert correspondence
\[
\RH\colon\mc{D}_{\mathrm{rh}}^b(X) \xrightarrow{\cong}    D_c^b(X, \C)
\] 
gives an equivalence between regular holonomic \emph{right} $\mc D$-modules $\mc{D}^b_{\mathrm{rh}}(X)$ and the bounded derived category $D_c^b (X, \C)$ of (Zariski) constructible sheaves, identifying regular holonomic $\mc D$-modules in degree 0 and  perverse sheaves. More specifically, we use the (shifted) de Rham functor for \emph{left} $\mc D$-modules and precompose it with the side-changing functor. In particular, writing $\omega_X^\Coh$ for the coherent dualizing complex, we send the right $\mc D$-module $\Omega_X^{\dim X}=\omega_X^{\Coh} [-\dim X]$ in degree 0 to $\bb C_X[\dim X]$. This equivalence commutes with the 6-functors, cf.~\cite{Paulin}. In particular, we can pass to smooth Artin stacks\footnote{See \cite{Vooys} for comparison of several approaches to equivariant constructible sheaves.} to obtain a fully faithful functor
\[
\RH^{-1}\colon \Ind D_c^b(X, \C) \to \Ind \mc{D}_c(X)
\]
whose image consists of regular holonomic ind-$\mc D$-modules, 
sending perverse sheaves on $X$ (which makes sense easily) to the heart of $\mc{D}_c(X)$, the bounded derived category of coherent $\mc D$-modules. Again, $\bb C_X[\dim X]$ is sent to $\omega_X^{\Coh} [-\dim X]$. 

\begin{lemma}
Let $Y$ be a smooth scheme with an action of an affine algebraic group $G$.
If there are only finitely many $G$-orbits, the functor $\RH^{-1}$ is an equivalence for $X=Y/G$. 
\end{lemma}

\begin{proof}
Regular holonomicity is automatic \cite[Theorem 11.6.1]{HottaTakeuchiTanisaki}. 
\end{proof}

\part{Moduli stack of \texorpdfstring{$L$}{L}-parameters}{\label{p: moduliofLparams}}

\section{\texorpdfstring{$L$}{L}-parameters and Vogan stacks}{\label{s: Lparamsandvoganstacks}}

For simplicity, let $G$ be a split reductive group in this paragraph. A canonical definition of $\Par_G$, defined over a certain stack $\mathrm{MG}_{\overline{\F}_q}$ over $\Z[1/p]$, is given in \cite[\S4.1]{ScholzeLLCmotives} in the context of the motivic geometrization of the local Langlands conjecture (we refer the reader there for the details). As discussed after Definition 4.3 of loc.~cit., a choice of isomorphism $\ov{\F}^{\times}_q \cong \Q / \Z[1/p]$ determines a map $\Spec(\Z[1/p]) \rightarrow \mathrm{MG}_{\overline{\F}_q}$. Upon a further choice of a uniformizer $\pi$ of $F$ and of roots $\pi^{1/n}$ for $n$ prime to $p$, the pullback of $\Par_G$ to $\Spec(\Z[1/p])$ naturally parametrizes $\wh{G}$-torsors on a classifying stack of a certain group scheme $\WD_F$ over $\Z[1/p]$ (but the choice of $\pi^{1/n}$ does not essentially change this parametrization).  In the proof of his Proposition 4.4, Scholze explains that $\WD_F$ when base-changed to $\Q$ is isomorphic to the classical Weil--Deligne group scheme $W_F \ltimes \Ga$. Since we work over the algebraically closed characteristic $0$ field $\Lambda$ throughout, this will be the context of the present paper, and motivates the following definitions (rather than the ones in \cite{DHKM, fargues--scholze, HansenMann}); the definitions below themselves do not involve any choices. Let us also note that Takaya and Lin \cite{takayalin} use a similar approach to the spectral side of the categorification of the relative local Langlands conjecture in the spirit of \cite{BZSV}. 

\subsection{The stack of Weil parameters}
We start by analyzing the stack of maps from the Weil group $W_F$, and will consider the $\Ga$-part later. First we recall a general definition: let $S$ be a classical stack and $H\to G$ a surjection of classical group stacks over $S$ that are schematic. We consider the stack
\[
\Map (BG, BH)\times_{\Map (BG, BG)}\{\id_{BG}\},
\]
which is the ``Weil restriction'' of $BH$ along the (non-representable) map $BG\to S$. 
It is the quotient of the stack (schematic over $S$) of sections
\[
\Sect (G, H)\coloneq \Hom (G, H)\times_{\Hom (G, G)}\{\id_G\}
\]
by the kernel of $H\to G$ (see Example \ref{ex: ExamplesDerivedStacks}, for instance). 

We regard the Weil group $W_F$ as a pro-(discrete group scheme); for any open subgroup $I\subset I_F$ that is normal in $W_F$, $W\coloneq W_F/I$ is finitely generated, in fact, an extension of $\Z$ by a finite group, and $W$ itself is a classical group scheme locally of finite type. 

Let $H$ be an affine algebraic group over $\Lambda$ with an action of $W_F$ that factors through a finite quotient $Q$. Then, we have surjections of groups $H\rtimes W_F\to W_F \to Q$. 

\begin{definition}
The stack of $H$-valued Weil parameters is 
\[
\LS^{W_F}_H\coloneq \Map (BW_F, B(H \rtimes W_F))\times_{\Map (BW_F, BW_F)}\{\id_{BW_F}\}. 
\]
It is the quotient of 
\[
\Sect (W_F, H)\coloneq \Sect (W_F, H\rtimes W_F)
\]
by $H$. We similarly define $\LS^W_H$, $\Sect (W, H)$ for all $W$ such that we have $W_F\to W=W_F/I\to Q$ for $I\subset I_F$ open and also normal in $W_F$ as above, and also define $\LS_H^{I_F}$, $\Sect(I_F, H)$.  
\end{definition}

\begin{lemma}
The stack $\LS_H^{I_F}$ is classical and a disjoint union of classifying stacks. Equivalently, $\Sect (I_F, H)\coloneq \Sect (I_F, H\rtimes I_F)$ is a disjoint union of $H$-orbits. 
\end{lemma}

\begin{proof}
There is a natural identification 
\[
\Sect(I_F, H)\simeq \Hom(I_F,H \rtimes Q)\times_{\Hom(I_F, Q)}\{p_I\}, 
\]
where $p_I\colon I_F\to Q$ is the canonical projection, since 
\[
H\rtimes I_F \simeq I_F\times_{Q} (H\rtimes Q).
\]
Any such map $I_F\to H\rtimes Q$, over some $\Lambda$-algebra, factors through a finite quotient of $I_F$. So, the claim follows from \cite[Proposition 2.28]{Zhu}, \cite[Lemma 3.1.8]{HansenMann} for finite groups.
\end{proof}

Fix a Weil parameter $\co{L}\lambda\colon W_F\to H\rtimes W_F$ with $\co{L}\lambda (w)=(\lambda (w), w)$. Set $\lambda_I\coloneq \lambda|_{I_F}$. By the previous lemma, $\co{L}\lambda_I$ determines an open and closed substack $\LS_{H, \lambda_I}^{W_F}$ (resp.~subscheme $\Sect(W_F, H)_{[\lambda_I]}$) of $\LS_H^{W_F}$ (resp.~$\Sect(W_F, H)$). This component can be described explicitly. When $H^\circ$ centralizes $\lambda_I$, we also write $\Sect (W_F, H)_{\lambda_I}$ for the open and closed subscheme where the restriction to $I_F$ is $\lambda_I$. 

Consider the normalizer of the image of ${}^L \lambda_I$ in $H \rtimes W_F$. It contains $\co{L}\lambda (W_F)$, and hence has the form $K_{\lambda}\co{L}\lambda(W_F)\cong K_{\lambda}\rtimes^\lambda W_F$, where $K_\lambda$ denotes $Z_{H}(\co{L}{\lambda}_I(I_F))$ with the action of $\frob^\Z\cong W_F/I_F$ via $\co{L}\lambda$. The following claims for the underlying classical scheme are clear:

\begin{lemma}
We have natural identifications
\[
\Sect(W_F, H)_{[\lambda_I], \cl}\cong \Sect (W_F, K_{\lambda})_{\lambda_I, \cl}\times^{K_\lambda}H, \quad 
\Sect(W_F, K_\lambda)_{\lambda_I,\cl}\cong K_\lambda\co{L}{\lambda}(\frob), 
\]
here $\Sect(W_F, K_\lambda)_{\lambda_I}$ denotes sections of $K_\lambda {}^L \lambda (W_F)\to W_F$ whose restriction to $I_F$ is ${}^L\lambda_I$. 
In particular, the connected component of $\LS_{H, \lambda_I, \cl}^{W_F}$ containing $\co{L}\lambda$ identifies with a connected component of
\[
K_\lambda \co{L}{\lambda}(\frob)/^\adj K_\lambda. 
\]
\end{lemma}

\begin{proof}
The second isomorphism is given by $\phi \in \Sect(W_F, K_{\lambda})_{\lambda_I, \cl} \mapsto \phi(\frob)$. Note that for $i \in I_F$ we have 
\begin{equation*}
    \phi(\frob) \co{L}{\lambda}_I(i)\phi(\frob)^{-1} = \co{L}{\lambda}_I(\frob \cdot i \cdot \frob^{-1}) = \co{L}{\lambda}(\frob) \co{L}{\lambda}_I(i) \co{L}{\lambda}(\frob)^{-1},
\end{equation*}
so we deduce that $\phi(\frob)^{-1}\co{L}{\lambda}(\frob) \in K_{\lambda}$.
\end{proof}

\begin{corollary}
\begin{enumerate}
    \item $\LS_H^{W_F}$ is the union of $\LS_H^W$ for quotients $W$ of $W_F$ with maps $W_F\to W\to Q$, and each $\LS_H^W$ is classical and of finite type. 
    \item $\LS_H^{W_F}$ is smooth, equidimensional of dimension $0$. 
    \item For any surjection $W'\to W$, $\LS_H^W\to \LS_H^{W'}$ is an open and closed immersion. 
\end{enumerate}    
\end{corollary}

\begin{proof}
The classicality follows since the underlying classical stack has the expected dimension by the previous lemma, cf.~\cite[Proposition 2.19]{Zhu}. The other assertions are clear. 
\end{proof}

By replacing $\co{L}{\lambda}(\frob)$ by an element of $K_{\lambda}^{\circ} \co{L}{\lambda}(\frob)$, we assume $\co{L}{\lambda}(\frob)$ is semisimple, staying in the same connected component. 

\begin{lemma}{\label{lem: frobfiniteorder}}
$Z_{K_\lambda}(\co{L}{\lambda}(\frob))^\circ \co{L}\lambda(\frob)$ contains an element $(t, \frob)$ whose image in $H\rtimes Q$ has finite order.
\end{lemma}

\begin{proof}
In this proof, we regard $\co{L}\lambda$ as a map $W\to H\rtimes Q$ for some $W$. We roughly follow the proof of \cite[Theorem 3.6]{DHKM}. 

The action of $\co{L}\lambda (\frob)$ on $K_\lambda^\circ$ is semisimple, and it preserves a Borel pair $(B_K, T_K)$ by Winter \cite{Winter} (see also Steinberg \cite{SteinbergEnd}). The normalizer of $(B_K, T_K)$ in $K_\lambda \co{L}\lambda(W)$ is of finite type and its identity component is $T_K$. So, there exists an integer $m\geq 1$ such that $\frob^m=e$ in $Q$ and  $\co{L}\lambda (\frob)^m \in T_K(\Lambda) \rtimes \{e\}$. On the other hand, there exists an integer $m'\geq 1$ such that $\co{L}\lambda (\frob)^{m'}$ commutes with $\co{L}\lambda (W)$ by the finiteness of the image of $I_F$ in $W$. So, we can replace $m$ so that $\co{L}\lambda (\frob)^m$ is in the centralizer of $\co{L}{\lambda}(W)$  in $T_K \rtimes \{e\}$. As $\co{L}{\lambda}(I_F)$ commutes with $T_K$, this centralizer agrees with $T_K^{\co{L}{\lambda}(\frob)}  \rtimes \{e\}$. By further replacing $m$, we may assume it is in $T_K^{\co{L}{\lambda}(\frob), \circ} \rtimes \{e\}$. It is a torus over $\Lambda$, and the $\Lambda$-points are $m$-divisible, so there exists $s\in T_K^{\co{L}{\lambda}(\frob), \circ}(\Lambda)$ such that $(s^m, e)=\co{L}\lambda (\frob)^m$. Then, $s^{-1}\co{L}{\lambda}(\frob)\in K_\lambda^\circ \co{L}{\lambda}(\frob)$ has finite order in $H\rtimes Q$. 
\end{proof}

Let $\lambda_1$ denote the Weil parameter corresponding to $(t, \frob)=s^{-1}\co{L}{\lambda}(\frob)$. We renormalize the identification by using $\lambda_1$. The advantage is that the action of $\frob^\Z\cong W_F/I_F$ factors through a finite quotient, so the notation $\LS_{K_{\lambda_1}, \triv}^{W_F}$ is allowed; namely, the restriction to the inertia is trivial. (Note that $K_\lambda=K_{\lambda_1}$ as groups.)

\begin{corollary}
$\LS_{H, \lambda_I}^{W_F} \simeq \LS_{K_{\lambda_1}, \triv}^{W_F}$. 
\end{corollary}

With the torus $T_K$ fixed, the element $t$ is not unique, and the ambiguities are elements of finite order in $(T_K^{\co{L}{\lambda}(\frob)})^\circ$. To normalize it, we make the following choice. 

Fix a splitting of $\Z$-modules $\Lambda^{\times} = \Lambda_{\tor}^{\times} \oplus \Lambda_{\tf}^{\times}$ where $\Lambda_{\tor}^{\times}$ is torsion and $\Lambda_{\tf}^{\times}$ is torsion-free, and assume further that $q \in \Lambda_{\tf}^{\times}$ for later use. For any torus $T$, the identification $T(\Lambda) = \Hom(X^{\ast}(T), \Lambda^{\times})$ and our fixed decomposition of $\Lambda^{\times}$ induce a decomposition 
\begin{equation*}
    T(\Lambda) = T(\Lambda)_{\tor} \oplus T(\Lambda)_{\tf}.
\end{equation*}
Note that this decomposition is functorial in $T$. Moreover, if $t \in T(\Lambda)_{\tf}$, then the Zariski closure of $t^\Z$ in $T(\Lambda)$ is a torus. Indeed, this is given by $A \subset T$ such that $X^{\ast}(A) = X^{\ast}(T) / \{ \chi : \chi(t) = 1\}$. Then $X^{\ast}(A)$ is torsion-free since $t$ is and by functoriality, $t \in A_{\tf}$. Hence for an element of $G(\Lambda)$ lying in some torus $T$, being in the torsion-free part is independent of the choice of torus, so we just call $t$ torsion-free. In particular, this implies that if $t \in T(\Lambda) \subset G(\Lambda)$ and $gtg^{-1}=t$, then $g$ preserves the decomposition $t=t_{\tor}t_{\tf}$ and hence commutes with each of $t_{\tor}$ and $t_{\tf}$.

\begin{lemma}{\label{lem: frobdecomp}}
The above choice determines a commuting decomposition $\co{L}{\lambda}(\frob)=(s,e)(t,\frob)$ such that $s$ is torsion-free, $(t, \frob)$ is of finite order in $H \rtimes Q$ and 
\[
Z_{K_\lambda} (\co{L}{\lambda}(\frob)) = Z_{Z_{K_\lambda} (t, \frob)}(s) \subset Z_{K_\lambda} (t, \frob). 
\]
Moreover, the $q$-eigenspace of $\co{L}{\lambda}(\frob)$ agrees with the $q$-eigenspace of $s$ in $\Lie K_\lambda^{(t,\frob)}$. 

Finally, such a decomposition is unique.
\end{lemma}

\begin{proof}
We first explain the equality at the level of Lie algebras. In the proof of Lemma \ref{lem: frobfiniteorder}, we now choose $s$ so that it is torsion-free under the above decomposition. Then, the action of $s$ has no nontrivial torsion eigenvalue on $\Lie K_\lambda$. On the other hand, as $(t, \frob)$ is of finite order, its eigenvalues are torsion. Since the actions of $s$ and $(t, \frob)$ are commuting, it follows that the simultaneous centralizer of $(t, \frob)$ and $(s,e)$ agrees with the centralizer of $\co{L}{\lambda}(\frob)$ (use that commuting semisimple operators are simultaneously diagonalizable and eigenvalues of the product are products of eigenvalues of the individual operators). The claim about the $q$-eigenspace follows by the same argument. 

Let $k\in Z_{K_\lambda} (\co{L}{\lambda}(\frob))(\Lambda)$. While $k$ may not preserve $T_K$, we can choose $m$ so that $\co{L}{\lambda}(\frob)^m$ generates a (connected) torus and $s$ can be constructed inside it. Then, by construction, $(ksk^{-1})^m=s^m$. By the current choice of $s$, this implies $ksk^{-1}=s$ as both are torsion-free. So, $k\in Z_{K_\lambda} (t, \frob)$. Now, the desired equality also follows. 

For uniqueness, if we have $\co{L}{\lambda}(\frob) = (s', e)(t', \frob)$, then choose a sufficiently large integer $N$ such that
\begin{equation*}
    (s', e)^N = \co{L}{\lambda}(\frob)^N = (s, e)^N.
\end{equation*}
Now, $s$ and $s'$ commute by the equality of centralizers in the first part, $s'/s$ is torsion-free and of finite order, hence $s=s'$.
\end{proof}

\begin{remark}
If we conjugate $\lambda$ by $h\in H(\Lambda)$, $s$ is also sent to the conjugate by $h$. In fact, this construction is functorial with respect to any map $H\to H'$ by the argument in the above proof. 
\end{remark}

Write $t_{\lambda}$ for the $t$ chosen above, and take $s_\lambda \in K_\lambda^{\co{L}{\lambda}(\frob), \circ}(\Lambda)$ with ${}^L \lambda(\frob)=(s_\lambda, e)(t_\lambda, \frob)=(t_\lambda, \frob)(s_\lambda, e)$. This is somewhat analogous to the hyperbolic-elliptic decomposition in \cite[\S 3.3]{CFMMX}. In the essentially discrete case, they indeed agree as follows. 

Our choice of $\Lambda^\times= \Lambda_{\tor}^{\times}\oplus \Lambda_{\tf}^{\times}$ determines a unique $q^{1/2}\in \Lambda_{\tf}^{\times}$. Suppose we are given a map $\co{L}{\phi}\colon W_F\times \SL_2\to H\rtimes W_F$ with the associated Weil parameter $\co{L}{\lambda}\colon W_F\to H \rtimes W_F$, 
\[
\co{L}{\lambda}(w)=\co{L}{\phi}(w, \diag(\lvert w\rvert^{1/2}, \lvert w\rvert^{-1/2}))
\]
determined by the choice of $q^{1/2}$. 

\begin{lemma}{\label{lem: tftorsdisclemma}}
If $\co{L}{\phi}(\frob, e)$ has finite order in $H\rtimes Q$, then $(t_\lambda, \frob)= \co{L}{\phi}(\frob, e)$ and $(s_\lambda, e)=\co{L}{\phi}(e, \diag(q^{1/2}, q^{-1/2}))$. 

If $H=\wh G$ and $\co{L}{\phi}$ is essentially discrete, then $\co{L}{\phi}(\frob,e) = c(t_{\lambda}, \frob)$ where $c \in Z(H)^{\Gamma, \circ}$ is torsion-free and $c^{-1}(s_\lambda, e)=\co{L}{\phi}(e, \diag(q^{1/2}, q^{-1/2}))$. In particular, if $H$ is semisimple, $c=1$.
\end{lemma}

\begin{proof}
For the first part, observe that $\co{L}{\phi}(e, \diag(q^{1/2}, q^{-1/2}))$ is contained in a torus $T$, and 
\[
\co{L}{\phi}(e, \diag(q^{1/2}, q^{-1/2}))\in T(\Lambda)_{\tf},
\]
by the choice of $q^{1/2}$. Moreover, $\co{L}{\phi}(\frob, e)$ is of finite order by assumption. It follows from uniqueness as in Lemma \ref{lem: frobdecomp} that $(s_\lambda, e)=\co{L}{\phi}(e, \diag(q^{1/2}, q^{-1/2}))$. 

For the second part, let $c$ be the torsion-free part of $\co{L}{\phi}(\frob, e)$ and hence $c \in S_{\phi}\coloneq Z_{\wh G}(\phi)$ (by Lemma \ref{lem: frobdecomp}). By the paragraph before Lemma \ref{lem: frobdecomp}, the Zariski closure of $c^\Z$ is a torus and hence $c \in S^{\circ}_{\phi} \subset Z(H)^{\Gamma, \circ}$. Then $c \co{L}{\phi}(e, \diag(q^{1/2}, q^{-1/2}))$ is torsion-free and $\co{L}{\lambda}(\frob)(c \co{L}{\phi}(e, \diag(q^{1/2}, q^{-1/2})))^{-1}$ is of finite order and commutes with $c \co{L}{\phi}(e, \diag(q^{1/2}, q^{-1/2}))$. The result then follows from the uniqueness in Lemma \ref{lem: frobdecomp}.
\end{proof}

We define 
\[
J_\lambda\coloneq Z_{K_\lambda}(t_\lambda, \frob) \supset  H_\lambda \coloneq Z_{J_\lambda} (s_\lambda)=Z_{K_\lambda} (\co{L}{\lambda}(\frob))
\]
following \cite[\S 5.3]{CFMMX}. Let $\overline{K_{\lambda_1}}$ denote the reductive quotient of $K_{\lambda_1}$.\footnote{Using the semisimplicity of $\lambda_1$, one can show that $\overline{K_{\lambda_1}}$ is the centralizer of $\lambda_1|_{I_F}$ in the reductive quotient $\overline{H}$ of $H$, cf.~the proof of Lemma \ref{lemma: NonreductiveLocalization}.}
We consider the map
\[
\LS_{H, \lambda_I}^{W_F} \simeq K_{\lambda_1} \frob /^\adj K_{\lambda_1} \to \overline{K_{\lambda_1}}\frob \sslash^\adj \overline{K_{\lambda_1}}, 
\]
and the formal completion of $\LS_{H, \lambda_I}^{W_F}$ along the fiber of the image $[\overline{\lambda}]$ of $\lambda$ is denoted by $(\LS_{H}^{W_F})^{\wedge}_\lambda$. (This construction does not depend on the choice of $\lambda_1$.)

\begin{proposition}\label{prop: ReductionToUnipotentWeilParameters}
With the above choices, there is a map
\[
\Map (B\Z, BJ_\lambda)\cong J_\lambda /^\adj J_\lambda \to \LS_{H, \lambda_I}^{W_F} 
\]
sending $1, s_\lambda\in J_\lambda^\circ (\Lambda)$ to $\lambda_1, \lambda$, and inducing a functorial identification
\[
\mc L^u (B H_\lambda)
\xrightarrow[\simeq]{sh_{s_\lambda}}
\mc L^u_{s_\lambda}(BH_\lambda)
\xrightarrow[\simeq]{l^u_{s_\lambda}}
\mc L^u_{s_\lambda}(BJ_\lambda)
\simeq
(\LS_{H}^{W_F})^{\wedge}_\lambda. 
\]
\end{proposition}

\begin{proof}
The essential surjectivity at the level of $\Lambda$-points follows from the next lemma, and then it is an equivalence by the Jordan decomposition. 

Therefore, it suffices to compare the cotangent complexes, or dually tangent complexes. Let $u\in \mc U_{H_\lambda}(\Lambda)$. Then, the tangent complex of the source at the point corresponding to $s_\lambda u \in J_\lambda (\Lambda)$ is
\[
\Lie J_\lambda \xrightarrow{1 - \Ad(s_\lambda u)} \Lie J_\lambda
\]
in degrees $-1, 0$, while the tangent complex at the image point is 
\[
\Lie K_{\lambda_1} \xrightarrow{1 - \Ad(s_\lambda u)\Ad (t_\lambda, \frob)} \Lie K_{\lambda_1}. 
\]
The latter decomposes according to eigenspaces of $\Ad (t_\lambda, \frob)$. On the eigenvalue-$1$ piece, it is identical to the former. It remains to observe that if the eigenvalue is not $1$, then it is a nontrivial root of unity. But $s_\lambda$ is torsion-free, and hence these pieces are acyclic. See also \cite[Propositions 2.16, 2.17]{YangZhu} (and \cite[Lemma 2.9, Proposition 2.44]{ZhuTame}). 
\end{proof}

\begin{lemma}\label{lem: FibersOfGITQuotientMap}
The (reduced) fiber of $[\overline{\lambda}]$ identifies with $\co{L}\lambda(\frob)\mc U_{H_\lambda}/^\adj H_\lambda$.
\end{lemma}

\begin{proof}
The case where $K_{\lambda_1}$ is reductive is standard and, for instance, follows from \cite[Proposition 3.18]{Mohrdieck}. The general case reduces to the reductive case by Lemma \ref{lemma: semisimpleconjugation} (applied to the disconnected group generated by $K_{\lambda_1}\frob$). 
\end{proof}

There is a slight reformulation in our main situation. 

\begin{definition}
Suppose $H=\wh G$ or $\wh P$ for a reductive group $G$ over $F$ or its semistandard parabolic subgroup $P$.
We define 
\[
X_G\coloneq \Sect(W_F, \wh G) \sslash \wh G, \quad X_P\coloneq \Sect(W_F, \wh P) \sslash \wh P
\]
to be the union of affine classical GIT quotients. They come with natural maps $\LS_{\wh G}^{W_F}\to X_G$, $\LS_{\wh P}^{W_F}\to X_P$.  
\end{definition}

\begin{lemma}\label{lem: X_P=X_M}
If $\wh P\to \wh M$ denotes the Levi quotient, then $X_P\cong X_M$ naturally. Moreover, every fiber of $\LS_{\wh P}^{W_F}\to X_P$ contains a unique closed point, which corresponds to the Frobenius-semisimple parameter.   
\end{lemma}

\begin{proof}
We work with one connected component containing $\lambda$, and we may assume $\wh M=Z_{\wh G}(S)$ for a nontrivial subtorus $S\subset \wh{T}^{\Gamma}\cap  K_{\lambda_1}$, and $\wh P$ is defined by a cocharacter $\chi$ of $S$. Take $mu\in K_{\lambda_1}$ with $m\in \wh M\cap K_{\lambda_1}$ and $u\in U_{\wh P}\cap K_{\lambda_1}$. It follows easily that the closure of the $\wh P \cap K_{\lambda_1}$-orbit of $mu\cdot \frob$ contains $m\cdot \frob$ (and hence its $\wh P \cap K_{\lambda_1}$-orbit). This immediately implies the desired claims. 
\end{proof}

\begin{corollary}\label{cor: XP to GG is finite}
The map $X_P\cong X_M\to X_G$ is finite. In particular, the inverse image of a connected component of $X_G$ is a union of finitely many connected components of $X_M$. 
\end{corollary}

\begin{proof}
The map $X_M\to X_G$ is also induced from the embedding $\wh{M}\to \wh{G}$. Then, the finiteness of $X_M\to X_G$ follows from, for instance, the twisted version of the Chevalley--Steinberg restriction theorem \cite[Theorem 1.1]{Mohrdieck} (or \cite[Proposition 6.6]{DHKM}) as in the proof of \cite[Theorem 2.3]{DHKMfiniteness}.  
\end{proof}

The following lemma will be used explicitly or implicitly throughout; our representative of $\lambda$ will be chosen using it. 

\begin{lemma}\label{lem: TJcontrol}
Assume $H=\wh G$. After replacing $\lambda$ by its $\wh G$-conjugate, all of the following hold: $\wh T\cap J^{\circ}_\lambda$ is a maximal (connected) torus of $J_\lambda$ and is contained in $\wh T^\Gamma$, and it contains $s_\lambda$; in particular, $s_\lambda \in \wh T^{\Gamma, \circ}$. Moreover, if $\lambda$ factors through a semistandard parabolic subgroup $\wh P$, then $\wh P$ is defined by some cocharacter factoring through $\wh T\cap J^{\circ}_\lambda$. 
\end{lemma}

\begin{proof}
Choose a maximal (connected) torus $\wh T_J$ of $J^{\circ}_\lambda$ containing $s_\lambda$. Then $Z_{{}^L G}(\wh T_J)=Z_{\wh G}(\wh T_J){}^L \lambda_1(W_F)$ is a Levi subgroup of ${}^L G$ in the sense of Borel. By \cite[Lemma 3.5]{Bor2}, after replacing $\lambda$ if necessary, this Levi subgroup becomes semistandard so that its connected center $\wh T_J (\ni s_\lambda)$ is contained in $\wh T^\Gamma$. Since $J_\lambda^\circ$ is reductive for $H=\wh G$, we have $Z_{J_\lambda^\circ}(\wh T_J)=\wh T_J$ and it follows that $\wh T\cap J_\lambda^\circ=\wh T_J$. 

For the second part, let $\overline{\lambda_P}$ denote the composite of $\lambda_P\coloneq \lambda$ with $\wh P\to \wh M$ for the Levi quotient $\wh M$. By Lemma \ref{lem: X_P=X_M}, ${}^L \lambda_P$ is $\wh P$-conjugate to ${}^L\overline{\lambda_P}$ post-composed with $\wh M\rtimes W_F\hookrightarrow \wh P\rtimes W_F$. Therefore, if $\nu$ is a cocharacter of $Z(\wh M)^\Gamma$ defining $\wh P$, then its corresponding $\wh P$-conjugate centralizes ${}^L \lambda (W_F)$, i.e., it is a cocharacter of $H_\lambda$, and it still defines $\wh P$. If necessary, we can take a further $H_\lambda^\circ \cap \wh P$-conjugate so that it is a cocharacter of $\wh T_J$.
\end{proof}

\subsubsection{Reduction to the graded unipotent case}
It is sometimes more convenient to make a further reduction as in \cite[\S 5.4]{CFMMX}, especially to apply \cite{LusztigGraded}. 

\begin{lemma}\label{lem: ReductionToGradedUnipotent}
Fix a maximal torus $\wh T_J$ of $J_\lambda$ containing $s_\lambda$. There exist a connected subgroup $I_\lambda \subset J_\lambda$ containing $\wh T_J$, a cocharacter $\chi\colon \Gm\to \wh T_J \subset I_\lambda$, and an integer $n\geq 1$ such that
\[
Z_{I_\lambda}(\chi)=H_\lambda^\circ, \quad \Lie I_\lambda (nr) = (\Lie J_\lambda)^{s_\lambda=q^r}, 
\]
where $r\in \frac{1}{n}\Z$, $\Lie I_\lambda (nr)$ is the $\chi$-weight $nr$ space of $\Lie I_\lambda$, and $q^r\in \Lambda_\tf^\times$, an element uniquely determined by $r$. In particular, $V_\lambda=\Lie I_\lambda (n)$ for the Vogan variety $V_\lambda$ in Definition \ref{def: VoganVariety}. 

If $G'$ is a reductive group, then we can make the construction functorial for all semistandard parabolic or Levi subgroups $G\subset G'$, i.e., we obtain $I_\lambda\to I_{\lambda'}$ with common $\chi$ and $n$, where $\lambda'$ is induced by $\lambda$. Moreover, we can make $I_{\lambda'}$ reductive and arrange that $\chi$ maps into $\wh T_J\cap I_{\lambda', \der}$. 
\end{lemma}

\begin{proof}
Take an integer $n\geq 1$ such that if $q^r\in \Lambda_\tf^\times$ is an eigenvalue of $s_\lambda$ in $\Lie J_\lambda$, then $nr\in \Z$. After replacing $n$ by a multiple if necessary, choose a cocharacter $\chi\colon \Gm\to \wh T_J $ such that $\Lie J_\lambda (nr) \supset (\Lie J_\lambda)^{s_\lambda=q^r}$ for any $r\in \frac{1}{n}\Z$. Then $I_\lambda=Z_{J_\lambda}(s_\lambda\chi (q^{1/n})^{-1})^{\circ}$ works. 

For the last part, note that we can choose $n$ satisfying the above condition for both $G$ and $G'$. Since $J_\lambda$ is a subgroup of $J_{\lambda'}$ containing a maximal torus of $J_{\lambda'}$, we can choose the same $\chi$. Then $I_\lambda\subset I_{\lambda'}$ is clear. 
\end{proof}

\begin{example}{\label{ex: Ilambda}}
Suppose $H=\wh G$ for a split $G$. If the parameter $\lambda$ arises from
\[
W_F\to \SL_2 \to \wh G;\quad w \mapsto \diag (\lvert w\rvert^{1/2}, \lvert w\rvert^{-1/2}) \mapsto \iota (\diag (\lvert w\rvert^{1/2}, \lvert w\rvert^{-1/2}))
\]
for some map $\iota\colon \SL_2\to \wh G$, then we can take $I_\lambda=\wh G$ and $n=2$.

Suppose now that $\phi$ is essentially discrete and $\tilde{\phi}\colon W_F \times \SL_2 \rightarrow \co{L}{G}$ is a lift and $\iota: \Gm \rightarrow \wh{G}$ is the restriction of $\tilde{\phi}$ to the maximal torus of $\SL_2$. Then we can take $I_{\lambda}$ to be $J^{\circ}_{\lambda}$ with $\chi = \iota$ and $n=2$. First note that by Lemma~\ref{lem: tftorsdisclemma} we have $s_{\lambda} = c\iota(q^{1/2})$. Thus, the eigenvalues of $s_{\lambda}$ under the adjoint action are half-integral powers of $q$. But then by Lemma \ref{lem: discequalsdistforunip} for $J^{\circ}_{\lambda}$ and the fact that distinguished unipotent orbits are even (\cite[Theorem 8.2.3]{CollingWoodMcGovern}), the eigenvalues of $s_{\lambda}$ are integral powers of $q$.
\end{example}

\subsection{The stack of Weil--Deligne parameters}
Let $\WD_F\coloneq W_F \ltimes \Ga$ denote the Weil--Deligne group, where $W_F$ acts on $\Ga$ via $W_F\to \frob^\Z$ and multiplication by powers of $q$. Recall the following form of local Tate duality:

\begin{lemma}[{\cite[Corollary 4.9]{takayalin}}]\label{lem: localTateDuality}
There is an isomorphism
\[
\Gamma (\WD_F, V)^{\vee} \simeq \Gamma(\WD_F, V^{\vee}(1)[2]),
\]
where $-(1)$ is the Tate twist, for any finite-dimensional representation $V$ of $\WD_F$. 
\end{lemma}

\begin{remark}
More precisely, the identification follows from the explicit description of group cohomology \cite[Proposition B.6]{takayalin}. 
\end{remark}

Let $H$ be an affine algebraic group over $\Lambda$ with an action of $W_F$ that factors through a finite quotient $Q$ as before.

\begin{definition}
The stack $\LS_H^{\WD_F}$ is defined as the fiber product
\[
\Map (B\WD_F, B(H\rtimes \WD_F))\times_{\Map (B\WD_F, B\WD_F)}\{\id_{B\WD_F} \}
\]
with a canonical map to $BH\simeq B(H\rtimes \WD_F) \times_{B\WD_F}*$. It is the quotient of $\Sect(\WD_F, H)$ by a natural action of $H$.

Furthermore, if $H$ is the dual group $\wh{G}$ of a reductive group $G$ over $F$, or its parabolic subgroup $\wh P$, we refer to $\Lambda$-points of $\Sect(\WD_F, H)$ as $L$-parameters, or Weil--Deligne ($L$-)parameters, of $G$, $P$ and write $\Par_G$, $\Par_P$ for $\LS_{\wh G}^{\WD_F}$, $\LS_{\wh P}^{\WD_F}$.
\end{definition} 

\begin{lemma}
For a point $(\lambda, N)\in \Sect(\WD_F, H)(\Lambda)$, 
\[
\bb L_{\LS^{\WD_F}_H, (\lambda, N)} \simeq
\Gamma (\WD_F, (\Lie H)_{(\lambda, N)}^{\vee}(1))[1], 
\]
where $(\Lie H)_{(\lambda, N)}$ denotes $\Lie H$ with the adjoint action induced by $(\lambda, N)$. 
\end{lemma}

\begin{remark}
Compare with \cite[Proposition 2.15, Remark 2.16]{Zhu}. See also \cite[Proposition 4.10]{takayalin}. 
\end{remark}

\begin{proof}
It follows from Lemma \ref{lem: MappingStackCotangentComplex}, Remark \ref{rem: LeftAdjointOfPullBack} and Example \ref{ex: GroupHomologyExample} (2) that the cotangent complex identifies with
\[
\pi_\natural (\Lie H)_{(\lambda, N)}^{\vee}[-1]
\simeq
\Gamma (\WD_F, (\Lie H)_{(\lambda, N)})^{\vee}[-1], 
\]
where $\pi\colon B\WD_F\to \Spec \Lambda$ denotes the structure map. Now, the claim follows from Lemma \ref{lem: localTateDuality}.
\end{proof}

Fix $\lambda\in \Sect(W_F, H)(\Lambda)$ such that $\lambda=\lambda_1$, and we define substacks and subschemes as before using the restriction to $I_F$. 

\begin{lemma}
The natural map
\[
\Sect (\WD_F, K_\lambda)_{\lambda_I} /K_\lambda\rightarrow \LS_{H, \lambda_I}^{\WD_F}   
\]
is an equivalence. 
\end{lemma}

\begin{proof}
It is an equivalence at the level of underlying classical stacks. The claim then follows from a comparison of cotangent complexes (see \S\ref{ss: DAGbasiclocaldefns}). 
We have
\begin{align*}
\Gamma (\WD_F, (\Lie H)_{(\lambda, N)})
&\simeq
\Gamma( \frob^\Z\ltimes \Ga, \Gamma (I_F, (\Lie H)_{(\lambda, N)})) \\
&\simeq 
\Gamma( \frob^\Z\ltimes \Ga, (\Lie K_{\lambda})_{(\lambda, N)})\\
& \simeq \Gamma(\WD_F, (\Lie K_{\lambda})_{(\lambda, N)}).
\end{align*} 
Note that the adjoint action of $\lambda_I$ is semisimple; for the second isomorphism, $I_F$ acts through a finite quotient, and hence the higher cohomology is torsion. But $(\Lie H)_{(\lambda, N)}$ is a $\Lambda$-vector space and hence this higher cohomology is trivial.
\end{proof}

Let $\mc U_H \subset H$, $\mc N_H\subset \Lie H$ denote the unipotent and nilpotent cones as reduced schemes; see \cite[\S4.2, \S4.3]{Cotner} for some properties stated in the case of reductive groups. We write $\wh{\mc U}_H$ for the formal completion along $\mc U_H$. 

\begin{proposition}\label{prop: LSDescriptions}
The connected components containing $\lambda$ of the following schemes are all equivalent:
\begin{enumerate}
    \item $\Sect(\WD_F, K_{\lambda})_{\lambda_I}$. 
    \item $(\Sect (W_F, K_{\lambda})_{\lambda_I}\times \wh{\mc U}_{K_{\lambda}})\times_{\wh{\mc U}_{K_{\lambda}}\times \wh{\mc U}_{K_{\lambda}}}\wh{\mc U}_{K_{\lambda}}$. 
    \item $(\Sect(W_F, K_{\lambda})_{\lambda_I} \times \mc U_{K_\lambda} )\times_{\mc U_{K_\lambda}\times \mc U_{K_\lambda}}\mc U_{K_\lambda}$. 
    \item $(\Sect(W_F, K_{\lambda})_{\lambda_I} \times \mc N_{K_\lambda} )\times_{\mc N_{K_\lambda}\times \mc N_{K_\lambda}}\mc N_{K_\lambda}$. 
    \item $(\Sect(W_F, K_{\lambda})_{\lambda_I} \times \Lie K_\lambda)\times_{\Lie K_\lambda\times \Lie K_\lambda}\Lie K_\lambda$. 
    \item $(\Sect(W_F, K_{\lambda})_{\lambda_I} \times \Lie K_\lambda)\times_{\Lie K_\lambda}\{0\}$.
\end{enumerate}
Here, the two maps $\Sect (W_F, K_{\lambda})_{\lambda_I}\times \wh{\mc U}_{K_{\lambda}} \to \wh{\mc U}_{K_{\lambda}}$ are induced by conjugation by $\frob$ and the $q$-th power on $\wh{\mc U}_{K_{\lambda}}$, while the two maps $\Sect(W_F, K_{\lambda})_{\lambda_I} \times \mc N_{K_\lambda} \to \mc N_{K_\lambda}$ are induced by the conjugation by $\frob$ and the multiplication by $q$ on $\mc N_{K_\lambda}$. Other left maps are similarly defined, except that the map $\Sect(W_F, K_{\lambda})_{\lambda_I} \times \Lie K_\lambda\to \Lie K_\lambda$ in (6) is the difference of the conjugation and the multiplication by $q$. The right maps are diagonal except the last one. 
\end{proposition}

\begin{remark}
If $H^\circ$ is reductive, the stack $\LS_H^{\WD_F}$ is classical as we will recall below, and it was first considered in \cite{BellovinGee, BooherPatrikis}, and then further discussed in \cite{DHKM, fargues--scholze, HellmannIwahori, Zhu}. The derived structure of (6) appears in \cite{HellmannIwahori}. 
\end{remark}

\begin{proof}
For (1) $\simeq$ (2), we may replace $\WD_F$, $W_F$ by $\Z\ltimes \Ga$, $\Z$. Then the claim follows from $\Map(B\Ga, BK_{\lambda})\simeq \wh{\mc U}_{K_{\lambda}}/^\adj K_\lambda$. (3) $\simeq$ (4) is via the exponential and logarithm. (5) $\simeq$ (6) is clear. 

Next, we prove (2) $\simeq$ (5). Since the exponential and logarithm induce an isomorphism $\wh{\mc U}_{K_{\lambda}}\cong \wh{\mc N}_{K_{\lambda}}$ as classical formal schemes\footnote{It is the correspondence between topologically unipotent matrices and topologically nilpotent matrices.}, it reduces to the claim that if $k\in \Lie K_{\lambda}$ and $k$, $qk$ are conjugate, then $k\in \mc N_{K_{\lambda}}$. It reduces to the case $K_{\lambda}=\GL_N$ for some integer $N$, in which case the given conditions imply that the characteristic polynomial of $k$ in $X$ is $X^N$ as we are working over $\Q$. 

It remains to show (4) $\simeq$ (5). Let $A\coloneq \Lie \overline{K_\lambda^\circ} \sslash \overline{K_\lambda^\circ}$, which is naturally a weighted affine space. Since the nilpotent cone is the (derived) fiber over $0$ of $\Lie K_\lambda \to A$ by Kostant's theorem (applied to $\overline{K_\lambda^\circ}$), (4) identifies with the base change of (5) along $\{0\}\to A \times_{A\times A}A$, where the left map in the fiber product is $(\theta, q^\bullet)$, with $\theta$ the weight-preserving automorphism of $A$ induced by $\Ad(t_\lambda, \frob)$ and $q^\bullet$ the map induced by multiplication by $q$ on $\Lie K_\lambda$, and the right map is the diagonal. Note that $A \times_{A\times A}A$ identifies with the scheme defined by the equation $\theta=q^\bullet$ inside $A$, whose underlying reduced scheme is $\{0\}$ as $\theta$ is of finite order and all weights are nonzero. As $\theta$ is weight-preserving, one can compute that the partial derivative of $\theta-q^{\bullet}$ at $0$ is the diagonal matrix $\zeta -q^\bullet$ for some roots of unity $\zeta$ (depending on coordinates). As all weights are nonzero, it indeed cuts out a smooth subscheme $\{0\}$. 
\end{proof}

\begin{remark}
While $\LS_{H}^{\WD_F}$ is completely canonical, the identification with other moduli involves choices. On the other hand, note that any element $a\in \Lambda^\times$ gives rise to an automorphism of $\LS_H^{\WD_F}$ by acting by multiplication on $\Lie K_\lambda$ in the above description. This accounts for the difference of identifications for different choices. Cf.~\cite[Proposition 3.12]{Zhu}. 
\end{remark}

\begin{proposition}\label{prop: properties_of_Sect}
Fix surjections $W_F \to W\to W'$ over $W_F\to Q$ with open kernels in $I_F$. 
\begin{enumerate}
    \item The scheme $\Sect (W\ltimes \Ga, H)$ is quasi-smooth, affine, almost of finite type, and its virtual dimension is $\dim H$. 
    \item The map $\Sect (W'\ltimes\Ga, H)\to \Sect (W\ltimes\Ga, H)$ is an open and closed immersion. 
    \item The stack $\Sect (W\ltimes \Ga, H)/H$ is an Artin stack almost of finite type, and it is perfect and QCA. In particular, our notation $\IndCoh(\Sect (W\ltimes \Ga, H)/H)$ makes sense, and it is compactly generated by coherent sheaves. 
    \item The stack $\Sect (W\ltimes \Ga, H)/H$ is quasi-smooth, in particular Gorenstein, and its dualizing complex is trivial in degree 0. 
\end{enumerate}
Suppose now $H^{\circ}$ is reductive. 
\begin{enumerate}
    \item[(5)] The scheme $\Sect (W\ltimes \Ga, H)$ is classical. In particular, it is a local complete intersection of pure dimension $\dim H$ and Cohen--Macaulay. 
    \item[(6)] The GIT quotient $\Sect (W\ltimes \Ga, H) \sslash H$ is an affine scheme of finite type. 
    \item[(7)] For any point $\lambda \in (\Sect(W_F, H) \sslash H)(\Lambda)$, the fiber \[
    \Sect (\WD_F, H)\times_{\Sect(W_F, H) \sslash H}\{\lambda\}\] contains a unique closed $H$-orbit. Any parameter in the closed orbit is Frobenius-semisimple with $N=0$ and, for any parameter $(\lambda', N)\in \Sect (\WD_F, H)(\Lambda)$ in the fiber, the Frobenius-semisimplification of $\lambda'$ (with $N=0$) is in the closed $H$-orbit. 
\end{enumerate}
\end{proposition}
\begin{proof}
(1)--(4): most claims follow from the previous discussion, and it remains to show the dualizing complex is trivial. This claim is essentially \cite[Propositions 2.35, 3.7]{Zhu}. In fact, it also follows from the descriptions we gave and Proposition \ref{prop: QuasismoothDualizingStack}. 

(5) It suffices to show the corresponding statement for the underlying classical stacks; then, the classicality follows from the quasi-smoothness and the dimension estimate. We can identify $\LS^{\WD_F}_H$ with some other definitions, cf.~\cite[Proposition 4.4]{ScholzeLLCmotives}. Then, this follows from \cite[Theorem 4.1]{DHKM}, \cite[Proposition 3.7]{Zhu}. Alternatively, one can argue as in \cite[Proposition 2.1]{HellmannIwahori}. (6) follows from Nagata's theorem on invariants. (7) follows from \cite{DHKM}. More precisely, \cite[Proposition 4.17]{DHKM} implies that the inclusion from the reduced locus $N=0$ induces an isomorphism of GIT quotients. Indeed, it is a surjective closed immersion into a reduced scheme\footnote{Recall the following argument: as $\Sect(\WD_F, H)$ is generically smooth \cite{BellovinGee, BooherPatrikis}, it is reduced by the Cohen--Macaulayness. And so is the GIT quotient.}, and it must be an isomorphism. So, we can consider instead the map
\[
\Sect (\WD_F, K_\lambda)_{\lambda_I} \to K_\lambda \frob \sslash K_\lambda
\]
and it satisfies the desired properties. 
\end{proof}

Assume now $W_F$ acts trivially on $H$ and $\lambda_I$ is trivial, and write $\LS_H^{\WD_F, \up}$,  $\Par_G^{\up}$, $\Par_P^{\up}$ for $\LS^{\WD_F}_{H, \lambda_I}$ in this case. As in \cite[\S 3.2]{BCHNconstructible}, we can identify $\LS_H^{\WD_F, \up}$ with a certain loop space, which they call $q$-twisted loops. 

\begin{lemma}\label{lem: loopInterpretationGlobally}
There are natural identifications
\begin{align*}
\LS^{\WD_F, \up}_{H} 
&\simeq
\mc L(\Lie H/(H\times \Gm))\times_{\mc L(B\Gm)}\{q\} \\
&\simeq
\mc L(\mc N_H/(H\times \Gm))\times_{\mc L(B\Gm)}\{q\},  
\end{align*}
where $\Gm$ acts on $\Lie H$ by the inverse of the scalar multiplication. \end{lemma}

\begin{proof}
    This follows from (4) and (5) of Proposition \ref{prop: LSDescriptions} using the presentation \eqref{eqn: chenfiberprod} and the identification
    \begin{equation*}
        \mc{L}(B\Gm) \simeq \Gm \times B\Gm.
    \end{equation*}
\end{proof}

\subsection{Vogan stacks}{\label{ss: voganstacks}}
The Frobenius-semisimple locus will play a key role. 
If $H=\wh G$, $\wh P$, recall that we write $X_G$, $X_P$ for the GIT quotients $\Sect (W_F, \overline{H}) \sslash \overline{H}$ respectively. 

\begin{definition}
A point $(\lambda, N)\in \Sect (\WD_F, H)(\Lambda)$ is Frobenius-semisimple (resp.~semisimple) if ${}^L\lambda(\frob)$ is semisimple (resp.~if moreover $N=0$).\footnote{Any point $(\lambda, N)$ projects to its Frobenius-semisimplification: for any Frobenius lift $\sigma$, replace ${}^L \lambda(\sigma)$ by its semisimple part. Together with $N$, it still gives a point of $\Sect(\WD_F, H)$, cf.~\cite[Lemma 2.21]{ZhuTame}.} The $H$-orbit of a semisimple point corresponds to a $\Lambda$-point of $X_H$. 

If $H$ is reductive, two equivalent notions of semisimplicity are (see \cite[Theorem 5.5, Corollary 5.13]{DatSurvey} and Proposition \ref{prop: properties_of_Sect})
\begin{itemize}
    \item $(\lambda, N)$ is in a closed $H$-orbit, 
    \item if $(\lambda, N)$ factors through a parabolic subgroup $P\rtimes W_F \subset H\rtimes W_F$, then $(\lambda, N)$ factors through some Levi subgroup of $P\rtimes W_F$.
\end{itemize}

For a semisimple $(\lambda, 0)$, we define $\LS^{\WD_F}_{H,\lambda}$ to be the fiber over the closed point $[\lambda]$ of the map
\[
\LS^{\WD_F}_H \to X_H=\Sect (W_F, \overline{H})\sslash \overline{H}. 
\]
If $H=\wh G, \wh P$, we also write $\Par_{G,\lambda}$, $\Par_{P,\lambda}$. 
\end{definition} 

\begin{definition}\label{def: VoganVariety}
    Fix a semisimple $\co{L}{\lambda}=(\lambda, 0) \in \Sect(\WD_F, H)(\Lambda)$. Define the \emph{Vogan variety}
    \begin{equation*}
        V_{\lambda} \coloneq \{ N \in \Lie H \mid \co{L}\lambda(I_F) \subset Z_{H\rtimes W_F}(N), \, \co{L}{\lambda}(\frob) \cdot N = qN \}. 
    \end{equation*}
    Let $H_{\lambda}$ denote the centralizer $Z_{H}(\co{L}{\lambda})$ and define the \emph{Vogan stack} $\V_\lambda$ to be the quotient $V_{\lambda} / H_{\lambda}$. 
\end{definition}

\begin{proposition}{\label{prop: cotangentdualizingcomplex}}
The Vogan stack $\V_{\lambda}$ is smooth, and the dualizing complex is given by the determinant line bundle $\det(\L_{\V_{\lambda}})[\rank \L_{\V_{\lambda}}]$ of the cotangent complex of $\V_{\lambda}$.
\end{proposition}

\begin{proof}
This is a special case of Proposition \ref{prop: QuasismoothDualizingStack}.
\end{proof}

\begin{lemma}
For each semisimple $\lambda$, we have a canonical embedding $\V_{\lambda} \hookrightarrow \LS^{\WD_F}_H$ identifying $\V_{\lambda}$ with $\LS^{\WD_F, \fss}_{H, \lambda}$, the Frobenius-semisimple reduced (and classical) closed locus of $\LS^{\WD_F}_{H, \lambda}$.
\end{lemma}
\begin{proof}
    We have a map $\V_\lambda\to \LS^{\WD_F, \fss}_{H, \lambda}$ induced by $N\mapsto (\lambda, N)$ and this induces the desired isomorphism.
\end{proof}

\begin{proposition}
Assume $H$ is reductive. The topological space $\V_{\lambda}$ has finitely many points. Equivalently, $V_{\lambda}$ is stratified by finitely many $H_{\lambda}$-orbits and there is a unique open orbit.
\end{proposition}

\begin{proof}
This is a variant of \cite[Proposition 5.6]{CFMMX}. 
By Lemma \ref{lem: ReductionToGradedUnipotent}, we have 
\[
V_\lambda=V_{\lambda^\ur}, \quad H^{\circ}_\lambda = H_{\lambda^\ur}, \quad \V_\lambda=V_\lambda/H_\lambda, \quad \V_{\lambda^\ur}=V_\lambda/ H_{\lambda}^\circ, 
\]
where $\lambda^{\ur}$ is a point of $\Sect (W_F, I_\lambda)^{\up}$ corresponding to $s_\lambda$. This case follows from \cite[Proposition 3.5, \S3.6]{LusztigGraded}. 
\end{proof}

On the other hand, the same is not true for $\Par_{G, \lambda}$.
\begin{example}
    In general, $\Par_{G, \lambda}$ can contain infinitely many $\wh{G}$-orbits. Consider the following example for $\GL_{21}$. We let $\lambda$ be the diagonal matrix with entries
    \begin{equation*}
          (1, q, q, q^{2}, q^{2}, q^{2}, ..., q^{5}).
    \end{equation*}
  Then $H_{\lambda} \cong \GL_1 \times \GL_2 \times \cdots\times \GL_6$. Let $N$ be the matrix such that for each diagonal $k \times k$ block of $H_{\lambda}$ giving a $\GL_k$ factor, $N$ is the identity matrix on the $k \times k$ block directly below the diagonal $k \times k$ block and then all other entries are $0$. Let $A \in Z_{H_{\lambda}}(N)$ and denote the $k \times k$ block of $A$ by $A_k$. Then the relation $AN=NA$ implies that the upper left $(k-1) \times (k-1)$ block of $A_k$ is $A_{k-1}$ and $A_k$ is upper triangular. Hence $Z_{H_{\lambda}}(N)$ can be identified with the upper triangular Borel $B \subset \GL_6$. In particular, for each unipotent $u \in B$, there is a unipotent element $A_u \in Z_{H_{\lambda}}(N)$ under this identification. Consider the Weil--Deligne parameters that map $\frob$ to $\lambda A_u$ for varying $u$ and have the fixed nilpotent $N$. Two such parameters are conjugate by $H_{\lambda}$ if and only if they are conjugate by $Z_{H_{\lambda}}(N)$, where under the above identification,  $Z_{H_{\lambda}}(N)$ acts on $A_u$ as conjugation by $B$ on $u$. It is well known (for instance, see \cite{HilleRohrle}) that the Borel subgroup of $\GL_6$ acts on its unipotent radical with infinitely many orbits, and this implies that the parameters $(\frob \mapsto \lambda A_u , N)$ fall into infinitely many $\GL_{21}$-orbits.
\end{example}

\begin{definition}
Suppose $H=\wh G$ for a reductive $G$. We consider a Weil--Deligne parameter $(\lambda, N)\in \Sect (\WD_F, \wh G)$.  
\begin{enumerate}
    \item It is essentially discrete if $Z_{\wh G}(\lambda, N)/Z(\wh G)^{\Gamma}$ is finite. 
    \item It is generic (or open) if $\Sect(\WD_F, \wh G)$ is smooth at $(\lambda, N)$. 
\end{enumerate}
\end{definition}

\begin{remark}
\begin{enumerate}
    \item If $(\lambda, N)$ is essentially discrete, then $\lambda$ is Frobenius-semisimple and $\dim \V_\lambda =-\dim Z(\wh G)^{\Gamma}$. We call such $\V_{\lambda}$ essentially discrete.
    \item The parameter $(\lambda, N)$ is generic if and only if $N$ is in the open orbit of the Vogan variety of the Frobenius-semisimplification of $\lambda$ \cite[Proposition 6.10]{DHKMllcinfamiliesbanal}.
\end{enumerate}
\end{remark}

\subsubsection{Relation with unipotent loops}{\label{sss: relationuniploops}}
Fix $\Lambda^{\times} = \Lambda_{\tor}^{\times} \oplus \Lambda_{\tf}^{\times}$ with $q\in \Lambda_{\tf}^{\times}$ as before. For any semisimple $\lambda$, we shall describe the completions of $\LS^{\WD_F}_H$ along $\LS_{H,\lambda}^{\WD_F}$ and $\V_\lambda$, for which we write
\[
\LS_{H, \lambda}^{\WD_F, \wedge}\simeq \Sect(\WD_F, K_\lambda)^{\wedge}_{\lambda_I, [\lambda]}/ K_\lambda, \quad \LS_{H, \V_\lambda}^{\WD_F, \wedge}\simeq \Sect(\WD_F, K_\lambda)^{\wedge}_{\lambda_I, V_\lambda\times^{H_\lambda} K_\lambda}/ K_\lambda,
\]
respectively, where $\lambda$ is also regarded as a semisimple parameter of $\Sect (\WD_F, K_\lambda)$. We also write $\Par_{G,\lambda}^{\wedge}$, $\Par^{\wedge}_{G,\V_\lambda}$, etc.

\begin{proposition}\label{prop: LoopSpaceIdentificationCompletions}
We have functorial identifications
\begin{equation*}
        \LS^{\WD_F, \wedge}_{H, \lambda} \simeq \mc{L}^{u}(\V_{\lambda}), \qquad \LS^{\WD_F, \wedge}_{H, \V_{\lambda}} \simeq \wh{\mc{L}}(\V_{\lambda}) \simeq \wh{\T}\V_{\lambda}[-1].
\end{equation*}    
\end{proposition}

We first check the unipotent case. 

\begin{lemma}
The statement of the proposition holds for $\LS^{\WD_F, \up}_H$ for any $H$ with trivial action of $W_F$.   
\end{lemma}

\begin{proof}
By Lemma \ref{lem: loopInterpretationGlobally}, 
\[
\LS_{H}^{\WD_F, \up} \simeq \mc L (\Lie H/ (H\times \Gm))\times_{\mc L(B\Gm)}\{q\},
\]
and, under this equivalence, we are looking at fibers of the map
\[
\mc L (\Lie H/ (H\times \Gm))\times_{\mc L(B\Gm)}\{q\} \to \mc L(BH) \to \overline{H} \sslash \overline{H}.
\]
As in (the proof of) Lemma \ref{lem: FibersOfGITQuotientMap}, 
\[
\mc L(BH)^{\wedge}_{\lambda}\simeq \lambda (\frob)\wh{\mc U}_{H_\lambda} /^\adj H_\lambda =\mc L_{\lambda(\frob)}^u(BH_\lambda) \hookleftarrow \wh{\mc L}_{\lambda(\frob)}(BH_\lambda). 
\]
Next observe that, since $\mc U_{H\times \Gm}=\mc U_H$, 
\begin{align*}
    \mc{L}^u_{(\lambda(\frob),q)}(B(H_{\lambda} \times \Gm)) & \simeq \wh{ \mc U}_{(\lambda(\frob), q)}/^\adj (H_{\lambda} \times \Gm) \simeq  (\wh{\mc U}_{\lambda(\frob)} /^\adj H_{\lambda}) \times (q\wh{\bb{G}}_\mathrm{m}/\Gm)  \\
     &\simeq  \mc{L}^u_{\lambda(\frob)}(BH_\lambda) \times q\wh{\bb G}_{\mathrm{m}}/ \Gm,
\end{align*}
by \cite[Definition 3.1.1, Propositions 2.1.25, 2.1.8]{Chenequivariant}, and where $\wh{ \mc U}_z$ is as in \cite[Definition 3.1.1]{Chenequivariant}. Similarly, 
\[
\wh{\mc{L}}_{(\lambda(\frob),q)}(B(H_{\lambda} \times \Gm))\simeq \wh{\mc{L}}_{\lambda(\frob)}(BH_\lambda)\times q\wh{\bb G}_{\mathrm{m}}/ \Gm,
\]
by \cite[Definition 3.1.1, Example 2.1.21, Proposition 2.1.8]{Chenequivariant}. Therefore, there is a map
\begin{align*}
\LS_{H, \lambda}^{\WD_F, \up, \wedge}
&\rightarrow
\mc L_{(\lambda(\frob), q)}^u(B(H\times \Gm))\times_{\mc L(B(H\times \Gm))} \mc L (\Lie H/ (H\times \Gm)) \\
&\simeq \mc L^u_{(\lambda (\frob), q)} (\Lie H / (H\times \Gm)),      
\end{align*}
identifying the source with the base change along $\{q\}\to q\wh{\bb G}_{\mathrm{m}}/\Gm$. 
If $H$ is reductive, using Theorem \ref{thm: LoopSpaceLocalization}, we deduce that
\[
\LS_{H, \lambda}^{\WD_F, \up, \wedge}
\simeq 
\mc L_{(\lambda(\frob), q)}^u (V_\lambda/ (H_\lambda\times \Gm))\times_{q\wh{\bb{G}}_{\mathrm{m}}/\Gm}\{q\}.
\]
Even if $H$ is not reductive, we still have a map
\[
\mc L_{(\lambda(\frob), q)}^u (V_\lambda/ (H_\lambda\times \Gm)) \to \mc L^u_{(\lambda (\frob), q)} (\Lie H / (H\times \Gm)),  
\]
and (\ref{eqn: chenfiberprod}) and Lemma \ref{lemma: semisimpleconjugation} imply that it is an equivalence on $\Lambda$-points. In fact, it is an equivalence of formal stacks by comparing (co)tangent complexes: at a point $(N, u)$, the tangent complexes are given respectively by $\Gamma (\Z,-)$ with respect to $(\lambda (\frob)u, q)$ applied to
\[
\Lie H_\lambda \times \Lie \Gm \xrightarrow{(h,c)\mapsto [h, N] -c N} V_\lambda
\]
and 
\[
\Lie H \times \Lie \Gm \xrightarrow{(h,c)\mapsto [h, N] -c N} \Lie H. 
\]
But we can take $(\lambda(\frob), q)$-invariants first, and then the latter identifies with the former.  
We similarly have
\[
\LS_{H, \V_\lambda}^{\WD_F, \up, \wedge} \simeq 
\wh{\mc L}_{(\lambda(\frob), q)} (V_\lambda/ (H_\lambda\times \Gm))\times_{q\wh{\bb{G}}_{\mathrm{m}}/\Gm}\{q\}.
\]

On the other hand, there is a canonical identification
\[
V_{\lambda}/ (H_{\lambda}\times \Gm)\simeq V_{\lambda}/ ((H_{\lambda}\times \Gm)/A)\times_{B((H_{\lambda}\times \Gm)/A)}B(H_{\lambda}\times \Gm), 
\]
giving a $z=(\lambda (\frob), q)$-trivialization, where $A=\langle z\rangle$, slightly generalizing a notion in \cite[Definition A.16]{BCHNconstructible}. Compare with \cite[Proposition A.17]{BCHNconstructible} together with \cite[Proposition 2.1.22]{Chenequivariant}. Now, since $z$ is central in $H_{\lambda}\times \Gm$ and we have the canonical $z$-trivialization above, we have the \emph{shift by $z$} maps as in the discussion before Corollary A.18 of \cite{BCHNconstructible} which give
\begin{equation*}
   \mc{L}^u_z(V_{\lambda} /( H_{\lambda} \times \Gm)) \simeq  \mc{L}^u(V_{\lambda} / (H_{\lambda} \times \Gm)),  \quad \wh{\mc{L}}_z(V_{\lambda} / (H_{\lambda} \times \Gm)) \simeq \wh{\mc{L}}(V_{\lambda} / (H_{\lambda} \times \Gm)).
\end{equation*}
Since loop spaces commute with fiber products, we have 
\[
\mc{L}^u(V_{\lambda} / (H_{\lambda} \times \Gm))\times_{\mc L^u(B\Gm)}\{e\}\simeq \mc{L}^u(\V_{\lambda}), \quad
\wh{\mc{L}}(V_{\lambda} / (H_{\lambda} \times \Gm))\times_{\wh{\mc L}(B\Gm)}\{e\}\simeq \wh{\mc{L}}(\V_{\lambda}),
\]
and this implies the result for $\LS^{\WD_F, \wedge}_{H, \lambda}$. For $\LS^{\WD_F, \wedge}_{H, \V_{\lambda}}$ we further use Theorem \ref{thm:exp} giving a functorial isomorphism 
\[
\exp\colon \wh{\bb{T}}\V_{\lambda}[-1] \simeq \wh{\mc{L}}(\V_{\lambda}).
\]
\end{proof}

\begin{proof}[Proof of Proposition \ref{prop: LoopSpaceIdentificationCompletions}]
By the previous lemma, it suffices to prove that
\[
\LS_{H, \lambda}^{\WD_F, \wedge} \simeq \LS_{J_\lambda, s_\lambda}^{\WD_F, \wedge}, \quad
\LS_{H, \V_\lambda}^{\WD_F, \wedge} \simeq \LS_{J_\lambda, \V_\lambda}^{\WD_F, \wedge}. 
\]
By Proposition \ref{prop: LSDescriptions}, we have
\[
\Sect(\WD_F, K_\lambda)^{\wedge}_{\lambda_I, [\lambda]} \simeq (\Sect (W_F, K_\lambda)_{\lambda_I, [\lambda]}^{\wedge}\times \Lie K_\lambda)\times_{\Lie K_\lambda}\{0\}. 
\]
Combined with Proposition \ref{prop: ReductionToUnipotentWeilParameters}, we have a functorial identification
\[
\Sect (\WD_F, K_\lambda)^{\wedge}_{\lambda_I, [\lambda]} \simeq ((s_\lambda \wh{\mc U}_{H_\lambda}t_\lambda\frob \times^{H_\lambda} K_\lambda) \times \Lie K_\lambda )\times_{\Lie K_\lambda } \{0\}. 
\]
There is a natural map from 
\[
((s_\lambda \wh{\mc U}_{H_\lambda}\times^{H_\lambda} J_\lambda)t_\lambda\frob \times \Lie J_\lambda )\times_{\Lie J_\lambda } \{0\}, 
\]
inducing
\[
\LS^{\WD_F, \wedge}_{J_\lambda, s_\lambda}\simeq\Sect (\WD_F, J_\lambda)^{\wedge}_{[s_\lambda]}/ J_\lambda \to
\LS^{\WD_F, \wedge}_{H, \lambda}\simeq\Sect(\WD_F, K_\lambda)^{\wedge}_{\lambda_I, [\lambda]}/K_\lambda. 
\]
We shall show that this map is an equivalence. At the level of $\Lambda$-points, this is clear since $(\Lie K_\lambda)^{(t_\lambda,\frob)}=\Lie J_\lambda$. 
It suffices to show that it induces an isomorphism of cotangent complexes at all $\Lambda$-points $(s_\lambda u, N)$ of $\Sect (\WD_F, J_\lambda)_{[s_\lambda]}$ with $u\in \mc{U}_{H_\lambda}$, $N\in \Lie J_\lambda$. Up to shifts, it is a natural map
\begin{multline*}
\Gamma (\WD_F, (\Lie J_{\lambda})_{(s_\lambda u, N)}^{\vee}(1))\simeq
\Gamma (\Z\ltimes \Ga, (\Lie J_{\lambda})_{(s_\lambda u, N)}^{\vee}(1))
\\
\leftarrow\Gamma (\Z\ltimes \Ga, (\Lie K_\lambda)_{(u\lambda(\frob), N)}^{\vee}(1)),     
\end{multline*}
and, by looking at the Hochschild--Serre spectral sequence, it suffices to show that
\begin{align*}
(\Lie J_{\lambda})_{(s_\lambda u, N)}^{\vee, N}(1)&\leftarrow
(\Lie K_{\lambda})_{(u\lambda (\frob), N)}^{\vee, N}(1), \\
(\Lie J_{\lambda})_{(s_\lambda u, N), N}^{\vee}&\leftarrow
(\Lie K_{\lambda})_{(u\lambda (\frob), N), N}^{\vee}    
\end{align*}
induce isomorphisms of the group cohomology with respect to $\Z$. As both $s_\lambda u$ and $N$ commute with $(t_\lambda,\frob)$ and $(t_\lambda,\frob)$ is semisimple, we may decompose $\Lie J_\lambda$, $\Lie K_\lambda$ into eigenspaces. Moreover, since an eigenvalue of $s_\lambda u$ that equals $q$ (resp.~$1$) modulo torsion must equal $q$ (resp.~$1$) itself, and the eigenvalues of $(t_\lambda,\frob)$ are torsion, we need only look at the summand $(\Lie K_\lambda)^{(t_\lambda, \frob)}$, which is exactly $\Lie J_\lambda$. 

The claim for the formal completion along the Vogan stack follows immediately. 
\end{proof}

\section{Eisenstein series and constant terms}{\label{s: EisandCT}}

We will introduce Eisenstein series in various settings, and relate them to each other. 

\subsection{Weight decompositions}{\label{ss: weightdecomposition}}
Let $\Delta$ be a diagonalizable algebraic group over $\Lambda$. If $\pi\colon \mc X\to S$ is a gerbe banded by $\Delta$ over an Artin stack $S$ almost of finite type, there is a decomposition
\[
\QCoh(\mc X) \simeq \prod_{\chi\in X^*(\Delta)}\QCoh(\mc X)^\chi,  \quad  \QCoh(S)\xrightarrow[\pi^*]{\simeq}\QCoh(\mc X)^0, 
\]
e.g.~\cite[Theorem 5.4, Proposition 5.7]{BerghSchnurer} if $S$ is classical, and \cite[Theorem 5.25]{BindaPorta}, \cite[Lemma 1.5.3]{HalpernLeistner} if $\Delta=\Gm$.
It is a $\Delta$-trivialization in the sense of \cite[Definition A.19]{BCHNconstructible}. In fact, we are only interested in the following situation:

\begin{lemma}
Suppose $\Delta\subset G$ is a central diagonalizable subgroup of an affine algebraic group $G$ and write $H\coloneq G/\Delta$. 
Suppose $G$ acts on a derived scheme $X$ almost of finite type. If $\Delta$ acts trivially on $X$ with an equivalence $X/\Delta\simeq X\times B\Delta$, then $X$ admits a $\Delta$-trivialization in the sense of \cite[Definition 4.1.1]{Chenequivariant}, \cite[Definition A.16]{BCHNconstructible}, i.e., there is an identification
\[
X/G \simeq X/H\times_{BH}BG. 
\]
In particular, 
\begin{align*}
\QCoh(X/ G)&\simeq \QCoh(X/H)\otimes_{\QCoh(BH)}\QCoh(B G) \\
&\simeq \prod_{\chi\in X^*(\Delta)}\QCoh(X/H)\otimes_{\QCoh(BH)} \QCoh(B G)^{\chi},     
\end{align*}
where $\QCoh(BG)^\chi$ is generated by representations of $G$ on which $\Delta$ acts through $\chi$. 
\end{lemma}

\begin{proof}
Under the assumption, the homotopy fixed points $X^\Delta$ are equivalent to $X$, so this is a special case of \cite[Proposition A.17]{BCHNconstructible}; the last part follows from Proposition \ref{prop:QCoh_satisfies_Kunneth}.  
\end{proof}

\begin{remark}
If $X$ is classical, the $\Delta$-trivialization $X/\Delta\simeq X\times B\Delta$ is canonical. In particular, if $X$ has the form $G\times^{G'}X'$ for some subgroup $G'\subset G$ containing $\Delta$, then the $\Delta$-trivialization is induced from that of $X'$ with its $G'$-action. 

If $X$ is not classical, but a closed subscheme of a classical $G$-scheme $X'$ on which $\Delta$ acts trivially, then the $\Delta$-trivialization of $\QCoh(X'/ G)$ induces that of $\QCoh(X/G)$ if $\mc O_{X/G}\in \QCoh(X'/G)^0$, since it is the category of $\mathcal{O}_{X/ G}$-modules in $\QCoh(X'/ G)$. If moreover $X$ identifies $ G$-equivariantly with the fiber product $X'\times_{X_1}X_2$ with classical $G$-schemes $X_1$, $X_2$ on which $\Delta$ acts trivially, this $\Delta$-trivialization agrees with the one induced from $X/\Delta\simeq X\times B\Delta$, which is in turn induced from $X\simeq X'\times_{X_1}X_2$. As in the previous paragraph, we can also pass to $G'\subset  G$ if all objects come from classical $G'$-schemes. 
\end{remark}

\begin{remark}
If $\Delta'\subset\Delta$ is a subgroup scheme, then a $\Delta'$-trivialization is also induced and it is compatible with the $\Delta$-trivialization in the obvious sense.  
\end{remark}

The above decomposition induces
\[
\Coh(X/ G) \simeq \bigoplus_{\chi\in X^*(\Delta)}\Coh(X/ G)^\chi, 
\]
and in turn
\[
\IndCoh(X/ G) \simeq \prod_{\chi\in X^*(\Delta)}\IndCoh(X/G)^\chi. 
\]
These decompositions are $\Delta$-trivializations in the sense of \cite[Definition A.19]{BCHNconstructible}. 

\begin{definition}
If a $\Delta$-trivialization of a $\QCoh(B G)$-linear stable $\infty$-category $\mc C$ is given, we write, for $\chi \in X^*(\Delta)$, 
\[
\mc C^\chi \coloneq\mc C^0 \otimes_{\QCoh(BH)}\QCoh( BG)^{\chi},\quad i^\chi\colon \mc C^\chi \to \mc C, \quad \pi^\chi\colon \mc C \to \mc C^\chi,
\]
which satisfy $\id \simeq \pi^\chi\circ i^\chi$. Note that $\pi^\chi$ is both left and right adjoint to $i^\chi$, and we shall identify $\mc C^\chi$ with a full subcategory of $\mc C$. 
\end{definition}

There are obvious functorialities:

\begin{lemma}
Let $ G\to G'$ be a map of affine groups of finite type such that it is injective on $\Delta$ and the image of $\Delta$ is central in $G'$. Let $H'$ denote $G'$ modulo the image of $\Delta$. If $f\colon X/ G\to Y/G'$ is a map and it lifts to a map $\widetilde{f}\colon X/H\to Y/H'$ together with compatible $\Delta$-trivializations, then $f^*$ (resp.~$f_*$) preserves full subcategories corresponding to $\chi\in X^*(\Delta)$ (resp.~if any base change of $X/G\to Y/G'$ to an affine scheme is QCA). Similar claims hold for $f_*^{\IndCoh}$, $f^!$ (as does $f^{*,\IndCoh}$ if $f$ is of finite Tor amplitude). 
\end{lemma}

\begin{proof}
Recall that $\QCoh(Y/G')\simeq \QCoh(Y/H')\otimes_{\QCoh(BH')}\QCoh(BG')$. 
Since $f^*, f_*$ are both $\QCoh(BG')$-linear under the assumption, the claim reduces to the case of trivial blocks, and it is clear for $f^*$ and follows from base change for $f_*$. 

The claim for $f_*^{\IndCoh}$ and $f^{*,\IndCoh}$ follows from the case of $\Coh\subset \QCoh$ as both functors are continuous. The case of $f^!$ follows from the $\QCoh(BG')$-linearity and base change as in the previous paragraph, since
\[
\IndCoh(Y/G')\simeq \IndCoh(Y/H')\otimes_{\QCoh(BH')}\QCoh(BG'). 
\]
\end{proof}

The weight decomposition extends, with similar functorial properties, to ind-algebraic stacks if $X$ is an ind-scheme almost of finite type and a $\Delta$-trivialization of $X/G$ is given. 
Similarly, if $X$ is a formal scheme, it further induces
\[
\wh{\Coh}{}^!(X/ G) \simeq \bigoplus_{\chi\in X^*(\Delta)}\wh{\Coh}{}^!(X/ G)^\chi 
\]
and a corresponding decomposition after ind-completions. 

\begin{example}
We use the notation in the main setup. Let $\wh M \subset \wh G$ be a semistandard Levi subgroup. There are natural $\Delta=Z(\wh M)^\Gamma$-trivializations in the following cases. 
\begin{enumerate}
    \item We may consider the decomposition of $\IndCoh(\Par_M)$ along $X^*(Z(\wh M)^\Gamma)$. Since $\omega_{\Par_M}\simeq \mc O_{\Par_M}$, Grothendieck--Serre duality exchanges $\Coh(\Par_M)^{\chi}$ and $\Coh(\Par_M)^{-\chi}$. 
    \item For any semisimple parameter $\lambda_M$, we may consider the decomposition of 
    \[ 
    \wh{\Coh}{}^!(\mc L^u(\V_{\lambda_M})), \quad \wh{\Coh}{}^!(\wh{\mc L}(\V_{\lambda_M})), \quad \wh{\Coh}{}^!(\widehat{\T}\V_{\lambda_M}[-1]), \quad \Coh(\T^* \V_{\lambda_M})
    \]
    along $X^*(Z(\wh M)^\Gamma)$; the first three decompositions preserve the subcategory $\Coh\subset \wh{\Coh}{}^!$, and all admit graded versions. The last one is induced from the $\Delta$-trivialization of 
    \[
    \Coh(\T^* V_{\lambda_M} /H_{\lambda_M})\subset \QCoh(\T^* V_{\lambda_M} /H_{\lambda_M})
    \]
    via the embedding $\T^*\V_{\lambda_M} \hookrightarrow (\T^*V_{\lambda_M})/H_{\lambda_M}$, or the interpretation of $\T^* \V_{\lambda_M}$ as $\mu^{-1}(0)/H_{\lambda_M}$ using the moment map $\mu$.  
\end{enumerate}
\end{example}

\begin{lemma}
Under the linear Koszul dualities
\begin{align*}
\IndCoh(\widehat{\T}\V_{\lambda_M}[-1]/\Gm) &\simeq \QCoh(\T^* \V_{\lambda_M}/\Gm), \\
\Ind\wh{\Coh}{}^!(\widehat{\T}\V_{\lambda_M}[-1]/\Gm) &\simeq \IndCoh(\T^* \V_{\lambda_M}/\Gm),     
\end{align*}
the full subcategories for $\chi\in X^*(\Delta)$ agree. 
\end{lemma}

\begin{proof}
It suffices to show the claim for the first equivalence, after forgetting the grading. As noted above, the $Z(\wh M)^\Gamma$-trivialization on $\QCoh(\T^* \V_{\lambda_M})$ is induced from that on $\QCoh(\T^* V_{\lambda_M}/H_{\lambda_M})$. In fact, it is induced from the $Z(\wh M)^\Gamma$-trivialization on $\QCoh(\V_{\lambda_M})$ via $*$-pushforward. By construction of the linear Koszul duality, this $*$-pushforward corresponds to 
\[
z^!\colon \IndCoh(\wh{\T}\V_{\lambda_M}[-1])\to \IndCoh(\V_{\lambda_M})=\QCoh(\V_{\lambda_M}), 
\]
and this functor preserves $Z(\wh M)^\Gamma$-trivializations. As $z^!$ is conservative (cf.~the proof of \cite[Theorem 2.1.19]{Chen2023categorical}), the functor $z^!$ determines the $Z(\wh M)^\Gamma$-decomposition of $\IndCoh(\wh{\T}\V_{\lambda_M}[-1])$. 
\end{proof}

\begin{remark}
As in the proof above, the weight decomposition of $\IndCoh(\widehat{\mathcal{L}}(\V_{\lambda_M}))$ is also determined by $z^!$ and the weight decomposition of $\QCoh(\V_{\lambda_M})$, and this implies that the identification $\exp$ in Theorem \ref{thm:exp} preserves the weight decompositions in this situation. 
\end{remark}

\subsection{Spectral Eisenstein series and constant terms}{\label{ss: spectraleisensteinconstantterm}}

Now let $G$ be a connected quasi-split reductive group and take semistandard Levi and parabolic subgroups $M\subset P \subset G$. Then the natural maps of $W_F$-stable groups $\wh{M}\leftarrow \wh{P}\rightarrow \wh{G}$ induce a diagram
\begin{equation*}
    \begin{tikzcd}
        \Par_M & \Par_P \arrow[l, swap,  "q"] \arrow[r, "p"] & \Par_G. 
    \end{tikzcd}
\end{equation*}

\begin{proposition}{\label{prop: pandqproperties}}
The map $p$ is proper and schematic, and $q$ is quasi-smooth with trivial relative dualizing complex in degree $0$. In particular, $q^! \simeq q^{\ast, \IndCoh}$ on $\IndCoh(\Par_M)$.      
\end{proposition}

\begin{proof}
The first claim is clear from the corresponding claim for $B(\wh{P}\rtimes \WD_F)\to B(\wh{G}\rtimes \WD_F)$. The quasi-smoothness of $q$ is proven in \cite[Lemma 3.21]{Zhu}, \cite[Lemma 3.5.3]{HansenMann}, by analyzing the cotangent complexes. Since the dualizing complexes of $\Par_M$, $\Par_P$ are trivial by Proposition \ref{prop: properties_of_Sect} (4), the relative one for $q$ is trivial as well. 
\end{proof}

As before, we have the following decomposition of $\IndCoh(\Par_M)$ as the product of full subcategories
\begin{equation*}
    \IndCoh(\Par_M) \simeq \prod\limits_{\chi \in X^*(Z(\wh{M})^{\Gamma})} \IndCoh(\Par_M)^{\chi}
\end{equation*}
with $i^\chi$ and $\pi^\chi$. 

\begin{definition}{\label{def: CTEis}}
Let $G$ be as in the previous paragraph. We define the \emph{spectral Eisenstein series} functor by
\begin{equation*}
    \Eis^{\spec}_P\colon \IndCoh(\Par_M) \rightarrow \IndCoh(\Par_G) \qquad \mc{F} \mapsto p_{\ast}q^{\ast, \IndCoh} \mc{F}\simeq p_* q^! \mc{F},
\end{equation*}
and the \emph{spectral constant term} functor
\begin{equation*}
    \CT^{\spec}_P\colon \IndCoh(\Par_G)\rightarrow \IndCoh(\Par_M) \qquad \mc{F} \mapsto q_{\ast}p^!\mc{F}.
\end{equation*}
We also have, for each $\chi \in X^*(Z(\wh{M})^{\Gamma})$, the compositions 
\begin{equation*}
    \CT^{\spec, \chi}_P \coloneq \pi^{\chi} \circ \CT^{\spec}_P \qquad \Eis^{\spec, \chi}_P \coloneq \Eis^{\spec}_P \circ i^{\chi}.
\end{equation*}
\end{definition}

\begin{remark}{\label{rem: CTEisproperties}}
\begin{enumerate}
    \item 
    We have that $\Eis^{\spec}_P$ is left adjoint to $\CT^{\spec}_P$ and $\Eis^{\spec, \chi}_P$ is left adjoint to $\CT^{\spec, \chi}_P$, and all four functors commute with colimits. 
    \item The functors $\Eis_P^\spec$ and $\CT^{\spec, \chi}_P$ preserve compact objects, i.e., coherent sheaves \cite[Corollary 3.5.4(ii), Theorem 3.7.1(i)]{HansenMann}.\footnote{The latter is much deeper than the former, but we will only use the former.} 
\item By Proposition \ref{prop: pandqproperties}, we have $\bb{D}_{\GS} \circ \Eis_P^\spec \simeq \Eis_P^\spec \circ \bb{D}_{\GS}$. 
\item Let $n \in N_{\wh{G}}(\wh T^{\Gamma,\circ})$ and consider the semistandard parabolic $\co{n}{\wh P} \coloneq \Int(n)(\wh P)$, which corresponds to the semistandard parabolic subgroup $\co{w}{P}$ of $G$ by \cite[Proposition 2.5.4(ii)]{SilbergerZink}, if $n$ maps to $w \in N_{\wh{G}}(\wh T^{\Gamma,\circ})/Z_{\wh{G}}(\wh T^{\Gamma,\circ}) \cong  W^\Gamma_{\wh{G}}$. Then we have a commutative diagram
    \begin{equation*}
        \begin{tikzcd}
            \Par_M \arrow[d, "n"] & \Par_P \arrow[l, swap,  "q"] \arrow[r, "p"] \arrow[d, "n"] & \Par_G \arrow[d, "n"] \\
            \Par_{\co{w}{M}}   & \Par_{\co{w}{P}}  \arrow[l, swap, "q"] \arrow[r, "p"] & \Par_G,
        \end{tikzcd}
    \end{equation*}
    where the vertical maps are isomorphisms induced by conjugation. In fact, the right vertical map induces a functor isomorphic to the identity on $\IndCoh(\Par_G)$ and we have
    \begin{equation*}
        \Int(n)^{\ast}\colon \IndCoh(\Par_{\co{w}{M}})^{w(\chi)} \rightarrow \IndCoh(\Par_M)^{\chi}.
    \end{equation*}
    Thus, $\Eis^{\spec, w(\chi)}_{\co{w}{P}} \simeq 
 \Eis^{\spec, \chi}_{P} \circ \Int(n)^{\ast}$. 
 
\end{enumerate}
\end{remark}

\begin{lemma}{\label{lem: eisctcompositions}}
Suppose that $P_1 \subset P_2 \subset G$ are semistandard parabolic subgroups with semistandard Levi subgroups $M_1\subset M_2$ respectively. Let $P'_1 = P_1 \cap M_2$ denote the semistandard parabolic subgroup of $M_2$. Then 
\begin{equation*}
    \Eis^{\spec}_{P_1} \simeq \Eis^{\spec}_{P_2} \circ \Eis^{\spec}_{P'_1}, \qquad \CT^{\spec}_{P_1} \simeq \CT^{\spec}_{P'_1} \circ \CT^{\spec}_{P_2}.
\end{equation*}
Additionally, if $\chi_1 \mapsto \chi_2$ under the restriction map $X^{\ast}(Z(\wh{M}_1)^{\Gamma}) \rightarrow X^{\ast}(Z(\wh{M}_2)^{\Gamma})$, then $\Eis^{\spec, \chi_1}_{P'_1}$ lands in $\IndCoh(\Par_{M_2})^{\chi_2}$.
\end{lemma}

\begin{proof}
This follows easily from $\Par_{P_2}\times_{\Par_{M_2}}\Par_{P'_1}\simeq \Par_{P_1}$, which follows from 
\[
B(\wh{P}_2\rtimes \WD_F)\times_{B(\wh{M}_2\rtimes \WD_F)} B(\wh{P'}_1\rtimes \WD_F) \simeq B(\wh{P}_1\rtimes \WD_F). 
\]
\end{proof}

\subsection{Eisenstein series for completed fibers}
Recall that $\Lambda$ is algebraically closed. Fix a semisimple parameter $\lambda\in X_G(\Lambda)$. We will consider a version of the spectral Eisenstein series for the completed fibers $\Par_{G,\lambda}^{\wedge}$. 

The fiber of $X_P\to X_G$ at $\lambda$ consists of finitely many points $\lambda_P \in X_P (\Lambda)$ by Corollary \ref{cor: XP to GG is finite}. Recall that there is a natural identification $X_P \cong X_M$ in Lemma \ref{lem: X_P=X_M}, and we write $\lambda_M$ for the image of $\lambda_P$. 

Note that the composite $\Par_P \to \Par_G \to X_G$ factors through $X_P$. Similarly, the composite $\Par_P\to \Par_M \to X_M\cong X_P$ agrees with the natural map $\Par_P\to X_P$. Therefore, we can consider the diagram
\begin{equation*}
    \begin{tikzcd}
                \coprod_{\lambda_M} \V_{\lambda_M} \arrow[d, "\coprod i_{\lambda_M}"] & \coprod_{\lambda_{P}} \V_{\lambda_{P}} \arrow[l, swap, "q_{\V}"] \arrow[r, "p_\V"] \arrow[d, "\coprod i_{\lambda_P}"] & \V_{\lambda} \arrow[d, " i_{\lambda}"] \\ 
                \coprod_{\lambda_M} \Par^{\wedge}_{M, \lambda_M} & \coprod_{\lambda_{P}} \Par^{\wedge}_{P, \lambda_{P}} \arrow[l, swap, "\hat{q}"] \arrow[r, "\hat{p}"]  & \Par^{\wedge}_{G, \lambda}, 
    \end{tikzcd}
\end{equation*}
where the disjoint unions range over all $\lambda_M$ and $\lambda_{P}$ in the fiber at $\lambda$. Note that these stacks are ind-(QCA algebraic). We adopt the notation that $p_{\lambda_P}$ and $q_{\lambda_P}$ are the restrictions of $p_{\V}$ and $q_{\V}$ to $\V_{\lambda_P}$. 
\begin{lemma}{\label{lem: propersmooth}}
    The map $q_{\lambda_P}$ is smooth and $p_{\lambda_P}$ is proper and schematic.
\end{lemma}
\begin{proof}
    To prove that $q_{\lambda_P}$ is smooth, we write it as
    \begin{equation*}
        \V_{\lambda_P} \rightarrow V_{\lambda_M}/H_{\lambda_P} \rightarrow \V_{\lambda_M},
    \end{equation*}
which is a composition of a vector bundle map and the base change of $BH_{\lambda_P} \rightarrow BH_{\lambda_M}$, both of which are smooth.

To prove $p_{\lambda_P}$ is proper and schematic, we write it as
\begin{equation*}
    \V_{\lambda_P} \rightarrow V_{\lambda}/H_{\lambda_P} \rightarrow \V_{\lambda},
\end{equation*}
which is a composition of a closed immersion and a base change of the map $BH_{\lambda_P} \rightarrow BH_{\lambda}$, both of which are proper.
\end{proof}

We record the following useful lemma.
\begin{lemma}{\label{lem: cartesianEisdiagram}}
    We have the following commutative diagram where both squares are cartesian.

\begin{equation*}
    \begin{tikzcd}
        \coprod_{\lambda_M} \Par^{\wedge}_{M, \lambda_M} \arrow[d, " \coprod\wh{i}_{\lambda_M}"]& \coprod_{\lambda_{P}} \Par^{\wedge}_{P, \lambda_{P}} \arrow[l, swap, "\hat{q}"] \arrow[r, "\hat{p}"] \arrow[d, "\coprod \wh{i}_{\lambda_{P}}"] & \Par^{\wedge}_{G, \lambda} \arrow[d, "\wh{i}_{\lambda}"] \\
        \Par_M & \Par_{P} \arrow[l, swap, "q"] \arrow[r, "p"] & \Par_G.
    \end{tikzcd}
\end{equation*} 
In particular, the map $\hat{p}$ is schematic and proper, and $\hat{q}$ is quasi-smooth. 
\end{lemma}

\begin{proof}
The diagram clearly commutes. To check that both squares are cartesian, note that formal completion commutes with pullback (Remark \ref{rem: formalcompletionpullback}), hence we just need to check that (for the right square and analogously for the left square) 
\begin{equation*}
    \coprod_{\lambda_P} (\Par_{P, \lambda_P})_{\red} = p^{-1}( \Par_{G, \lambda})_{\red},
\end{equation*}
which is clear by definition.
\end{proof}

We define, on ind-coherent sheaves, $\Eis^{\spec, \wedge}_{P, \lambda}$ by $\hat{p}_{\ast}\hat{q}^!$ and $\CT^{\spec, \wedge}_{P, \lambda}$ by $\hat{q}_{\ast} \hat{p}^!$, and we can also define $\chi$-graded versions of these functors for $\chi \in X^*(Z(\wh{M})^{\Gamma})$. 

\begin{proposition}\label{prop: EisComparisonCompletedFiber}
There are commutative diagrams
\[
    \begin{tikzcd}
        \IndCoh(\Par_M) \arrow[d, "(\coprod\wh{i}_{\lambda_M})^!"] \arrow[r, "\Eis_P^\spec"] & \IndCoh(\Par_G) \arrow[d, "\wh{i}_\lambda^!"] \\    
        \IndCoh (\coprod_{\lambda_M} \Par^{\wedge}_{M, \lambda_M}) \arrow[r, "\Eis_{P, \lambda}^{\spec,\wedge}"]  & \IndCoh(\Par^{\wedge}_{G, \lambda}), 
    \end{tikzcd}
\]
\[
      \begin{tikzcd}
        \IndCoh(\Par_M)  \arrow[r, "\Eis_P^\spec"] & \IndCoh(\Par_G)  \\    
        \IndCoh (\coprod_{\lambda_M} \Par^{\wedge}_{M, \lambda_M}) \arrow[u, "(\coprod\wh{i}_{\lambda_M})_*"] \arrow[r, "\Eis_{P, \lambda}^{\spec,\wedge}"]  & \IndCoh(\Par^{\wedge}_{G, \lambda}) \arrow[u, "\wh{i}_{\lambda*}"],
    \end{tikzcd}
\]
and the same holds for the $\chi$-graded versions.    
\end{proposition}

\begin{proof}
The proposition follows from Lemma \ref{lem: cartesianEisdiagram}. The first commutative diagram uses base change on the right square, and the second uses base change on the left.
\end{proof}

The functor $\Eis^{\spec, \wedge}_{P,\lambda}$ preserves continuous ind-coherent sheaves $\wh{\Coh}{}^!$ (e.g.~one can use Theorems \ref{thm: completeness_along_fiber}, \ref{thm: completeness_along_fiber_parabolic} to reduce it to the non-formal setting), and by passing to the ind-completion, we obtain a functor defined on $\Ind\wh{\Coh}{}^!$. 
\begin{corollary}\label{cor: EisComparisonCompletedFiber}
There are commutative diagrams
\[
    \begin{tikzcd}
        \Coh(\Par_M) \arrow[d, "(\coprod\wh{i}_{\lambda_M})^!"] \arrow[r, "\Eis_P^\spec"] & \Coh(\Par_G) \arrow[d, "\wh{i}_\lambda^!"] \\    
        \wh{\Coh}{}^! (\coprod_{\lambda_M} \Par^{\wedge}_{M, \lambda_M}) \arrow[r, "\Eis_{P,\lambda}^{\spec,\wedge}"]  & \wh{\Coh}{}^!(\Par^{\wedge}_{G, \lambda}),
    \end{tikzcd} 
    \]
    \[
    \begin{tikzcd}
        \IndCoh(\Par_M) \arrow[r, "\Eis_P^\spec"] & \IndCoh(\Par_G) \\    
        \wh{\Coh}{}^! (\coprod_{\lambda_M} \Par^{\wedge}_{M, \lambda_M}) \arrow[u, "(\coprod\wh{i}_{\lambda_M})_*"]  \arrow[r, "\Eis_{P,\lambda}^{\spec,\wedge}"]  & \wh{\Coh}{}^!(\Par^{\wedge}_{G, \lambda}) \arrow[u, "\wh{i}_{\lambda*}"] ,
    \end{tikzcd}
    \]
and the same holds for their ind-completions and $\chi$-graded versions. 
\end{corollary}

For each $\lambda_P$, we write
$\Eis_{P, \lambda_P}^{\spec,\wedge}$ for the restriction of $\Eis^{\spec, \wedge}_{P,\lambda}$ to $\Par^{\wedge}_{M, \lambda_M}$. 

Let $\lambda_P, \lambda'_P$ be two points lifting $\lambda$. At the level of representatives, there exists $g\in \wh{G}(\Lambda)$ such that ${}^L\lambda'_P=g{}^L\lambda_P g^{-1}$. There exists $w\in W_{\wh M}\backslash W_{\wh G}/W_{\wh M}$ such that $g\in \wh{P}(\Lambda)w \wh P(\Lambda)$; one can take a minimal representative, for instance. It is clear that $\lambda_P=\lambda'_{P}$ in $X_P(\Lambda)$ if and only if $g\in Z_{\wh G}({}^L \lambda'_P)\wh P Z_{\wh G}({}^L \lambda_P)$ if and only if $g\in\wh P Z_{\wh G}({}^L \lambda_P) $. 

Fix a semisimple lift $\lambda_P\colon W_F\to \wh P$. By Lemma \ref{lem: X_P=X_M}, choose a $\wh P$-conjugate $\lambda_M\colon W_F\to \wh M$. We can further assume $S\coloneq (Z_{\wh M}({}^L \lambda_M)\cap \wh T^{\Gamma})^\circ (\supset Z(\wh M)^{\Gamma,\circ})$ is a maximal torus of $Z_{\wh M}({}^L \lambda_M)$ by Lemma \ref{lem: TJcontrol}. Note that it is also maximal in $Z_{\wh G}({}^L \lambda_M)$ as $Z_{\wh G}(S)\subset \wh M$. 

\begin{lemma}\label{lem: ParPParGFibers}
Any class $\lambda'_P $ in the fiber of $\lambda$ is the class of $w{}^L\lambda_M w^{-1}$ for some $w \in (W_{\wh M}\backslash W_{\wh G})^\Gamma\cong W_{\wh M}^\Gamma \backslash W_{\wh G}^\Gamma$. 
\end{lemma}

\begin{proof}
As above, any lift ${}^L \lambda'_P=g{}^L \lambda_P g^{-1}$ is $\wh P$-conjugate to some $\lambda'_M\colon W_F \to \wh M$ such that $S'\coloneq (Z_{\wh M}({}^L \lambda'_M)\cap \wh T^\Gamma)^\circ$ is a maximal torus of $Z_{\wh M}({}^L \lambda'_M)$. Then, by $Z_{\wh G}({}^L \lambda_M)$-conjugation, we may assume that $g$ conjugates $S$ to $S'$. Since $\wh T$ centralizes $S, S'\subset \wh T^\Gamma$, we can write $g=mn\in Z_{\wh G}(S')N_{\wh G}(\wh T)$ such that ${}^L \lambda'_M=mn{}^L \lambda_M n^{-1}m^{-1}$, $nSn^{-1}=S'$, and we may assume $m=e$. 

Write $w\in W_{\wh G}$ for the class of $n\in N_{\wh G}(\wh T)$. For any $\gamma \in W_F$, we have $\lambda_M (\gamma ) \in Z_{\wh M}(S), n\lambda_{M}(\gamma) \gamma(n)^{-1}\in Z_{\wh M}(S')$. Therefore, 
\[
\gamma (n)\in Z_{\wh M}(S') n \subset \wh M n,
\]
and this implies that 
\[
w \in (W_{\wh M}\backslash W_{\wh G})^\Gamma\cong W_{\wh M}^\Gamma \backslash W_{\wh G}^\Gamma
\]
and hence we can replace $w$ by an element of $W_{\wh G}^\Gamma$. 
\end{proof}

Fix a representative $\lambda_P=\lambda_M$ as above, and suppose $\lambda'_P=\lambda'_M$ is attached to $w\in (W_{\wh M}\backslash W_{\wh G})^\Gamma$. There exists a commutative diagram
\[
\begin{tikzcd}
    \V_{\wh M,\lambda'_M} \arrow[d, "w^{-1}", "\simeq"'] 
    & \V_{\wh P, \lambda'_P} \arrow[l, "q_{\lambda'_P}"'] \arrow[r, "p_{\lambda'_P}"] \arrow[d, "w^{-1}", "\simeq"'] 
    & \V_\lambda \arrow[d, "w^{-1}", "\simeq"'] \\
    \V_{{}^{w^{-1}}\wh M,\lambda_M} 
    & \V_{{}^{w^{-1}}\wh P,\lambda_P} \arrow[l, "q_{\lambda_P}"'] \arrow[r, "p_{\lambda_P}"]
    & \V_\lambda, 
\end{tikzcd}
\]
and $\Eis^{\spec, \wedge}_{P, \lambda'_P}$ identifies with $\Eis^{\spec, \wedge}_{{}^{w^{-1}} P, \lambda_P}$, while the $\chi$-part corresponds to the $\chi \circ \mathrm{Ad}(w)$-part. 

\subsection{Eisenstein series for formal loop spaces}
Recall that the completed fibers identify with the unipotent loop spaces of Vogan stacks by Proposition \ref{prop: LoopSpaceIdentificationCompletions}, and we have the following commutative diagram
\begin{equation*}
\begin{tikzcd}
\coprod_{\lambda_M}\wh{\mc{L}}(\V_{\lambda_M}) \arrow[d, "\coprod \wh{i}_{\lambda_M ,\mathrm{Fss}}"]  
&\coprod_{\lambda_P}\wh{\mc{L}}(\V_{\lambda_P}) \arrow[l, "\wh{\mc{L}}(q_{\V})"'] \arrow[d, "\coprod \wh{i}_{\lambda_P,\mathrm{Fss}}"] \arrow[r, "\wh{\mc{L}}(p_{\V})"]
&\wh{\mc{L}}(\V_\lambda) \arrow[d, "\wh{i}_{\lambda, \mathrm{Fss}}"]
\\ 
\coprod_{\lambda_M}\Par_{M,\lambda_M}^\wedge \arrow[d, "\coprod \wh{i}_{\lambda_M}"]
&\coprod_{\lambda_P}\Par_{P,\lambda_P}^\wedge \arrow[l, "\hat{q}"'] \arrow[r, "\hat{p}"] \arrow[d, "\coprod \wh{i}_{\lambda_P}"]
&\Par_{G,\lambda}^\wedge \arrow[d, "\wh{i}_{\lambda}"]\\
\Par_M & \Par_P \arrow[r, "p"] \arrow[l, swap, "q"] & \Par_G.
\end{tikzcd}
\end{equation*}
These are all ind-(QCA algebraic) stacks. 
The upper right square is cartesian by the $2$ out of $3$ lemma for pullback squares since we have already shown the bottom right square is cartesian and the right rectangle is cartesian by the same logic. Hence $\wh{\mc{L}}(p_\V)$ is proper and schematic. Beware that the upper left square may not be cartesian. 

Next, we consider the Eisenstein series for the formal loop spaces defined as
\[
\Eis_{P}^{\wh{\mc{L}}} \coloneq \wh{\mc{L}}(p_\V)_* \circ \wh{\mc{L}}(q_\V)^!, \quad
\Eis_{P,\lambda_P}^{\wh{\mc{L}}} \coloneq \wh{\mc{L}}(p_{\lambda_P})_* \circ \wh{\mc{L}}(q_{\lambda_P})^!,
\]
which are defined on $\IndCoh(\coprod_{\lambda_M}\wh{\mc{L}}(\V_{\lambda_M}))$ and $\IndCoh(\wh{\mc{L}}(\V_{\lambda_M}))$ respectively. One can also consider $\chi$-versions.

\begin{proposition}{\label{prop: unipvsformal}}
The functors
\begin{gather}\label{eq: unip vs formal}
\wh{i}_{\lambda, \mathrm{Fss}}^{!}\colon \wh{\Coh}{}^! (\Par_{G,\lambda}^\wedge) \xrightarrow{}
\wh{\Coh}{}^{!} (\wh{\mc{L}}(\V_\lambda)), \\ 
\wh{i}_{\lambda_P, \mathrm{Fss}}^{!}\colon\wh{\Coh}{}^! (\Par_{P,\lambda_P}^\wedge) \xrightarrow{}
\wh{\Coh}{}^! (\wh{\mc{L}}(\V_{\lambda_P})), \quad
\wh{i}_{\lambda_M,\mathrm{Fss}}^{!}\colon\wh{\Coh}{}^! (\Par_{M,\lambda_M}^\wedge) \xrightarrow{}
\wh{\Coh}{}^! (\wh{\mc{L}}(\V_{\lambda_M}))
\notag
\end{gather}
are all equivalences. 
\end{proposition}

\begin{proof}
We explain the proof of the first equivalence, the others being analogous (using Theorem \ref{thm: completeness_along_fiber_parabolic} in the parabolic case). We apply Theorem \ref{thm: completeness_along_fiber} with $X$ being the base change of $\Sect(\WD_F, \wh{G})$ along the map $\Spec \mc O_{X_G,\lambda}^\wedge\to X_G$. Then Corollary \ref{cor: cohhatcompletion} implies we have an equivalence
\begin{equation*}
    \Coh(X/\wh{G}) \simeq \wh{\Coh}{}^!(\Par^{\wedge}_{G, \lambda}).
\end{equation*}
Also, by Theorem \ref{thm: completeness_along_closedpoint} and the same corollary, we have an equivalence
\begin{equation*}
      \Coh(X / \wh{G}) \simeq \wh{\Coh}{}^!(\wh{\mc{L}}(\V_{\lambda})),
\end{equation*}
induced by $\wh{i}^!_{\lambda, \mathrm{Fss}}$. The composition of these is the desired equivalence.

\end{proof}

\begin{lemma}
Under the equivalences \eqref{eq: unip vs formal}, $\hat{q}^!$ identifies with $\wh{\mc{L}}(q_\V)^!$ and $\hat{p}_*$ identifies with $\wh{\mc{L}}(p_\V)_*$. After the ind-completion, the right adjoint $\hat{p}^!$ of $\hat{p}_*$ identifies with the right adjoint $\wh{\mc{L}}(p_\V)^!$ of $\wh{\mc{L}}(p_\V)_*$.
\end{lemma}

\begin{proof}
The claim for $\hat{q}^!$ is clear. The claim for $*$-pushforward follows from base change, and the claim for $\hat{p}^!$ follows by passing to the right adjoints.
\end{proof}

In particular, $\Eis_{P}^{\wh{\mc{L}}}$ preserves continuous ind-coherent sheaves. This also follows from \cite[Theorem 2.3.6]{Chen2023categorical} by using the identification with the completed odd tangent stacks. 

\begin{corollary}{\label{cor: unipvsformalEis}}
There is a commutative diagram
\[
\begin{tikzcd}
\wh{\Coh}{}^! (\coprod_{\lambda_M} \Par^\wedge_{M,\lambda_M}) \arrow[r, "\Eis_{P,\lambda}^{\spec, \wedge}"] \arrow[d, "\simeq"] 
& \wh{\Coh}{}^! (\Par^\wedge_{G, \lambda}) \arrow[d, "\simeq"] \\
\wh{\Coh}{}^! (\coprod_{\lambda_M }\wh{\mc{L}}(\V_{\lambda_M})) \arrow[r, "\Eis_{P}^{\wh{\mc{L}}}"]
& \wh{\Coh}{}^! (\wh{\mc{L}}(\V_\lambda)). 
\end{tikzcd}
\] 
The same claim holds for the $\lambda_P$-components and the $\chi$-versions. 
\end{corollary}

Recall that we have functorial identifications between the formal loop spaces and the completed odd tangent stacks in Proposition \ref{prop: LoopSpaceIdentificationCompletions}. Namely, we have the following diagram
\[
\begin{tikzcd}
\coprod_{\lambda_M} \wh{\T}\V_{\lambda_M}[-1]
& \coprod_{\lambda_P} \wh{\T}\V_{\lambda_P}[-1]
\arrow[l, "{\wh{\T}(q_\V)[-1]}"']
\arrow[r, "{\wh{\T}(p_\V)[-1]}"]
& \wh{\T}\V_{\lambda}[-1], 
\end{tikzcd}
\]
and define
\[
\Eis_{P}^{\wh{\T}[-1]}\coloneq \wh{\T}(p_\V)[-1]_*  \circ \wh{\T}(q_\V)[-1]^!, \quad
\Eis_{P, \lambda_P}^{\wh{\T}[-1]}\coloneq \wh{\T}(p_{\lambda_P})[-1]_*  \circ \wh{\T}(q_{\lambda_P})[-1]^!.
\]
As before, they preserve continuous ind-coherent sheaves, and we clearly have the following.

\begin{proposition}{\label{prop: loopvstangentEis}}
On continuous ind-coherent sheaves, $\Eis_{P}^{\wh{\T}[-1]}$ identifies with $\Eis_{P}^{ \wh{\mc{L}}}$, and the same holds for the $\lambda_P$-components and the $\chi$-versions.
\end{proposition}

There is also a $\Gm$-equivariant version. Since both maps are $\Gm$-equivariant, $\Eis^{\wh{\T}[-1]}_{P}$, $\Eis_{P, \lambda_P}^{\wh{\T}[-1]}$ lift to this setting so that they are compatible with forgetting $\Gm$-equivariant structures.

\subsection{Eisenstein series for cotangent stacks}{\label{ss: eisensteinseriesforcotangents}}
We shall pass to the Koszul dual side. Define 
\[
\Eis_P^{\T^*}\coloneq (\breve{d}p_\V)_*  \circ (\breve{d}q_\V)^{!}, \quad
\Eis_{P,\lambda_P}^{\T^*}\coloneq (\breve{d}p_{\lambda_P})_*  \circ (\breve{d}q_{\lambda_P})^{!}, 
\]
where $(\breve{d}p_\V)_*$, $(\breve{d}q_\V)^{!}$, $(\breve{d}p_{\lambda_P})_*$, and $(\breve{d}q_{\lambda_P})^{!}$ are the proper $*$-pushforward and the smooth $!$-pullback on $\IndCoh(\T^* (-))$ preserving $\Coh(\T^*(-))$, defined in \cite[\S 2.3]{Chen2023categorical} by using correspondences among the cotangent stacks. In particular, in Chen's notation, we have
\begin{equation*}
     (\breve{d}p_{\lambda_P})_* \coloneq \breve{\pi}_{p_{\lambda_P} \ast} \circ \breve{\delta}^{\ast}_{p_{\lambda_P}} \qquad (\breve{d}q_{\lambda_P})^{!} \coloneq \breve{\delta}_{q_{\lambda_P} \ast} \circ \breve{\pi}^!_{q_{\lambda_P}}, 
\end{equation*}
for correspondences 
\begin{equation}{\label{eqn: eiscorrespondences}}
    \begin{tikzcd}
        & \T^{\ast}\V_{\lambda_M}|_{\V_{\lambda_P}} \arrow[dl, swap, "\breve{\pi}_{q_{\lambda_P}}"] \arrow[dr, " \breve{\delta}_{q_{\lambda_P}}"]& \\
       \T^{\ast}\V_{\lambda_M} &  & \T^{\ast}\V_{\lambda_P} ,
    \end{tikzcd} \qquad 
      \begin{tikzcd}
        & \T^{\ast}\V_\lambda |_{\V_{\lambda_P}} \arrow[dr, "\breve{\pi}_{p_{\lambda_P}}"] \arrow[dl, swap, " \breve{\delta}_{p_{\lambda_P}}"]& \\
       \T^{\ast}\V_{\lambda_P} &  & \T^{\ast}\V_{\lambda}.
    \end{tikzcd}
\end{equation}
The transitivity of these functors is proven in Corollary \ref{cor: eiscompositioncotan}. These functors lift to $\Gm$-equivariant sheaves. 

\begin{proposition}{\label{prop: eiscotangentcompatibility}}
We have the following commutative diagram
\[
\begin{tikzcd}
\wh{\Coh}{}^! (\coprod_{\lambda_M} \wh{\T}\V_{\lambda_M}[-1])^\Gm \arrow[d, "\prod_{\lambda_M} \unshear\circ \kappa", "\simeq"'] \arrow[r, "{\Eis_{P}^{\wh{\T}[-1]}}"] & \wh{\Coh}{}^! (\wh{\T}\V_{\lambda}[-1])^\Gm \arrow[d, "\unshear \circ \kappa", "\simeq"'] \\
\Coh(\coprod_{\lambda_M} \T^* \V_{\lambda_M})^\Gm \arrow[r, "\Eis_{P}^{\T^*}"]
& \Coh(\T^* \V_{\lambda})^\Gm. 
\end{tikzcd}
\]
The same holds for the $\lambda_P$-components and the $\chi$-versions. 
\end{proposition}

\begin{proof}
The vertical equivalences come from Theorem \ref{thm: renormalizedcontinuousKoszul} and the diagram commutes because of the functorialities in \cite[Theorem 2.2.10, Proposition 2.3.3]{Chen2023categorical} (see also the proof of Theorem 2.3.6 of loc.~cit.).
\end{proof}

There is another natural Eisenstein series functor defined as follows. Let the composition of the correspondences in \eqref{eqn: eiscorrespondences} be
\begin{equation}{\label{eqn: Ycorrcomp}}
    \begin{tikzcd}
        & Y^{\geq 0} \arrow[dl, swap, "\breve{\pi}_{q_{\lambda_P}} \circ \breve{\delta}_{p_{\lambda_P}}"] \arrow[dr, " \breve{\pi}_{p_{\lambda_P}} \circ \breve{\delta}_{q_{\lambda_P}}"]& \\
       \T^{\ast}\V_{\lambda_M} &  & \T^{\ast}\V_{\lambda} ,
    \end{tikzcd}
\end{equation}
where
\begin{equation*}
    Y^{\geq 0} = \T^{\ast}\V_{\lambda_M}|_{\V_{\lambda_P}} \times_{\T^{\ast} \V_{\lambda_P}} \T^{\ast} \V_\lambda |_{\V_{\lambda_P}}.
\end{equation*}
Then we define
\begin{equation}{\label{eq: stareis}}
    \Eis^{\ast}_{P,\lambda_P} \coloneq (\breve{\pi}_{p_{\lambda_P}} \circ \breve{\delta}_{q_{\lambda_P}})_{\ast} (\breve{\pi}_{q_{\lambda_P}} \circ \breve{\delta}_{p_{\lambda_P}})^{\ast}.
\end{equation}
\begin{proposition}{\label{prop: cotaneiscompatibility}}
    We have
    \begin{equation*}
       \Eis^{\T^{\ast}}_{P,\lambda_P} \simeq \Eis^{\ast}_{P,\lambda_P} \circ  (\omega_{\V_{U}} \otimes -), 
    \end{equation*}
    where $\omega_{\V_{U}} \coloneq \omega_{\V_{\lambda_P}} \otimes \omega^{-1}_{\V_{\lambda_M}}$ and we note we are implicitly descending $\omega_{\V_{\lambda_P}}$ to $\V_{\lambda_M}$ which is valid since $\omega_{\V_{\lambda_P}}$ is a character of $H_{\lambda_P}$ (up to a shift) and these come via pullback from $H_{\lambda_M}$.
\end{proposition}
\begin{proof}
    By base change, the only difference between $\Eis^{\T^{\ast}}_{P,\lambda_P}$ and $\Eis^{\ast}_{P,\lambda_P}$ is between $\breve{\pi}^!_{q_{\lambda_P}}$ and $\breve{\pi}^{\ast}_{q_{\lambda_P}}$.  The map
    \begin{equation*}
\breve{\pi}_{q_{\lambda_P}}\colon \T^{\ast}\V_{\lambda_M}|_{\V_{\lambda_P}} \rightarrow \T^{\ast}\V_{\lambda_M}
    \end{equation*}
is given by base change along the derived $0$-fiber of the map
\begin{equation*}
    (V^{\vee}_{\lambda_M} \times V_{\lambda_P}) / H_{\lambda_P} \rightarrow (V^{\vee}_{\lambda_M} \times V_{\lambda_M})/ H_{\lambda_M}.
\end{equation*}
Hence it suffices to compute the relative dualizing complex of this latter map and this can be done using the formulas in Proposition \ref{prop: QuasismoothDualizingStack}.
\end{proof}
We record the following computation.
\begin{lemma}{\label{lem: TVgeqdualizingcomputation}}
The dualizing complex of $Y^{\geq 0}$ identifies with 
\[
\omega_{\V_{\lambda_P}}^2 \omega_{\V_{\lambda}}^{-1}\omega_{\V_{\lambda_M}}^{-1}[2(\dim \V_{\lambda} + \dim \V_{\lambda_M} - \dim \V_{\lambda_P})].
\]
\end{lemma}
\begin{proof}
   To compute the dualizing complex of $Y^{\geq 0}$, we need to compute those of $\T^{\ast}\V_{\lambda_M}|_{\V_{\lambda_P}}$, $\T^{\ast}\V_{\lambda}|_{\V_{\lambda_P}}$ and $\T^{\ast}\V_{\lambda_P}$.

   We compute $\T^{\ast}\V_{\lambda}|_{\V_{\lambda_P}}$ by observing the relative dualizing complex of $\T^{\ast}\V_{\lambda}|_{\V_{\lambda_P}}$ over $\V_{\lambda_P}$ equals by base change the relative dualizing complex of $\T^{\ast}\V_{\lambda}$ over $\V_{\lambda}$. Then the dualizing complex of $\T^{\ast}\V_{\lambda}$ is computed by base change along the moment map of $BH_\lambda \to (\Lie H_\lambda)^\vee /H_\lambda$ which has dualizing complex $\det(\Lie H_{\lambda})^{\vee}[-\dim H_{\lambda}]$. This implies $\T^{\ast}\V_{\lambda}$ has dualizing complex $\mc{O}_{\T^{\ast}\V_{\lambda}}[2 \dim \V_{\lambda}]$, so $\T^{\ast}\V_{\lambda}|_{\V_{\lambda_P}}$ has dualizing complex $\omega_{\V_{\lambda_P}}\omega^{-1}_{\V_{\lambda}}[2 \dim \V_{\lambda}]$.

   The same logic shows that the dualizing complex of $\T^{\ast}\V_{\lambda_M}|_{\V_{\lambda_P}}$ is $\omega_{\V_{\lambda_P}}\omega^{-1}_{\V_{\lambda_M}}[2 \dim \V_{\lambda_M}]$ and that $\T^{\ast}\V_{\lambda_P}$ has dualizing complex $\mc{O}_{\T^{\ast}\V_{\lambda_P}}[2 \dim \V_{\lambda_P}]$. This gives the desired result.
\end{proof}
\begin{definition}{\label{defn: twistedGS}}
    We define the twisted Grothendieck--Serre dual on $\Coh(\T^{\ast}\V_{\lambda})$ by the formula 
    \begin{equation*}
\D'_{\GS}\mc F \coloneq \omega_{\V_\lambda}\otimes \D_{\GS}\mc F, 
    \end{equation*}
    where $\omega_{\V_\lambda}$ is pulled back along the structure map $\T^* \V_{\lambda}\to \V_\lambda$. 
\end{definition}
\begin{proposition}{\label{prop: twistedGSEiscompatibility}}
    We have the compatibility
    \begin{equation*}
        \D'_{\GS} \Eis^{\T^{\ast}}_{P, \lambda_P} \simeq \Eis^{\T^{\ast}}_{P, \lambda_P} \D'_{\GS} [2( \dim \V_{\lambda} - \dim \V_{\lambda_P})].
    \end{equation*}
\end{proposition}
\begin{proof}
    Let $s = \breve{\pi}_{q_{\lambda_P}} \circ \breve{\delta}_{p_{\lambda_P}}$ and $t = \breve{\pi}_{p_{\lambda_P}} \circ \breve{\delta}_{q_{\lambda_P}}$ as in \eqref{eqn: Ycorrcomp}. Then $t$ is proper since it is the composite of a base change of $\breve{\delta}_{q_{\lambda_P}}$ which is a closed immersion, and a base change of $p_{\lambda_P}$ which is proper by Lemma \ref{lem: propersmooth}. Similarly, $s$ is quasi-smooth. Hence, using Proposition \ref{prop: cotaneiscompatibility}, we have
    \begin{equation*}
        \D'_{\GS} \Eis^{\T^{\ast}}_{P, \lambda_P} \simeq \omega_{\V} \otimes \D_{\GS} t_{\ast}s^{\ast} ( \omega_{\V_U}\otimes-) \simeq \omega_{\V} \otimes t_{\ast}s^{!} \D_{\GS} ( \omega_{\V_U}\otimes-).
    \end{equation*}
    Now we use Lemma \ref{lem: TVgeqdualizingcomputation} to get the above equals
    \begin{equation*}
      t_{\ast} s^{\ast} (\omega_{\V} \otimes \omega^2_{\V_P}\omega^{-1}_{\V}\omega^{-1}_{\V_M} \otimes \omega^{-1}_{\V_U} \omega^{-1}_{\V_M} \otimes  \D'_{\GS} -)[2( \dim \V_{\lambda} - \dim \V_{\lambda_P})]
    \end{equation*}
    where we are using $\omega_{\T^*\V_{\lambda_M}}\simeq \mc O_{\T^*\V_{\lambda_M}}[2\dim \V_{\lambda_M}]$. Finally, we observe that the twist reduces simply to $\omega_{\V_U}$ as desired.
\end{proof}

\subsection{Eisenstein series of \texorpdfstring{$\mathcal{D}$}{D}-modules}
We further consider Eisenstein series for \emph{right} $\mathcal{D}$-modules on Vogan stacks. Define
\[
\Eis_P^{\mc{D},!}\colon \mc{D}(\coprod_{\lambda_M} \V_{\lambda_M}) \xrightarrow{q_\V^!} \mc{D} (\coprod_{\lambda_P} \V_{\lambda_P}) \xrightarrow{p_{\V*}} \mc{D}(\V_\lambda). 
\]
It lifts to filtered $\mc{D}$-modules:
\[
\Eis_P^{F\mc{D}, !}\colon F\mc{D}(\coprod_{\lambda_M} \V_{\lambda_M}) \xrightarrow{q_\V^!} F\mc{D} (\coprod_{\lambda_P} \V_{\lambda_P}) \xrightarrow{p_{\V*}} F\mc{D}(\V_\lambda), 
\]
and it preserves the coherent subcategories $F\mathcal{D}_c$. We define $\Eis_{P,\lambda_P}^{\mc{D}, !}$ and $\Eis_{P,\lambda_P}^{F\mc{D}, !}$ similarly. 
We define $\Eis_P^{\mc{D},*}$, $\Eis_P^{F\mc{D},*}$, $\Eis_{P, \lambda_P}^{\mc{D},*}$, $\Eis_{P, \lambda_P}^{F\mc{D},*}$ by replacing $!$-pullback by $*$-pullback.  

\begin{proposition}
We have the following commutative diagram:
\[
\begin{tikzcd}
\mc{D} (\coprod_{\lambda_M} \V_{\lambda_M}) 
\arrow[d, "\Eis_P^{\mc{D}, !}"]
& 
F\mc{D}_c(\coprod_{\lambda_M} \V_{\lambda_M}) 
\arrow[l, "\mathrm{un}"'] 
\arrow[r, "\mathrm{gr}"]
\arrow[d, "\Eis_P^{F\mc{D}, !}"] 
& 
\Coh(\coprod_{\lambda_M} \T^{\ast}\V_{\lambda_M})^\Gm 
\arrow[d, "\Eis_P^{\T^*}"] 
\\
\mc{D}(\V_\lambda) 
& 
F\mc{D}_c(\V_\lambda) 
\arrow[l, "\mathrm{un}"'] 
\arrow[r, "\mathrm{gr}"]
& 
\Coh(\T^{\ast}\V_\lambda)^\Gm. 
\end{tikzcd}
\]
\end{proposition}

\begin{proof}
    This follows from Remark \ref{rem: Dmodulefunctoriality}.
\end{proof}
\begin{remark}
    One expects that for a smooth Artin stack $X$, there is an equivalence
    \begin{equation*}
        \gr \circ \D_{\mathrm{Ver}} \simeq \D'_{\GS} \circ \gr(-2\dim X) [- 2 \dim X],
    \end{equation*}
    where $\D_{\mathrm{Ver}}$ is the Verdier duality functor for (filtered) $\mc{D}$-modules, $\D'_{\GS}$ is the generalization of Definition \ref{defn: twistedGS} to $\T^{\ast}X$ in the obvious way, and $\gr\colon F\mc{D}_c(X) \rightarrow \Coh(\T^{\ast}X)^{\Gm}$. At least when $X$ is a smooth scheme, this is essentially proven in \cite[Proposition 2.3]{SaitoSchnellDuality}, and the general case should follow by smooth descent as in \cite[Proposition 3.3.18]{Chen2023categorical}. 
    \end{remark}
\subsection{Eisenstein series of constructible sheaves}
Finally, we consider constructible sheaves on Vogan stacks, assuming $\Lambda=\C$.\footnote{Everything here except RH also makes sense for any algebraically closed $\Lambda$ with the coefficient field $\C$ replaced by $\overline{\Q}_\ell$ for some $\ell$.} 
We define
\[
\Eis_P^{\mathrm{cons},!}\coloneq p_{\V*}\circ q_{\V}^!\colon D_c^b(\coprod_{\lambda_M} \V_{\lambda_M}) \xrightarrow{q_\V^!} D_c^b (\coprod_{\lambda_P} \V_{\lambda_P}) \xrightarrow{p_{\V*}} D_c^b(\V_\lambda),  
\]
and $\Eis_{P,\lambda_P}^{\mathrm{cons}, !}$, $\Eis_P^{\mathrm{cons},*}$, $\Eis_{P,\lambda_P}^{\mathrm{cons},*}$ similarly. We have $ \Eis_P^{\mathrm{cons}, !} \circ \RH\simeq \RH \circ \Eis_P^{\mathcal{D},!} $ and so on. The transitivity of $\Eis_{P,\lambda_P}^{\mathrm{cons}, !}$ follows from proper base change. 

A \emph{left} adjoint of $\Eis_{P, \lambda_P}^{\mathrm{cons},!}$ is given by $\CT_{P, \lambda_P}^{\mathrm{cons}}\coloneq q_{\lambda_P !}\circ p_{\lambda_P}^*$, while a \emph{right} adjoint is given by $q_{\lambda_P *}(p_{\lambda_P}^!(-) \otimes \omega_{\V_{\lambda_P}/\V_{\lambda_M}}^{\mathrm{cons,-1}})$ as $q_{\lambda_P}$ is smooth, where the superscript $\mathrm{cons}$ indicates that it is the dualizing complex in the constructible setting. 

Moreover, the left adjoint of $\Eis_{P, \lambda_P}^{\mathrm{cons},!}$ and the right adjoint of $\Eis_{P^\op, \lambda_{P^\op}}^{\mathrm{cons},!}$ can be related via Braden's hyperbolic localization \cite{Braden} since all constructible sheaves here are $\Gm$-monodromic, cf.~\cite[Proposition 1.8]{LiuTsaiVilonenXue}. 

\section{Cuspidal constructible sheaves and \texorpdfstring{$\mc D$}{D}-modules on Vogan stacks}{\label{s: DMod}}
For pure inner forms, Vogan's form of the local Langlands correspondence \cite{VoganLLC} has a different geometric realization from the categorical conjecture: irreducible representations correspond to simple equivariant perverse sheaves on Vogan varieties, and such sheaves are intermediate extensions of local systems. In fact, the study of such equivariant perverse sheaves goes back to the work of Lusztig. Based on the works of \cite{Chen2023categorical} and \cite{BCHNconstructible}, we can think of the constructible picture as the Koszul dual of the coherent picture up to deformation, and we use such a relationship later in our more generalized setup involving cuspidal sheaves to classify cuspidal coherent sheaves. We record here some material on cuspidal local systems and their relation with (filtered) $\mc D$-modules in our favorite generality. 

\subsection{Results of Lusztig in the graded setting}{\label{ss: lusztiggraded}}
We briefly recall results of Lusztig \cite{LusztigGraded} (see also \cite{Chatterjee}). His results are true for any algebraically closed $\Lambda$ and $\overline{\Q}_{\ell}$-sheaves for any $\ell$ (recall that $\Lambda$ always has characteristic $0$). Moreover, the same holds true for $\Lambda=\C$ and $\C$-constructible sheaves as all monodromy involved is finite.

The setting is as follows: $G$ is split and the unramified semisimple parameter $\lambda\colon \frob^\Z\to \wh G$ arises from a cocharacter $\iota$ of $\wh G$ that factors through the fixed maximal torus $\wh T$, $H_{\lambda}=(\wh G)^\iota$, and $\lambda(\frob)=\iota(q^{1/n})$ for some integer $n\geq 1$ so that $V_\lambda=(\Lie \wh{G})(n)$. 

Lusztig started from the ungraded setting. Let $\mc L$ be an irreducible $\wh G$-equivariant local system on a nilpotent orbit in $\Lie \wh {G}$. Lusztig defined a notion of \emph{cuspidal} local system in \cite[\S 2.2]{LusztigCuspidalI}. As in \cite[p.~179]{LusztigFourier} (and \cite{LusztigCharacterI, LusztigCharacterII-III, LusztigCharacterV}), these cuspidal local systems are closely related to various other notions in \cite{LusztigIntersection, LusztigCharacterI, LusztigCharacterII-III}, and in particular, they correspond to cuspidal local systems on the unipotent variety of $\wh G$ defined in \cite[Definition 2.4]{LusztigIntersection}. The latter are classified in the same paper. By definition, the class of cuspidal local systems is preserved by any automorphism of $\wh G$. 

We now move to the graded setting. Let $\mc L$ be an irreducible $H_{\lambda}$-equivariant local system on a single $H_{\lambda}$-orbit in $V_{\lambda}$. It is \emph{cuspidal} \cite[\S 4.4]{LusztigGraded} if it is the $*$-restriction of a cuspidal local system in the previous paragraph. If $\mc L$ is cuspidal, then, by the argument in \cite[\S 4.3]{LusztigIntersection}, after replacing $\iota$ and $n$ (without changing $H_{\lambda}$ and $V_{\lambda}$), we may assume that $n=2$, and $\iota$ determines the semisimplification of the $L$-parameter $\phi$ arising from a map $\tilde{\phi}\colon \SL_2 \to \wh G$ whose associated nilpotent orbit is distinguished in the sense of \cite{BalaCarterI} (equivalently, the centralizer $S_\phi$ of $\phi$, equivalently that of $\tilde{\phi}$, is finite modulo $Z(\wh{G})^{\Gamma}$, as proven below), and the $H_{\lambda}$-orbit associated to $\mc{L}$ is open (and dense) in $V_{\lambda}$. 
(Moreover, the cuspidal local system in the previous paragraph that restricts to $\mc L$ is unique.)

For convenience, we record the following well-known lemma.
\begin{lemma}{\label{lem: discequalsdistforunip}}
Let $G$ be a split connected reductive group over $F$ and $\phi=(\lambda, N)$ a unipotent $L$-parameter of $G$ such that $\lambda(\frob) = s_{\lambda}$, the torsion-free part, i.e., $t_\lambda$ is trivial. Then $\phi$ is essentially discrete if and only if the associated nilpotent element is distinguished (as defined in \cite[\S8.2]{CollingWoodMcGovern}, for instance).
\end{lemma}
\begin{proof}
If $G$ is semisimple and $\phi$ is discrete, then our assumptions imply that a lift $\tilde{\phi}\colon W_F \times \SL_2 \rightarrow \co{L}{G}$ will have trivial restriction to $W_F$ (Lemma \ref{lem: tftorsdisclemma}). Hence the centralizer of $\tilde{\phi}$ is just the centralizer of the restriction to $\SL_2$. This is a maximal reductive subgroup of the centralizer of $\phi$ by \cite[Proposition 3.7]{BMImaiYoucis} and hence its identity component is trivial by the discreteness assumption. But the connected centralizer of the map $\SL_2 \rightarrow \wh{G}$ equals the connected reductive part of the centralizer of the nilpositive element of the corresponding $\SL_2$-triple (for instance \cite[Lemma 3.7.3]{CollingWoodMcGovern}) which is $N$. By definition, this nilpotent element is distinguished if the connected reductive part of its centralizer is trivial in $\wh{G}$, and this is therefore implied by discreteness of $\phi$. Conversely, if $N$ is distinguished, then the centralizer of $\phi$ does not contain a nontrivial torus, meaning that $\phi$ cannot factor through a proper Levi subgroup of $\wh{G}$ and hence is discrete by \cite[Corollary 3.16]{BMImaiYoucis}.

For non-semisimple $G$, we argue as above with $\wh{G_{\der}}$.
\end{proof}

Let us recall some results from \cite{LusztigVanishing, LusztigGraded}.

\begin{proposition}[{\cite[\S 4.4(b)]{LusztigGraded}}]{\label{prop: constrcleanness}}
If $\mc L$ is cuspidal, then the intermediate extension $\mc{IC} (\mc L)$ of $\mc L$ is clean, i.e., it agrees with both the $!$- and $*$-extensions. 
\end{proposition}

\begin{remark}
Together with the classification of cuspidal local systems, it follows that the Verdier dual of $\mc{IC} (\mc L)$ also has the form $\mc{IC} (\mc L')$ for some cuspidal $\mc L'$ (up to a shift). This can be deduced from Corollary \ref{cor: LusztigCuspidal_is_equivalent_to_CT_vanishing} below too.   
\end{remark}

\begin{corollary}\label{cor: Hom_between_cuspidal_perverse_sheaves}
Let $\mc{IC}(\mc L)$, $\mc{IC}(\mc L')$ be two simple perverse sheaves (up to shifts so that they are the intermediate extensions of local systems in degree 0 and possibly coming from different orbits in $V_{\lambda}$). If one of them is cuspidal and $\mc{IC}(\mc L)\not\simeq \mc{IC}(\mc L')$, then $\mathrm{Hom}(\mc{IC}(\mc L), \mc{IC}(\mc L'))\simeq 0$. If both are cuspidal, then
\[
\mathrm{Hom}(\mc{IC}(\mc L), \mc{IC}(\mc L'))\simeq \mathrm{Hom}(\mc L, \mc L') \simeq \mathrm{Hom}_{\pi_0(S_{\phi})}(\rho, \rho')\otimes \Sym(\Lie Z(\wh G)^\vee[-2]), 
\]
where $\rho$, $\rho'$ are irreducible representations of the finite group $\pi_0(S_{\phi})$ corresponding to $\mc L$, $\mc L'$.
\end{corollary}

\begin{proof}
Suppose $\mc{IC}(\mc L)$ is cuspidal. Then, by using the adjunction for $!$-extension, 
\[
\mathrm{Hom}(\mc{IC}(\mc L), \mc{IC}(\mc L'))\simeq \mathrm{Hom}_{D_c^b(BS_\phi)}(\rho, \rho'),
\]
where $\rho'$ corresponds to the $*$-restriction of $\mc{IC}(\mc L')$ to the open orbit of $V_{\lambda}$; if the $*$-restriction vanishes, we take $\rho'=0$. If the $*$-restriction does not vanish, then $\mc L'$ is a simple local system on the open orbit of $V_{\lambda}$, and $\rho'$ is an irreducible representation of $\pi_0(S_{\phi})$, and the desired statements follow, with the standard computation of equivariant cohomology. 

The case where $\mc{IC}(\mc L')$ is cuspidal is similar, by using the adjunction for $*$-extension instead.
\end{proof}

In particular, all the cuspidal $\mc{IC} (\mc L)$ generate a semisimple $\C$-linear abelian category and a stable full $\infty$-subcategory of $D_c^b(\V_\lambda)$ (which is semisimple if $G$ is semisimple). We define a \emph{cuspidal constructible sheaf} to be an object of the latter. It follows that the left- and the right-orthogonal to cuspidal constructible sheaves in $D_c^b (\V_\lambda)$ are both generated by all non-cuspidal simple perverse sheaves. 

\begin{proposition}[{\cite[\S 4.4(c)]{LusztigGraded}}]\label{prop: CT_cuspidal_constructible_sheaves=0}
Let $\wh M\subset \wh P\subset \wh G$ be proper semistandard Levi and parabolic subgroups of $\wh{G}$ containing $\wh T$. For any lift $\lambda_P$ of $\lambda$, 
$\CT_{P, \lambda_P}^{\mathrm{cons}}\mc{IC} (\mc L)\simeq 0$ for any cuspidal $\mc L$. 
\end{proposition}

\begin{proof}
After passing to a conjugate of $(P, \lambda_P)$, we may assume the representative $\lambda$ is the composite of a representative $\lambda_P$ and $\wh P\to \wh G$. The vanishing of $\CT_{P, \lambda_P}^{\mathrm{cons}}\mc{IC} (\mc L)$ can be checked on stalks at closed points of $V_{\lambda_M}$, and the cohomology groups of these stalks are, by proper base change and the cleanness of $\mc{IC}(\mc L)$, exactly those appearing in \cite[\S 4.4(c)]{LusztigGraded}, which are shown to vanish. 
\end{proof}

\begin{corollary}
The left- and the right-orthogonal to cuspidal constructible sheaves both contain the essential image of $\Eis_{P,\lambda_P}^{\mathrm{cons}, !}$.  
\end{corollary}

\begin{proof}
By the previous proposition and the adjunction, the essential image of $\Eis_{P, \lambda_P}^{\mathrm{cons}, !}$ is contained in the right-orthogonal. By the description of left- and right-orthogonal above, it is also contained in the left-orthogonal. 
\end{proof}

In fact, the following is proved in \cite{LusztigGraded}:

\begin{proposition}\label{prop: OrthogonalCuspidalConstructibleSheaves}
The left- and the right-orthogonal to cuspidal constructible sheaves are both generated, as an idempotent-complete stable $\infty$-category, by the images of $\Eis_{P,\lambda_P}^{\mathrm{cons}, !}$ for proper $P$ as above. 
\end{proposition}

\begin{remark}
Simple perverse sheaves whose support in $V_\lambda$ is contained in the complement of the open orbit lie in both the left and right orthogonal complement of the category of cuspidal constructible sheaves, and the claim that they are contained in the category generated by $\Eis_{P,\lambda_P}^{\mathrm{cons}, !}$ follows from considering a $\Theta$-stratification of $\V_\lambda$ as in \S \ref{sss:ThetaStratificationCotangent}. For the $!$-extension of a local system on the open orbit, the essential claim is that cuspidality is equivalent to the vanishing of certain compactly supported cohomology groups as in \cite[\S 4.4(c)]{LusztigGraded}; see below.  
\end{remark}

\begin{proof}
It suffices to show that any non-cuspidal simple perverse sheaf $\mc{IC}(\mc L)$ is, up to a shift, a direct summand of some $\Eis_{P,\lambda_P}^{\mathrm{cons}, !}\mc F$ with $P$ proper. By \cite[Proposition 7.5]{LusztigGraded}, $\mc{IC}(\mc L)$ is normal in the sense of \cite[\S 7.1]{LusztigGraded} meaning that some shift of $\mc{IC}(\mc L)$ is a direct summand of Lusztig's induction of a cuspidal pair \cite[\S 6]{LusztigGraded}. By the following lemmas, Lusztig's induction agrees with $\Eis_{P,\lambda_P}^{\mathrm{cons}, !}$ up to a shift.  
\end{proof}

\begin{lemma}
Assume $\lambda$ is the composite of $\lambda_P$ and $\wh P \to \wh G$. 
Lusztig's induction \cite[\S 6]{LusztigGraded} agrees with $\Eis_{P, \lambda_P}^{\mathrm{cons},*}$. 
\end{lemma}

\begin{proof}
As in \cite[Lemma 4.3]{Chatterjee}, Lusztig's induction can be written using the Bernstein--Lunts induction equivalence applied to the diagram
\[
V_{\lambda_M} \leftarrow V_{\lambda_P} \leftarrow H_{\lambda}\times V_{\lambda_P} \to H_{\lambda}\times^{H_{\lambda_P}} V_{\lambda_P} \rightarrow V_{\lambda}. 
\]
This translates to $\Eis_{P, \lambda_P}^{\mathrm{cons},*}$ in the stacky formulation.\footnote{\cite[Theorem 4.1]{Chatterjee} does not hold in the stacky formulation as the usual normalization of dualizing complexes is different.}
\end{proof}

\begin{lemma}
The functors $\Eis_{P, \lambda_P}^{\mathrm{cons},*}$ and $\Eis_{P, \lambda_P}^{\mathrm{cons},!}$ agree up to a shift.
\end{lemma}

\begin{proof}
It suffices to show that the relative dualizing complex $\omega_{\V_{\lambda_P} /\V_{\lambda_M}}$ is a constant sheaf up to a shift. As the map $\V_{\lambda_P} \to \V_{\lambda_M}$ is smooth, it is a local system in one degree. As $H_{\lambda_P}$ is connected in the current setting, any $H_{\lambda_P}$-equivariant local system on $V_{\lambda_P}$ is constant as desired. 
\end{proof}

We conclude the following alternative characterization of cuspidality:

\begin{corollary}\label{cor: LusztigCuspidal_is_equivalent_to_CT_vanishing}
A constructible sheaf $\mc F \in D_c^b (\V_\lambda)$ is cuspidal if and only if $\CT_{P, \lambda_P}^{\mathrm{cons}}\mc F\simeq 0$ for all proper $P$ and all $\lambda_P$. 
\end{corollary}

\begin{proof}
One direction is Proposition \ref{prop: CT_cuspidal_constructible_sheaves=0}. For the other direction, we see that the assumption implies that $\mc F$ is left-orthogonal to all non-cuspidal simple perverse sheaves by (the proof of) the previous proposition. In particular, it is the $!$-extension of a sheaf on the open orbit. Moreover, it is left-orthogonal to the intermediate extension of any non-cuspidal local system. It follows that $\mc F$ is a finite colimit of cuspidal simple perverse sheaves. 
\end{proof}

\begin{remark}
In the $\Z/m\Z$-graded setting, the above condition is called supercuspidal in \cite[Definition 1.12]{LiuTsaiVilonenXue}. They prove that a clean extension from a distinguished open orbit is supercuspidal \cite[Proposition 1.15]{LiuTsaiVilonenXue}, and the same argument would work in the $\Z$-graded setting here. 

As already mentioned, the left adjoint of $\Eis_{P, \lambda_P}^{\mathrm{cons},!}$ and the right adjoint of $\Eis_{P^\op, \lambda_{P^\op}}^{\mathrm{cons},!}$ can be related via Braden's hyperbolic localization \cite{Braden}. This means one can directly show that the supercuspidal sheaves are also \emph{right}-orthogonal to the essential image of $\Eis_{P, \lambda_P}^{\mathrm{cons},!}$.  

\end{remark}

\subsection{Reduction to the graded setting for cuspidal constructible sheaves}{\label{ss: CCSreductiontograded}}

Now we let $\lambda\in X_G (\Lambda)$ be a general semisimple $L$-parameter. In this subsection we prove results relating cuspidal local systems, or rather constructible sheaves, on $\V_{\lambda} \hookrightarrow \Par_G$ with those of $\V_{\lambda^{\ur}} \hookrightarrow \LS^{\up}_{I_{\lambda}}\coloneq \LS_{I_\lambda}^{\WD_F, \up}$ from the previous setting; we recall the notation below. This is closely related to works of Aubert--Moussaoui--Solleveld \cite{AMS1, AMS2} and Solleveld \cite{Solleveld}. 

\begin{definition}\label{def: CuspidalLocalSystem}
A constructible sheaf $\mc F \in D^b_c(\V_\lambda)$ is \emph{cuspidal} if $\CT_{P, \lambda_P}^{\mathrm{cons}} \mc F=0$ for all proper semistandard $P$ and $\lambda_P$. Equivalently, it is \emph{left}-orthogonal to $\Eis_{P,\lambda_P}^{\mathrm{cons},!}$ for all proper semistandard $P$ and $\lambda_P$.  
\end{definition}

\begin{remark}
The $\infty$-category of cuspidal constructible sheaves is an idempotent-complete stable $\infty$-category. 
\end{remark}

Choose a representative $\lambda\colon W_F\to \wh G$; we use the same notation. Recall that, fixing an isomorphism of $\Z$-modules $\Lambda^{\times} = \Lambda_{\tor}^{\times} \oplus \Lambda_{\tf}^{\times}$ with $q^{1/2} \in \Lambda_{\tf}^{\times}$, we have the decomposition ${}^L \lambda(\frob)=(s_\lambda, e)(t_\lambda, \frob)=(t_\lambda, \frob)(s_\lambda, e)$, and
\[
K_\lambda=Z_{\wh G}({}^L \lambda_{I}), \quad J_\lambda\coloneq Z_{K_\lambda}(t_\lambda, \frob) \supset  H_\lambda \coloneq Z_{J_\lambda} (s_\lambda)=Z_{K_\lambda} (\co{L}{\lambda}(\frob)). 
\]

Choose $\lambda$ as in Lemma \ref{lem: TJcontrol} and write $\wh T_J\coloneq \wh T\cap J_\lambda^\circ$. 
Taking a connected reductive subgroup $I_\lambda \subset J_\lambda$, a cocharacter \[
\chi\colon \Gm\to \wh T_I\coloneq \wh T\cap I_\lambda=\wh T_J,
\]
and an integer $n\geq 1$ as in Lemma \ref{lem: ReductionToGradedUnipotent}, we have
\[
Z_{I_\lambda}(s_\lambda)=Z_{I_\lambda}(\chi)=H_\lambda^\circ, \quad V_\lambda=\Lie I_\lambda (n). 
\]
This is the setting in which \cite{LusztigGraded} is available. There is a canonical map
\[
\V_{\lambda^\ur} \coloneq V_\lambda /H_\lambda^\circ  \to \V_\lambda.  
\]

We first observe that

\begin{lemma}\label{lem: reduction_to_LusztigGraded}
If a constructible sheaf $\mc F\in D_c^b (\V_\lambda)$ is cuspidal, then its pullback to $\V_{\lambda^\ur}$ is cuspidal in the sense of Lusztig (depending on $I_\lambda$ and $\chi$).
\end{lemma}

\begin{proof}
Let $\mc F^{\ur}$ denote the pullback of $\mc F$. By Corollary \ref{cor: LusztigCuspidal_is_equivalent_to_CT_vanishing}, it suffices to show that $\CT^{\mathrm{cons}}_{P_{I},\lambda_{P}^{\ur}} \mc F^{\ur}\simeq 0$ for all proper semistandard parabolic subgroups $\wh P_I$ of $I_\lambda$ containing $\wh T_I$ and all $\lambda_P^{\ur}$ lifting $\lambda^\ur$ such that the representative $\lambda^\ur$ is the composite of a representative $\lambda^\ur_P$ with $\wh P_I\to I_\lambda$. 

Fix $\wh P_{I}$ with a semistandard Levi subgroup $\wh M_{I}$ and suppose it arises from a cocharacter $\nu\in X_*(Z(\wh M_{I})^\circ)\subset X_*(\wh T_I)$ via the dynamic method. Then $\nu$ determines a proper semistandard parabolic subgroup $\wh P\rtimes W_F$ of $\wh G\rtimes W_F$ with a semistandard Levi subgroup $\wh M=Z_{\wh G}(\nu)$. 

Let $\lambda_{P}^\ur\colon \frob^\Z\to \wh P_{I}$ be the lifting of $\lambda^\ur$ as above, and let $\lambda_{M}^\ur$ denote the resulting map $\frob^\Z\to \wh M_{I}$. Since $\nu$ centralizes ${}^L\lambda (W_F)$, ${}^L \lambda$ factors through $\wh M\rtimes W_F$, and we write $\lambda_M$, $\lambda_P$ for liftings. We have the following commutative diagram
\[
\begin{tikzcd}
    \V_{\lambda_{M}^{\ur}} \arrow[d] & \arrow[l] \V_{\lambda_{P}^{\ur}} \arrow[d] \arrow[r] & \V_{\lambda^{\ur}} \arrow[d] \\
    \V_{\lambda_M}  & \arrow[l] \V_{\lambda_P} \arrow[r] & \V_{\lambda}. 
\end{tikzcd}
\]
The left square is cartesian: $\pi_0 (H_{\lambda_P})\cong \pi_0 (H_{\lambda_M})$ (just because the kernel of the surjection $H_{\lambda_P}\to H_{\lambda_M}$ is contained in $U_{\wh P}$, and hence connected). Therefore, the vanishing $\CT^{\mathrm{cons}}_{P_I,\lambda_{P}^{\ur}} \mc F^{\ur}\simeq 0$ follows from the cuspidality of $\mc F$ and base change. 
\end{proof}

\begin{lemma}{\label{lem: GIcenters}}
The Vogan stack $\V_{\lambda}$ is essentially discrete if and only if $\V_{\lambda^{\ur}}$ is essentially discrete and $Z(\wh{G})^{\Gamma, \circ} = Z(I_{\lambda})^{\circ}$.
\end{lemma}
\begin{proof}
Let $\phi\colon W_F \ltimes\Ga\to \wh G$ and $\phi^\ur\colon \frob^\Z\ltimes \Ga\to I_\lambda$ be corresponding representatives of the open orbits of $V_\lambda=V_{\lambda^\ur}$. Then the connected centralizer of $\phi^\ur$ in $I_{\lambda}$ equals the connected centralizer of $\phi$ in $\wh{G}$. On the other hand, $\phi$, $\phi^\ur$ are essentially discrete if and only if their connected centralizers are $Z(\wh{G})^{\Gamma, \circ}$, $Z(I_{\lambda})^{\circ}$ respectively. We always have $Z(\wh{G})^{\Gamma, \circ} \subset Z(I_{\lambda})^{\circ}$, which means essential discreteness of $\phi$ implies essential discreteness of $\phi^\ur$ and the equality of connected centers. Conversely, if we have an equality of connected centers, then essential discreteness of $\phi^\ur$ implies that of $\phi$. 
\end{proof}

\begin{corollary}
If a nonzero cuspidal constructible sheaf exists, then $\lambda$ is the semisimplification of an essentially discrete $L$-parameter $\phi$.
\end{corollary}

\begin{proof}
By Lemmas \ref{lem: reduction_to_LusztigGraded}, \ref{lem: GIcenters}, and \cite{LusztigGraded} as we recalled above, it suffices to show $Z(\wh{G})^{\Gamma, \circ} = Z(I_{\lambda})^{\circ}$. Suppose $Z(I_{\lambda})^{\circ}(\subset \wh T_I)$ is strictly larger than $Z(\wh{G})^{\Gamma, \circ}$. Then there is a proper semistandard Levi subgroup $\wh M$ containing $I_\lambda$, and $\lambda$ lifts to $\lambda_M$. For any semistandard parabolic subgroup $\wh P=\wh M \wh U$ with $\lambda_P=\lambda_M$, it is easy to see that $\CT_{P,\lambda_P}^{\mathrm{cons}} \mc F\neq 0$ if $\mc F\neq 0$, since $*$-pullback along the \'etale cover $\V_{\lambda^\ur}\to \V_\lambda$ is conservative. 
\end{proof}

Assuming $\lambda$ is the semisimplification of an essentially discrete $L$-parameter $\phi$, we want to discuss the converse direction next. 

\begin{lemma}{\label{lem: lambdaIparabolic}}
Suppose that $Z(\wh{G})^{\Gamma, \circ}= Z(I_{\lambda})^{\circ}$ and $\wh M\subset \wh P \subset \wh{G}$ are semistandard Levi and parabolic subgroups with Levi decomposition $\wh P=\wh M \wh U$. If $\lambda$ lifts to $\wh P$ and $I_{\lambda} \subset \wh P$, then $\wh M =\wh P= \wh{G}$. 
\end{lemma}

\begin{proof}
Since $\lambda$ lifts to $\wh{P}$ and is assumed semisimple, it in fact factors through $\wh{M}$ up to $\wh{P}$-conjugacy by Lemma \ref{lem: X_P=X_M}. This implies (after possibly replacing $\lambda$ with a $\wh P$-conjugate) that $Z(\wh{M})^{\Gamma, \circ} \subset I_{\lambda}$. 

Now let $\ov{\wh M} = \wh P/\wh U$. As $I_{\lambda} \subset \wh P$, consider the projection $\ov{I_{\lambda}} \subset \ov{\wh M}$. Since $I_{\lambda}$ is reductive, this projection induces an isomorphism $I_{\lambda} \cong \ov{I_{\lambda}}$. Now, 
\[
Z(\wh M)^{\Gamma, \circ} \cong Z(\ov{\wh M})^{\Gamma, \circ} \subset Z(\ov{I_{\lambda}})^{\circ} \cong Z(I_{\lambda})^{\circ} = Z(\wh{G})^{\Gamma, \circ},
\]
where the last part is the assumption. Hence, $\wh M= \wh{G}$ as desired. 
\end{proof}

\begin{lemma}
Suppose $\lambda$ is the semisimplification of an essentially discrete $L$-parameter $\phi$. If $\mc F^{\ur}$ is cuspidal, then $\mc F$ is cuspidal. 
\end{lemma}

\begin{proof}
By assumption and Lemma \ref{lem: GIcenters}, $Z(I_\lambda)^\circ = Z(\wh G)^{\Gamma, \circ}$. After possibly conjugating $(P, \lambda_P)$, we need to show $\CT_{P,\lambda_P}^{\mathrm{cons}} \mc F\simeq 0$ for all proper semistandard $\wh P$ and $\lambda_P$ such that the representative $\lambda$ factors through $\wh P$, and $\lambda_P=\lambda$. By (the proof of) Lemma \ref{lem: TJcontrol}, $\lambda$ can in fact be taken to factor through $\wh{M}$ and $\wh P_I\coloneq\wh P \cap I_{\lambda}$ is a proper semistandard parabolic subgroup of $I_\lambda$ since $I_\lambda \not\subset \wh P$ by Lemma \ref{lem: lambdaIparabolic}.

As in the proof of Lemma \ref{lem: reduction_to_LusztigGraded}, we have a lift $\lambda_P^\ur\colon \frob^\Z\to \wh P_I$ of the representative $\lambda^\ur$ corresponding to the lift $\lambda_P\colon W_F\to \wh P$ of $\lambda$, and we have the commutative diagram
\[
\begin{tikzcd}
    \V_{\lambda_{M}^{\ur}} \arrow[d] & \arrow[l] \V_{\lambda_{P}^{\ur}} \arrow[d] \arrow[r] & \V_{\lambda^{\ur}} \arrow[d] \\
    \V_{\lambda_M}  & \arrow[l] \V_{\lambda_P} \arrow[r] & \V_{\lambda}
\end{tikzcd}
\]
with the left square being cartesian. Since the $*$-pullback along $\V_{\lambda^\ur_M}\to \V_{\lambda_M}$ is conservative, $\CT_{P,\lambda_P}^\mathrm{cons} \mc F\simeq 0$ reduces to $\CT_{P_I,\lambda_P^\ur}^\mathrm{cons} \mc F^\ur\simeq 0$ by base change. 
\end{proof}

\subsection{Classification of cuspidal simple perverse sheaves}

\begin{definition}[{cf.~\cite[Definition 6.9]{AMS1}}]\label{def:purecuspidal}
Suppose $\lambda$ is the semisimplification of an essentially discrete $L$-parameter $\phi$. 
A representation $\rho$ of $\pi_0(S_\phi)$ is \emph{cuspidal} if every irreducible summand of $\rho$ as a semisimple representation of $\pi_0 (Z_{H_\lambda^\circ}(\phi))$ via the natural map
\[
\pi_0 (Z_{H_\lambda^\circ}(\phi)) \to \pi_0 (Z_{H_\lambda}(\phi))\cong \pi_0(Z_{\wh G}(\phi))
\]
is cuspidal in the sense of Lusztig (applied to the group $I_\lambda$). 
\end{definition}

\begin{theorem}\label{thm: ClassificationCuspidalPerverseSheaves}
Suppose $\lambda$ is the semisimplification of an essentially discrete $L$-parameter $\phi$.  
\begin{enumerate}
    \item The intermediate extension from the open orbit of the local system corresponding to a cuspidal irreducible representation is cuspidal, and it is a clean extension.\footnote{Cuspidality would also be characterized by saying $\V_\lambda$ is essentially discrete and $\mc F$ is a clean extension from the open orbit.}
    \item Any cuspidal simple perverse sheaf is the intermediate extension of the local system corresponding to a cuspidal irreducible representation from the open orbit. 
    \item The category of cuspidal constructible sheaves is semisimple if $G$ is semisimple, and the abelian category of cuspidal perverse sheaves is semisimple for any $G$.
    \item The $\infty$-category of cuspidal constructible sheaves is closed under Verdier duality. 
\end{enumerate}
\end{theorem}

\begin{proof}
We have proven that $\mc F$ is cuspidal if and only if $\mc F^\ur$ is cuspidal, and (4) follows from the graded case as the map $\V_{\lambda^\ur}\to \V_\lambda$ is finite \'etale. Again by the finiteness, $\mc F$ must be a direct summand of the pushforward of $\mc F^\ur$. If $\mc F$ is simple, $\mc F^\ur$ is a direct sum of shifts of clean extensions of cuspidal local systems by the classification of Lusztig. This implies that $\mc F$ is also a direct sum of shifts of clean extensions of certain simple local systems $\mc L$; if $\mc L$ corresponds to an irreducible representation $\rho$ of $\pi_0(S_\phi)$, then $\rho$ is cuspidal in the sense above. In fact, any such cuspidal $\rho$ must appear in the pushforward of some cuspidal constructible sheaf on $\V_{\lambda^\ur}$. This proves (1) and (2). Finally, (3) follows immediately as in the proof of Corollary \ref{cor: Hom_between_cuspidal_perverse_sheaves}.
\end{proof}

\begin{remark}
For the most part, one can argue more directly: using the $\Theta$-stratification of $\V_\lambda$ in the graded case from \S \ref{sss:ThetaStratificationCotangent}, one can show that the category generated by $\Eis_{P, \lambda_P}^{\mathrm{cons}, !}\mc F$ contains all sheaves supported on the complement of the open orbit. It follows that any cuspidal $\mc F$ is a $!$-extension. By using that the cuspidal sheaves are also right-orthogonal to the essential image of $\Eis_{P, \lambda_P}^{\mathrm{cons}, !}$, it also follows that $\mc F$ is the $*$-extension, i.e., a clean extension. In particular, this directly shows that the \emph{abelian} category of cuspidal perverse sheaves is semisimple and generated by such clean extensions $\mc{IC} (\mc L)$. This can be used to describe the orthogonal to cuspidal constructible sheaves; see below. It would not be difficult to show that the open orbit is distinguished as well.  

The most important result we want to make use of is the explicit classification done earlier in \cite{LusztigIntersection}, which we can access only via the definition of cuspidality in \cite{LusztigGraded}. 
\end{remark}

\begin{corollary}\label{cor: ConstructibleCuspidalComplement}
Both the left- and right-orthogonal to the cuspidal constructible sheaves are generated by $\Eis_{P, \lambda_P}^{\mathrm{cons},!}\mc F$ with proper $P$ and cuspidal $\mc F$, or equivalently, by all non-cuspidal simple perverse sheaves.  

The stable $\infty$-category $D^b_c (\V_{\lambda})$ is generated under finite colimits and retracts by $\Eis_{P, \lambda_P}^{\mathrm{cons}, !}\mc F$, where $\mc F\in D^b_c(\V_{\lambda_M})$ is cuspidal and $P$ could possibly equal $G$. 
\end{corollary}

\begin{proof}
Using that cuspidal simple perverse sheaves are clean extensions, it is easy to describe both the left- and right-orthogonal: the restriction to the open orbit is in the subcategory generated by local systems corresponding to non-cuspidal irreducible representations of $\pi_0(S_\phi)$. It is generated by all the non-cuspidal simple perverse sheaves. The idempotent-complete stable $\infty$-category generated by $\Eis_{P, \lambda_P}^{\mathrm{cons},!}\mc F$ is also generated by all non-cuspidal simple perverse sheaves; this reduces to the case of $\V_{\lambda^\ur}$ and it is discussed in the proof of Proposition \ref{prop: OrthogonalCuspidalConstructibleSheaves}.

The last claim is now clear from the transitivity of Eisenstein series. 
\end{proof}

\subsection{Miscellany}
If $\Lambda=\C$, a $\mc D$-module on $\V_\lambda$ is called cuspidal if the corresponding constructible sheaf under the Riemann--Hilbert correspondence is so. It has characterizations using Eisenstein series or constant terms, and it makes sense for a general $\Lambda$. We will use the following criterion for the underlying $\mc D$-module of a filtered $\mc D$-module to be cuspidal. 

\begin{lemma}{\label{lem: filteredcusplem}}
Let $\widetilde{\mc M}$ be a filtered $\mc D$-module on $\V_{\lambda}$ that is compact. Suppose that for every proper semistandard Levi subgroup $\wh M$, every lift $\lambda_M$ and every cuspidal simple perverse sheaf $\mc N$ on $\V_{\lambda_M}$, a filtered lift $\widetilde{\mc N}$ of the $\mc D$-module corresponding to $\mc{N}$ is given, and 
\[
\Hom (\widetilde{\mc M}, \Eis_{P,\lambda_P}^{F\mc D, !} \widetilde{\mc N})\simeq 0
\]
holds as an object of $\QCoh (\A^1 /\Gm)$ for any proper $P$, $\lambda_P$. Then, the underlying $\mc D$-module of $\widetilde{\mc M}$ is cuspidal. 
\end{lemma}

\begin{proof}
This follows from \cite[Proposition 3.3.16]{Chen2023categorical}, the definition of cuspidal constructible sheaves and Corollary \ref{cor: ConstructibleCuspidalComplement}. 
\end{proof}

We will make use of the following lemma as well.

\begin{lemma}{\label{lem: Dmodulehomequality}}
Let $\mc V$ be a vector bundle on $\V_{\lambda}$ with trivial filtration. For any $\mc M\in F\mc{D}_c(\V_{\lambda})$, 
\[
\Hom_{F\mc{D}(\V_{\lambda})} (\ind \mc V, \mc M)\simeq \Hom_{\QCoh(\V_{\lambda})^{\textnormal{Fil}}} (\mc V, \res \mc M)\simeq \Gamma (\V_{\lambda}, \mc V^{\vee}\otimes \res \mc M),
\]
as a filtered complex, has a separated filtration. In particular, if its derived fiber at $B\Gm$ vanishes, then it vanishes as an object of $\QCoh (\A^1 /\Gm)$. 
\end{lemma}
\begin{proof}
The first identification is the adjunction, and the second uses that $\underline{\mathrm{Hom}} (\mc V, \res \mc M)\simeq \mc V^\vee \otimes \res \mc M$ as filtered quasi-coherent sheaves. The filtration on $\mc V^{\vee}\otimes \res \mc M$, and hence on $\Gamma (\V_{\lambda}, \mc V^{\vee}\otimes \res \mc M)$, is separated because the filtration of $\res \mc M$ is separated by the definition of $F\mc{D}_c(\V_{\lambda})$. It follows that every piece of the filtration vanishes if all the associated grades vanish. 
\end{proof}

\section{Cuspidal enhanced \texorpdfstring{$L$}{L}-parameters}{\label{s: amscomparison}}
In this section, we review some conjectures in the literature on the classical Langlands correspondence and enhanced $L$-parameters corresponding to supercuspidal representations.

Let $\phi=(\lambda, N)$ be a Frobenius-semisimple Weil--Deligne parameter for a quasi-split reductive group $G$ over $F$ and let $\tilde{\phi}\colon W_F \times \SL_2 \rightarrow \co{L}{G}$ be a lift to a parameter of $W_F \times \SL_2$.  Define $S^1_{\tilde{\phi}, \Sc}$ to be the pre-image in $\wh{G}_{\Sc}$ of the image of $S_{\tilde{\phi}}$ in $\wh{G}_{\ad}$. There is an obvious inclusion
\begin{equation*}
    S^1_{\tilde{\phi}, \Sc} \hookrightarrow Z_{S^1_{\tilde{\phi}|_{W_F}, \Sc}}(u),
\end{equation*}
inducing an isomorphism
\begin{equation}{\label{eqn: AMScomponentgroupiso}}
    \pi_0(S^1_{\tilde{\phi}, \Sc}) \cong \pi_0(Z_{S^1_{\tilde{\phi}|_{W_F}, \Sc}}(u)),
\end{equation}
where $u$ is the unipotent element of $\wh{G}$ (equivalently $\wh{G}_{\Sc}$) corresponding to the nilpotent $N$ in the definition of $\phi$. Further, we consider the group $\pi_0(Z_{(S^1_{\tilde{\phi}|_{W_F}, \Sc})^{\circ}}(u))$.

Now for each character $\chi \in X^{\ast}(Z(\wh{G}_{\Sc})^{\Gamma})$, choose a lift $\tilde{\chi}$ to $X^{\ast}(Z(\wh{G}_{\Sc}))$. An enhanced $L$-parameter in the sense of \cite{AMS1} (following \cite{ArthurLpacketnote}) is then a pair $(\tilde{\phi}, \rho)$ such that $\rho \in \Irr(\pi_0(S^1_{\tilde{\phi}, \Sc}), \tilde{\chi})$, meaning that $\rho \in \Irr(\pi_0(S^1_{\tilde{\phi}, \Sc}))$ and has central character $\tilde{\chi}$. Arthur conjectures a bijection
\begin{equation*}
    \coprod\limits_{z \in H^1(F, G_{\ad})} \Irr(G_z) \longleftrightarrow \coprod\limits_{\tilde{\phi} / \sim } \quad  \coprod\limits_{z \in H^1(F, G_{\ad})} \Irr(\pi_0(S^1_{\tilde{\phi}, \Sc}), \tilde{\chi}_z),
\end{equation*}
where $G_z$ ranges over inner twists of $G$ and we can recover $G_z$ from $\rho$ by restricting to $Z(\wh{G}_{\Sc})^{\Gamma}$ and applying the isomorphism
\begin{equation*}
    X^{\ast}(Z(\wh{G}_{\Sc})^{\Gamma}) \cong H^1(F, G_{\ad})
\end{equation*}
of Kottwitz \cite[Proposition 6.4]{Kot5}.

\begin{conjecture}[{\cite[Conjecture 6.10]{AMS1}}]{\label{conj: ams}}

An enhanced parameter $(\tilde{\phi}, \rho)$ corresponds to a supercuspidal representation on the automorphic side if and only if $\tilde{\phi}$ is essentially discrete and $\rho$, thought of as a representation of $\pi_0(Z_{S^1_{\tilde{\phi}|_{W_F}, \Sc}}(u))$ via \eqref{eqn: AMScomponentgroupiso}, satisfies that each irreducible factor of its restriction to $\pi_0(Z_{(S^1_{\tilde{\phi}|_{W_F}, \Sc})^{\circ}}(u))$ corresponds to a cuspidal local system on the unipotent orbit of $u$ inside $(S^1_{\tilde{\phi}|_{W_F}, \Sc})^{\circ}$ (cf.~\S\ref{ss: lusztiggraded}). We call $(\tilde{\phi}, \rho)$ a cuspidal enhanced parameter.
\end{conjecture}

Note that Conjecture \ref{conj: ams} has been verified in a number of cases. In \cite[Proposition 3.6]{Moussaouicusp}, Moussaoui shows that Conjecture \ref{conj: ams} is compatible with the description of supercuspidals in terms of alternating representations of essentially discrete parameters without gaps of Moeglin \cite[Theorem 2.5.1]{Moeglinmultiplicityone} for split $\Sp_{2n}, \SO_n$ and the description for $\GL_n$ of \cite{HT1} and \cite{Henniart}. In \cite[Example 6.11]{AMS1}, Conjecture \ref{conj: ams} is proven for inner twists of general linear groups. In \cite[Example 6.13]{AMS1}, Aubert--Moussaoui--Solleveld verify the conjecture for quasi-split unitary groups using the description of Moeglin \cite[Theorem 8.4.4]{Moeglindiscunitary} analogous to the split classical case.

\subsection{The conjecture for pure inner forms}

Conjecture \ref{conj: ams} is stated for all inner forms of $G$, whereas the results of this paper are stated for extended pure inner forms classified by $B(G)_{\bas}$, or more generally, for groups arising from the Kottwitz set $B(G)$. An intermediary notion is the set of pure inner forms classified by $H^1(F,G)$. When $G$ is semisimple, we have that $H^1(F,G)$ and $B(G)_{\bas}$ agree. When $G$ is adjoint, $H^1(F,G)$ and $H^1(F, G_{\ad})$ agree. To compare our results to Conjecture \ref{conj: ams}, we rephrase it in terms of pure inner forms. 

The classical local Langlands correspondence for pure inner forms is expected (see \cite{KalethaLLCnonQS}) to yield a bijection
\begin{equation*}
    \coprod\limits_{z \in H^1(F,G)} \Irr(G_z) \longleftrightarrow \coprod\limits_{\tilde{\phi} / \sim } \quad \coprod\limits_{z \in H^1(F,G)} \Irr(\pi_0(S_{\tilde{\phi}}), \chi_z).
\end{equation*}
\begin{conjecture}{\label{conj: amspure}}
    A pair $(\tilde{\phi}, \rho)$ with $\rho \in \Irr(\pi_0(S_{\tilde{\phi}}), \chi_z)$ corresponds to a supercuspidal representation on the automorphic side if and only if $\tilde{\phi}$ is essentially discrete and $(u, \rho')$ is a cuspidal pair in the sense of Lusztig, for each irreducible summand $\rho'$ in the restriction of $\rho$ to $\pi_0(Z_{S^{\circ}_{\tilde{\phi}|_{W_F}}}(u))$, where $\rho$ is thought of as a representation of $\pi_0(Z_{S_{\tilde{\phi}|_{W_F}}}(u))$. 
\end{conjecture}
Notice that $Z_{S_{\tilde{\phi}|_{W_F}}}(u)$ is in general different from $S_{\phi}$ (the centralizer in $\wh{G}$ of the Weil--Deligne parameter $\phi$). However, there is an isomorphism of component groups since we have 
\begin{equation*}
    S_{\tilde{\phi}} = Z_{S_{\tilde{\phi}|_{W_F}}}(\tilde{\phi}|_{\SL_2}),
\end{equation*}
and in any (potentially disconnected) reductive group, the component group of the centralizer of an $\SL_2$ homomorphism agrees with that of the centralizer of the associated nilpotent. We also use that $S_{\tilde{\phi}}$ is the reductive quotient of $S_{\phi}$ and hence they have the same component groups (for instance, \cite[Proposition 3.7]{BMImaiYoucis}). 

When $G$ is adjoint, clearly Conjectures \ref{conj: ams} and \ref{conj: amspure} are the same. We now argue that in a suitable sense, Conjecture \ref{conj: ams} restricted to pure inner forms identifies with Conjecture \ref{conj: amspure}.

Let $\chi \in X^{\ast}(\pi_0(Z(\wh{G})^{\Gamma}))$ and consider $\chi_{\Sc}$ the pullback to $Z(\wh{G}_{\Sc})^{\Gamma}$. Choose a lift $\tilde{\chi}_{\Sc}$ to $Z(\wh{G}_{\Sc})$.
\begin{lemma}
    There is a bijection
    \begin{equation*}
        \Irr(\pi_0(S^1_{\tilde{\phi}, \Sc}), \tilde{\chi}_{\Sc}) \cong \Irr(\pi_0(S_{\tilde{\phi}}), \chi),
    \end{equation*}
    that preserves the notion of cuspidality for each set.
\end{lemma}
\begin{proof}
There is no natural map relating $S^1_{\tilde{\phi}, \Sc}$ and $S_{\tilde{\phi}}$. So we compare representations of these groups through the group $S^+_{\tilde{\phi}}$ defined in \cite[\S4.2]{KalethaLLCnonQS} as the pre-image in $\wh{\bar{G}}$ of $S_{\tilde{\phi}}$ where $\bar{G} = G/Z(G_{\der})$. Recall we have $\bar{G} \cong G_{\ad} \times G^{\ab}$ and hence 
\begin{equation*}
    \wh{\bar{G}} \cong \wh{G}_{\Sc} \times Z(\wh{G})^{\circ}.
\end{equation*}
In \S 4.6 of loc.~cit.~Kaletha denotes by $Z(\wh{\bar{G}})^+$ the pre-image of $Z(\wh{G})^{\Gamma}$ in $\wh{\bar{G}}$. Kaletha then proves the product map induces an isomorphism
\begin{equation*}
    S^+_{\tilde{\phi}} \times^{Z(\wh{\bar{G}})^+} Z(\wh{G}_{\Sc}) \rightarrow S^1_{\tilde{\phi}, \Sc}.
\end{equation*}
Define $\chi_{\rig} \in X^{\ast}(\pi_0(Z(\wh{\bar{G}})^+))$ by pulling back $\tilde{\chi}_{\Sc}$ along the projection 
\begin{equation*}
    Z(\wh{\bar{G}})^+ \rightarrow Z(\wh{G}_{\Sc}).
\end{equation*}
Then given $\rho \in \Irr(\pi_0(S^+_{\tilde{\phi}}), \chi_{\rig})$, the representation $\rho \otimes \tilde{\chi}_{\Sc}$ of $S^+_{\tilde{\phi}} \times Z(\wh{G}_{\Sc})$ descends to a representation of $S^1_{\tilde{\phi}, \Sc}$ and this construction induces a bijection 
\begin{equation}{\label{eqn: Amsvsrig}}
    \Irr(\pi_0(S^+_{\tilde{\phi}}), \chi_{\rig}) \cong \Irr(\pi_0(S^1_{\tilde{\phi}, \Sc}), \tilde{\chi}_{\Sc}).
\end{equation}

Next we define $\chi'_{\rig}$ to be the pullback of $\chi$ along 
\begin{equation*}
    \pi_0(Z(\wh{\bar{G}})^+) \rightarrow \pi_0(Z(\wh{G})^{\Gamma}),
\end{equation*}
and relate $\Irr(\pi_0(S^+_{\tilde{\phi}}), \chi'_{\rig})$ and $\Irr(\pi_0(S_{\tilde{\phi}}), \chi)$. There is a natural surjection
\begin{equation*}
    \Pi\colon S^+_{\tilde{\phi}} \rightarrow S_{\tilde{\phi}}.
\end{equation*}
 Hence any $\rho \in \Irr(\pi_0(S_{\tilde{\phi}}), \chi)$ pulls back to an element of $\Irr(\pi_0(S^+_{\tilde{\phi}}), \chi'_{\rig})$. For the converse, a $\rho \in \Irr(\pi_0(S^+_{\tilde{\phi}}), \chi'_{\rig})$ has a central character that is trivial on the kernel of $\Pi$ and hence descends to an element of $\Irr(\pi_0(S_{\tilde{\phi}}), \chi)$. This gives a bijection
 \begin{equation}{\label{eqn: purevsrig}}
     \Irr(\pi_0(S^+_{\tilde{\phi}}), \chi'_{\rig}) \cong \Irr(\pi_0(S_{\tilde{\phi}}), \chi).
 \end{equation}

Next we construct a bijection
\begin{equation}{\label{eqn: relrigbij}}
    \Irr(\pi_0(S^+_{\tilde{\phi}}), \chi_{\rig}) \cong \Irr(\pi_0(S^+_{\tilde{\phi}}), \chi'_{\rig}).
\end{equation}
We have that $\chi_{\rig}$ and $\chi'_{\rig}$ agree when pulled back to $\pi_0(Z(\wh{G}_{\Sc})^{\Gamma})$ and hence correspond to the same underlying inner twist of $G$. As in \cite[\S6.3]{KallocalBGvsRig}, there exists a character $x_{\rel}$ of $\pi_0(S^+_{\tilde{\phi}})$ such that tensoring with $x_{\rel}$ induces \eqref{eqn: relrigbij}. Indeed, we take the groups $Z_n \subset Z(G)$ defined to be the pre-image of $(Z(G)/Z(G_{\der}))[n]$. We let $S^+_{\tilde{\phi}, n}$ be the pre-image of $S^+_{\tilde{\phi}}$ in $\wh{G / Z_n}$ and it is harmless to pass to this group since the characters in \eqref{eqn: relrigbij} are trivial on the kernel of $\wh{G/Z_n} \rightarrow \wh{\bar{G}}$. Then for sufficiently large $n$, Kaletha argues there is some $x \in Z^1(W, Z_n)$ ($W$ is Kaletha's rigid gerbe) such that $z'_{\rig} = xz_{\rig}$, where $z_{\rig}, z'_{\rig}$ are rigid inner twists associated to $\chi_{\rig}, \chi'_{\rig}$. Then Kaletha defines $x_{\rel}$ as the pullback of a corresponding character of $Z^1(\Gamma, \wh{Z_n})$ ($\wh{Z_n}$ signifies Cartier dual and this character exists by the duality in \cite[\S6.2]{KallocalBGvsRig}) via the map 
\begin{equation*}
    (-d)\colon \pi_0(S^+_{\tilde{\phi}, n}) \rightarrow Z^1(\Gamma, \wh{Z_n}), \qquad s \mapsto \tilde{\phi}(w_{\sigma})s^{-1}\tilde{\phi}(w_{\sigma})^{-1}s.
\end{equation*}
Note that the definition of $(-d)$ does not use the $\SL_2$ part of $\tilde{\phi}$ and hence also gives a map from $S^+_{\tilde{\phi}|_{W_F}, n}$ that is trivial on its identity component. It follows that tensoring with $x_{\rel}$ preserves cuspidality, since it is trivial on $Z_{(S^+_{\tilde{\phi}|_{W_F}, n})^{\circ}}(u)$.

It remains to argue that the bijections \eqref{eqn: Amsvsrig}, \eqref{eqn: purevsrig} preserve the notions of cuspidality. The point is that 
\begin{equation*}
    S^1_{\tilde{\phi}|_{W_F}, \Sc}, \qquad S^+_{\tilde{\phi}|_{W_F}}, \qquad S_{\tilde{\phi}|_{W_F}},
\end{equation*}
all agree up to central factors and component groups, and Lusztig's notion of cuspidal pairs is insensitive to this. More precisely, the comparisons above induce maps
\begin{equation*}
    S^{+, \circ}_{\tilde{\phi}|_{W_F}} \rightarrow S^{1,  \circ}_{\tilde{\phi}|_{W_F}, \Sc}, \qquad   S^{+, \circ}_{\tilde{\phi}|_{W_F}} \rightarrow S^{  \circ}_{\tilde{\phi}|_{W_F}},
\end{equation*}
which are surjective with central kernel. Pullback preserves unipotent varieties and parabolics and hence the definition of cuspidality in \cite[Definition 2.4]{LusztigIntersection}.
\end{proof}

\part{Cuspidal Coherent Sheaves}{\label{p: cuspidalsheaves}}

\section{First definitions and properties}{\label{s: cuspsheavesdef}}

We are now ready to define categories of cuspidal and basic sheaves. These are the key objects of study in this paper.

\subsection{Cuspidal and basic sheaves}{\label{ss: defncuspsheaves}}

\begin{definition}{\label{def: cuspdef}}
    We say a sheaf $\mc{F} \in \IndCoh(\Par_G)$ is \emph{right-cuspidal} (resp.~\emph{right-basic}) if 
    \begin{equation*}
       \Hom(\Eis^{\spec}_{P} \mc{F}', \mc{F})\simeq 0,
    \end{equation*}
    for each proper standard parabolic $P \subset G$ and each $\mc{F}' \in \IndCoh(\Par_M)^\chi$, for each $P$-dominant (resp.~strictly $P$-dominant) $\chi \in X^*(Z(\wh{M})^{\Gamma})$ relative to our fixed pinning. 

   One might expect an entirely analogous definition of left-cuspidal sheaves, but we found that definition to be weaker than what was needed. We instead give the following more restrictive definition. Let $\Eis_- \subset \IndCoh(\Par_G)$ (resp.~$\Eis^{\circ}_- \subset \IndCoh(\Par_G)$)  be the category generated under all colimits by $\Eis^{\spec}_{P^{\op}} \mc{F}'$ for each proper standard parabolic $P \subset G$ and each $\mc{F}' \in \IndCoh(\Par_M)^{\chi}$, for some $P$-dominant (resp.~strictly $P$-dominant) $\chi \in X^*(Z(\wh{M})^{\Gamma})$. Then we say a sheaf $\mc{F} \in \IndCoh(\Par_G)$ is \emph{left-cuspidal} (resp.~\emph{left-basic}) if 
    \begin{equation*}
        \Hom(\mc{F}, \Eis_-) \simeq 0 \qquad \left(\text{resp.~}\Hom(\mc{F}, \Eis^{\circ}_-) \simeq 0\right).
    \end{equation*}
    
 We denote these full subcategories of $\IndCoh(\Par_G)$ by $\LCusp(\Par_G)$ (resp.~$\LBas(\Par_G)$) and $\RCusp(\Par_G)$ (resp.~$\RBas(\Par_G)$). 
 We also write
 \begin{align*}
 \LCusp^{\Coh}(\Par_G)&=\LCusp(\Par_G) \cap \Coh(\Par_G), \\ 
\RCusp^{\Coh}(\Par_G)&=\RCusp(\Par_G) \cap \Coh(\Par_G), \\
\LBas^{\Coh}(\Par_G)&=\LBas(\Par_G) \cap \Coh(\Par_G), \\
 \RBas^{\Coh}(\Par_G)&=\RBas(\Par_G) \cap \Coh(\Par_G).     
 \end{align*}
\end{definition}

\begin{remark}{\label{rem: BasCohaltdef}}
 Note that this subtlety of definitions does not appear for $\LCusp^{\Coh}(\Par_G)$ and $\LBas^{\Coh}(\Par_G)$ since then we can move the colimit out of the $\Hom$.
 
 Indeed, recall from Remark \ref{rem: CTEisproperties} that the Eisenstein functor preserves coherent sheaves and commutes with colimits. Using this, we see the following: an $\mc F \in \Coh(\Par_G)$ is left-cuspidal/basic if and only if
\[
\Hom(\mc F, \Eis^\spec_{P^\op} \mc F')\simeq 0
\]
for $\mc F'\in \Coh(\Par_M)^\chi$. This is clearly necessary, and sufficient since any $\mc G \in \Eis_-$ is the colimit of some $\Eis^\spec_{P_i^{\op}}\mc G_i, \mc G_i\in \Coh(\Par_{M_i})^{\chi_i}$, and the vanishing for $\Eis^\spec_{P_i^\op}\mc G_i$ implies the vanishing for $\mc G$:
\[
\Hom(\mc F,  \colim \Eis^\spec_{P_i^\op} \mc  G_i)\simeq \colim \Hom(\mc F, \Eis^\spec_{P_i^\op} \mc  G_i)\simeq 0. 
\]
If $\mc F\in \IndCoh(\Par_G)$, then a similar claim holds for being right-cuspidal/basic as we have
\[
\Hom(\Eis^\spec_{P} \colim\mc  G_i, \mc F)\simeq  \Hom(\colim\Eis^\spec_{P}\mc  G_i, \mc F)\simeq \lim \Hom(\Eis^\spec_{P} \mc  G_i, \mc F). 
\]
\end{remark}

\begin{lemma}{\label{lem: cuspconstanttermdef}}
    An $\mc{F} \in \IndCoh(\Par_G)$ is right-cuspidal (resp.~right-basic) if and only if 
    \[
    \CT^{\spec, \chi}_{P} \mc{F} \simeq 0
    \]
    for each $P$-dominant (resp.~strictly $P$-dominant) $\chi \in X^*(Z(\wh{M})^{\Gamma})$ and proper standard parabolic subgroup $P \subset G$. 
\end{lemma}
\begin{proof}
By adjointness, we have that $\mc{F}$ is right-cuspidal if and only if \[
\Hom( \mc{F}', \CT^{\spec, \chi}_{P} \mc{F})\simeq 0\] for each proper standard parabolic subgroup $P \subset G$, each $P$-dominant $\chi \in X^*(Z(\wh{M})^{\Gamma})$, and $\mc{F}' \in \IndCoh(\Par_M)^{\chi}$. Equivalently, $\CT^{\spec, \chi}_{P} \mc{F} \simeq 0$. The case of right-basic sheaves is the same. 
\end{proof}
\begin{remark}{\label{rem: cuspbasicproperties}}
\begin{enumerate}
 \item  It appears that Definition \ref{def: cuspdef} depends on our choice of $\wh{B}$, but in fact it does not. We could alternatively say $\mc{F} \in \IndCoh(\Par_G)$ is right-cuspidal (resp.~right-basic) if \begin{equation*}
        \Hom(\Eis^{\spec, \chi}_{P} \mc{F}', \mc{F}) \simeq 0,
    \end{equation*}
    for every proper semistandard $P$ and any $P$-dominant  (resp.~strictly $P$-dominant) $\chi$. This certainly implies the usual definition of right-cuspidal (resp.~right-basic) but they are equivalent by conjugating. 
    
    Similarly, we can further work with more general parabolic subgroups, but we will only use those containing $\wh T$. 
\item We have $\mc{F} \in \RCusp(\Par_G)$ if and only if $\CT^{\spec, \chi}_P\mc{F}\simeq 0$ for each $P$-dominant $\chi$ and \emph{maximal} standard parabolic $P$. Indeed, one direction is clear and the other follows from Lemma \ref{lem: eisctcompositions}. A similar statement holds for right-basic sheaves. 
\end{enumerate}
\end{remark}

\begin{lemma}{\label{lem: ParGdomswitch}}
    A sheaf $\mc{F} \in \IndCoh(\Par_G)$ is right-basic if and only if for each proper semistandard $P \subset G$, and for each defining cocharacter $\nu \in X_{\ast}(Z(\wh{M})^{\Gamma} \cap \wh{G}_{\der})$ of $\wh{P}$, one has
    \begin{equation*}
        \Hom(\Eis^\spec_P \mc{F}', \mc{F}) = 0 
    \end{equation*}
    for each $\mc{F}'$ of $Z(\wh M)^{\Gamma}$-weight $\chi$ such that $\langle \nu, \chi \rangle > 0$.  A similar statement holds for right-cuspidal sheaves. For the ``if'' statement, one need only check a single defining $\nu$ for each $P$.
\end{lemma}
\begin{proof}
    By conjugating, we may assume $P$ is standard. Then this condition implies right-basic since $\chi$ strictly $P$-dominant implies $\langle \nu, \chi \rangle > 0$ for any such defining $\nu$. On the other hand, suppose that $\mc{F}$ is right-basic and we have fixed such a defining $\nu$ and $\mc{F}'$ of $Z(\wh M)^{\Gamma}$-weight $\chi$ such that $\langle \nu, \chi \rangle > 0$. Then we can write
    \begin{equation*}
        \nu = \sum\limits_Q c_Q \nu_Q,
    \end{equation*}
    where $Q \supset P$ runs over maximal standard parabolic subgroups of $G$, $\nu_Q$ is the sum of the $\Gamma$-orbit of fundamental coweights for the omitted simple roots of the Levi $\wh M_Q$, and $c_Q$ is a positive rational number. Then $\langle \nu, \chi \rangle > 0$ implies that $\langle \nu_Q, \chi\rangle > 0$ for some $Q$.
    
    Now let $P_{M_Q} = P \cap M_Q$ be the corresponding parabolic subgroup of $M_Q$. We have
    \begin{equation*}
        \Hom(\Eis^\spec_P \mc{F}', \mc{F}) \simeq \Hom(\Eis^\spec_Q \Eis^\spec_{P_{M_Q}} \mc{F}', \mc{F}). 
    \end{equation*}
To prove the vanishing, it suffices to show that the $Z(\wh{M}_Q)^\Gamma$-weight $\chi_Q = \chi|_{Z(\wh{M}_Q)^{\Gamma}}$ of $\Eis^\spec_{P_{M_Q}} \mc{F}'$ is strictly $Q$-dominant. Since $Q$ is maximal, this is implied by
\begin{equation*}
    0 < \langle \nu_Q, \chi \rangle = \langle \nu_Q, \chi_Q \rangle.
\end{equation*}
\end{proof}

\begin{lemma}{\label{lem: leftrightcuspduality}}
    Let $\mc F\in \Coh(\Par_G) $. We have $\mc{F} \in \LCusp^{\Coh}(\Par_G)$ if and only if $\D_{\GS}\mc{F} 
    \in \RCusp^{\Coh}(\Par_G)$ if and only if $c^* \D_{\GS}\mc{F} 
    \in \RCusp^{\Coh}(\Par_G)$, for $c$ the Chevalley involution associated with the fixed pinning. 
\end{lemma}
\begin{proof}
For the first part, we have
    \begin{align*}
        \Hom(\mc{F}, \Eis^{\spec, \chi}_{P^{\op}} \mc{F}') &\simeq \Hom(\D_{\GS}(\Eis^{\spec, \chi}_{P^{\op}} \mc{F}'), \D_{\GS}\mc{F}) \\ & \simeq \Hom(\Eis^{\spec, -\chi}_{P^{\op}} \D_{\GS}\mc{F}', \D_{\GS}\mc{F})\\
        & \simeq \Hom(\Eis^{\spec, \chi'}_{P'} (\D_{\GS}\mc{F}')', \D_{\GS}\mc{F}),
    \end{align*}
    where the last equivalence is by conjugation, $P'$ is some standard parabolic, and $(\D_{\GS}\mc{F}')'$ is the conjugate of $\D_{\GS}\mc{F}'$. The second part then follows by applying $c^{\ast}$.
\end{proof}

\begin{lemma}\label{lem:RCusp(ParG)Presentable}
\begin{enumerate}
    \item The stable $\infty$-categories $\RCusp(\Par_G)$, $\RBas(\Par_G)$ are presentable, compactly generated, and closed under limits and colimits in $\IndCoh(\Par_G)$. The subcategories $\RCusp^\Coh(\Par_G)$, $\RBas^\Coh(\Par_G)$ are closed under finite colimits and retracts. 
    \item The inclusion $\RCusp(\Par_G)\to \IndCoh(\Par_G)$ admits a left adjoint $\mc F\mapsto \mc F^\rc$, and the fiber of the unit $\mc F \to \mc F^{\rc}$
    belongs to the subcategory generated under colimits by $\Eis_{P}^{\spec} \mc{G}$ with $P$-dominant $\mc{G}$ and proper standard parabolic $P$. 
    \item The functor $\mc{F} \mapsto \mc F^\rc$ preserves compact objects. 
\end{enumerate}
\end{lemma}

\begin{proof}
First, the subcategory $\RCusp(\Par_G)$ is closed under limits in $\IndCoh(\Par_G)$ and also under filtered colimits since right-cuspidality can be tested on coherent sheaves as in Remark \ref{rem: BasCohaltdef}. It follows that it is also closed under colimits as we are working with stable $\infty$-categories. The same works for $\RBas(\Par_G)$. The claim for coherent sheaves in (1) follows immediately. 

The presentable stable $\infty$-category $\IndCoh(\Par_G)$ is compactly generated by $\Coh(\Par_G)$ (and the same holds for Levi subgroups $M$). It follows that $\Eis_{P}^{\spec} \mc{G}$ with $P$-dominant $\mc{G}$ and proper standard parabolic $P$ generates under colimits a presentable $\infty$-subcategory of $\IndCoh(\Par_G)$. Therefore, its right orthogonal $\RCusp(\Par_G)$ satisfies the desired properties. More precisely, by \cite[Proposition 1.4.4.11]{LurieHA}, there exists a $t$-structure on $\IndCoh(\Par_G)$ such that the presentable subcategory generated by $\Eis_{P}^{\spec} \mc{G}$ with $P$-dominant $\mc{G}$ corresponds to $\IndCoh(\Par_G)_{\geq 0}$.  Its right-orthogonal $\IndCoh(\Par_G)_{<0}$ is $\RCusp(\Par_G)$. This proves (2), the truncation $\mc F\mapsto \mc F_{<0}$ being the left adjoint. Since $\RCusp$ is closed under filtered colimits in $\IndCoh(\Par_G)$, the left adjoint $(-)^\rc$ preserves compact objects as in (3). 
\end{proof}

In particular, since the left adjoint commutes with colimits and any sheaf in $\IndCoh(\Par_G)$ is a colimit of coherent sheaves, it follows that $\RCusp(\Par_G)$ is generated under colimits by sheaves of the form $\mc{F}^{\rc}$ for $\mc{F} \in \Coh(\Par_G)$. It is a priori not clear that $\mc F^\rc$ is still coherent. This is equivalent to the semiorthogonal decomposition in the proof inducing one on $\Coh(\Par_G)$. 

We eventually prove that $\LCusp(\Par_G) = \RCusp(\Par_G)$ (Theorem \ref{thm:LCusp=RCusp}), and they are indeed compactly generated by cuspidal \emph{coherent} sheaves (Theorem \ref{thm: Generators}). For instance, it seems far from clear from the definition that $\LCusp(\Par_G)$ is closed under limits.

\subsection{Reduction to completed fibers and completed Vogan stacks}{\label{ss: reductiontocompletedfiber}}
In this subsection we give techniques for reducing the study of cuspidal sheaves to completed fibers at semisimple $L$-parameters.

Fix a semisimple $L$-parameter $\lambda\colon W_F\to \wh G$; its associated point of $X_G(\Lambda)$ is also denoted by $\lambda$. 

\begin{definition}{\label{def: cuspcompletion}}
We say that $\mc F\in \Ind\wh{\Coh}{}^!(\Par^{\wedge}_{G, \lambda})$ (resp.~$\mc F\in \IndCoh(\Par^{\wedge}_{G, \lambda})$) is right-cuspidal if we have, for all proper \emph{semistandard} $P\subset G$ and every lift $\lambda_P$ of $\lambda$,  
\[
\Hom(\Eis^{\spec, \wedge}_{P, \lambda_P} \mc{G}, \mc{F})\simeq 0
\]
for each $\mc{G} \in \Ind\wh{\Coh}{}^{!}(\Par^{\wedge}_{M, \lambda_M})^\chi$ (resp.~$\mc{G} \in \IndCoh(\Par^{\wedge}_{M, \lambda_M})^\chi$), or equivalently for each $\mc{G} \in \wh{\Coh}{}^{!}(\Par^{\wedge}_{M, \lambda_M})^\chi$ (resp.~$\mc{G} \in \Coh(\Par^{\wedge}_{M, \lambda_M})^\chi$), for $P$-dominant $\chi \in X^*(Z(\wh{M})^{\Gamma})$. 

Similarly, we say that $\mc F\in \Ind\wh{\Coh}{}^!(\wh{\mc L}(\V_\lambda))$ is right-cuspidal if we have, for all proper semistandard $P\subset G$ and every lift $\lambda_P$ of $\lambda$,  
\[
\Hom(\Eis^{\wh{\mc L}}_{P, \lambda_P} \mc{G}, \mc{F})\simeq 0
\]
for each $\mc{G} \in \Ind\wh{\Coh}{}^{!}(\wh{\mc L}(\V_{\lambda_M}))^\chi$, or equivalently for each $\mc{G} \in \wh{\Coh}{}^{!}(\wh{\mc L}(\V_{\lambda_M}))^\chi$, for $P$-dominant  $\chi \in X^*(Z(\wh{M})^{\Gamma})$. We define right-cuspidal sheaves $\mc F\in \Ind\wh{\Coh}{}^!(\wh{\T}\V_\lambda[-1])$ by using $\Eis^{\wh{\T}[-1]}_{P,\lambda_P}$.  

We also make the analogous definitions of right-basic sheaves in all these situations.
\end{definition}

\begin{remark}    
Note that $\mc F\in \wh{\Coh}{}^! (\Par_{G,\lambda}^\wedge)\subset \IndCoh(\Par_{G,\lambda}^\wedge)$ is right-cuspidal if and only if it is so as $\mc F\in \Ind\wh{\Coh}{}^!(\Par_{G,\lambda}^\wedge)$. 
\end{remark}

\begin{lemma}{\label{lem: ihatpreservescusp}}
\begin{enumerate} 
    \item The $!$-restriction and its right adjoint identify the right-cuspidal and right-basic objects of $\Ind\wh{\Coh}{}^! (\Par^\wedge_{G,\lambda})$ and those of $\Ind\wh{\Coh}{}^! (\wh{\mc L}(\V_\lambda))$. 
    \item The functors
\begin{align*}
\wh{i}_{\lambda}^!\colon \Coh(\Par_G) &\to \wh{\Coh}{}^! (\Par^\wedge_{G,\lambda}), \quad \IndCoh(\Par_G) \to \IndCoh (\Par^\wedge_{G,\lambda}) \\
\wh{i}_{\lambda, \Fss}^!\colon \Coh(\Par_G) &\to \wh{\Coh}{}^! (\wh{\mc L}(\V_\lambda))
\end{align*}
preserve right-cuspidal and right-basic objects. 
    \item The right adjoint functors of $\wh{i}_{\lambda}^!$ and $\wh{i}_{\lambda, \Fss}^!$ on $\Ind\wh{\Coh}{}^!$ preserve right-cuspidal and right-basic objects. In particular, if the $*$-pushforward of $\mc F\in \Coh(\V_\lambda)$ is right-cuspidal/basic on $\Par^\wedge_{G,\lambda}$, equivalently on $\wh{\mc L}(\V_\lambda)$, then so is $i_{\lambda*}\mc F\in \Coh(\Par_G)$.   
    \item The right adjoint functor of 
    \[
    \Ind\wh{\Coh}{}^!(\Par^\wedge_{G,\lambda}) \to \IndCoh(\Par^\wedge_{G,\lambda})
    \]
    preserves right-cuspidal and right-basic objects. 
    \item The equivalence 
    \[
    \Ind\wh{\Coh}{}^! (\wh{\mc L}(\V_\lambda))\simeq
    \Ind\wh{\Coh}{}^! (\wh{\T}\V_\lambda[-1]) 
    \]
    preserves the right-cuspidal and right-basic objects. 
\end{enumerate}
\end{lemma}

\begin{remark}
Part (2) is strengthened in Remark \ref{rem: ihatpreservescuspIndCoh}, but it uses the main classification results.     
\end{remark}

\begin{proof}
The statement (1) follows from Proposition \ref{prop: unipvsformal} and Corollary \ref{cor: unipvsformalEis}; (2) and (3) follow from Proposition \ref{prop: EisComparisonCompletedFiber}, Corollary \ref{cor: EisComparisonCompletedFiber}, and Lemma \ref{lem:right_adjoint_from_hatCoh}; (4) holds since right-cuspidality/basicness can be tested on maps from Eisenstein series of $\wh{\Coh}{}^!(\Par^\wedge_{M,\lambda_M})$; and (5) follows from Proposition \ref{prop: loopvstangentEis}. 
\end{proof}

\subsection{Reduction to cotangent stacks}\label{ss: reductiontocotangent}
Next, we study the relation with linear Koszul duality. 
We define right-cuspidal and right-basic ind-coherent sheaves on $\T^* \V_\lambda$ analogously. 
In the graded setting, we define right-cuspidal and right-basic sheaves in the same way using the graded versions of Eisenstein series. 

\begin{lemma}\label{lem: HatCohForgettingGm}
There are natural equivalences
\begin{align*}
\wh{\Coh}{}^! (\wh{\T}\V_\lambda[-1]) &\simeq \Coh(\T\V_\lambda[-1] \times_{\Lie H_\lambda \sslash H_\lambda} \Spec \Lambda'), \\
\wh{\Coh}{}^! (\wh{\T}\V_\lambda[-1] /\Gm) &\simeq \Coh(\T\V_\lambda[-1] /\Gm),
\end{align*}
where $\Lambda'$ denotes the completed local ring of $\Lie H_\lambda \sslash H_\lambda$ at the origin, and the $*$-pullback corresponds to $*$-pullback, i.e., forgetting grading and base change further to $\Lambda'$. In particular, the image of the $*$-pullback generates  $\Ind\wh{\Coh}{}^! (\wh{\T}\V_\lambda[-1])$ under colimits, and the $*$-pushforward, the right adjoint of $*$-pullback, is conservative. Moreover, if $\mc F_1\in \wh{\Coh}{}^!(\wh{\T}\V_\lambda[-1]/\Gm)$ and $\mc F_2\in \Ind\wh{\Coh}{}^!(\wh{\T}\V_\lambda[-1]/\Gm)$, then
\[
\Hom_{\wh{\T}\V_\lambda[-1]}(\mc F_1, \mc F_2)\simeq (\bigoplus_{w\in \Z} \Hom(\mc F_1, \mc F_2 \langle w \rangle))\otimes_{\mc O(\Lie H_\lambda \sslash H_\lambda)}\Lambda';
\]
if $\mc F_2\in \wh{\Coh}{}^!(\wh{\T}\V_\lambda[-1]/\Gm)$, it is also equivalent to
\[
\prod_{w\in \Z} \Hom(\mc F_1, \mc F_2 \langle w \rangle).
\]
\end{lemma}

\begin{proof}
This follows from Corollary \ref{cor: cohhatcompletion} and Theorem \ref{thm: completeness_along_closedpoint} (which is applicable because the underlying classical stack of $\T\V_\lambda[-1]$ and $\Lie H_\lambda /H_\lambda$ have the same GIT quotients by Lemma \ref{lem: ReductionToGradedUnipotent}). The generation claim holds, since $\Spec \Lambda'\to \Lie H_\lambda \sslash H_\lambda$ is regular; it is not of finite type, but it does not cause any issue. The last part reduces to the case where $\mc F_2\in \wh{\Coh}{}^!(\wh{\T}\V_\lambda[-1]/\Gm)$ as both sides commute with filtered colimits. It holds, since $\oplus_w\Hom(\mc F_1, \mc F_2\langle w\rangle)$ has finitely generated cohomology over $\mc O(\Lie H_\lambda \sslash H_\lambda)$, and $\Lambda'$ has the same description in terms of $\mc O(\Lie H_\lambda \sslash H_\lambda)$. 
\end{proof}

\begin{lemma}{\label{lem: cusppullpush}}
A sheaf $\mc F\in \Ind\wh{\Coh}{}^!(\wh{\T}\V_\lambda[-1] /\Gm)$ is right-cuspidal (resp.~right-basic) if and only if its $*$-pullback to $\wh{\T}\V_\lambda[-1]$ is right-cuspidal (resp.~right-basic). If $\mc F^{\mathrm{ungr}}\in \Ind\wh{\Coh}{}^! (\wh{\T}\V_\lambda[-1] )$ is right-cuspidal (resp.~right-basic), then so is its $*$-pushforward to $\Ind\wh{\Coh}{}^! (\wh{\T}\V_\lambda[-1] /\Gm)$. The corresponding claims hold for $\IndCoh(\T^* \V_\lambda/\Gm)$.      
\end{lemma}

\begin{proof}
If $\mc F^{\mathrm{ungr}}$ is right-cuspidal, then the right-cuspidality of the $*$-pushforward of $\mc F^{\mathrm{ungr}}$ follows from the adjunction as Eisenstein series are compatible with forgetting gradings. 

For the first claim, note that the $*$-push-pull of $\mc F\in \wh{\Coh}{}^! (\wh{\T}\V_\lambda[-1] /\Gm)$ is equivalent to $\prod_{w\in \Z}\mc F\langle w\rangle $ by Lemma \ref{lem: HatCohForgettingGm}, and $\mc F$ is still a retract of the $*$-push-pull in general as it commutes with filtered colimits. If the $*$-pullback of $\mc F$ is right-cuspidal, then the right-cuspidality of $\mc F$ follows from the second claim proven in the previous paragraph. If $\mc F$ is right-cuspidal, the $*$-push-pull is right-cuspidal by Lemma \ref{lem: HatCohForgettingGm}. Since $\Ind\wh{\Coh}{}^! (\wh{\T}\V_{\lambda_M}[-1])$ is generated under colimits by the image of the $*$-pullback of graded sheaves by Lemma \ref{lem: HatCohForgettingGm}, it suffices to use the Eisenstein series of such sheaves to check the right-cuspidality, and the right-cuspidality of the $*$-pullback of $\mc F$ follows again from the adjunction. 

Similar arguments work for the right-basic sheaves and the case of $\T^*\V_\lambda$. 
\end{proof}

\begin{lemma}{\label{lem: koszulcusp}}
Under linear Koszul duality, 
\[
\Ind\wh{\Coh}{}^!(\wh{\T}\V_\lambda[-1]/\Gm) \simeq
\IndCoh(\T^*\V_\lambda/\Gm), 
\]
right-cuspidal and right-basic objects are preserved. 
\end{lemma}

\begin{proof}
    This follows from Proposition \ref{prop: eiscotangentcompatibility}.
\end{proof}

Therefore, it would suffice to classify the right-cuspidal ind-coherent sheaves $\RCusp(\T^*\V_\lambda)$ on $\T^*\V_\lambda$, which is a presentable stable $\infty$-category.

\subsection{Sheaves on cotangent stacks}

Let us record the following explicitly:

\begin{lemma}\label{lem:RCusp(TV)Presentable}
\begin{enumerate}
    \item The stable $\infty$-categories $\RCusp(\T^*\V_\lambda)$, $\RBas(\T^*\V_\lambda)$ are presentable, compactly generated, and closed under limits and colimits in $\IndCoh(\T^*\V_\lambda)$. The subcategories $\RCusp^\Coh(\T^*\V_\lambda)$, $\RBas^\Coh(\T^*\V_\lambda)$ are closed under finite colimits and retracts. 
    \item The inclusion $\RCusp(\T^*\V_\lambda)\to \IndCoh(\T^*\V_\lambda)$ admits a left adjoint $\mc F\mapsto \mc F^\rc$, and the fiber of the unit $\mc F\to \mc F^\rc$
    belongs to the subcategory generated under colimits by $\Eis_{P}^{\T^*} \mc{G}$ with $P$-dominant $\mc{G}$ and proper semistandard parabolic $P$. 
    \item The functor $\mc{F} \mapsto \mc F^\rc$ preserves compact objects. 
\end{enumerate}
\end{lemma}

\begin{proof}
The proof of Lemma \ref{lem:RCusp(ParG)Presentable} works. 
\end{proof}

We also define left-cuspidal and left-basic ind-coherent sheaves on $\T^*\V_\lambda$, and the categories $\LCusp(\T^*\V_\lambda)$ and $\LBas(\T^*\V_\lambda)$ by adapting Definition \ref{def: cuspdef}.

\begin{lemma}\label{lem: Cotangentleftrightcuspidality}
Let $\mc F\in \Coh(\T^* \V_{\lambda})$. Then $\mc F$ is right-cuspidal (resp.~right-basic) if and only if $\D'_{\GS}\mc F$ is left-cuspidal (resp.~left-basic). 
\end{lemma}

\begin{proof}
    This follows from Proposition \ref{prop: twistedGSEiscompatibility}.
\end{proof}

\begin{lemma}{\label{lem: cotangentdomswitch}}
    A sheaf $\mc{F} \in \Coh(\T^{\ast}\V_{\lambda})$ is right-basic if and only if for each proper semistandard $P \subset G$, and for each defining cocharacter $\nu \in X_{\ast}(Z(\wh{M})^{\Gamma} \cap \wh{G}_{\der})$ of $\wh{P}$, one has
    \begin{equation*}
        \Hom(\Eis^{\T^{\ast}}_P \mc{F}', \mc{F}) = 0 
    \end{equation*}
    for each $\mc{F}' \in \IndCoh(\T^{\ast} \V_{\lambda_M})$ of $Z(\wh M)^{\Gamma}$-weight $\chi$ such that $\langle \nu, \chi \rangle > 0$.  A similar statement holds for each combination of left/right and cuspidal/basic sheaves. For the ``if'' statement, one need only check a single defining $\nu$ for each $P$.
\end{lemma}

\begin{proof}
    The proof is the same as that of Lemma \ref{lem: ParGdomswitch} using the analogue of Lemma \ref{lem: eisctcompositions} (proven in Corollary \ref{cor: eiscompositioncotan}). The proofs of the left versions are analogous since we are restricting to $\mc{F}$ coherent, where the definitions of left/right cuspidality are analogous.
\end{proof}

\section{Statements of main classification results}{\label{s: statementsmain}}

In this section we formulate the classification of cuspidal sheaves.

\subsection{Classification: statements}

\begin{theorem}\label{thm:LCusp=RCusp}
If $\mc F\in \Coh(\Par_G)$ is right-cuspidal, then $\D_{\GS} \mc F$ is also right-cuspidal. In particular, $\mc F$ is also left-cuspidal, i.e., 
\[
\LCusp^{\Coh}(\Par_G)=\RCusp^{\Coh}(\Par_G).
\]
In fact, 
\[
\LCusp(\Par_G)=\RCusp(\Par_G).
\]
\end{theorem}

We call such sheaves \emph{cuspidal}, and for simplicity formulate our results below for them. However, many of these results are in practice proven first for right-cuspidal sheaves and Theorem \ref{thm:LCusp=RCusp} is deduced as a corollary; see Corollary \ref{cor: LCusp=RCusp}.  

Cuspidal sheaves are supported on essentially discrete Vogan varieties.

\begin{theorem}\label{thm: discretesupport}
If $\mc F\in \Ind\wh{\Coh}{}^!(\wh{\T}\V_\lambda[-1])\simeq \Ind\wh{\Coh}{}^!(\mc L^u(\V_\lambda))$ is cuspidal and nonzero, then $\V_\lambda$ is essentially discrete. The same holds for $\mc F\in \IndCoh(\mc L^u(\V_\lambda))$.
\end{theorem}
\begin{proof}
    The second part follows from the first part and Lemma \ref{lem: ihatpreservescusp}, since the right adjoint of the essentially surjective functor
    \[
    \Ind\wh\Coh{}^!(\Par^\wedge_{G,\lambda})\to \IndCoh(\Par^\wedge_{G,\lambda})
    \]
    is conservative. We reduce the first part to Lusztig's graded setting in \S\ref{ss: reductiondiscsuppcvbclassgraded} and the simply-connected case in \S\ref{ss: reductiontosimplyconnecteddiscsuppandCVBclass}. The result then follows via Koszul duality from the result for cotangent stacks (Lemmas \ref{lem: HatCohForgettingGm}, \ref{lem: cusppullpush}, and Theorem \ref{thm: ClassificationCuspidalCotangent}) and the left-cuspidal version proven in Proposition \ref{prop: cuspreductiontofinite}.
\end{proof}

\begin{definition}{\label{def: cuspgeneral}}
Let $G$ be reductive and $\phi$ an essentially discrete $L$-parameter with semisimplification $\lambda$. A representation $\rho$ of $S_{\phi}$ is \emph{cuspidal} if it decomposes into Lusztig's cuspidal representations (as in \S \ref{ss: lusztiggraded}) as a representation of $S_\phi \cap I_{\lambda,\der}$. 
\end{definition}

\begin{remark}
\begin{enumerate}
    \item If $\rho$ is irreducible, the central character of $\rho$ restricted to $Z(\wh G)^\Gamma \subset S_\phi$ determines an element $b_\rho\in X^*(Z(\wh G)^\Gamma)$. 
    \item If $\rho$ is inflated from an irreducible representation of $\pi_0 (S_\phi)$, i.e., in the case $b_\rho$ gives rise to a pure inner form of $G$, this recovers the previous notion of cuspidality in Definition \ref{def:purecuspidal}. 
    \item By Lusztig's classification, the restriction of a cuspidal $\rho$ to the finite group $S_\phi \cap I_{\lambda,\der}$ is determined by the restriction to $Z(I_{\lambda,\der})$. Therefore, the latter also determines, together with $b_\rho$, the restriction to $Z(\wh G)^{\Gamma} (S_{\phi}\cap I_{\lambda, \der})$, and in particular to $S_{\phi}\cap I_\lambda$. 
\end{enumerate}
\end{remark}

The main classification theorem is the following

\begin{theorem}\label{thm: CVBclassification}
\begin{enumerate}
    \item For any essentially discrete $L$-parameter $\phi$ with semisimplification $\lambda$ and a cuspidal irreducible representation $\rho$ of the centralizer $S_\phi$ of $\phi$, there exists a unique vector bundle $\mc V_\rho$ on $\V_\lambda$ whose restriction to $BS_\phi$ corresponds to $\rho$ such that its $*$-pushforward to $\Par_G$ is cuspidal. 
    \item Suppose $\rho, \rho'$ are irreducible cuspidal representations of $S_{\phi}, S_{\phi'}$ respectively. 
    \begin{enumerate}
    \item $\Hom_{\Par_G} (\mathcal{V}_{\rho}, \mathcal{V}_{\rho'})\simeq 0$ unless $\phi$ and $\phi'$ are $\wh G$-conjugate and, under the induced isomorphism $S_{\phi}\cong S_{\phi'}$ by some, equivalently any, conjugation, $\rho$ corresponds to $\rho'$. 
    \item The $\Lambda$-algebra $\mathrm{Ext}_{\Par_G} (\mc V_{\rho}, \mc{V}_{\rho})$ identifies with the exterior algebra generated by $\Lie Z(\wh G)^\Gamma$ in degree 1. 
    \end{enumerate}
    More precisely, the restriction from $\wh{\T}\V_{\lambda} [-1]$ to $\wh {\T}BS_\phi [-1]$ is fully faithful on cuspidal sheaves in $\wh{\Coh}{}^!(\wh{\T}\V_{\lambda} [-1])$. 
    \item The $\ast$-pushforwards of $\mc{V}_{\rho}$ associated to $\lambda$ generate the category of cuspidal sheaves set-theoretically supported on $\Par_{G, \lambda}$ under colimits. Moreover, if $G$ is semisimple, then the category of cuspidal ind-coherent sheaves is semisimple with simple objects given up to shift by $\ast$-pushforwards of the $\mc V_{\rho}$.
\end{enumerate}
\end{theorem}
\begin{proof}
    We reduce it to Lusztig's graded setting by \S \ref{ss: reductiondiscsuppcvbclassgraded} and further to the simply-connected case in \S\ref{ss: reductiontosimplyconnecteddiscsuppandCVBclass}. 
    The first statement of (3) in the simply-connected (or semisimple) case is proven in Corollary \ref{cor: CohCompletionSOD} (whose main input is Theorem \ref{thm: ClassificationCuspidalCotangent}). 
    By Lemma \ref{lem: ihatpreservescusp}, the cuspidality of a vector bundle $\mc V$ on $\V_\lambda$ is independent of the setting we consider. 
    Then (1) and (2) for cuspidal vector bundles follow by Koszul duality from the analogous results for cotangent stacks (Theorem \ref{thm: ClassificationCuspidalCotangent}) in the simply-connected graded case, since $\Ind\wh{\Coh}{}^!(\wh{\T}\V_\lambda[-1]
    )$ corresponds to $\IndCoh(\T^*\V_\lambda)$ after passing to $\Gm$-equivariant sheaves on both sides and the Koszul duals of $\mc V_\rho$ are vector bundles on $\T^*\V_\lambda$ (in general, they are vector bundles on the underlying classical stack, but here the derived structure is trivial). Now the general fully faithful statement in (2) follows from the case of cuspidal vector bundles, the first part of (3), and Corollary \ref{cor: SODhatCoh}. The final statement of (3) follows because for any cuspidal $\mc{F} \in \IndCoh(\Par_G)$, the $!$-restriction $\wh{i}^!_{\lambda} \mc{F} \in \Ind\wh{\Coh}{}^!(\mc{L}^u(\V_\lambda))$ is zero unless $\V_\lambda$ is discrete by Theorem \ref{thm: discretesupport} and by Corollary \ref{cor: SODhatCoh} and Remark \ref{rem: ihatpreservescuspIndCoh}. Therefore $\mc F$ vanishes outside the union of closed substacks $\Par_{G,\lambda}$ for $\V_\lambda$ being discrete, and $\mc F$ is set-theoretically supported on the union.\footnote{This argument works only for $\Ind\wh{\Coh}{}^!$ and fails for $\IndCoh$.} Then, again by Corollary \ref{cor: SODhatCoh} and Remark \ref{rem: ihatpreservescuspIndCoh}, it follows that $\wh{i}^!_{\lambda} \mc{F}$ is a sum of pushforwards of $\mc{V}_{\rho}$. Such pushforwards for different $\lambda$ have no nontrivial maps on $\Par_G$, so $\mc F$ must be a direct sum of such $*$-pushforwards on $\Par_G$. 
\end{proof}

\begin{remark}
In (2), if $\phi$ and $\phi'$ are not $\wh G$-conjugate, then the vanishing $\Hom (\mathcal{V}_{\rho}, \mathcal{V}_{\rho'})\simeq 0$ clearly holds, and the full faithfulness claim implies the rest. The main point of (3) for semisimple $G$ is that cuspidal vector bundles generate cuspidal ind-coherent sheaves, and the semisimplicity follows from (2).  
\end{remark}

We identify $\mc V_\rho$ with its $*$-pushforward below if there is no confusion, and call such sheaves \emph{cuspidal vector bundles}. As a consequence of Theorem \ref{thm: CVBclassification} we deduce the following properties of cuspidal vector bundles.

\begin{corollary}\label{cor: CVBproperty}
\begin{enumerate}
    \item $\bb D_\GS \mc V_{\rho}\cong \mc V_{\rho^{\vee}}[-\dim Z(\wh G)^{\Gamma}]$. 
    \item The set of cuspidal vector bundles is preserved by any $\Gamma_F$-equivariant automorphism of $\wh G$ preserving $\wh T$.     
    \item Let $c$ denote the Chevalley involution with respect to $(\wh B, \wh T)$ (and its pinning). Then $c^* \D_\GS \mc V_{\rho}$ is cuspidal and isomorphic to $\mc V_{c^* \rho^{\vee}}[-\dim Z(\wh G)^{\Gamma}]$ under the induced isomorphism $c\colon S_{c\circ\phi}\cong S_{\phi}$. 
    \item For $z\in Z(\wh G)^{\Gamma,\circ}\subset I_\lambda$, we consider the translation $z\lambda$, $z\phi$. Under the identification $\V_\lambda\cong \V_{z\lambda}$, the cuspidal vector bundle $\mc V_{\rho}$ corresponds to the cuspidal vector bundle $\mc V_{\rho_z}$ where $\rho_z$ corresponds to $\rho$ under the identification $S_\phi\cong S_{z\phi}$.   
\end{enumerate}
\end{corollary}

\begin{proof}
This follows from Theorem \ref{thm: CVBclassification} together with the fact that the various operations in the statement preserve cuspidal vector bundles. For instance, the stability under $\D_\GS$ reduces to the (simply-connected) graded setting and it follows from the characterization (or the construction) of cuspidal vector bundles in \S \ref{s: CVBs}.  
\end{proof}

We can also describe generators of cuspidal sheaves. Consider the $!$-pullback of $\mc O_{\Lie Z(\wh G)^\Gamma}\boxtimes \mc V_\rho$ to the completion along the zero section of $\T\V_\lambda[-1]$, living in degree 0 with respect to the $t$-structure.\footnote{More canonically, $!$-pullback along $\T\V_\lambda[-1]\to \T\V_\lambda[-1] / \Lie Z(\wh G)^{\Gamma}$ of $\mc V_\rho$ on the zero section, shifted by $[-\dim Z(\wh G)^\Gamma]$. We apply linear Koszul duality in this way, cf.~Example \ref{ex: derivedVB}.} It is cuspidal, and, by restricting to $(\Lie Z(\wh G)^\Gamma)^\wedge_0 \times BS_\phi$, it corresponds to $\mathcal{O}_{(\Lie Z(\wh G)^\Gamma)^\wedge_0}\boxtimes \rho$ and, in particular, its endomorphism algebra is the ring of formal power series $\mc O(\Lie Z(\wh G)^\Gamma)^{\wedge}_0$, concentrated in degree 0 as above. We write $\wh{\mc V}_\rho$ for this generator. The claim that they form generators is proven on the Koszul dual side (Theorem \ref{thm: ClassificationCuspidalCotangent}). 

\subsection{Comparison with the work of Aubert--Moussaoui--Solleveld}
We now describe how our work is compatible with \cite{AMS1}. In \S\ref{s: amscomparison}, we briefly recalled the statement and context of their \cite[Conjecture 6.10]{AMS1} which gives a conjectural description of which enhanced $L$-parameters in the classical local Langlands correspondence correspond to supercuspidal representations on the automorphic side.

We wish to relate Conjecture \ref{conj: amspure} and Theorem \ref{thm: CVBclassification}. These statements deal with essentially discrete $L$-parameters, so we assume this in the ensuing discussion.

Firstly, we remark that by Lemma \ref{lem: tftorsdisclemma}, we have that
\begin{equation*}
   S_{\tilde{\phi}|_{W_F}} = J_{\lambda}, 
\end{equation*}
where $\lambda$ is the semisimplification of $\tilde{\phi}$. Then by Example \ref{ex: Ilambda}, we may take $I_{\lambda}$ to be $J^{\circ}_{\lambda}$. In the context of Definition \ref{def: cuspgeneral}, we say that a representation $\rho$ of $S_{\phi}
$ is \emph{pure} if it factors through $\pi_0(S_{\phi})$. Then there is an identification of pure, cuspidal $\rho$ and the cuspidal enhanced $L$-parameters of Conjecture \ref{conj: amspure}, where the difference between graded and ungraded cuspidality is discussed in \S \ref{ss: lusztiggraded}. Hence for pure $\rho \in \Irr(S_{\phi})$, a cuspidal vector bundle $\mc{V}_{\rho}$ exists precisely when $(\tilde{\phi}, \rho)$ is a cuspidal enhanced $L$-parameter in the sense of Conjecture \ref{conj: amspure}.
\subsection{Generators of cuspidal sheaves}
When $G$ is not semisimple, the final statement of part (3) of Theorem \ref{thm: CVBclassification} does not hold. In this subsection we explain the general situation.

Given two semisimple parameters $\lambda, \lambda'\colon W_F \to \wh G$, we write $\lambda \sim \lambda'$ if they differ by a map $W_F \to \frob^\Z \to Z(\widehat G)^{\Gamma, \circ}$ up to conjugation. Let $\bb{V}_{[\lambda]}$ denote the universal family of Vogan stacks defined as 
\[
\{ (\lambda', N), \lambda' \colon W_F\to \wh G, N\in V_{\lambda' } \mid \lambda \sim \lambda' \} /\wh G. 
\]
There exists a finite \'etale cover
\[
Z(\wh G)^{\Gamma, \circ}\times \V_\lambda\to \bb{V}_{[\lambda]}, 
\]
which is a torsor under the finite group $\{z \in Z(\wh G)^{\Gamma, \circ}(\Lambda) \mid z\cdot\lambda\,\text{is conjugate to}\,\lambda\}$. For any cuspidal $\rho$, we define $\widetilde{\mc V}_{\rho}$ to be the pushforward of $\mc O_{Z(\wh G)^{\Gamma, \circ}}\boxtimes \mc V_{\rho}$, which is a vector bundle on $\V_{[\lambda]}$. We identify $\widetilde{\mc V}_\rho$ with its $*$-pushforward to $\Par_G$. By the construction, $\wh{i}^!_\lambda \widetilde{\mc V}_\rho$ is a direct sum of $\wh{\mc V}_{\rho'}$ for certain conjugates $\rho'$ of $\rho$ with multiplicities. More precisely, the $\rho'$ that occur are exactly those for which there exists $z\in Z(\wh G)^{\Gamma, \circ}$ and $g\in \wh G$ such that $(\lambda, N)=g(z\cdot \lambda, N)g^{-1}$ (at the level of fixed representatives, $(\lambda, N)$ being essentially discrete) and $ \rho\cong \Ad(g)_*\rho'$ under the resulting identification $\Ad(g)\colon S_{\phi}\xrightarrow{\cong} S_\phi$. 

It follows from Corollary \ref{cor: CVBproperty} that

\begin{proposition}{\label{prop: cuspproperties}}
    \begin{enumerate}
    \item The set $\{\widetilde{\mc{V}}_{\rho}\}$ is preserved by any $\Gamma_F$-equivariant automorphism of $\wh G$ preserving $\wh T$. 
    \item $\D_\GS \widetilde{\mc V}_{\rho}\cong \widetilde{\mc V}_{\rho^{\vee}}$, $c^*\D_\GS\widetilde{\mc V}_{\rho}\cong \widetilde{\mc V}_{c^*\rho^{\vee}}$.      
    \end{enumerate}
\end{proposition}

The following result is analogous to the structure of cuspidal Bernstein blocks on the geometric side. 

\begin{theorem}\label{thm: Generators}
\begin{enumerate}
    \item The stable $\infty$-category of cuspidal ind-coherent sheaves is generated under colimits by $\widetilde{\mc V}_{\rho}$ with $[\lambda]$ and irreducible cuspidal $\rho$ varying.  
    \item Suppose $\rho$, $\rho'$ are irreducible cuspidal representations of $S_{\phi}, S_{\phi'}$ respectively. Then $\Hom_{\Par_G}(\widetilde{\mc V}_{\rho}, \widetilde{\mc V}_{\rho'})=0$ if and only if $\lambda \sim \lambda'$ and, under some induced isomorphism $S_{\phi}\cong S_{\phi'}$, $\rho$ corresponds to $\rho'$. And $\mathrm{Ext}^i_{\Par_G} (\widetilde{\mc V}_{\rho}, \mc F)=0$ for $i>0$ and any cuspidal coherent sheaf $\mc F$ in degree 0. 
    \item If $\rho$ is irreducible and cuspidal, then $\End (\widetilde{\mc V}_{\rho})$ is an Azumaya algebra over its center, and the center is the ring of regular functions on the following torus
    \[
    Z(\wh G)^{\Gamma, \circ}/\{z\in Z(\wh G)^{\Gamma, \circ}\mid \exists g, g\lambda g^{-1}=z\lambda,  \mathrm{Ad}(g)^* \mc V_\rho\cong \mc V_\rho\}. 
    \] 
    In particular, simple modules correspond to closed points of this space.
    \end{enumerate}
\end{theorem}
\begin{proof}
    We reduce to the simply-connected graded case in \S\ref{ss: reductiongenerators}. In the simply-connected case, the theorem follows from Theorem \ref{thm: CVBclassification}.
\end{proof}

\begin{remark}\label{rem: Azumaya}
An explicit description of the Azumaya algebra in (3) is given in \S\ref{ss: reductiongenerators}. It is a crossed product of the form
\[
\mc O(Z(\wh G)^{\Gamma,\circ})\rtimes_{\mu_\rho}
\{z\in Z(\wh G)^{\Gamma, \circ}\mid \exists g, g\lambda g^{-1}=z\lambda,  \mathrm{Ad}(g)^* \mc V_\rho\cong \mc V_\rho\}
\]
for a 2-cocycle $\mu_\rho$ whose class in the second group cohomology class of the right factor with values in $\Lambda^\times$ is the one associated with $\rho$ and the extension 
\[
\{g\in \wh G\mid \exists z\in Z(\wh G)^{\Gamma,\circ}, g\phi g^{-1}=z\phi, {}^g \rho\cong \rho\}\supset S_\phi, 
\]
cf.~\cite[Appendix A]{AbeOhara}. 
This is related to \cite{AbeOhara}: one difference is that there exists a surjection
\[
Z(\wh G)^{\Gamma,\circ} \to X_{\mathrm{unr}}({}^L G)\coloneq (Z(\wh G)^{I_F})_{\frob}^\circ, 
\]
and they work with the image under this map. Our class is pulled back from the class obtained in \cite[Proposition 4.1.6]{AbeOhara}, at least up to sign. 
\end{remark}
Although we state everything for cuspidal sheaves, we will prove them for right-cuspidal sheaves, because of the following implication:

\begin{corollary}\label{cor: LCusp=RCusp}
Suppose Theorem \ref{thm: Generators} holds for right-cuspidal sheaves. Then
\[
\LCusp(\Par_G)=\RCusp(\Par_G)
\]
and it is compactly generated by
\[
\LCusp^{\Coh}(\Par_G)=\RCusp^{\Coh}(\Par_G). 
\]
\end{corollary}

We will simply write $\Cusp^{\Coh}(\Par_G)\subset \Cusp(\Par_G)$ for these categories. 

\begin{proof}

It follows from Theorem \ref{thm: Generators} for right-cuspidal sheaves and Proposition \ref{prop: cuspproperties} (2) that
\[
\LCusp^{\Coh}(\Par_G)=\RCusp^{\Coh}(\Par_G),
\]
and $\RCusp(\Par_G)$ is compactly generated by $\RCusp^{\Coh}(\Par_G)$. Therefore, we obtain a semiorthogonal decomposition of $\Coh(\Par_G)$ using $\RCusp^{\Coh}(\Par_G)$ and Eisenstein series. By passing to the Grothendieck--Serre dual, we obtain another semiorthogonal decomposition of $\Coh(\Par_G)$ using $\LCusp^{\Coh}(\Par_G)$ and Eisenstein series. By passing to the ind-completion, we obtain a semiorthogonal decomposition of $\IndCoh(\Par_G)$, and it implies that $\LCusp(\Par_G)$ is the ind-completion of $\LCusp^{\Coh}(\Par_G)$. Since $\LCusp(\Par_G)$ and $\RCusp(\Par_G)$ are ind-completions in $\IndCoh(\Par_G)$ of the same category, they are equal. 
\end{proof}

\begin{conjecture}
We denote a fixed Whittaker datum by $\mf w$ and choose $q^{1/2}\in \overline{\Q}_{\ell}$. 
Recall the unconditionally defined functors with coefficients in $\overline{\Q}_{\ell}$ using the spectral action \cite{fargues--scholze}
\[
a_{\mf w}\colon \QCoh(\Par_G)\to D(\Bun_G), \quad c_{\mf w}\colon D(\Bun_G) \to \QCoh(\Par_G);
\]
$c_{\mf w}$ is expected to send compact objects to coherent sheaves and the resulting ind-completion
\[
\mathbf{L}_{\mf w}\colon D(\Bun_G) \rightarrow \IndCoh(\Par_G) 
\]
is expected to be an equivalence, cf.~\cite{HansenMann}. 

\begin{enumerate}
    \item The functor $c_{\mf w}\circ i_{b!}$, $b\in B(G)_\bas$, is $t$-exact on cuspidal Bernstein blocks. If $\pi$ is an irreducible supercuspidal representation of $G_b (F)$, then $c_{\mf w}\circ i_{b!}\pi$ is isomorphic to some $\mc V_{\rho}$, $\rho$ cuspidal. 
    \item The coherent sheaf $\widetilde{\mc V}_{\rho}$ is the image of the $!$-extension of a compact projective generator of a cuspidal Bernstein block of $\Rep (G_b(F))^{\heartsuit}$. 
\end{enumerate}
\end{conjecture}

(2) is related to \cite[Conjecture 4.2.3]{AbeOhara}.

\begin{definition}
Let $M \subset P$ be standard Levi and parabolic subgroups of $G$. For $\widetilde{\mc V}_{\rho_M} \in \Coh(\Par_M)^{\kappa_M(b_M)}$ with irreducible $\rho_M$ as above and $b_M \in B(M)_{\bas}$ with image $b \in B(G)$, we call 
\[
\Eis_P^{\spec} (\widetilde{\mc V}_{\rho_M}) \in \Coh(\Par_G)^{\kappa(b)},
\]
a dominant/antidominant \emph{generalized coherent Springer sheaf} if $b_M$ is $P$-dominant/antidominant. 
\end{definition}

\begin{corollary}{\label{cor: ParGgenbycohspringer}}
The presentable stable $\infty$-category $\IndCoh (\Par_G)$ is generated under colimits by dominant/antidominant generalized coherent Springer sheaves. In particular, $\Coh (\Par_G)$ is generated by them under finite colimits and retracts.
\end{corollary}
\begin{proof}
    This follows from the semiorthogonal decomposition in Theorem \ref{thm: cuspSOD} together with the description of generators of cuspidal sheaves in Theorem \ref{thm: Generators}. 
\end{proof}

\subsection{Classification of cuspidal sheaves on cotangent stacks}

Theorem \ref{thm: CVBclassification} is deduced via a local analysis of strata of $\Par_G$. More precisely, we study the completion of $\Par_G$ at $\V_{\lambda}$, and this identifies with the completed odd tangent stack $\wh{\T}\V_{\lambda}[-1]$ by Proposition \ref{prop: LoopSpaceIdentificationCompletions}. Then by Koszul duality, we reduce Theorem \ref{thm: CVBclassification} to analogous results on $\T^{\ast} \V_{\lambda}$. In the course of the proof of these results, we also reduce to Lusztig's graded setting which is roughly the unramified setting of \cite[\S5]{CFMMX}. For simplicity, we omit discussion of these reductions in the present subsection. 

In this subsection, we formulate our results for cotangent stacks. We believe these are of independent interest; indeed such spaces and categories of sheaves are studied intensively in the literature (e.g.~\cite{PTDolbeault, BuPadurariuToda}), and our proof is largely inspired by them. 

Fix a semisimple $L$-parameter $\lambda \colon W_F\to \wh G$, and also $\phi=(\lambda, N)$ for $N\in V_\lambda$ in the open orbit. 

\begin{theorem}\label{thm: ClassificationCuspidalCotangent}
\begin{enumerate}
    \item If $\mc F\in \Coh(\T^* \V_\lambda)$ is right-cuspidal, then $\D'_{\GS}\mc F\coloneq\omega_{\V_\lambda}\otimes\D_{\GS}\mc F$ is right-cuspidal and
    \begin{equation*}
            \LCusp^{\Coh}(\T^*\V_\lambda)=\RCusp^{\Coh}(\T^*\V_\lambda). 
    \end{equation*}
    Moreover, $\RCusp(\T^{\ast}\V_{\lambda})$ is compactly generated by $\RCusp^{\Coh}(\T^{\ast}\V_{\lambda})$ and
    \[
    \LCusp(\T^*\V_\lambda)=\RCusp(\T^*\V_\lambda). 
    \]
    We simply call these sheaves cuspidal. 
    \item[(1)'] $\Coh(\T^*\V_\lambda)$ admits a semiorthogonal decomposition $\langle \Cusp^{\Coh}, \Eis^{\T^{\ast}}_+ \rangle$, where $\Eis^{\T^{\ast}}_+$ is generated by  the image of $\Eis^{\T^{\ast}}_P$ with $P$-dominant weights for all proper semistandard $P$. An analogous statement holds with $P$-dominance replaced by $P$-antidominance. 
    \item The restriction
\[
\Coh(\T^* \V_\lambda) \to \Coh(\T^* BS_\phi)
\]
is fully faithful on cuspidal sheaves. 
    \item The essential image of cuspidal sheaves under the previous functor is generated under finite colimits by $*$-pushforwards of cuspidal representations of $S_\phi$ along the zero section. 
    \item The coherent sheaf corresponding to a cuspidal representation $\rho$ by the previous two parts is a vector bundle on the underlying classical stack of $\T^* \V_\lambda$. 
    \item The $\infty$-category $\Coh(\T^* \V_\lambda)$ is generated under finite colimits and retracts by cuspidal sheaves and $\Eis^{\T^*,\chi}_P$ of cuspidal sheaves with $P$-dominant $\chi$.
\end{enumerate}
In particular, the cuspidal ind-coherent sheaves are compactly generated by cuspidal coherent sheaves, and there exists a nonzero cuspidal sheaf only if $\phi$ is essentially discrete. 
\end{theorem}
\begin{proof}
    (1)' is an immediate consequence of (1). After reduction to Lusztig's graded setting in \S\ref{ss: reductionurclassificationcuspidalcotangent}, the first two statements in (1) follow for perfect complexes from Corollary \ref{cor: leftrightbasequal} and Lemma \ref{lem: Cotangentleftrightcuspidality}; for coherent sheaves they are proved at the end of \S \ref{ss: basgen}, together with (2), which uses Proposition \ref{prop: cuspextopen}. The proof of (4) is by an explicit classification of cuspidal vector bundles carried out in \S \ref{s: CVBs}. After reducing to the essentially discrete case in \S \ref{ss: reductiondiscClassificationCuspidalCotangent} (proving the last claim in Proposition \ref{prop: cuspreductiontofinite}), the last part of (1) and (3) are proven in \S \ref{ss: basgen}. Finally, (5) is proven in \S \ref{ss: SODcotangent}. 
\end{proof}

\begin{remark}
The vector bundle on the underlying classical stack corresponding to $\rho$ is Koszul-dual to $\wh{\mc V}_\rho [\dim Z(\wh G)^\Gamma]$ by Example \ref{ex: derivedVB}.  
\end{remark}

We have the following results for basic sheaves:

\begin{theorem}\label{thm: BasicCotangentMain}
\begin{enumerate}
    \item If $\mc F\in \Coh(\T^* \V_\lambda)$ is right-basic, then $\D'_{\GS}\mc F$ is right-basic. Moreover, 
    \[
    \LBas^{\Coh}(\T^*\V_\lambda)=\RBas^{\Coh}(\T^*\V_\lambda), 
    \]
    and
     \[
    \LBas(\T^*\V_\lambda)=\RBas(\T^*\V_\lambda). 
    \]
      We simply call these sheaves basic, and denote the categories by $\Bas^{\Coh}(\T^{\ast}\V_{\lambda})$ and $\Bas(\T^{\ast}\V_{\lambda})$. Moreover, $\Bas(\T^{\ast}\V_{\lambda})$ is compactly generated by $\Bas^{\Coh}(\T^{\ast}\V_{\lambda})$.
   
    \item $\Eis^{\T^*}_{P,\lambda_P}\mc F$ is basic if $\mc F$ is basic and the $Z(\wh M)^{\Gamma}$-weights of $\mc F$ are basic in $B(G)$. 
    \item The basic coherent sheaves are generated under finite colimits and retracts by cuspidal coherent sheaves and Eisenstein series of basic coherent sheaves as above. 
\end{enumerate}
\end{theorem}
\begin{proof}
In the graded setting, the proof of $(1)$ is Corollary \ref{cor: LBas=RBasTV} and $(3)$ is Theorem \ref{thm: cotanbasgenbycuspeis}. Given (1) for Levi subgroups, $(2)$ is Proposition \ref{prop: EisBasisBas}. The general case of $(1)$ and $(3)$ is Proposition \ref{prop: refinedcotanSOD}(4).\footnote{Alternatively, it would be possible to reduce all statements to the graded case.} 
\end{proof}

This result is related to limit categories in \cite{PTDolbeault} and semiorthogonal decompositions discussed in \cite{BuPadurariuToda}, while we use their results and work out the precise relation only in the graded setting. 

\section{Reduction to Lusztig's graded setting}{\label{s: cohreductiontograded}}
This section is analogous to \S \ref{ss: CCSreductiontograded}. We will reduce Theorems \ref{thm: discretesupport}, \ref{thm: CVBclassification}, \ref{thm: Generators}, and \ref{thm: ClassificationCuspidalCotangent} to $I_\lambda$, $I_{\lambda, \der}$, or its simply-connected cover. Recall that it suffices to consider right-cuspidal sheaves (see Corollary \ref{cor: LCusp=RCusp}), and we shall do so. 

\subsection{Reduction to the unramified split case}
Let $\lambda\colon W_F \to \wh G$ be a semisimple $L$-parameter. By Lemma \ref{lem: TJcontrol}, we may assume $\wh T_I\coloneq\wh {T}\cap I_\lambda$ is a maximal torus of $I_\lambda$ contained in $\wh T^\Gamma$ and $s_\lambda\in \wh T_I$ (see \S\ref{ss: CCSreductiontograded}). Let us consider the finite \'etale cover $\pi\colon \V_{\lambda^\ur}=V_\lambda/H_\lambda^\circ\to \V_\lambda$ together with the induced finite \'etale cover
\[
\pi^{\wh{\T}[-1]}\colon \wh{\T}\V_{\lambda^\ur}[-1] \to \wh{\T}\V_{\lambda}[-1]. 
\]

\begin{lemma}\label{lem: CuspILambdaCenter}
If a nonzero right-cuspidal sheaf $\mc F\in \Ind\wh{\Coh}{}^! (\wh{\T}\V_\lambda[-1])$ exists, then $Z(\wh G)^{\Gamma,\circ}=Z(I_\lambda)^{\circ}$.     
\end{lemma}

\begin{proof}
Suppose $Z(I_\lambda)^\circ$ is strictly larger than $Z(\wh G)^{\Gamma, \circ}$. Then, there is a proper semistandard Levi subgroup $\wh M$ such that $I_\lambda\subset \wh M$, and ${}^L \lambda$ factors through $\wh M\rtimes W_F$. Then, for any semistandard parabolic subgroup $\wh P\supset \wh M$ with $\wh M$ being a Levi subgroup, $V_{\lambda}=V_{\lambda_P}=V_{\lambda_M}$, $H_\lambda^\circ=H_{\lambda_P}^\circ=H_{\lambda_M}^\circ$, and $H_{\lambda_P}=H_{\lambda_M}$ (where the last equality follows from the others and from $\pi_0(H_{\lambda_P}) = \pi_0(H_{\lambda_M})$ as in the proof of Lemma \ref{lem: reduction_to_LusztigGraded}). Therefore, $\Eis_{P,\lambda_P}^{\wh \T [-1]}$ is just the $*$-pushforward and admits a right adjoint $\CT_{P,\lambda_P}^{\wh \T [-1]}$ given by the $!$-pullback to $\wh{\T}\V_{\lambda_M}[-1]$. It is nonzero if $\mc F$ is nonzero, as $V_\lambda/H_{\lambda_M}\to V_\lambda/H_\lambda$ is a finite \'etale cover. Let $\chi\in X^* (Z(\wh M)^{\Gamma})$ be a $Z(\wh M)^\Gamma$-weight appearing in the $!$-pullback of $\mc F$ to $\wh{\T}\V_{\lambda_M}[-1]$. Then there exists some $\wh P$ such that $\chi$ is $P$-dominant, contradicting the right-cuspidality of $\mc F$. 
\end{proof}

\begin{lemma}{\label{lem: maxparabolicreslem}}
Suppose $Z(\wh{G})^{\Gamma, \circ}= Z(I_{\lambda})^{\circ}$ and let $\wh P_I$ be a maximal semistandard parabolic subgroup of  $I_{\lambda}$ with semistandard Levi subgroup $\wh M_I$. Then there exists a maximal semistandard parabolic subgroup $P \subset G$, with Levi subgroup $M$ such that $\lambda$ factors through $\wh M$, $\wh P \cap I_{\lambda}= \wh P_I$, and $\wh{M} \cap I_{\lambda} = \wh{M}_I$. 
\end{lemma}
\begin{proof}
Let us consider the groups $\wh M_1\subset \wh P_1 \subset \wh{G}$ constructed as in the second paragraph of the proof of Lemma \ref{lem: reduction_to_LusztigGraded} (we are using a different notation here). If $\wh P_1$ is maximal, we are done. Otherwise, let $\wh P\subset \wh G$ be a maximal semistandard parabolic subgroup containing $\wh P_1$ (recall from Lemma \ref{lem: TJcontrol} that $\wh P$ is defined by a cocharacter factoring through $\wh T_I$ whose centralizer is written as $\wh M$). By construction $\lambda$ factors through $\wh{M}_1$ hence through $\wh M\supset \wh M_1$. Then by Lemma \ref{lem: lambdaIparabolic}, we cannot have $I_{\lambda} \subset \wh P$, and hence we must have $\wh P \cap I_{\lambda} = \wh P_I$, as the intersection is a parabolic subgroup, and $\wh M \cap I_{\lambda}=\wh M_I$. 
\end{proof}

\begin{proposition}\label{prop: PullbackIlambdaCuspidal}
If a sheaf $\mc F\in \Ind\wh{\Coh}{}^! (\wh{\T}\V_\lambda[-1])$ is right-cuspidal, then $(\pi^{\wh{\T}[-1]})^! \mc F$ is right-cuspidal with respect to $I_\lambda$. It also holds for ind-coherent sheaves and unipotent loop spaces. 
\end{proposition}

\begin{proof}
Pick a maximal semistandard parabolic subgroup $\wh P_I$ of $I_\lambda$ through which $\lambda^{\ur}$ factors. We need to show that $(\pi^{\wh{\T}[-1]})^! \mc F$ is right-orthogonal to $\Eis_{P_I, \lambda^\ur_{P_I}}^{\wh \T[-1]} \mc G$ such that $\mc G$ has $P_I$-dominant $Z(\wh M_I)$-weights. We may assume $\mc F\neq 0$ and hence $Z(\wh G)^{\Gamma,\circ}=Z(I_\lambda)^\circ$ by Lemma \ref{lem: CuspILambdaCenter}. Then, we take a maximal semistandard parabolic subgroup $P\subset G$ such that $\wh P \cap I_\lambda=\wh P_I$ and $\wh M \cap I_\lambda=\wh M_I$ as in Lemma \ref{lem: maxparabolicreslem}. As in the proof of Lemma \ref{lem: reduction_to_LusztigGraded}, we have the commutative diagram
\[
\begin{tikzcd}
    \V_{\lambda_{M_I}^{\ur}} \arrow[d] & \arrow[l] \V_{\lambda_{P_I}^{\ur}} \arrow[d] \arrow[r] & \V_{\lambda^{\ur}} \arrow[d] \\
    \V_{\lambda_M}  & \arrow[l] \V_{\lambda_P} \arrow[r] & \V_{\lambda} 
\end{tikzcd}
\]
with the left square being cartesian. Since the vertical maps are \'etale, the induced commutative diagram
\begin{equation}{\label{eqn: urcompletedtandiagram}}
  \begin{tikzcd}
    \wh{\T}\V_{\lambda_{M_I}^{\ur}}[-1] \arrow[d, "\pi^{\wh{\T}[-1]}"] & \arrow[l] \wh{\T}\V_{\lambda_{P_I}^{\ur}}[-1] \arrow[d] \arrow[r] & \wh{\T}\V_{\lambda^{\ur}}[-1] \arrow[d, "\pi^{\wh{\T}[-1]}"] \\
    \wh{\T}\V_{\lambda_M}[-1]  & \arrow[l] \wh{\T}\V_{\lambda_P}[-1] \arrow[r] & \wh{\T}\V_{\lambda}[-1]
\end{tikzcd}  
\end{equation}
is the base change of the previous one, and the left square is cartesian. By base change, $(\pi^{\wh{\T}[-1]})_*\Eis_{P_I, \lambda^\ur_{P_I}}^{\wh \T[-1]}\mc G$ identifies with $\Eis_{P, \lambda_P}^{\wh \T[-1]}(\pi^{\wh{\T}[-1]})_* \mc G$. 

To conclude the right-cuspidality of $(\pi^{\wh{\T}[-1]})^! \mc F$, by the right-cuspidality of $\mc F$, it suffices to show that $(\pi^{\wh{\T}[-1]})_* \mc G$ has $P$-dominant $Z(\wh M)^{\Gamma}$-weights. The dominance can be checked by restricting to $Z(\wh M)^{\Gamma,\circ}=Z(\wh M_I)^\circ$ (equality holds since both parabolics are maximal and the left is contained in the right), and such weights are the same $Z(\wh M_I)^\circ$-weights appearing in $\mc G$. Therefore, they are $P_I$-dominant. Finally, we observe the following: a $Z(\wh M_I)^\circ$-weight (or $Z(\wh M_I)$-weight) is $P_I$-dominant if and only if it is $P$-dominant. Indeed, since $P$ is assumed maximal, the $P$-dominance reduces to checking the coefficient of the omitted simple root $\sigma$ restricted to the maximal torus $\wh T_I$ of $H_{\lambda}$. As $\wh P_I$ is maximal in $I_{\lambda}$, we again just have to test the coefficient of the omitted simple root $\sigma_I$, which will be equal to the restriction of a root (not necessarily simple) appearing in the unipotent radical of $\wh{P}$. The coefficient of $\sigma$ is non-negative if and only if the coefficient of $\sigma_I$ is non-negative (since a character of $Z(\wh{M})^{\Gamma, \circ}=Z(\wh{M}_I)^{\circ}$ vanishes on all other simple roots). Therefore these dominance conditions are equivalent. 
\end{proof}

\begin{proposition}\label{prop: PushforwardIlambdaCuspidal}
Suppose $Z(\wh G)^{\Gamma,\circ}=Z(I_\lambda)^{\circ}$. If $\mc F\in \Ind\wh{\Coh}{}^! (\wh{\T}\V_{\lambda^\ur}[-1])$ is right-cuspidal, then $(\pi^{\wh{\T}[-1]})_* \mc F$ is right-cuspidal. It also holds for ind-coherent sheaves and unipotent loop spaces. 
\end{proposition}

\begin{proof}
By Lemma \ref{lem: ParPParGFibers}, after conjugating $(P, \lambda_P)$, it suffices to consider a maximal semistandard parabolic subgroup $\wh P$ through which the representative $\lambda$ factors, with $\lambda_P=\lambda$. Then $\lambda$ factors through $\wh M$ by Lemma \ref{lem: TJcontrol}. By Lemma \ref{lem: lambdaIparabolic}, $\wh P_I\coloneq \wh P\cap I_\lambda$ is a proper semistandard parabolic subgroup of $I_\lambda$. 

Note that the right square of the diagram \eqref{eqn: urcompletedtandiagram} need not be cartesian because the injective map $\pi_0(H_{\lambda_P}) \rightarrow \pi_0(H_{\lambda})$ need not be surjective. We have the following identification
\begin{equation*}
    \wh{\T}\V_{\lambda_P}[-1] \times_{\wh{\T}\V_{\lambda}[-1]} \wh{\T}\V_{\lambda^{\ur}}[-1]  \simeq \coprod\limits_{[x] \in \pi_0(H_{\lambda}) / \pi_0(H_{\lambda_P})} \wh{\T}\V_{\lambda_{\co{x}{P}_I}^{\ur}}[-1],
\end{equation*}
where $\co{x}{P}_I$ is the conjugate of $P_I$ by some representative $x$ of $[x] \in \pi_0(H_{\lambda}) / \pi_0(H_{\lambda_P})$ normalizing $I_\lambda$ and $\wh{T}_I$.

Now for each $[x] \in \pi_0(H_{\lambda}) / \pi_0(H_{\lambda_P})$, we have the commutative diagram
\begin{equation}{\label{eqn: contrivedfiberprod}}
    \begin{tikzcd}
    \wh{\T}\V_{\lambda_{\co{x}{M}_I}^{\ur}}[-1] \arrow[d, "\co{x}{\pi}^{\wh{\T}[-1]}"] & \wh{\T}\V_{\lambda_{\co{x}{P}_I}^{\ur}}[-1] \arrow[r, "\co{x}{p}_I"] \arrow[d, "\co{x}{\pi}^{\wh{\T}[-1]}"] \arrow[l, swap, "\co{x}{q}_I "] & \wh{\T}\V_{\lambda^{\ur}}[-1] \arrow[d, "\pi^{\wh{\T}[-1]}"] \\
    \wh{\T}\V_{\lambda_M}[-1]  & \arrow[l, swap, "q"] \wh{\T}\V_{\lambda_P}[-1] \arrow[r, "p"] & \wh{\T}\V_{\lambda}[-1].
\end{tikzcd}
\end{equation}
 We have $Z(\wh{M})^{\Gamma, \circ} \subset Z(\wh{M}_I)^{\circ}$. Let $\nu$ be a defining cocharacter of $\wh{P}$ which we assume factors through $\wh{G}_{\der}$ and note that we may as well assume $\nu$ factors through $Z(\wh{M})^{\Gamma, \circ}$ implying that $\nu$ factors through $Z(\wh{M}_I)^{\circ} \subset I_{\lambda}$ and hence that $\nu$ is a defining cocharacter for $\wh{P}_I$ in $I_{\lambda}$. Then $\co{x}{\nu}$ is a defining cocharacter of $\co{x}{P}_I$. Let $\mc{G} \in \Ind\wh{\Coh}{}^!(\wh{\T}\V_{\lambda_M}[-1])$ and suppose $\mc{G}$ has $P$-dominant central weights. Now let $\chi \in X^{\ast}(Z(\co{x}{\wh M}_I))$ be a central weight of $(\co{x}{\pi}^{\wh{\T}[-1]})^!\mc{G}$. By base-change of the right square of (\ref{eqn: contrivedfiberprod}), we have 
\begin{equation*}
     (\pi^{\wh{\T}[-1]})^!\Eis^{\wh{\T}[-1]}_{P,\lambda_P}\mc{G} \simeq \bigoplus\limits_{[x]} \Eis^{\wh{\T}[-1]}_{\co{x}{P}_I, \lambda^\ur_{\co{x}{P}_I}}(\co{x}{\pi}^{\wh{\T}[-1]})^!\mc{G}.
\end{equation*}
Then we claim that $\langle \co{x}{\nu}, \chi\rangle \geq 0$. If we can show this, then the right-cuspidality of $(\pi^{\wh{\T}[-1]})_* \mc F$ follows from that of $\mc F$ by adjunction and the above equation.

We now discuss the claim, which we prove using (an analogue of) Lemma \ref{lem: ParGdomswitch}. We first handle the case where $[x]$ is trivial. We can write $\nu$ uniquely as a sum $z_I + \sigma_I$ (resp. $z + \sigma$) where $z_I$ (resp. $z$) factors through $Z(I_{\lambda})^{\circ}$ (resp. $Z(\wh{G})^{\Gamma, \circ}$) and $\sigma_I$ (resp. $\sigma$) lies in the non-negative rational span of positive coroots of $I_{\lambda}$ (resp. $\wh{G}$). Since $Z(\wh{G})^{\Gamma, \circ} = Z(I_{\lambda})^{\circ}$ and the coroots of $I_{\lambda}$ lie in the span of coroots of $\wh{G}$, we deduce $z_I = z$ by uniqueness of decomposition. In fact, since $\nu$ factors through $\wh G_\der$, we have $z_I=z=0$.

Then let $\chi_I$ be a central weight of $(\pi^{\wh{\T}[-1]})^! \mc G$ and $\chi$ its restriction to $Z(\wh{M})^{\Gamma, \circ}$, which is then $P$-dominant by assumption. We have 
\begin{equation*}
    0  \leq \langle \sigma, \chi \rangle = \langle \sigma_I , \chi_I \rangle,
\end{equation*}
as desired.

 Now we reduce the case of nontrivial $[x]$ to the trivial one. The diagrams for $\Eis^{\wh{\T}[-1]}_{\co{x}{P}_I, \lambda^\ur_{\co{x}{P}_I}}$ and $\Eis^{\wh{\T}[-1]}_{P_I, \lambda^\ur_{P_I}}$ differ by twisting by $x$. This means we replace $\nu$ with $\co{x}{\nu}$ and $\chi$ with $\co{x}{\chi}$. Then we have
 \begin{equation*}
     \langle \nu, \chi \rangle = \langle \co{x}{\nu}, \co{x}{\chi} \rangle,
 \end{equation*}
 showing that $\langle \co{x}{\nu}, \co{x}{\chi} \rangle \geq 0$ is equivalent to $\langle \nu, \chi \rangle \geq 0$. This completes the reduction.
\end{proof}
\subsection{Reduction of Theorems \ref{thm: discretesupport} and \ref{thm: CVBclassification} to $I_\lambda$}{\label{ss: reductiondiscsuppcvbclassgraded}}

We may assume a nonzero right-cuspidal sheaf $\mc F\in \Ind\wh{\Coh}{}^! (\wh{\T}\V_\lambda[-1])$ exists. By Lemma \ref{lem: CuspILambdaCenter}, $Z(\wh{G})^{\Gamma, \circ} = Z(I_{\lambda})^{\circ}$. 

\begin{proposition}{\label{prop: cuspurpushpull}}
Assume $Z(\wh{G})^{\Gamma, \circ} = Z(I_{\lambda})^{\circ}$ and let $\mc V$ and $\mc{V}^{\ur}$ be vector bundles on $\V_\lambda$ and $\V_{\lambda^{\ur}}$ respectively. The right-cuspidality of $\mc V$ is equivalent to the right-cuspidality of $\pi^* \mc V$ with respect to $I_\lambda$. Moreover, $\pi_* \mc V^\ur$ is right-cuspidal if $\mc V^\ur$ is right-cuspidal. 
\end{proposition}

\begin{proof}
This follows from Propositions \ref{prop: PullbackIlambdaCuspidal} and \ref{prop: PushforwardIlambdaCuspidal} because $\pi^* \mc V\cong \pi^! \mc V$ and $\mc V$ is a direct summand of $\pi_* \pi^! \mc V\cong \pi_* \pi^* \mc V$.       
\end{proof}

Since a vector bundle $\mc V$ on $\V_{\lambda}$ is a summand of the vector bundle $\pi_* \pi^* \mc V$, any right-cuspidal vector bundle on $\V_\lambda$ is a summand of $\pi_* \mc V^\ur$ for some right-cuspidal vector bundle $\mc V^\ur$ on $\V_{\lambda^\ur}$. 

Assume Theorems \ref{thm: discretesupport} and \ref{thm: CVBclassification} hold for $I_\lambda$. Theorem \ref{thm: discretesupport} follows from the case of $I_\lambda$, Lemma \ref{lem: GIcenters}, and Proposition \ref{prop: PullbackIlambdaCuspidal}. We shall discuss Theorem \ref{thm: CVBclassification} (1) and (2) for cuspidal vector bundles. If $\mc V^\ur=\mc V_{\rho^\ur}$ for a cuspidal representation $\rho^\ur$, then
\[
(\pi_* \mc V_{\rho^\ur})|_{BS_\phi} \cong (\pi|_{BS_{\phi^\ur}})_* \rho^\ur \cong \Ind_{S_{\phi^\ur}}^{S_\phi} \rho^\ur,
\]
as a sheaf on $BS_{\phi}$ or a representation of $S_\phi$. The induction $\Ind_{S_{\phi^\ur}}^{S_\phi} \rho^\ur$ is the direct sum of finitely many cuspidal irreducible representations $\rho_i$ of $S_\phi$ (recall that an automorphism of $I_\lambda$ preserves the set of cuspidal representations by \S\ref{ss: lusztiggraded}). Note also that any cuspidal representation $\rho$ of $S_\phi$ appears in $\Ind_{S_{\phi^\ur}}^{S_\phi} \rho^\ur$ for some cuspidal $\rho^\ur$ by definition. On the other hand, since $\pi^*\pi_* \mc V_{\rho^\ur}\cong \pi^! \pi_* \mc V_{\rho^\ur}$ is a cuspidal vector bundle, $\pi^*\pi_* \mc V_{\rho^\ur}$ is a finite direct sum of $\mc V_{\rho^\ur_j}$ where $\rho^\ur_j$ is a conjugate of $\rho^\ur$. 

\begin{lemma}\label{lem: OddTangentRestriction}
The restriction to the open substack $\wh{\T}BS_\phi [-1]$ 
\[
\Hom_{\wh{\T}\V_{\lambda}[-1]}(\pi_* \mc V^\ur_1, \pi_* \mc V^\ur_2) \to \Hom_{\wh{\T}BS_\phi[-1]} (\pi_* \mc V^\ur_1|_{\wh{\T}BS_\phi[-1]}, \pi_* \mc V^\ur_2|_{\wh{\T}BS_\phi[-1]})
\]
is an equivalence.   
\end{lemma}

\begin{proof}
This follows from the adjunction, the above discussion, and Theorem \ref{thm: CVBclassification} (2) for $I_\lambda$. 
\end{proof}

All in all, these imply that the additive category of summands of $\pi_* \mc V^\ur$ with $\mc V^\ur$ varying is embedded into the semisimple category of representations of $S_\phi$, and $\mc V_{\rho}$ exists for any cuspidal $\rho$ and it is characterized as the one mapping to $\rho$. Moreover, the restriction to the open substack is fully faithful on the shifts of $\mc V_\rho$ as desired. 

For part (3), after using Lemma \ref{lem: ihatpreservescusp}, we only need to consider $\IndCoh(\Par_{G,\lambda}^\wedge)$. It reduces to the graded setting by Propositions \ref{prop: PullbackIlambdaCuspidal} and \ref{prop: PushforwardIlambdaCuspidal} together with the above discussion on cuspidal vector bundles, since the relevant operations preserve ind-coherent sheaves.  

\subsection{Reduction of Theorems \ref{thm: discretesupport} and \ref{thm: CVBclassification} to $I_{\lambda, \der}$ and its simply-connected cover}{\label{ss: reductiontosimplyconnecteddiscsuppandCVBclass}}

Assume $Z(\wh{G})^{\Gamma, \circ} = Z(I_{\lambda})^{\circ}$. Let us consider the simply-connected cover $I_{\lambda, \Sc}\to I_{\lambda, \der}$. After replacing $\chi, n$ in Lemma \ref{lem: ReductionToGradedUnipotent} if necessary, we may assume $\chi$ factors through $I_{\lambda, \der}$ and lifts to $I_{\lambda, \Sc}$, and we can talk about $H_{\lambda_\Sc}$, the centralizer of $\chi$ in $I_{\lambda, \Sc}$ (we also use the notation $H_{\lambda_{\der}}$ for the centralizer of $\chi$ in $I_{\lambda, \der}$). Let $\phi^{\ur}_{\der}$ be a central twist of the $L$-parameter $\phi^{\ur}$ that factors through $I_{\lambda,\der}$, and let $\phi^{\ur}_{\Sc}$ be a lift to $I_{\lambda, \Sc}$. Recall that $H_\lambda \cap I_{\lambda}=H_\lambda^\circ$ is the centralizer of $\chi$ in $I_\lambda$.
Let us consider maps
\[
BZ(\wh G)^{\Gamma,\circ}\times V_{\lambda} / H_{\lambda_\Sc}\to BZ(\wh G)^{\Gamma,\circ}\times V_{\lambda} / H_{\lambda_\der}  \to \V_{\lambda^\ur}; 
\]
these are gerbes banded by finite groups
\[
\ker(I_{\lambda, \Sc}\to I_{\lambda,\der}), \quad \ker(Z(\wh G)^{\Gamma,\circ}\times H_{\lambda_\der} \to H_\lambda^\circ),
\]
respectively. The maps of the open substacks
\[
BZ(\wh G)^{\Gamma,\circ}\times BS_{\phi^\ur_\Sc}\to BZ(\wh G)^{\Gamma,\circ}\times BS_{\phi^\ur_\der} \to BS_{\phi^\ur}
\]
are gerbes banded by the same groups. 
Moreover, the induced maps
\[
\wh{\T}(BZ(\wh G)^{\Gamma,\circ}\times V_{\lambda} / H_{\lambda_\Sc})[-1]\to \wh{\T}(BZ(\wh G)^{\Gamma,\circ}\times  V_{\lambda} / H_{\lambda_\der}) [-1] \to \wh{\T}\V_{\lambda^\ur}[-1]
\]
identify with the base change of the previous maps. We set $\V_{\lambda_{\der}^\ur}\coloneq V_\lambda /H_{\lambda_{\der}}$, $\V_{\lambda_{\Sc}^\ur}\coloneq V_\lambda /H_{\lambda_{\Sc}}$. 

Consider the map
\[
\mathrm{pr}\colon
\wh{\T}(BZ(\wh G)^{\Gamma,\circ}\times \V_{\lambda_{\Sc}^\ur} )[-1] \to
BZ(\wh G)^{\Gamma,\circ}\times \wh{\T}\V_{\lambda_\Sc^\ur}[-1] 
\]
over $\wh{\T}\V_{\lambda_\Sc^\ur}[-1]$. On the target, we can consider right-cuspidality: for each $Z(\wh G)^{\Gamma,\circ}$-weight, the corresponding summand, regarded as an object of $\Ind\wh{\Coh}{}^!(\wh{\T}\V_{\lambda_\Sc^\ur}[-1])$, is right-cuspidal. Note that the source is the completion of $Z(\wh G)^{\Gamma,\circ}\times B Z(\wh G)^{\Gamma,\circ}\times \wh{\T}\V_{\lambda_\Sc^\ur}[-1]$ along $\{e\}\times BZ(\wh G)^{\Gamma,\circ}\times \V_{\lambda_\Sc^\ur}$, and we may apply Theorem \ref{thm: completeness_along_closedpoint} and Corollary \ref{cor: cohhatcompletion} (cf.~proof of Proposition \ref{prop: unipvsformal}). Namely, we may work with
\[
\Spec ((\Lambda'\otimes_\Lambda \mc{O}(Z(\wh G)^{\Gamma,\circ}))^\wedge_{z=e} \otimes_{\Lambda'} \Lambda'')/ (Z(\wh G)^{\Gamma,\circ}\times H_{\lambda_\Sc}) \to
BZ(\wh G)^{\Gamma,\circ}\times (\Spec \Lambda'')/ H_{\lambda_\Sc} \to
\Spec \Lambda'
\]
for a suitable complete local ring $\Lambda'$ and a $\Lambda'$-algebra $\Lambda''$ with an action of $H_{\lambda_\Sc}$, and the ind-coherent sheaves there. The first map is not of finite type but flat with regular fibers, and everything here works fine. In particular, $\mathrm{pr}^!$ admits a conservative right adjoint $\mathrm{pr}_?$ and the image of $\mathrm{pr}^!$ generates under colimits all sheaves. 

\begin{lemma}\label{lem: reductiontosimplyconnectedcover}
A sheaf $\mc F^\ur \in \Ind\wh{\Coh}{}^! (\wh{\T}\V_{\lambda^\ur}[-1])$ is right-cuspidal if and only if the $!$-pullback of $\mc F^\ur$ to $\wh{\T}(BZ(\wh G)^{\Gamma,\circ}\times \V_{\lambda_{\der}^\ur} )[-1]$ is right-cuspidal with respect to $Z(\wh G)^{\Gamma,\circ}\times I_{\lambda,\der}$ if and only if the $!$-pullback to  $\wh{\T}(BZ(\wh G)^{\Gamma, \circ}\times \V_{\lambda_{\Sc}^\ur})[-1]$ is right-cuspidal with respect to $Z(\wh G)^{\Gamma,\circ}\times I_{\lambda,\Sc}$.\footnote{This equivalence also holds for ind-coherent sheaves and unipotent loop spaces.} Moreover, these conditions are also equivalent to the right-cuspidality with respect to $ I_{\lambda,\Sc}$ of $\mathrm{pr}_?$ applied to the $!$-pullback of $\mc F^\ur$. Analogous statements hold for right-basic sheaves. 

If these equivalent conditions for right-cuspidal sheaves hold and Theorem \ref{thm: CVBclassification} holds for $I_{\lambda, \Sc}$, then the $!$-pullback of $\mc F^\ur$ to $\wh{\T}(BZ(\wh G)^{\Gamma,\circ}\times \V_{\lambda_\Sc^\ur})[-1]$ is scheme-theoretically supported on $\wh{\T}BZ(\wh G)^{\Gamma,\circ}[-1]\times \V_{\lambda_\Sc^\ur}$ (i.e., the $!$-pullback is the pushforward from the latter closed substack). 
\end{lemma}

\begin{proof}
For any $(\wh P_I, \lambda^\ur_{P_I})$, we have the commutative diagram 
\[
\begin{tikzcd}
    \wh{\T}(BZ(\wh G)^{\Gamma,\circ}\times \V_{\lambda_{M_I, \Sc}^{\ur}})[-1] \arrow[d] & \arrow[l] \wh{\T}(BZ(\wh G)^{\Gamma, \circ}\times \V_{\lambda_{P_I, \Sc}^{\ur}})[-1] \arrow[d] \arrow[r] & \wh{\T}(BZ(\wh G)^{\Gamma,\circ}\times \V_{\lambda_{\Sc}^{\ur}})[-1] \arrow[d] \\
    \wh{\T}(BZ(\wh G)^{\Gamma,\circ}\times \V_{\lambda_{M_I, \der}^{\ur}})[-1] \arrow[d] & \arrow[l] \wh{\T}(BZ(\wh G)^{\Gamma, \circ}\times \V_{\lambda_{P_I, \der}^{\ur}})[-1] \arrow[d] \arrow[r] & \wh{\T}(BZ(\wh G)^{\Gamma,\circ}\times \V_{\lambda_{\der}^{\ur}})[-1] \arrow[d] \\
    \wh{\T}\V_{\lambda^\ur_{M_I}}[-1]  & \arrow[l] \wh{\T}\V_{\lambda^\ur_{P_I}}[-1] \arrow[r] & \wh{\T}\V_{\lambda^\ur}[-1]
\end{tikzcd}
\]
with obvious notation. All squares are cartesian, and all vertical maps are gerbes banded by the finite groups above. Since the $P_I$-dominance only depends on $\Lie (I_{\lambda, \der}\cap \wh P_{I})$, the first part follows easily. The right-cuspidality can be checked after applying $\mathrm{pr}_?$: $!$-pullback commutes with Eisenstein series and, for Levi subgroups, the image of $\mathrm{pr}^!$ generates under colimits all sheaves. This completes the proof of the first paragraph, and the second paragraph follows from the above identification, the same result for $\V_{\lambda^{\ur}_{\Sc}}$ (Theorem \ref{thm: CVBclassification} (3)), and Lemma \ref{lem: ihatpreservescusp} by the Barr--Beck--Lurie argument. 
\end{proof}

This reduces Theorem \ref{thm: discretesupport} to the case of $I_{\lambda, \Sc}$, since $\mathrm{pr}_?$ is conservative as in the above proof. It also follows that any right-cuspidal vector bundle, $\mc{V}$, on $\V_{\lambda^\ur}$ is a sum of descents of vector bundles $\chi \boxtimes \mc V_{\Sc}$ on $BZ(\wh G)^{\Gamma,\circ}\times \V_{\lambda_\Sc}$ with $\mc V_{\Sc}$ right-cuspidal  with respect to $I_{\lambda, \Sc}$ (indeed the pullback of $\mc{V}$ to $BZ(\wh G)^{\Gamma,\circ}\times \V_{\lambda_\Sc}$ decomposes as such a sum and the summands descend). And $\chi \boxtimes \mc V_{\rho_\Sc}$ for a character $\chi$ of $Z(\wh G)^{\Gamma, \circ}$ and a right-cuspidal vector bundle $\mc V_{\rho_{\Sc}}$ on $\V_{\lambda_\Sc}$ descends to $\V_{\lambda^\ur}$ if and only if $\chi\boxtimes \rho_{\Sc}$ is trivial on
\[
\ker (Z(\wh G)^{\Gamma, \circ}\times H_{\lambda_\Sc} \to H_{\lambda}^\circ) \subset Z(\wh G)^{\Gamma, \circ}\times S_{\phi_\Sc^\ur}, 
\]
equivalently, it descends to a representation of $S_{\phi^\ur}$. We have reduced Theorem \ref{thm: CVBclassification} (1) to the case of $I_{\lambda,\Sc}$. 

It also follows easily that
\[
\Hom (\chi_1 \boxtimes \mc V_{\rho_{\Sc,1}}, \chi_2 \boxtimes \mc V_{\rho_{\Sc,2}}) \simeq 0
\]
if $\chi_1 \not\cong \chi_2$. If $\chi_1\cong \chi_2 \cong \chi$ for some $\chi$, then
\[
\Hom (\chi \boxtimes \mc V_{\rho_{\Sc,1}}, \chi \boxtimes \mc V_{\rho_{\Sc,2}}) \simeq
\Hom (\mathbf{1} \boxtimes \mc V_{\rho_{\Sc,1}}, \mathbf{1} \boxtimes \mc V_{\rho_{\Sc,2}})
\]
for the trivial character $\mathbf{1}$. By Theorem \ref{thm: CVBclassification} (2) for $I_{\lambda, \Sc}$, the restriction is fully faithful. Now, it suffices to use linear Koszul duality and the following lemma (but beware of the degrees) to deduce part (2) for $I_\lambda$. 

\begin{lemma}
For any irreducible vector bundle $\chi \boxtimes \rho_{\Sc}$ on 
\[
\T^* (BZ(\wh G)^{\Gamma,\circ}\times BS_{\phi_{\Sc}^\ur})\simeq \T^* BZ(\wh G)^{\Gamma,\circ}\times BS_{\phi_{\Sc}^\ur},
\]
its $\mathrm{Ext}$ algebra is the exterior algebra of $\Lie Z(\wh G)^{\Gamma}$ in degree $-1$. 
\end{lemma}

\begin{proof}
The stack $\T^* BZ(\wh G)^{\Gamma,\circ}$ is the quotient by $Z(\wh G)^{\Gamma,\circ}$ of the spectrum of the dg algebra generated by $\Lie Z(\wh G)^{\Gamma}$ in degree $-1$. Then
\begin{equation*}
    \mathrm{End}(\chi \boxtimes \rho_{\Sc}) \simeq \mc{O}(\T^* BZ(\wh G)^{\Gamma,\circ}) \otimes \mathrm{End}_{S_{{\phi}^{\ur}_{\Sc}}}(\rho_{\Sc}) \simeq \bigwedge\nolimits^{\bullet}(\Lie Z(\wh G)^{\Gamma}), 
\end{equation*}
since $\rho_{\Sc}$ is irreducible, and with $\Lie Z(\wh G)^{\Gamma}$ in degree $-1$.
\end{proof}

For part (3), note that, by the footnote to Lemma \ref{lem: reductiontosimplyconnectedcover}, if $\mc F^\ur\in \IndCoh(\LS^{\up}_{I_{\lambda}})$ is right-cuspidal and set-theoretically supported on $\mc L^u(\V_{\lambda^\ur})$, then its $!$-pullback to 
\[
Z(\wh G)^{\Gamma,\circ} \times B Z(\wh G)^{\Gamma,\circ}\times \mc L^u(\V_{\lambda^\ur_\Sc})
\]
is right-cuspidal and set-theoretically supported on $BZ(\wh G)^{\Gamma,\circ}\times \mc L^u(\V_{\lambda^\ur_\Sc})$. Arguing as before, we see that the $*$-pushforward to 
\[
\IndCoh(BZ(\wh G)^{\Gamma,\circ}\times\mc L^u(\V_{\lambda^\ur_\Sc}))
\]
preserves right-cuspidality, and hence the $*$-pushforward is a sum of cuspidal vector bundles on $\V_{\lambda^\ur_\Sc}$ twisted by characters of $Z(\wh G)^{\Gamma,\circ}$. Therefore, the $!$-pullback of $\mc F^\ur$ is a colimit of such twisted cuspidal vector bundles. By the above discussion on cuspidal vector bundles, it follows that $\mc F^\ur$ is also a colimit of cuspidal vector bundles as desired.  

\subsection{Reduction of Theorem \ref{thm: ClassificationCuspidalCotangent}}{\label{ss: reductionurclassificationcuspidalcotangent}}
Part (5) of Theorem \ref{thm: ClassificationCuspidalCotangent} will be deduced from part (1) in \S \ref{ss: SODcotangent}. On the other hand, in essentially the same way as above (or by working Koszul-dually), we can first reduce Theorem \ref{thm: ClassificationCuspidalCotangent} except part (5) to the setting $\T^* \V_{\lambda^\ur}$ for the connected reductive group $I_\lambda$. We can then reduce it further to $I_{\lambda, \der}$, or $I_{\lambda, \Sc}$, noting the following points: for part (1), use that twisting by the dualizing sheaf $\omega_{BZ(\wh G)^{\Gamma,\circ}}$ preserves the cuspidality to compare the effects of the twisted Grothendieck--Serre duality. Part (2) is analogous to the arguments in previous sections. For parts (3) and (4), use that $\T^* BZ(\wh G)^{\Gamma,\circ}$ has the underlying classical stack $BZ(\wh G)^{\Gamma,\circ}$. 

Let us note, nevertheless, that if one is interested only in Theorems \ref{thm: discretesupport}, \ref{thm: CVBclassification}, \ref{thm: Generators} for the unipotent connected component $\LS^\up_{I_{\lambda, \Sc}}$, then it suffices to prove Theorem \ref{thm: ClassificationCuspidalCotangent} only for parameters in the unipotent connected component $\LS^\up_{I_{\lambda, \Sc}}$ by the discussion on linear Koszul duality (\S \ref{ss: koszuloddtangent}, Lemma \ref{lem: cusppullpush}, Lemma \ref{lem: koszulcusp}). 

\subsection{Reduction of Theorem \ref{thm: Generators}}{\label{ss: reductiongenerators}} 

We assume Theorems \ref{thm: discretesupport}, \ref{thm: CVBclassification}, and also Theorem \ref{thm: Generators} for the unipotent connected component $\LS^\up_{I_{\lambda, \Sc}}$. We aim to deduce Theorem \ref{thm: Generators} in the general case. The argument below is not formal, and we use Corollary \ref{cor: SODhatCoh} and Remark \ref{rem: ihatpreservescuspIndCoh} which are consequences of Theorem \ref{thm: ClassificationCuspidalCotangent}. 

Suppose $\mc F\in \IndCoh(\Par_G)$ is right-cuspidal. By Theorem \ref{thm: discretesupport} and Remark \ref{rem: ihatpreservescuspIndCoh}, for any semisimple $L$-parameter $\lambda$ such that $\wh{i}^!_\lambda\mc F\not\simeq 0$ in $\Ind\wh{\Coh}{}^! (\Par^\wedge_{G,\lambda})$, the Vogan stack $\V_\lambda$ is essentially discrete. Recall that the closure of the open substack of essentially discrete $L$-parameters is a disjoint union of irreducible components (essentially \cite[Proposition 2.32]{ZhuTame}). Pick one such irreducible component $\V_{[\lambda]}$ which is the closure of unramified twists of $\phi$ with semisimplification $\lambda$, and we may assume that $\mc F$ vanishes outside the closed substack of parameters whose semisimple part lies in the component $\V_{[\lambda]}$. 

We first prove Theorem \ref{thm: Generators} (1). Consider the translation map
\[
Z(\wh G)^{\Gamma, \circ}\times BZ(\wh G)^{\Gamma, \circ} \times \V_{\lambda_\Sc}  \to \V_{[\lambda]}, 
\]
and define $\widetilde{\mc V}_{\chi\boxtimes \rho_\Sc}$ for a cuspidal representation $\rho_\Sc$ of $S_{\phi^\ur_\Sc}$ and $\chi\in X^*(Z(\wh G)^{\Gamma,\circ})$ to be the pushforward of $\mc O_{Z(\wh G)^{\Gamma, \circ}}(\chi)\boxtimes \mc V_{\rho_\Sc}$. It is the finite direct sum, possibly zero, of several $\widetilde{\mc V}_{\rho}$, and $\widetilde{\mc V}_{\rho}$ for any cuspidal $\rho$ appears in some $\widetilde{\mc V}_{\chi\boxtimes\rho_\Sc}$. Therefore, it suffices to show that these sheaves are right-cuspidal and generate right-cuspidal sheaves. 

Now consider the following ``product'' map \'etale around $Z(\wh G)^{\Gamma, \circ}\times BZ(\wh G)^{\Gamma, \circ}\times \mc L^u(\V_{\lambda_\Sc})$:
\[
Z(\wh G)^{\Gamma, \circ}\times BZ(\wh G)^{\Gamma, \circ} \times \LS_{I_{\lambda,\Sc}}^\up \to \Par_{G, \lambda_I}\subset \Par_G, 
\]
where $\Par_{G, \lambda_I}\coloneq \LS_{\wh G,\lambda_I}^{\WD_F}$ and  $\LS_{I_{\lambda,\Sc}}^\up$ denotes the connected component of unipotent $L$-parameters (use Proposition \ref{prop: LoopSpaceIdentificationCompletions} to check the \'etaleness). 

\begin{lemma}\label{lem: PassingToSimplyConnectedCover}
If $\mc F\in \IndCoh (\Par_G)$ is right-cuspidal, then its $!$-pullback to $Z(\wh G)^{\Gamma, \circ}\times BZ(\wh G)^{\Gamma, \circ}\times \LS_{I_{\lambda,\Sc}}^\up$ is right-cuspidal on the completion along $Z(\wh G)^{\Gamma, \circ}\times BZ(\wh G)^{\Gamma, \circ}\times \mc L^u(\V_{\lambda_\Sc})$. If $\mc F$ vanishes outside the closed substack of parameters whose semisimple part lies in the component $\V_{[\lambda]}$, the converse holds. Under this support assumption, if moreover the $!$-pullback to $Z(\wh G)^{\Gamma, \circ}\times BZ(\wh G)^{\Gamma, \circ}\times \LS_{I_{\lambda,\Sc}}^\up$ vanishes, then $\mc F\simeq  0$. 
\end{lemma}

\begin{proof}
Let $\wh M_{I}\subset \wh P_{I}$ be maximal semistandard Levi and parabolic subgroups of $I_{\lambda}$, and let $\wh M_{I, \Sc}\subset \wh P_{I, \Sc}$ be the resulting ones for $I_{\lambda, \Sc}$. Using Lemma \ref{lem: maxparabolicreslem}, we take a maximal semistandard parabolic subgroup $P\subset G$ with Levi subgroup $M$ such that $\wh P\cap I_\lambda=\wh P_I$ and $\wh M\cap I_\lambda=\wh M_I$; since $\lambda$ factors through $\wh M$, any unramified twist of $\lambda$ factors through $\wh M$. We have the commutative diagram
\[
\begin{tikzcd}
    Z(\wh G)^{\Gamma, \circ}\times BZ(\wh G)^{\Gamma, \circ}\times \LS_{M_{I, \Sc}}^\up \arrow[d] &Z(\wh G)^{\Gamma, \circ}\times BZ(\wh G)^{\Gamma, \circ}\times \LS_{P_{I, \Sc}}^\up \arrow[l] \arrow[d] \arrow[r] & Z(\wh G)^{\Gamma, \circ}\times BZ(\wh G)^{\Gamma, \circ}\times \LS_{I_{\lambda, \Sc}}^\up \arrow[d] \\
    \Par_{M} & \Par_P \arrow[l] \arrow[r] & \Par_{G}
\end{tikzcd}
\]
with the left square being cartesian around $Z(\wh G)^{\Gamma, \circ}\times BZ(\wh G)^{\Gamma, \circ}\times \mc L^u (V_\lambda /H_{\lambda_\Sc})$. The right-cuspidality of the $!$-restriction of $\mc F$ follows from base change by an argument analogous to Proposition \ref{prop: PullbackIlambdaCuspidal} with equality of $Z(\wh{G})^{\Gamma, \circ}$ and $Z(I_{\lambda})^{\circ}$ following from Lemma \ref{lem: CuspILambdaCenter}. For the converse, it suffices to show that the $*$-pushforward preserves right-cuspidality, and it is analogous to the proof of Proposition \ref{prop: PushforwardIlambdaCuspidal}. 

Since the right vertical map is an \'etale cover around $Z(\wh G)^{\Gamma, \circ}\times BZ(\wh G)^{\Gamma, \circ}\times \mc L^u(V_\lambda / H_{\lambda_\Sc})$, the $!$-restriction is conservative on ind-coherent sheaves set-theoretically supported on $\V_{[\lambda]}$. 
\end{proof}

The right-cuspidal ind-coherent sheaves on $Z(\wh G)^{\Gamma, \circ}\times BZ(\wh G)^{\Gamma, \circ} \times \LS_{I_{\lambda, \Sc}}^\up$ set-theoretically supported on $Z(\wh G)^{\Gamma, \circ}\times BZ(\wh G)^{\Gamma, \circ}\times \mc L^u(V_\lambda /H_{\lambda_\Sc})$ are generated by sheaves $\mc O_{Z(\wh G)^{\Gamma, \circ}\times BZ(\wh G)^{\Gamma, \circ}}(\chi)\boxtimes \mc V_{\rho_\Sc}$: we have the equivalence
\[
\IndCoh(Z(\wh G)^{\Gamma,\circ}\times BZ(\wh G)^{\Gamma,\circ} \times \LS^\up_{I_{\lambda,\Sc}}) \simeq
\QCoh(Z(\wh G)^{\Gamma,\circ}\times BZ(\wh G)^{\Gamma,\circ}) \otimes \IndCoh(\LS^\up_{I_{\lambda,\Sc}}), 
\]
and Eisenstein series are $\QCoh(Z(\wh G)^{\Gamma,\circ}\times BZ(\wh G)^{\Gamma,\circ})$-linear, and the support condition implies the second factor is generated by $\mc V_{\rho_\Sc}$. In particular, by Lemma \ref{lem: PassingToSimplyConnectedCover}, the sheaves $\widetilde{\mc V}_{\chi\boxtimes \rho_\Sc}$ are right-cuspidal. Since the stable $\infty$-category of right-cuspidal ind-coherent sheaves on $\Par_G$ set-theoretically supported on $\V_{[\lambda]}$ is presentable and the $!$-pullback is conservative on this subcategory, it follows that the sheaves $\widetilde{\mc V}_{\chi\boxtimes \rho_\Sc}$ generate this category under colimits. By their coherence, they also generate the idempotent-complete stable $\infty$-category of right-cuspidal coherent sheaves supported on $\V_{[\lambda]}$. We have proved Theorem \ref{thm: Generators} (1).

To prove (2), first observe that it suffices to consider the case $\lambda\sim \lambda'$, and we may assume $\phi=\phi'$. Consider the map of completions
\[
Z(\wh G)^{\Gamma,\circ}\times BZ(\wh G)^{\Gamma,\circ} \times (\LS^\up_{I_{\lambda,\Sc}})^\wedge_{\V_{\lambda^\ur_\Sc}} \to \Par^\wedge_{G,\V_{[\lambda]}},
\]
which is (not necessarily representable) finite \'etale. To compute the relevant $\Hom$, we may work with the above map, since $\widetilde{\mc V}_\rho$, $\widetilde{\mc V}_{\rho'}$ are set-theoretically supported on $\V_{[\lambda]}$. Pick an irreducible summand $\chi\boxtimes \rho_{\Sc}$ of $\rho$ restricted to $Z(\wh G)^{\Gamma,\circ}\times S_{\phi^\ur_\Sc}$. The $*$-pullback of $\widetilde{\mc V}_{\rho}$ to $Z(\wh G)^{\Gamma, \circ}\times BZ(\wh G)^{\Gamma, \circ}\times \V_{\lambda^\ur_\Sc}$ is a finite direct sum of $\mc O_{Z(\wh G)^{\Gamma,\circ}}\boxtimes \chi_i\boxtimes \mc V_{\rho_i}$, where $\chi_i\boxtimes\rho_i$ is the conjugate of $\chi\boxtimes\rho_\Sc$ by some $g\in \wh G$ such that $g\lambda g^{-1}$ is an unramified twist of $\lambda$. By passing to $\IndCoh((\LS^\up_{I_{\lambda,\Sc}})^\wedge_{\V_{\lambda^\ur_\Sc}})$, it follows from Theorem \ref{thm: CVBclassification} (2) for $I_{\lambda,\Sc}$ that there is no nonzero map
\[
\mc O_{Z(\wh G)^{\Gamma,\circ}}\boxtimes \chi_i \boxtimes \mc V_{\rho_i} \to \mc O_{Z(\wh G)^{\Gamma,\circ}}\boxtimes \chi'_j\boxtimes \mc V_{\rho'_j}
\]
in nonzero degrees, where $\chi'_j$, $\rho'_j$ are defined similarly for $\rho'$. Moreover, in degree 0, we can compute over $\V_{[\lambda]}$. Then the adjunction for $Z(\wh G)^{\Gamma,\circ}\times \V_\lambda \to \V_{[\lambda]}$ and Theorem \ref{thm: CVBclassification} (2) imply that $\Hom$ is nonzero if and only if the condition in the first claim holds. The second claim also follows since $\mc F$ is a retract of a finite colimit of sheaves $\widetilde{\mc V}_{\rho'}$ by (1), and we may work with a single $\rho'$ by the first claim of (2), and then it is clear by adjunction. 

Finally, we prove (3). By the proof of (2), 
\[
\End_{\Par_G} (\widetilde{\mc V}_\rho) \cong \End_{\V_{[\lambda]}} (\widetilde{\mc V}_\rho). 
\]
Now, the right-hand side identifies with the endomorphism ring of the pushforward $\mc E$ of $\mc O_{Z(\wh G)^{\Gamma,\circ}}\boxtimes \mc V_\rho$ to 
\[
(Z(\wh G)^{\Gamma, \circ}\times \V_\lambda )/\{z\in Z(\wh G)^{\Gamma, \circ}\mid \exists g, g\lambda g^{-1}=z\lambda,  \mathrm{Ad}(g)^* \mc V_\rho\cong \mc V_\rho\}. 
\] 
The closed substack determined by the condition $N=0$ identifies with
\[
(Z(\wh G)^{\Gamma, \circ} \times BH_\lambda)/\{z\in Z(\wh G)^{\Gamma, \circ}\mid \exists g, g\lambda g^{-1}=z\lambda,  \mathrm{Ad}(g)^* \mc V_\rho\cong \mc V_\rho\}. 
\] 
Therefore, the canonical map
\[
\End(\mc E) \to \End (\mc E|_{N=0})
\]
has a section. Since $\End_{\V_\lambda}(\mc V_{\rho})\cong \Lambda$, it is in fact an isomorphism. Using the adjunction, the algebra $\End (\mc E|_{N=0})$ can be described as a crossed product
\[
\mc O(Z(\wh G)^{\Gamma,\circ}) \rtimes_{\mu_\rho} \{z\in Z(\wh G)^{\Gamma, \circ}\mid \exists g, g\lambda g^{-1}=z\lambda,  \mathrm{Ad}(g)^* \mc V_\rho\cong \mc V_\rho\}
\]
as in Remark \ref{rem: Azumaya}, and we conclude that $\End(\widetilde{\mc V}_\rho)$ is an Azumaya algebra over its center (e.g.~\cite[Corollary 3.2.7]{AbeOhara}), which has the desired form. 

\subsection{Summary}
If Theorem \ref{thm: ClassificationCuspidalCotangent} holds for split semisimple groups, or even simply-connected ones, and unramified parameters arising from maps $\Gm\to \wh G$, then we obtain all classification theorems both on $\wh{\T}\V_\lambda[-1]$ and $\T^* \V_\lambda$. 

\section{Cuspidal sheaves and the geometry of the cotangent stack}{\label{s: cotangentgeometry}}

We start to prove Theorems \ref{thm: ClassificationCuspidalCotangent} and \ref{thm: BasicCotangentMain}, working in the graded setting. 

We assume throughout the following simplified situation as in \S \ref{ss: lusztiggraded}. Let $G$ be a split reductive group with a cocharacter $\iota$ of $\wh{G}$. For simplicity, we may assume $\iota$ factors through $\wh{G}_{\der}$. Fix Borel pairs $(B, T)$ and $(\wh{B}, \wh{T})$, and we assume $\iota$ factors through $\wh{T}$. Let $H$ denote the centralizer of $\iota$, $V$ the weight $n$ subspace of $\Lie \wh{G}$ for an integer $n > 0$. There are finitely many $H$-orbits in $V$ \cite[Proposition 3.5]{LusztigGraded}, and it is a prehomogeneous vector space. We denote $V/H$ by $\V$.

In this section, all Eisenstein series are taken only for $\lambda_{\wh{P}}=\lambda$, so we will write $\Eis^{\T^*}_{\wh P}$ or $\Eis^{\T^*}_{P}$. 

\begin{remark}
This section relies crucially on Halpern-Leistner's theory developed in \cite{HalpernLeistner} and his main results are based on \cite[Appendix A]{HalpernLeistner}. Unfortunately, some statements and proofs in \cite[Appendix A]{HalpernLeistner} are problematic, as pointed out in \cite[Remark 2.5.13(i)]{HansenMann}. In particular, the proof of \cite[Theorem A.0.1]{HalpernLeistner} is incorrect and we think Part (3) is false. The parts except (3) and (4) follow from \cite[Proposition 2.5.12]{HansenMann}. Although the proof of \cite[Corollary A.0.1.1(1)]{HalpernLeistner} mentions Part (3), it actually follows from \cite[Corollary A.0.1.1(2)]{HalpernLeistner} whose proof itself works thanks to \cite[Proposition 2.5.12]{HansenMann}, since the cotangent complex of a closed immersion vanishes in non-negative degrees. Part (3) is not used in any other place in our understanding. 

In fact, Part (4) at least for QCA stacks $X$ also follows from \cite{HansenMann}: it reduces to the truncated case (the first half). Recall that $\Psi_X\colon \IndCoh(X)\to \QCoh(X)$ induces an equivalence on the bounded below subcategories. The corresponding statement under this equivalence is proven in \cite[Theorem 2.10.7]{HansenMann}.
\end{remark}

\subsection{$\Theta$-stratifications of cotangent stacks}

\subsubsection{Weak $\Theta$-actions and baric structures}{\label{sss: derivedtheta strat}}
We first recall the theory of weak $\Theta$-actions, baric structures, and (derived) $\Theta$-stratifications as developed in \cite{HalpernLeistner}. We denote $\Theta \coloneq \A^1/ \Gm$, where $\A^1$ carries the weight $1$ action of $\Gm$ and we study actions of the monoid $\Theta$ on a stack $X$. In loc.~cit., this is known as a \emph{weak} $\Theta$-\emph{action} (\cite[Definition 1.1.1]{HalpernLeistner}). A weak $\Theta$-action induces a \emph{baric structure} on $\QCoh(X)$ which is a family of semiorthogonal decompositions indexed by $w \in \Z$,
\begin{equation*}
    \QCoh(X) = \langle \QCoh(X)^{<w}, \QCoh(X)^{\geq w}\rangle,
\end{equation*}
with truncation functors $\beta^{\geq w}, \beta^{<w}$ such that any object $\mc{F} \in \QCoh(X)$ sits in compatible fiber sequences
\begin{equation*}
    \beta^{\geq w}(\mc{F}) \rightarrow \mc{F} \rightarrow \beta^{<w}(\mc{F}).
\end{equation*}
The truncation functors induce a baric structure on $\Aperf^*(X)$ \cite[Proposition 1.1.2]{HalpernLeistner}.

\begin{lemma}{\label{lem: baricnaturality}}
Baric structures satisfy natural compatibilities: namely, if $f: X \to Y$ is $\Theta$-equivariant, then (\cite[Proposition 1.1.2]{HalpernLeistner})
    \begin{itemize}
    \item $f^{\ast}(\QCoh(Y)^{\geq w}) \subset \QCoh(X)^{\geq w}$,
    \item $f^{\ast}(\QCoh(Y)^{< w}) \subset \QCoh(X)^{< w}$,
    \item $f_{\ast}(\QCoh(X)^{< w}) \subset \QCoh(Y)^{< w}$.
\end{itemize}
\end{lemma}

For simplicity, assume $X$ is the quotient of a quasi-affine scheme almost of finite type by an affine algebraic group. Consider the stack of graded (resp.~filtered) points of $X$ given by $\Map(B\Gm, X)$ (resp.~$\Map(\Theta, X)$). Evaluation at the open and closed points of $\Theta$ gives a diagram
\begin{equation*}
    \begin{tikzcd}
        \Map(B\Gm, X) \arrow[r, swap, bend right=15, "{\sigma}"] & \Map( \Theta, X) \arrow[r, "{\ev_1}"] \arrow[l,swap, bend right=15, "{\ev_0}"] & X,
    \end{tikzcd}
\end{equation*}
such that $\sigma$ is induced by the projection $\Theta \rightarrow B\Gm$ and satisfies $\ev_0 \circ \sigma \simeq \id$. Note also that $\Map(\Theta, X)$ carries a canonical weak $\Theta$-action by pre-composition. 
\begin{example}{\label{ex: TVnugeq0}}
   The following examples are ubiquitous in this section. By \cite[Theorem 1.4.8]{HLInstability}
    \begin{equation*}
       \Map(\Theta, \V) \simeq \bigsqcup\limits_{\nu \in X_{\ast}(\wh T)/ W_H} \V^{\nu \geq 0}    \end{equation*}
    where $\V^{\nu \geq 0}$ is simply the quotient of the (smooth) subspace $V^{\nu \geq 0}$ of $V$ with non-negative $\nu$-weights by $H_{P_\nu}\coloneq H\cap P_\nu$ (cf.~\cite[Propositions 1.4.1, 1.4.5]{HLInstability}). Similarly, using the 0-fiber description of $\T^*\V$ via the moment map
    \[
    \mu\colon \T^* V=V\times V^\vee \to (\Lie H)^\vee, \qquad \T^*\V\simeq \mu^{-1}(0)/H,
    \]
    we can also deduce from \cite[Theorem 1.4.8]{HLInstability} that 
    \begin{equation*}
\Map(\Theta, \T^{\ast}\V) \simeq \bigsqcup\limits_{\nu \in X_{\ast}(\wh T)/ W_H} \T^{\ast}\V^{\nu \geq 0}, 
    \end{equation*}
    where $\T^{\ast}\V^{\nu \geq 0}$ is the base change of the map induced by $\mu$, $(\T^*V)^{\nu \geq 0}/H_{P_\nu}\to (\Lie H)^{*, \nu \geq 0}/H_{P_\nu}$ along $BH_{P_\nu}\xrightarrow{0} (\Lie H)^{*, \nu \geq 0}/H_{P_\nu}$. Note that these are quotients of quasi-smooth affine schemes by affine algebraic groups, and $\T^{\ast}\V^{\nu \geq 0}$ (resp.~$\V^{\nu \geq 0}$) carries a natural weak $\Theta$-action and hence its quasi-coherent sheaves have a baric structure and we have a diagram
    \begin{equation*}
        \begin{tikzcd}
\T^{\ast}\V^{\nu = 0} \arrow[r, swap, bend right=15, "{\sigma}"] & \T^{\ast}\V^{\nu \geq 0} \arrow[l,swap, bend right=15, "{\ev_0}"] \arrow[r, "\ev_1"]& \T^{\ast}\V
        \end{tikzcd},
    \end{equation*}
    and similarly for $\V^{\nu \geq 0}$. The map $\ev_0$ is quasi-smooth and we will take $*$-pullback of coherent sheaves freely. The map $\ev_1$ is proper schematic by \cite[Lemma 5.1.7]{HLInstability}. In the situation of (derived) $\Theta$-strata that we discuss later, $\sigma$ is quasi-smooth, cf.~\cite[Lemma 3.2.7, Example 3.3.2]{HalpernLeistner}. If $\nu$ determines a parabolic subgroup $P$, we define $\Eis^*_P\coloneq \ev_{1*}\circ\ev_0^*$ or just write $\Eis^*$ if the context is clear. Note that $\Eis^*$ preserves coherent sheaves.  

    In this setting, we have the following useful compatibilities of baric structures.
    \begin{lemma}{\label{lem: HLlemmas}}
    Let $w\in \Z$.
        \begin{enumerate}
            \item~(\cite[Lemma 1.5.4]{HalpernLeistner})~A complex $\mc{F} \in \Aperf^*(\T^{\ast}\V^{\nu \geq 0})$ lies in $\Aperf^*(\T^{\ast}\V^{\nu \geq 0})^{\geq w}$ (resp.~$\Aperf^*(\T^{\ast}\V^{\nu \geq 0})^{< w}$) if and only if $\sigma^{\ast} \mc{F} \in \Aperf^*(\T^{\ast}\V^{\nu = 0})^{\geq w}$ (resp.~$\sigma^{\ast} \mc{F} \in \Aperf^*(\T^{\ast}\V^{\nu = 0})^{< w}$).
            \item~(\cite[Lemma 1.5.6]{HalpernLeistner})~Define \[\Aperf^*(\T^{\ast}\V^{\nu \geq 0})^{w} \coloneq \Aperf^*(\T^{\ast}\V^{\nu \geq 0})^{\geq w} \cap \Aperf^*(\T^{\ast}\V^{\nu \geq 0})^{<w +1}\] and analogously for $\Aperf^*(\T^{\ast}\V^{\nu = 0})^{w}$. Then $\ev^{\ast}_0$ induces an equivalence
\begin{equation*}
    \Aperf^*(\T^{\ast}\V^{\nu = 0})^{w} \simeq \Aperf^*(\T^{\ast}\V^{\nu \geq 0})^{w},
\end{equation*}
with inverse given by either $\beta^{\geq w} \circ \ev_{0 \ast}$ or $\sigma^*$. The same holds for quasi-coherent sheaves.
    \item $\QCoh(\T^*\V^{\nu\geq 0})^{\geq w}\subset \QCoh(\T^*\V^{\nu\geq 0})$ is generated under colimits by $\ev_0^* \mc F$ for $\mc F\in \Perf(\T^*\V^{\nu=0})^{w'}$ with $w' \geq w$.  
        \end{enumerate}
    \end{lemma}
    \begin{proof}
        For (1) and (2), the proofs go through from loc.~cit. We note that Halpern-Leistner works with a $\Theta$-stratum but this is not necessary for his arguments applied to the present situation.

        Let us prove (3). First note that $\QCoh(\T^*\V^{\nu\geq 0})$ is generated under colimits by vector bundles $\mc V$. By (1), any such $\mc V$ has bounded $\nu$-weights, and it follows by applying (2) to truncations that it is a successive extension of vector bundles pulled back from $\T^*\V^{\nu=0}$, and hence attached to representations of $H_{M_\nu}\coloneq H\cap M_\nu$. It follows that $\QCoh(\T^*\V^{\nu\geq 0})^{\geq w}$ is generated under colimits by $\ev_0^{*} \mc V$ for vector bundles $\mc V$ with $\nu$-weights $\geq w$ as desired. 
    \end{proof}

    On the space $\T^{\ast} \V^{\nu = 0}$, we have by \cite[Lemma 1.5.3]{HalpernLeistner} the following compatibility with the weight decompositions of \S \ref{ss: weightdecomposition} (with notation analogous to the above)
    \begin{equation}{\label{eqn: bariccentralweightcompatibility}}
\QCoh(\T^{\ast}\V^{\nu=0})^w \simeq  \prod_{\substack{\chi \in X^{\ast}(Z(M_{\nu})^{\circ});\\ \langle \nu, \chi \rangle = w}}\QCoh(\T^{\ast}\V^{\nu=0})^{\chi}.
    \end{equation}
    Suppose that $\nu' = z + \nu$ where $z \in X_{\ast}(Z(\wh{G}))$. Then we have 
    \begin{equation*}
        \T^{\ast}\V^{\nu \geq 0} \simeq \T^{\ast}\V^{\nu' \geq 0}.
    \end{equation*}
    If $\mc{F} \in \Aperf^*(\T^{\ast}\V^{\nu' \geq 0})^{\leq w}$ with $Z(\wh G)^\circ$-weight $\chi$, then by part (1) of Lemma \ref{lem: HLlemmas} and the equivalence \eqref{eqn: bariccentralweightcompatibility},
    \begin{equation}{\label{eqn: basiccentralshift}}
        \mc{F} \in \Aperf^*(\T^{\ast}\V^{\nu \geq 0})^{\leq w - \langle z, \chi \rangle},
    \end{equation}
    and analogously for $\mc{F} \in \Aperf^*(\T^{\ast}\V^{\nu' \geq 0})^{\geq w}$.
\end{example}
\begin{lemma}{\label{lem: eiscotanfiberprodcomparison}}
    We continue with the notation of Example \ref{ex: TVnugeq0}. Then the stack $\T^{\ast}\V^{\nu \geq 0}$ identifies with $Y^{\geq 0}$ defined in \S\ref{ss: eisensteinseriesforcotangents}, where the relevant parabolic is $P_\nu$. 
\end{lemma}
\begin{proof}
Recall that the space $Y^{\geq 0}$ is defined via a fiber product
\begin{equation*}
    \begin{tikzcd}
        & Y^{\geq 0}  \arrow[dl] \arrow[dr]& \\
        \T^{\ast}\V_{M_\nu}|_{\V_{P_\nu}} \arrow[dr] & & \T^{\ast}\V|_{\V_{P_\nu}}\arrow[dl] \\
        & \T^{\ast} \V_{P_\nu} &
    \end{tikzcd}.
\end{equation*}
Each object in the above diagram is the quotient by $H_{P_\nu}$ of the $0$-fiber of the moment map applied to some object, and bottom maps lift before taking the $0$-fiber. Since fiber products commute, we can compute the fiber product defining $Y^{\geq 0}$ before taking $0$-fiber.

Note that $(\Lie H_{M_\nu})^\vee \times_{(\Lie H_{P_\nu})^\vee} (\Lie H)^\vee \simeq (\Lie H_{P_\nu^{\op}})^\vee$ and $V^\vee_{M_\nu} \times_{ V^\vee_{P_\nu}} V^\vee \simeq V^\vee_{P_\nu^{\op}}$. Rewriting, $Y^{ \geq 0}$ becomes
\begin{equation*}
    [\{0\} \times_{(\Lie H_{P_\nu^{\op}})^\vee} (V_{P_\nu} \times V^\vee_{P_\nu^{\op}})] / H_{P_\nu}.
\end{equation*}

We now compare this with $\T^{\ast} \V^{\nu \geq 0}$. As in Example \ref{ex: TVnugeq0}, it is given by the quotient by $H_{P_\nu}$ of the fiber $(\mu^{\geq 0})^{-1}(0)$, where $\mu^{\geq 0}: (\T^{\ast}V)^{\nu \geq 0} \rightarrow (\Lie H)^{\vee, \nu \geq 0}$ is the restriction of $\mu$, which identifies with the above (note $(\Lie H)^{\vee, \nu \geq 0} = (\Lie H^{\nu \leq 0})^\vee$).
\end{proof}
\begin{corollary}{\label{cor: eiscompositioncotan}}
Let $\wh{P}_1 \subset \wh{P}_2 \subset \wh{G}$ be semistandard  parabolic subgroups with Levi decompositions
\begin{equation*}
    \wh{P}_1 = \wh{M}_1\wh{U}_1, \qquad \wh{P}_2 = \wh{M}_2\wh{U}_2,
\end{equation*}
so that $\wh{M}_1 \subset \wh{M}_2$. Let $\wh{P}'_1 = \wh{P}_1 \cap \wh{M}_2$. Then we have
\begin{equation*}
    \Eis^{\ast}_{P_1} \simeq \Eis^{\ast}_{P_2} \circ \Eis^{\ast}_{P'_1}, \qquad   \Eis^{\T^{\ast}}_{P_1} \simeq \Eis^{\T^{\ast}}_{P_2} \circ \Eis^{\T^{\ast}}_{P'_1}.
\end{equation*}
\end{corollary}
\begin{proof}
    The result for $\Eis^{\ast}$ follows from the assertion that
    \begin{equation*}
        Y^{\geq 0}_{\wh{P}'_1} \times_{\T^{\ast}\V_{\wh{M}_2}} Y^{\geq 0}_{\wh{P}_2} \simeq Y^{\geq 0}_{\wh{P}_1}.
    \end{equation*}
    By the computations in Lemma \ref{lem: eiscotanfiberprodcomparison}, we can compute this fiber product before taking $0$-fiber. Then the result follows from the fact that the following diagram is cartesian
    \begin{equation*}
        \begin{tikzcd}
            & V_{\wh{P}_1} \times V^{\vee}_{\wh{P}^{\op}_1} \arrow[dl] \arrow[dr] & \\
         V_{\wh{P}'_1} \times V^{\vee}_{\wh{P}^{' \op}_1}  \arrow[dr] & &  V_{\wh{P}_2} \times V^{\vee}_{\wh{P}^{\op}_2} \arrow[dl] \\
            & V_{\wh{M}_2} \times V^{\vee}_{\wh{M}_2} & 
        \end{tikzcd}
    \end{equation*}
For $\Eis^{\T^{\ast}}$, we use the result for $\Eis^{\ast}$ plus projection formula to combine the twists, using that $\omega_{\V_{U_1}} = \omega_{\V_{U'_1}}\omega_{\V_{U_2}}$.
\end{proof}

Now let $U/H = \U \subset \V$ denote an open substack and let $j: \T^{\ast}\U \rightarrow \T^{\ast}\V$ be the induced map on cotangent stacks. Let $\nu$ be a cocharacter of $\wh T$. This induces a weak $\Theta$-action on $(\T^{\ast}\V)^{\nu \geq 0}$ and we can consider the commutative diagram 
\begin{equation}{\label{eqn: naturalitydiagram}}
    \begin{tikzcd}
\T^{\ast}\U^{\nu = 0} \arrow[d, "j^{=0}"] & \T^{\ast}\U^{\nu \geq 0} \arrow[r, "p_\U"] \arrow[d, "j^{\geq 0}"] \arrow[l, swap, "q_\U"]&  \T^{\ast}\U \arrow[d, "j"] \\
 \T^{\ast}\V^{\nu = 0} & \T^{\ast}\V^{\nu\geq 0} \arrow[r, "p"] \arrow[l, swap, "q"] & \T^{\ast}\V .
    \end{tikzcd}
\end{equation}
Note that $p_\U$ may not be proper and it is so if and only if $\U^{\nu\geq 0}\to \U$ is proper if and only if $U^{\nu\geq 0}$ is closed in $U$. Here, $U^{\nu \geq 0}$ consists of the points of $U$ whose contracting limit under $\nu$ lies in $U_{\wh{M}} \coloneq U \cap \Lie \wh{M}$, and could be smaller than $V^{\nu\geq 0}\cap U$. 

\begin{lemma}[{cf.~\cite[Lemma 6.10]{PTDolbeault}}]{\label{lem: eisfiberprod}}
    \begin{enumerate}
\item The left square is cartesian. If $\U^{\nu = 0}$ is nonempty and $p_\U$ is proper, then the right square is cartesian.
\item Let $\Eis^{\ast}$ be $p_{\ast}q^{\ast}$ and $\Eis^{\ast}_\U$ be $p_{\ast \U}q^{\ast}_\U$. Then 
\begin{equation*}
    j_{\ast}\Eis^{\ast}_\U \simeq \Eis^{\ast} j^{=0}_{\ast},
\end{equation*}
and if $\U^{\nu = 0}$ is nonempty and $p_\U$ is proper, then
\begin{equation*}
        j^{\ast}\Eis^{\ast} \simeq \Eis^{\ast}_\U (j^{=0})^{\ast}.
\end{equation*}

    \end{enumerate}
\end{lemma}
\begin{proof}
To prove the left square is cartesian, we may carry out the computation before taking $0$-fiber. Let $\wh{P}$ (resp.~$\wh{M}$) be the parabolic (resp.~Levi) subgroup determined by $\nu$. Using computations from the proof of Lemma \ref{lem: eiscotanfiberprodcomparison}, we need to check the following diagram is cartesian
\begin{equation*}
    \begin{tikzcd}
        U_{\wh{M}} \times V^{\vee}_{\wh{M}} \arrow[d] & U^{\nu \geq 0} \times V^{\vee}_{\wh{P}^{\op}}  \arrow[l] \arrow[d] \\
        V_{\wh{M}} \times V^{\vee}_{\wh{M}} & V_{\wh{P}} \times V^{\vee}_{\wh{P}^{\op}} \arrow[l].
    \end{tikzcd}
\end{equation*}
Then the left square is cartesian since $U$ is open. 

Now we consider the right square of \eqref{eqn: naturalitydiagram}. We need to show \begin{equation}{\label{eqn: cotanmap}}
    \T^{\ast}\U^{\nu \geq 0} \rightarrow \T^{\ast}\V^{\nu \geq 0} \times_{\T^{\ast}\V} \T^{\ast}\U,
\end{equation} 
is an equivalence. Note that $\T^{\ast}\U^{\nu \geq 0}$ is nonempty if $\U^{\nu=0}$ is nonempty. The map is an open immersion since both the source and target are open substacks of $\T^{\ast}\V^{\nu \geq 0}$. On the other hand, $\T^{\ast}\V^{\nu \geq 0} \times_{\T^{\ast}\V} \T^{\ast}\U$ is proper over $\T^{\ast}\U$ as it is the base change of a proper map. Hence, \eqref{eqn: cotanmap} is proper if $p_{\U}$ is proper. Finally, the target is connected as it contracts to $\V^{\nu \geq 0} \cap \U$, which is a quotient of a nonempty Zariski open subset of an affine space. This implies \eqref{eqn: cotanmap} is an isomorphism under the stated assumptions. 

For compatibility of $\Eis^{\ast}$ with $j_{\ast}$ (resp.~$j^{\ast}$) we use base-change of the left (resp.~right) square.
\end{proof}

\begin{remark}
Suppose $\U=Hc^\circ /H$ is the open orbit and $c^\circ \in \Lie \wh M$. When using twists, note also that $\omega_{\V_U}$ restricts to $\omega_{B(H_{c^\circ} \cap \wh P)}\otimes \omega_{B(H_{c^\circ} \cap \wh M)}^{-1}$. 
\end{remark}

\subsubsection{$\Theta$-stratification of the cotangent stack}\label{sss:ThetaStratificationCotangent}
For simplicity, assume $X$ is the quotient of a quasi-affine scheme almost of finite type by an affine algebraic group. 
A (derived) $\Theta$-\emph{stratum} is a (nonempty) union $\mc{S}$ of connected components of $\Map(\Theta, X)$ such that $\ev_1|_{\mc{S}}: \mc{S} \rightarrow X$ is a closed immersion \cite[Definition 1.2.2]{HalpernLeistner}. The \emph{center} of a $\Theta$-stratum $\mc{S}$ is the open and closed substack of $\Map(B\Gm, X)$ given by $\mc{Z} \coloneq \sigma^{-1}(\mc{S})$ and $\sigma$ and $\ev_0$ induce a diagram 
\begin{equation*}
    \begin{tikzcd}
        \mc{Z} \arrow[r, swap, bend right=15, "{\sigma}"] & \mc{S} \arrow[r, "{\ev_1}"] \arrow[l,swap, bend right=15, "{\ev_0}"] & X
    \end{tikzcd}.
\end{equation*}

A (derived) $\Theta$-\emph{stratification} of $X$ \cite[Definition 2.3.1]{HalpernLeistner} is a totally ordered set $I$, an increasing collection of open substacks $X_{\leq \alpha} \subset X$ for $\alpha \in I$, and for each $X_{\leq \alpha}$, a $\Theta$-stratum $\mc{S}_{\alpha}$ such that 
\begin{equation*}
   X_{\leq \alpha} \setminus \ev_1(|\mc{S}_{\alpha}|) = \bigcup\limits_{\alpha' < \alpha} X_{\leq \alpha'}. 
\end{equation*}

We now return to the previous setting. Let $\mathcal{C}$ denote the set of $H$-orbits in $V$ with a partial order induced from the closure relation; $\{0\}$ is the unique maximal orbit and the open orbit of $V$ is the unique minimal one. We fix a total order of $\mathcal{C}$ respecting the partial order.

We first claim that $\mathcal{C}$ induces a $\Theta$-stratification of $\V$ with $\Theta$-strata given by the points of $\V$ which correspond to $\mathcal{C}$. We formulate a more precise result as follows. Pick $c \in C \in \mathcal{C}$, and let $\iota_c$ denote the cocharacter of $\wh{G}$ associated with $c$ via the Jacobson--Morozov theorem; if $c=0$, we take $\iota_0$ to be the trivial cocharacter. Taking an $H$-conjugate if necessary, we may assume $\iota_c$ factors through $\wh T$ and we let $P_c=M_c U_c$ denote the semistandard parabolic subgroup of $\wh{G}$ associated with $\nu_c \coloneq n\iota_c - 2\iota$; its Lie algebra is precisely the one considered in \cite[\S 5]{LusztigGraded}. 

\begin{proposition}{\label{prop: nuthetastrat}}
We write $H_{P_c}=H\cap P_c$, $H_{M_c}=H\cap M_c$, $V_{P_c}=V\cap \Lie P_c$, $V_{M_c}=V\cap \Lie M_c$. 
\begin{enumerate}
    \item $c\in V_{M_c} \subset V_{P_c}$ and is fixed by the adjoint action of $\nu_c$. 
    \item $H_c \coloneq Z_{H}(c)\subset P_c$ and $H_c=(H_c \cap M_c)(H_c \cap U_c)$. 
    \item The identity component $H_{M,c}^{\circ}$ of $H_{M,c}\coloneq H_c \cap M_c$ is reductive and $H_{M,c}^{\circ}$ is the reductive quotient of $H^{\circ}_c$. 
    \item $BH_c=BZ_{H_{P_c}}(c)\to V_{P_c}/H_{P_c}$ is an open dense substack. 
    \item Under the natural map $V_{P_c}/H_{P_c}\to V/H$, the inverse image of $\bigcup_{C'> C}C'/H$ is the complement of $BZ_{H_{P_c}}(c)\subset \V_{P_c}$. In particular, $C/H\simeq BH_c$ is a (smooth) $\Theta$-stratum there and its center is $BZ_{H_{M_c}}(c)=BH_{M,c}$.
\end{enumerate}
\end{proposition}

\begin{proof}
This is proven in \cite{LusztigGraded} and we shall cite results there. 
(1) is clear. The first part of (2) is Proposition 5.4 and the second part follows from (1) since $\nu_c$ factors through $H_c$ and hence induces a Levi decomposition on $H_c$. Part (3) follows from 5.3(a), Proposition 4.2 and (2). Part (4) is Proposition 5.9. The map $V_{P_c}/H_{P_c}\to V/H$ is proper, and its image is closed. As the image contains $C/H$ and is contained in the closure $\overline{C}/H$ of $C/H$ by (4), the image agrees with $\overline{C}/H$. If $hc\in V_{P_c}$ for some $h\in H$, then $h\in H_{P_c}$ by Lemma 5.7. This implies (5):
\[
H_{P_c}c=
V_{P_c} \cap (V\setminus \bigcup_{C'> C}C')\supset
(V\setminus \bigcup_{C'> C}C')^{\nu_c \geq 0}
\]
is an equality as the set of limit points of $H_{P_c}c$ with respect to $\nu_c$ is $H_{M_c}c\subset C$. Note that the first part of (2) implies $H_{P_c}c/H_{P_c} = C/H \subset (V \setminus \bigcup\limits_{C' > C} C')/H$ is a closed immersion.
\end{proof}

Recall that the zero fiber $\mu^{-1}(0)$ is the union of the conormal bundles $\T_C^* V$ of $C\in \mathcal{C}$. In particular, the closure of $\T_C^* V$ is contained in $\bigcup_{C'\geq C} \T_{C'}^*V$. We claim that it induces a (derived) $\Theta$-stratification of $\mu^{-1}(0)/H$ as proven below.

We canonically identify $\T^*V\cong V\times V^\vee$ with a subspace of $\Lie \wh{G}$. 

\begin{proposition}{\label{prop: TVthetastrat}}
\begin{enumerate}
    \item Any $d$ such that $(c,d) \in \T_{C, c}^* V \cap (\Lie P_c \times \Lie P_c)$ is contracting with respect to $\nu_c$, and $Z_{H}(c,d)\subset P_c$. 
    \item Under the natural map (at the level of reduced classical stacks) \[(\mu^{-1}(0)\cap (\Lie P_c\times \Lie P_c))/H_{P_c}\to \mu^{-1}(0)/H,\]  the inverse image of $(\bigcup_{C'> C}\T_{C'}^*V)/H$ is the complement of $(\T_{C}^* V\cap (\Lie P_c\times \Lie P_c))/H_{P_c}$. Therefore, by \cite[Lemma 1.2.3]{HalpernLeistner}, $\T_C^* V/H$ inside $(\mu^{-1}(0)\setminus \bigcup_{C'> C}\T_{C'}^*V)/H$ upgrades uniquely to a $\Theta$-stratum. Its center is the open substack $\T^* BH_{M,c}$ of $\T^* \V_{M_c}$ determined by $(c,0)$.  
\end{enumerate}
\end{proposition}

\begin{proof}
(1) Since $d\in (\Lie \wh{G})^c$, the $\Gm$-action on $d$ via $\iota_c$ has a non-negative weight \cite[Lemma 3.4.3]{CollingWoodMcGovern}, and hence the $\nu_c$-weight is positive as $d\in V^\vee\subset \Lie \wh G$. The rest follows as $H_c \subset P_c$ by (2) of Proposition \ref{prop: nuthetastrat}.

(2)  From the previous results, we have, at the level of reduced closed subschemes, 
\[
(\T_C^* V \cap (H_{P_c} c\times V^\vee))=\mu^{-1}(0)^{\nu_c\geq 0}\setminus \bigcup_{C'> C}\T_{C'}^*V
\hookleftarrow 
(\mu^{-1}(0)\setminus \bigcup_{C'> C}\T_{C'}^*V)^{\nu_c\geq 0},
\]
and it is an isomorphism as the limit points in the leftmost scheme with respect to $\nu_c$ are in $H_{M_c} c \times \{0\}$. This means $(\mu^{-1}(0)^{\nu_c\geq 0}\setminus \bigcup_{C'> C}\T_{C'}^*V)/ H_{P_c}$ with the derived structure is a $\Theta$-stratum (cf.~\cite[Lemma 1.2.3]{HalpernLeistner}). The claim for the center is now clear. 
\end{proof}

Therefore, there are unique derived structures on these strata coming from cocharacters $\nu_c$ so that they form a $\Theta$-stratification.

\subsection{A characterization of basic sheaves}
We briefly remark on the various dualizing complexes. We have computed the dualizing complex $\omega_{\V}$ in Proposition \ref{prop: cotangentdualizingcomplex}. It is a line bundle up to a shift and corresponds to a character of $H$. Its total ($\wh T$-)weight is given by
\begin{equation*}
   \wt(\omega_{\V}) \coloneq  -\sum\limits_{v \in \wt(V)} v + \sum\limits_{w \in \wt(\Lie H)} w.
\end{equation*}
For a subgroup $G' \subset \wh{G}$ normalized by $\wh T$ inducing an inclusion $\Lie G' \subset \Lie \wh{G}$, we denote by $\wt(\V_{G'})$ the following weight 
\begin{equation*}
    \wt(\V_{G'}) \coloneq    -\sum\limits_{v \in \wt(V \cap \Lie G')} v + \sum\limits_{w \in \wt(\Lie H \cap \Lie G')} w.
\end{equation*}
Thus, $\wt(\V_{\wh{U}}) = \wt(\omega_{\V_{\wh{U}}})$, for instance. We write $\wt(V_{G'})$ for the negative of the first term in the above equation (so it is compatible with the weights of a representation) and $\wt(\Lie H \cap \Lie G')$ for just the second.

We have the following characterization of basic and cuspidal coherent sheaves.

\begin{proposition}{\label{prop: basicweightcharacterization}}
Let $b \in B(G)_{\bas}$ and $\mc F \in \Coh(\T^*\V)^{\kappa (b)}$. Then $\mc{F}$ is right-basic if and only if for all cocharacters $\nu$ defining proper parabolic subgroups $\wh P_{\nu} \subset \wh G$, 
\[
p^{!,\QCoh} \mc F \in \QCoh(\T^* \V^{\nu\geq 0})^{\leq \langle \nu, \wt(\V_{ U}) \rangle + \langle \nu, \nu_G(b)\rangle}, 
\]
where we make sense of the pairing $\langle \nu, \nu_G(b) \rangle$ by considering $\nu \in X^{\ast}(T)$ and $U$ is the unipotent radical of $\wh{P}_{\nu}$. For the right-basicness, it suffices to check that, for each maximal $\wh{P}$, there exists a defining $\nu$ such that $p^{!,\QCoh}\mc{F}$ satisfies the weight inequality relative to $\nu$.

Suppose moreover $\mc F \in \Perf(\T^*\V)^{\kappa (b)}$. Then $\mc{F}$ is right-basic if and only if for all cocharacters $\nu$ above, 
\[
p^! \mc F \in \QCoh(\T^* \V^{\nu\geq 0})^{\leq \langle \nu, \wt(\V_{ U}) \rangle + \langle \nu, \nu_G(b)\rangle},
\]
and by the twisted Grothendieck--Serre duality $\D'_{\GS}\coloneq \D_\GS\otimes \omega_{\V}$, left-basicness is equivalent to the condition that, for all $\nu$,
\[
p^* \mc F \in \QCoh(\T^{\ast} \V^{\nu\geq 0})^{\geq \langle \nu, \wt(\V_{U})\rangle+ \langle \nu, \nu_G(b)\rangle}. 
\]

The statement for a left/right-cuspidal $\mc{F}$ is analogous but the above weight inequalities are strict. 
\end{proposition}
Note that since $b$ is basic, we also have
\begin{equation*}
    \langle \nu, \nu_G(b) \rangle = \langle \nu, \kappa(b) \rangle,
\end{equation*}
where we view $\kappa(b)$ as an element of $X_{\ast}(T)_{\Q}$ via
\begin{equation*}
    X^{\ast}(Z(\wh{G}))_{\Q} \simeq X_{\ast}(A_G)_{\Q} \rightarrow X_{\ast}(T)_{\Q},
\end{equation*}
for $A_G$ the maximal split torus in the center of $G$.
\begin{proof}
We shall prove the claim for $\mc F\in \Coh(\T^*\V)^{\kappa (b)}$ and the claims for perfect complexes follow immediately. 

First, by scaling by an appropriate central cocharacter as in \eqref{eqn: basiccentralshift} and subtracting the central part, it suffices to prove the statement for $\nu$ mapping into $\wh{G}_{\der}$ and in this case the $\langle \nu, \nu_G(b)\rangle$ term is $0$. Secondly, by Corollary \ref{cor: eiscompositioncotan}, it suffices to check basicness/cuspidality for maximal parabolics. Finally, we will show that the $\Eis$-vanishing as in the definition of basic or cuspidal sheaves relative to a fixed parabolic $\wh{P}$ is equivalent to the weight condition as in the statement for any such cocharacter $\nu$ such that $\wh{P} = \wh{P}_{\nu}$. This proves that to show basicness or cuspidality, it suffices to check the weight condition for a single $\nu$ for each $\wh{P}$.

By an analogue of Remark \ref{rem: BasCohaltdef} for cotangent stacks, the sheaf $\mc{F}$ is right-basic if and only if for all strictly $\wh{P}$-dominant $\mc{F'} \in \Coh(\T^{\ast}\V^{\nu = 0})$, we have 
 \begin{equation}{\label{eqn: eishomvanishing}}
   \Hom_{\QCoh}(q^{\ast} (\mc{F}' \otimes \omega_{\V_U}), p^{!,\QCoh}\mc{F})\simeq \Hom_{\Coh}(\Eis^{\T^*}_{\wh P} \mc{F}' , \mc{F}) \simeq 0.
 \end{equation}
By Lemma \ref{lem: cotangentdomswitch}, right-basic is also equivalent to the same vanishing but for $\mc{F'}$ with central character $\chi$ satisfying $\langle \nu, \chi \rangle > 0$. 

To finish the proof, we need to show the vanishing in \eqref{eqn: eishomvanishing} is equivalent to the statement that the $\nu$-weight of $p^{!,\QCoh}\mc{F}$ is strictly less than the $\nu$-weight of $q^{\ast} (\mc{F}' \otimes \omega_{\V_U})$. If $p^{!,\QCoh}\mc{F} \in \QCoh(\T^{\ast}\V^{\nu \geq 0})^{ \leq \langle \nu, \wt(\V_U)\rangle}$, then the desired vanishing follows from semiorthogonality and functoriality of the baric structure along $q^{\ast}$.

Conversely, suppose that the vanishing in \eqref{eqn: eishomvanishing} holds for each such $\mc{F}'$. By Lemma \ref{lem: HLlemmas} (3), it implies that $p^{!,\QCoh}\mc{F}$ is right-orthogonal to $\QCoh(\T^*\V^{\nu \geq 0})^{> \langle \nu, \wt(\V_U)\rangle}$, and hence belongs to $\QCoh(\T^*\V^{\nu\geq 0})^{\leq \langle \nu, \wt(\V_U)\rangle}$. 
\end{proof}

Now the first part of Theorem \ref{thm: ClassificationCuspidalCotangent} (1) for perfect complexes in the graded setting follows from Lemma \ref{lem: Cotangentleftrightcuspidality} and

\begin{corollary}{\label{cor: leftrightbasequal}}
    Let $b \in B(G)_{\bas}$ and $\mathcal{F}\in \Perf(\T^{\ast}\V)^{\kappa(b)}$. Then $\mathcal{F}$ is left-basic (resp.~left-cuspidal) if and only if it is right-basic (resp.~right-cuspidal).
\end{corollary}

\begin{proof}
It suffices to prove it for basic sheaves and the argument for cuspidal sheaves is completely analogous. 

By Proposition \ref{prop: basicweightcharacterization}, we have that $\mc{F}\in \Perf(\T^*\V)$ is right-basic if and only if for each cocharacter $\nu$ of $\wh{T}$ defining a proper parabolic subgroup, we have 
\begin{equation*}
    p^!\mc{F} \in \QCoh(\T^{\ast}\V^{\nu \geq 0})^{ \leq \langle \nu, \wt(\V_U) \rangle + \langle \nu,  \nu_G(b)\rangle}
\end{equation*}
and left-basic if and only if
\begin{equation*}
    p^{\ast}\mc{F} \in \QCoh(\T^{\ast}\V^{\nu \geq 0})^{ \geq \langle \nu, \wt(\V_U) \rangle + \langle \nu,  \nu_G(b)\rangle}.
\end{equation*}
We have the following diagram
\begin{equation*}
    \begin{tikzcd}
        & \T^{\ast}\V^{\nu \geq 0} \arrow[ld, swap, "q"] \arrow[rd, "p"] & \\
        \T^{\ast}\V^{\nu = 0} \arrow[ru, bend left = 30, "\sigma"] \arrow[rd, bend right = 30, "\sigma^{\op}"] & & \T^{\ast}\V \\
        &  \T^{\ast}\V^{-\nu \geq 0} \arrow[lu, swap, "q^{\op}"] \arrow[ru, "p^{\op}"]  &
    \end{tikzcd}
\end{equation*}
such that $p \sigma  \simeq p^{\op} \sigma^{\op}$ and $q \sigma \simeq \id \simeq q^{\op} \sigma^{\op}$.

By part (1) of Lemma \ref{lem: HLlemmas} and Lemma \ref{lem: TVgeqdualizingcomputation}, $\mc{F}$ is right-basic if and only if
\begin{equation*}
   \QCoh(\T^{\ast}\V^{\nu = 0})^{ \leq \langle \nu, \wt(\V_U) \rangle + \langle \nu, \nu_G(b)\rangle}  \ni \sigma^{\ast} p^! \mc{F} \simeq \sigma^{\ast} p^{\ast} \mc{F} \otimes \omega^2_{\V_{\wh{P}}}\omega_{\V}^{-1} \omega_{\V_{\wh{M}}}^{-1}[2(\dim \V_M - \dim \V_P)].
\end{equation*}
Equivalently,
\begin{equation*}
    \sigma^{\op \ast}p^{\op \ast} \mc{F} \simeq \sigma^{\ast}p^{\ast} \mc{F} \in \QCoh(\T^{\ast}\V^{\nu = 0})^{ \leq \langle \nu, \wt(\V_{U^{\op}}) \rangle + \langle \nu, \nu_G(b)\rangle},
\end{equation*}
and hence
\begin{equation*}
    \sigma^{\op \ast}p^{\op \ast} \mc{F} \in \QCoh(\T^{\ast}\V^{-\nu = 0})^{ \geq \langle -\nu, \wt(\V_{U^{\op}}) \rangle + \langle -\nu, \nu_G(b)\rangle},
\end{equation*}
using the description in the equivalence \eqref{eqn: bariccentralweightcompatibility}.

But this is exactly (again by part (1) of Lemma \ref{lem: HLlemmas}) the condition for $\mc{F}$ to be left-basic (for $P^{\op}$ determined by $- \nu$).
\end{proof}
On account of Corollary \ref{cor: leftrightbasequal}, we may now simply refer to basic and cuspidal perfect complexes on $\T^{\ast} \V$.

\subsection{Window categories}{\label{ss: windowcat}}

This is a review of the main result of \cite{HalpernLeistner}, restricted to our setting above. For $c \in C\in \mathcal{C}$, we write $U_C=V\setminus \bigcup_{C'> C} C'$, $\U_C=U_C/H$, $\T^{\ast}\U_C$ with the $\Theta$-stratum $i\colon S \hookrightarrow \T^{\ast}\U_C$ with center $\pi\colon S\to Z\hookrightarrow \overline{Z}=\T^{\ast}\V_{M_c}$ arising from the construction in \S \ref{sss:ThetaStratificationCotangent}. For integers $u,v$, we define
\begin{align*}
\mathcal{G}^{\geq u}&=\{F \in \Coh(\T^{\ast}\U_C) \mid i^* F\in \Aperf^*(S)^{\geq u}\subset \QCoh(S)\}, \\
\mathcal{G}^{< v}&=\{F \in \Coh(\T^{\ast}\U_C) \mid i^{!} F\in \QCoh(S)^{< v}\}, \\
\mathcal{G}^{[u,v)}&=\mathcal{G}^{\geq u} \cap \mathcal{G}^{ < v}, \quad
\mathcal{W}^{u}=\mathcal{G}^{[u,u)}=\mathcal{G}^{\geq u} \cap \mathcal{G}^{ < u}. 
\end{align*}

We write $\Eis^*_{P_c}=i_*\pi^*\colon \Coh(Z)\to \Coh(\T^{\ast}\U_C)$ for the $*$-Eisenstein series \eqref{eq: stareis}.

\begin{theorem}[{\cite[Theorems 2.2.2, 3.3.1, Example 3.3.2]{HalpernLeistner}}]{\label{thm: HL3.3.1}}

\begin{enumerate}
    \item The restriction functor $\mathcal{W}^w \to \Coh(\T^{\ast}\U_C \setminus S)$ is an equivalence for any integer $w$. 
    \item $\Eis^{\ast}_{P_c}$ is fully faithful on $\Coh(Z)^w$, and its essential image is the objects of $\mathcal{G}^{[w,w+1)}$ that are set-theoretically supported on $S$. 
    \item For integers $u \leq w \leq v$, the category $\mathcal{G}^{[u,v)}$ admits a semiorthogonal decomposition
    \begin{equation*}
    \mathcal{G}^{[u,v)}=\langle \Eis^*_{P_c} \Coh(Z)^u, \dots, \Eis^*_{P_c} \Coh(Z)^{w-1}, \mathcal{W}^w, \Eis^*_{P_c} \Coh(Z)^{w}, \dots, \Eis^*_{P_c} \Coh(Z)^{v-1} \rangle. 
    \end{equation*}
\end{enumerate}    
\end{theorem}
It is convenient for us to rewrite the conditions. Let $\sigma\colon Z\hookrightarrow S$ denote the natural inclusion. 
\begin{lemma}\label{lem: HLEquivalenceWeightConditions}
Write $a$ for the $\nu_c$-weight of $\sigma^* \omega_S$. Let $\mc F_1, \mc F_2 \in \Coh(\T^{\ast}\U_C)$. The following are equivalent:
\begin{enumerate}
    \item $\mc F_1, \D_{\GS} \mc F_2\in \mc{G}^{\geq w}$,
    \item $\D_{\GS} \mc F_1,\mc F_2\in \mc{G}^{< a+1 -w}$, 
    \item $\sigma^* i^* \mc F_1, \sigma^* i^* \D_{\GS} \mc F_2\in \Aperf^*(Z)^{\geq w}$,
    \item $\pi_* i^! \D_{\GS} \mc F_1, \pi_* i^! \mc F_2 \in \QCoh(Z)^{< a+1 - w}$.
\end{enumerate}
\end{lemma}

\begin{proof}
The equivalences among (1), (2), and (3) are \cite[Proposition 3.2.6, Lemma 3.2.8]{HalpernLeistner}. (2) $\Rightarrow$ (4) is clear. For the converse, assume $i^!\D_\GS\mc F_1\neq 0$ and take the minimal $w'$ such that $\D_{\GS} \mc F_1\in \mc G^{< w' +1}$, which exists by coherence. It suffices to show $a+1-w\geq w'+1$. In the triangle
\[
\beta^{\geq w'}i^! \D_{\GS}\mc F_1 \to i^! \D_{\GS} \mc F_1 \to \beta^{< w'} i^! \D_{\GS} \mc F_1,
\]
$0\not\simeq \beta^{\geq w'}i^! \D_{\GS}\mc F_1\in \QCoh(S)^{w'}$ and hence 
$0\not\simeq \beta^{\geq w'}\pi_* \beta^{\geq w'}i^! \D_{\GS}\mc F_1 \in \QCoh(Z)^{\leq w'}$ by Lemma \ref{lem: HLlemmas}. Since $\pi_* \beta^{< w'} i^! \D_{\GS} \mc F_1 \in \QCoh(Z)^{< w'}$, it follows that
\[
\pi_* i^! \D_{\GS} \mc F_1 \in \QCoh(Z)^{< w' +1}, \qquad 
\pi_* i^! \D_{\GS} \mc F_1 \not\in \QCoh(Z)^{< w'}. 
\]
Therefore, $a+1-w \geq w'+1$.
\end{proof}

It is more natural to work with the twisted Grothendieck--Serre duality $\D'_{\GS}\coloneq\D_\GS\otimes \omega_{\U_C}$.

\begin{lemma}
The sum of $a$ and $\langle \nu_c, \wt(\V) \rangle$ equals $2\langle \nu_c, \wt(\V_{P_c})\rangle$.  
\end{lemma}
\begin{proof}
Follows from the formula $a = \langle \nu_c, \wt(\omega^2_{\V_{P_c}}\omega^{-1}_{\V}\omega^{-1}_{\V_{M_c}})\rangle$.
\end{proof}

Set $w_C=\langle \nu_c, \wt(\V_{P_c})\rangle$, which is independent of $c\in C$. As in Proposition \ref{prop: cotaneiscompatibility}, we have
\begin{equation*}
    \Eis^{\T^{\ast}}_{P_c}\simeq\Eis^*_{P_c} (\omega_{\V_{U_c}}\otimes-),
\end{equation*}
where $\omega_{\V_{U_c}}$ denotes $\omega_{\V_{P_c}} \omega^{-1}_{\V_{M_c}}$.
We are mainly interested in
\[
\mathcal{G}^{[w_C,w_C+1)}=\langle \Eis^*_{P_c} \Coh(Z)^{w_C}, \mathcal{W}_C\rangle=\langle \Eis^{\T^{\ast}}_{P_c} \Coh(Z)^{0}, \mathcal{W}_C\rangle,
\]
where $\mathcal{W}_C=\mathcal{W}^{w_C+1}$. 

\begin{corollary}{\label{cor: Gstableundertwistedduality}}
The category $\mathcal{G}^{[w_C, w_C+1)}$ is stable under the twisted Grothendieck--Serre duality $\D'_\GS$. 
\end{corollary}

We write $j_C \colon \T^*  (( U_C\setminus C)/H) \hookrightarrow \T^* \U_C$ and $j_{C \ush} \colon \Coh(\T^* \U_C\setminus S) \xrightarrow[\simeq]{j_C^{*-1}} \mathcal{W}_C\hookrightarrow \mathcal{G}^{[w_C, w_C + 1)}$. It is the left adjoint of $j^*_C \colon \mathcal{G}^{[w_C, w_C + 1)}\to \Coh(\T^* \U_C\setminus S)$.

By composing $j_{C'\ush}$ for all $C'> C$, we also obtain  $j_{\ush} \colon \Coh(\T^{\ast}\U_C)\to \Coh (\T^{\ast}\V)$ for $j\colon \T^* \U_C\hookrightarrow \T^{\ast}\V$ which satisfies an adjunction relation 
\begin{equation*}
    \Hom(j_{\ush} \mc{G} , \mc{F}) \simeq \Hom(\mc{G}, j^{\ast} \mc{F}),
\end{equation*}
for all sheaves $\mc{F} \in \Coh(\T^{\ast} \V)$ such that $\mc{F}|_{\T^* \U_{C'}} \in \mc{G}^{[w_{C'}, w_{C'}+1)}_{C'}$ for each $C' > C$. Note that if $C'>C$ are consecutive entries in the total order, then $\T^{\ast}\U_{C'} \setminus S = \T^{\ast}\U_C$. 

Let $C^{\circ} \subset V$ be the open orbit, corresponding to the initial orbit in the total order. We write
\begin{equation*}
    j^{\circ}: \T^{\ast}\U_{C^{\circ}} \rightarrow \T^{\ast}\V.
\end{equation*}

\begin{remark}
The choice of total order is unimportant. Suppose $C_1 > C_2 > C_3$ are consecutive and swapping $C_1$ and $C_2$ still gives a total order refining the closure order; then $j_{\ush}\mc F$ makes sense for $\mc F\in \Coh(\T^{\ast}\U_{C_3})$ in both cases. Our supposition implies
\[
\overline{C_1}\cap C_2 =\varnothing, \qquad C_1 \cap \overline{C_2} =\varnothing,
\]
and hence $\T^*_{C_1}\V$, $\T^*_{C_2}\V$ are closed in the union $\T^*_{C_1}\V \cup \T^*_{C_2}\V \cup \T^{\ast}\U_{C_3}$ and have empty intersection. Then the restriction of $j_{C_1 \ush}j_{C_2 \ush}\mc F$ to the complement of $\T^*_{C_2}\V$ agrees with $j_{C_1 \ush}^{\text{second}} \mc F$, where we are using the order exchanging $C_1$ and $C_2$. Namely, $j_{C_1 \ush}j_{C_2 \ush}\mc F$ is simply gluing $j_{C_2\ush}\mc F$ and $j_{C_1 \ush}^{\text{second}} \mc F$, and this description is symmetric with respect to $C_1$ and $C_2$. In particular, for the open stratum $C^\circ$,
\[
j^{\circ}_{\ush} \colon \Coh(\T^{\ast}\U_{C^\circ})\to \Coh (\T^{\ast}\V)
\]
is independent of the choice of total order. 
\end{remark}

\subsection{Limit category}{\label{ss: limcat}}
We will show (Proposition \ref{prop: lim=bas}; see also Theorem \ref{thm: cotanbasgenbycuspeis}) that basic coherent sheaves agree with the limit category of P\u{a}durariu--Toda \cite{PTDolbeault} if $G$ is semisimple, and here we only show that the limit category is contained in $\LBas^\Coh(\T^*\V)$.

The limit category is related under Koszul equivalence to the graded matrix factorizations of the magic category of (the dual of) $\mc{V} \coloneq (\T^*V \times (\Lie H)^\vee)/H$.

\begin{definition}[{\cite[Corollary 3.35]{PTqBPSsymquiverswithpotential}, \cite[Proposition 6.5]{PTDolbeault}}]{\label{defn: magicat}}
    The \emph{magic category} $\M(\mc{V}^{\vee}) \subset \Coh(\mc{V}^{\vee})$ is the full subcategory generated by vector bundles $B_{\chi}$ coming from representations of $H$ of highest weight $\chi$, such that 
    \begin{equation*}
        \chi +\rho - \frac{1}{2}\delta \in \left( \bigoplus\limits_{\beta \in \wt (\T^{\ast} V \times \Lie H)} [0, \frac{1}{2}\beta]\right) \eqqcolon \mc{P}_{\magic},
    \end{equation*}
    where $\delta = \wt(\omega_{\V})$ and $\rho$ is the half-sum of positive roots of $H$. 
\end{definition}

We have the following characterization of all weights of vector bundles in $\M(\mc{V}^{\vee})$ as opposed to just the dominant ones.
\begin{lemma}[{\cite[Lemma 2.34]{PTqBPSsymquiverswithpotential}, \cite[Lemma 2.9]{HalpernLeistnerSam}}]{\label{lem: magicvbweightcharacterization}}
A vector bundle $B \in \Coh(\mc{V}^{\vee})$ is in $\M(\mc{V}^{\vee})$ if and only if the $\wh T$-weights $\chi$ of the representation of $H$ satisfy that for each cocharacter $\nu \in X_{\ast}(\wh T)$, we have
\begin{equation}\label{eq: MagicAllWeights}
    \langle \nu, \chi \rangle \in \frac{1}{2}\left[ -\langle\nu , \det \L^{\nu > 0}_{\mc{V}^{\vee}} \rangle,   \langle \nu, \det \L^{\nu > 0}_{\mc{V}^{\vee}} \rangle \right] + \langle \nu, \frac{1}{2}\delta \rangle,
\end{equation}
where $\delta =  \wt(\omega_{\V})$. 
\end{lemma}

\begin{example}\label{ex: MagicCategoryExample}
Suppose $V/H=BH$, in particular $\delta=0$. If $B_\chi\in \M(\mc V^\vee)$, then $B_\chi$ is trivial. 
\end{example}

Eisenstein series for $\mc{V}^{\vee}$ are defined analogously to the $\T^{\ast}\V$ case. Namely, given a cocharacter $\nu$ of $\wh{T}$, there is an associated connected component $\mc{V}^{\vee, \nu \geq 0}$ (resp.~$\mc{V}^{\vee, \nu = 0}$) of $\Map(\Theta, \mc{V}^{\vee})$ (resp.~$\Map(B\Gm, \mc{V}^{\vee})$) and a diagram of quotient stacks
    \begin{equation*}
    \begin{tikzcd}
        \mc{V}^{\vee, \nu = 0} \arrow[r, swap, bend right=15, "{\sigma}"] & \mc{V}^{\vee, \nu \geq 0} \arrow[r, "{\ev_1}"] \arrow[l,swap, bend right=15, "{\ev_0}"] & \mc{V}^{\vee}
    \end{tikzcd}.
\end{equation*}
We then define $\Eis^{\mc{V}}_{\nu} \coloneq \ev_{1 \ast}\ev^{\ast}_0$ (or we write $\Eis^{\mc V}_{P_\nu}$ instead). In \cite{PTDolbeault} and elsewhere in the literature, this functor is called a \emph{(categorical) Hall product/induction}.

\begin{lemma}[{cf.~\cite[Lemma 10.1]{PTDolbeault}}]{\label{lem: PT10.1lem}}
Let $T' \subset \wh T \subset H$ be the largest sub-torus of $\wh T$ that acts trivially on $\T^{\ast}V \times \Lie H$ and let $w$ be the restriction of $\frac{1}{2}\delta$ to $T'$. Then 
\begin{enumerate}
\item  $\Coh(\mc{V}^{\vee})^w$ is generated under finite colimits by $\M(\mc{V}^{\vee})$ and $\Eis^{\mc{V}}_{\nu} \mc{F}$ where $\mc{F} \in \Coh(\mc{V}^{\vee, \nu=0})^w$ is in the subcategory generated by vector bundles of highest weight $\chi''$ such that
\begin{equation*}
    \langle \nu, \chi'' \rangle < -\frac{1}{2} \langle \nu, \det \L^{\nu > 0}_{\mc{V}^{\vee}} \rangle + \langle \nu, \frac{1}{2}\delta \rangle,
\end{equation*}
for each cocharacter $\nu$ giving a face of $\mc{P}_{\magic}$,
\item  we have 
    \begin{equation*}
        \Hom(\M(\mc{V}^{\vee}), \Eis^{\mc{V}}_{\nu} \mc{F}) \simeq 0,
    \end{equation*}
    for each $\mc{F}$ as in (1) and general $\nu$.
\end{enumerate}
Analogous statements hold for any other $T'$-weight. 
\end{lemma}
\begin{proof}
The proof of (1) is \cite[Lemma 3.16]{PTqBPSsymquiverswithpotential}, keeping in mind that we take $\ell=0$ in that lemma so that their $\res$ is the identity functor. 

The proof of (2) is as in the proof of \cite[Proposition 3.22]{PTqBPSsymquiverswithpotential}, using the weight characterization of Lemma \ref{lem: magicvbweightcharacterization}. 
\end{proof}
\begin{corollary}[{\cite[Definition 3.12]{PTDolbeault}}]
The category $\M(\mc{V}^{\vee}) \subset \Coh(\mc{V}^{\vee})$ can be alternatively defined as the full subcategory whose objects $\mc{E}$ satisfy that for each map $\nu: B\Gm \rightarrow \mc{V}^{\vee}$, we have
\begin{equation*}
    \wt(\nu^{\ast} \mc{E}) \subset \frac{1}{2}\left[ - c_1(\nu^{\ast}\L^{>0}_{\mc{V}^{\vee}}) , c_1(\nu^{\ast}\L^{>0}_{\mc{V}^{\vee}})\right] + \langle \nu, \frac{1}{2}\delta \rangle.
\end{equation*}
Our notation follows \cite{PTDolbeault} so $c_1: K(B\Gm)_{\Q} \rightarrow \Q$ takes $A=\oplus_m A_m \in \Coh(B\Gm)$ to $\sum_m m \cdot \dim A_m$.
\end{corollary}

\begin{proof}
    This is essentially \cite[Corollary 3.35]{PTqBPSsymquiverswithpotential}. Rather, the definition of $\M(\mc{V}^{\vee})$ in \cite{PTDolbeault} is the one given here, and \cite[Proposition 6.5]{PTDolbeault} proves Definition \ref{defn: magicat} is equivalent to this one. In \cite{PTqBPSsymquiverswithpotential}, both definitions are given and their equivalence is \cite[Corollary 3.35]{PTqBPSsymquiverswithpotential}.
\end{proof}

\begin{corollary}{\label{cor: magicSODs}}
We have a semiorthogonal decomposition of $\Coh(\mc{V}^{\vee})$ given by
\begin{equation*}
\langle \Eis_{-}^{\mc{V}, \circ} , \M(\mc{V}^{\vee}) \rangle,
\end{equation*}
   where $\Eis^{\mc{V}, \circ}_-$ is the category generated under finite colimits by the $\Eis^{\mc{V}}_{\nu} \mc{F}$ varying over all cocharacters $\nu$ with $\mc{F}$ having central weights satisfying the inequality of part (1) of Lemma \ref{lem: PT10.1lem} relative to $\nu$ (but we are not assuming $\nu$ satisfies the face condition).
\end{corollary}
\begin{proof}
    We have to check the categories generate $\Coh(\mc{V}^{\vee})$ and the requisite $\Hom$-vanishing. The generation is part (1) of Lemma \ref{lem: PT10.1lem} and the vanishing is part (2).
\end{proof}

We now introduce the category of graded matrix factorizations. We give the vector bundle $\mc{V}^{\vee} \rightarrow (\T^{\ast}V)/ H$ the weight $2$ $\Gm$-action.  The moment map
\begin{equation*}
    \mu: \T^{\ast}V / H \rightarrow (\Lie H)^\vee /H,
\end{equation*}
induces a function $f: \mc{V}^{\vee} \rightarrow \A^1$ given by
\begin{equation*}
    f: (y, v) \mapsto \langle \mu(y), v \rangle.
\end{equation*}
\begin{definition}[{\cite[\S2.7]{PTDolbeault}, \cite[\S2.6.2]{PTC3I}}]
    The category of \emph{graded matrix factorizations}, $\MF^{\gr}(\mc{V}^{\vee}, f)$ has objects $(P, d_P)$ where $P \in \Coh_{\Gm}(\mc{V}^{\vee})$ and $d_P: P \rightarrow P(1)$ is a $\Gm$-equivariant morphism of weight $1$ such that $d_P \circ d_P = f$. 

    For $\bb{M} \subset \Coh(\mc{V}^{\vee})$ a full subcategory, we denote $\bb{M}_{\Gm} \subset \Coh_{\Gm}(\mc{V}^{\vee})$ the full subcategory with objects in $\bb{M}$ after forgetting the $\Gm$-action and $\MF^{\gr}(\bb{M}, f)$ the matrix factorizations with factors in $\bb{M}_{\Gm}$.
\end{definition}
\begin{corollary}{\label{cor: MFSODs}}
    We have the following semiorthogonal decomposition of $\MF^{\gr}(\mc{V}^{\vee}, f)$.
    \begin{equation*}
  \langle \MF^{\gr}(\Eis_{-}^{\mc{V}, \circ}, f) , \MF^{\gr}(\M(\mc{V}^{\vee}), f) \rangle.
    \end{equation*}
\end{corollary}
\begin{proof}
By \cite[Proposition 2.5(b)]{PTC3I}, the semiorthogonal decomposition of $\Coh(\mc{V}^{\vee})$ in Corollary \ref{cor: magicSODs} lifts to a semiorthogonal decomposition of $\MF^{\gr}(\mc{V}^{\vee}, f)$.
\end{proof}

\begin{theorem}[{\cite[Theorem 2.3.3]{TodaBook}. See also \cite{Shipman}, \cite{Isik}, \cite{Hirano}, \cite{oblomkovrozansky}}]{\label{thm: MFKoszulequiv}}
There is a Koszul equivalence 
    \begin{equation*}
        \Psi: \Coh(\T^{\ast}\V) \xrightarrow{\sim} \MF^{\gr}(\mc{V}^{\vee}, f),
    \end{equation*}
    induced by tensoring with the Koszul factorization
    \begin{equation*}
        (\mc{O}_{\mc{V}^{\vee}} \otimes_{\mc{O}_{\T^{\ast}V/H}} \mc{O}_{\T^{\ast}\V}, d_K).
    \end{equation*}
\end{theorem}
\begin{lemma}[{\cite[Proposition 2.27]{PTqBPSsymquiverswithpotential}}]{\label{lem: MFcotaneiscompatibility}}
The Koszul equivalence of Theorem \ref{thm: MFKoszulequiv} satisfies the compatibility
\begin{equation*}
 \Psi \circ \Eis^{\ast}_{\nu} \simeq \Eis^{\mc{V}}_{\nu}\circ ( \det(\mc{V}^{\nu > 0})^{\vee}[\rank \mc{V}^{\nu >0}]\otimes   \Psi).
\end{equation*}
In particular, $\Eis^{\mc{V}}_{\nu}$ and $\Eis^{\T^{\ast}}_{P_\nu}$ differ under $\Psi$ by a twist of
\begin{equation*}
    \det(\mc{V}^{\nu > 0})^{\vee}[\rank \mc{V}^{\nu >0}] \otimes \omega_{\V_U} .
\end{equation*}
We call this twist $\omega$. 
\end{lemma}

Let $i: \T^{\ast}\V \hookrightarrow \T^{\ast}V/H$ denote the inclusion.
\begin{definition}[{\cite[Lemma 3.7]{PTDolbeault}}]{\label{def: limcat}}
The \emph{limit category} $\mrL(\T^{\ast}\V) \subset \Coh(\T^{\ast}\V)$ \footnote{in \cite{PTDolbeault} this is denoted $\mrL(\T^{\ast}\V)_{\omega^{1/2}_{\V}}$ or $\mrL(\T^{\ast}\V)_{\frac{1}{2}}$ but we drop the subscript as this is the only limit category we consider.} is defined to be the full subcategory of objects $\mc{E}$ such that for each
\begin{equation*}
    \nu: B\Gm \rightarrow \T^{\ast}\V,
\end{equation*}
we have
\begin{equation*}
    \wt(\nu^{\ast}i^{\ast}i_{\ast}\mc{E}) \subset \frac{1}{2}[c_1(\nu^{\ast}\L^{<0}_{\mc{V}}), c_1(\nu^{\ast}\L^{>0}_{\mc{V}})] + \frac{1}{2}c_1(\nu^{\ast} \omega_{\V}).
\end{equation*}
In \cite[Definition 3.3]{PTDolbeault}, limit categories are defined for more general stacks and satisfy the functorialities of the following proposition. In particular, they are defined for the stacks $\T^{\ast}\V_P$. In this paper, we primarily use the definition above as it allows us to avoid discussing regularizations of $\nu$. 
\end{definition}
\begin{proposition}[{\cite[\S5.2, \S5.4]{PTDolbeault}}]{\label{prop: Limcatfunctorialities}}
Let $f: X \rightarrow Y$ be a map of smooth QCA stacks. Then we consider the following functors and their analogues for $\Coh$ (see \S \ref{ss: eisensteinseriesforcotangents}, \cite[Definition 2.1.23]{Chen2023categorical})
\begin{equation*}
    (\breve{d}f)_{\ast}: \IndCoh(\T^{\ast}X) \rightarrow \IndCoh(\T^{\ast}Y), \quad   \left(\text{resp.~}(\breve{d}f)^!: \IndCoh(\T^{\ast}Y) \rightarrow \IndCoh(\T^{\ast}X)\right),
\end{equation*}
when $f$ is proper schematic (resp.~smooth).

Then if $f$ is projective (resp.~smooth), we have
\begin{equation*}
(\breve{d}f)_{\ast}\mrL(\T^{\ast}X) \subset \mrL(\T^{\ast}Y), \qquad \left(\text{resp.~} (\breve{d}f)^{!}\mrL(\T^{\ast}Y) \subset \mrL(\T^{\ast}X)\right).
\end{equation*}

In particular, $\Eis^{\T^{\ast}}_P \mrL(\T^{\ast}\V_M) \subset \mrL(\T^{\ast}\V)$. Note that making sense of $\mrL(\T^{\ast}\V_P)$ requires the second part of Definition \ref{def: limcat}.    
\end{proposition}
\begin{proposition}[{\cite[Proposition 3.14]{PTDolbeault}}]{\label{prop: limkoszulequiv}}
    The Koszul equivalence of Theorem \ref{thm: MFKoszulequiv} induces an equivalence
\begin{equation*}
    \MF^{\gr}(\M(\mc{V}^{\vee}), f) \simeq \mrL(\T^{\ast}\V).
\end{equation*}
\end{proposition}
\begin{proof}
    The argument in \cite[Proposition 3.14]{PTDolbeault} uses \cite[Lemma 2.23]{PTqBPSsymquiverswithpotential}. This describes, for a subcategory $S \subset \MF^{\gr}(\mc{V}^{\vee}, f)$ whose underlying category is generated by direct sums of vector bundles (as is the case for $\M(\mc{V}^{\vee})$), the subcategory of $\Coh(\T^{\ast}\V)$ that corresponds to $S$ under the above Koszul equivalence.
\end{proof}

\begin{lemma}{\label{lem: limitselfdual}}
The duality $\D'_{\GS}$ induces an equivalence
\begin{equation*}
    \D'_{\GS}: \mrL(\T^{\ast}\V) \simeq \mrL(\T^{\ast}\V)^{\op}.
\end{equation*}
\end{lemma}
\begin{proof}
    This follows from \cite[Lemma 3.8, Lemma 3.9]{PTDolbeault} and their proofs.
\end{proof}

\begin{proposition}{\label{prop: magickoszuleiscompatibility}}
Here, the superscript $0$ denotes $Z(\wh G)^\circ$-weight. We have an equivalence
\begin{equation*}
    \Psi: \Eis_-^0 \simeq \Eis^{\mc{V}}_-(\MF^{\gr})^0,
\end{equation*}
and analogously for the $\circ$- and $+$- versions of $\Eis^0\coloneq\Eis^{\T^*}\cap \Coh(\T^*\V)^0$, where $\Eis_{-}^{0}$ (resp.~$\Eis^{0 \circ}_-$) is the category generated under finite colimits by Eisenstein series $\Eis^{\T^{\ast}}_{P_\nu} \mc{F}$ such that for every central weight $\chi$ of $\mc{F}\in \Coh(\T^*\V^{\nu=0})^0$, then $\langle \nu, \chi \rangle \leq 0$ (resp.~$\langle \nu, \chi \rangle < 0$) and $\Eis^{\mc{V}}_-(\MF^{\gr})$ is the category generated by $\Eis^{\mc{V}}_\nu \mc{F}$ where $\mc{F} \in \MF^{\gr}(\mc{V}^{\vee,\nu=0}, f)^{ \leq -\frac{1}{2} c_1(\nu^* \L^{\nu > 0}_{\mc{V}^{\vee}} ) + \langle \nu, \frac{1}{2}\delta \rangle}$.
\end{proposition}
\begin{proof}  
Suppose that $\Eis^{\T^{\ast}}_{P_\nu} \mc{F} \in \Eis^0_-$ such that $\mc{F}$ has central weight $\chi$ and $\langle \nu, \chi \rangle  \leq 0$. It is easy to see that $\Psi(\mc{F})$ has central weight $\chi$ as well, so $\Psi(\mc{F}) \otimes \omega$ has central weight $\chi + \wt(\omega)$. 

We explicitly compute the twist $\omega$ of Lemma \ref{lem: MFcotaneiscompatibility}, where the sums are over weights in the indicated spaces.  We have
\begin{equation*}
\langle \nu, \wt(\omega_{\V_U}) \rangle = \langle \nu, -\sum\limits_{V^{\nu > 0}} + \sum\limits_{(\Lie H)^{\nu > 0}} \rangle.
\end{equation*}
Then
\begin{equation*}
    \langle \nu, \det(\mc{V}^{\nu > 0})^{\vee} \rangle = \langle \nu,  - \sum\limits_{(\Lie H)^{\nu \geq 0}}  \rangle,
\end{equation*}
so
\begin{equation}{\label{eqn: omegaformula}}
    \langle \nu, \wt(\omega) \rangle =\langle \nu,   - \sum\limits_{V^{\nu > 0}} \rangle.
\end{equation}
We also compute
\begin{align*}
    -\frac{1}{2}c_1(\nu^{\ast} \L^{\nu>0}_{\mc{V}^{\vee}}) + \langle \nu, \frac{1}{2}\delta \rangle &= -\frac{1}{2}\langle \nu , \sum\limits_{V^{\nu > 0}} - \sum\limits_{V^{\nu<0}} \rangle - \langle \nu, \frac{1}{2}\sum\limits_V \rangle\\
    &= \langle \nu, - \sum\limits_{V^{\nu>0}} \rangle.
\end{align*}
Thus,
\begin{equation*}
    \langle \nu, \chi \rangle \leq 0 \text{ iff }  \langle \nu, \chi + \wt(\omega) \rangle \leq \langle \nu, \wt(\omega) \rangle = -\frac{1}{2}c_1(\nu^{\ast} \L^{\nu>0}_{\mc{V}^{\vee}}) + \langle \nu, \frac{1}{2}\delta \rangle,
\end{equation*}
as desired.
\end{proof}

\begin{corollary}{\label{cor: LimsubsetBas}}
\[
\mrL(\T^*\V) \subset \LBas^{\Coh}(\T^{\ast}\V)^0.
\]
\end{corollary}
\begin{proof}
 Note that $\mrL(\T^*\V)\subset\Coh(\T^*\V)^0$. Passing to the Koszul dual side and considering left orthogonal complements,  by Proposition \ref{prop: magickoszuleiscompatibility}, it suffices to show that $\Eis_{-}^{\mc{V}, \circ}(\MF^{\gr}) \subset \MF^{\gr}(\Eis_{-}^{\mc{V}, \circ}, f)$ and this is clear by definition.
\end{proof}

It will also be useful to control $\mrL  (\T^*\U_C)$ or $\mrL (Z)$. The notation below is different from the main text. 

\begin{proposition}\label{prop: magicategorySOD}
Let $G$ be a (possibly disconnected) reductive group with a maximal torus $T$. Set $\mc V^\vee=\Lie G /G$ over $BG$. For a semistandard Levi subgroup $M$ associated with a cocharacter, write $\M'(\Lie M/M) \subset \Coh(\Lie M/M)$ for the category generated under finite colimits by vector bundles associated with representations of $M/([M, M]^\circ)$. 

We have a semiorthogonal decomposition
\[
\Coh(\mc V^\vee)=\langle \Eis_{P(\chi)}^{\mc V}\M'(\Lie M(\chi)/M(\chi))^{-\chi|_{Z(M(\chi))^\circ}}, \chi \in X^* (T)_\Q / W_G\rangle,
\]
where $M(\chi)$ is the stabilizer of $\chi\in X^*(T)_\Q \subset (\Lie G)^\vee$ and $P(\chi)\supset M(\chi)$ is associated with $\chi$ so that $\chi$ is strictly $P(\chi)$-dominant and $W_G\coloneq N_G (T)/Z_G(T)$. Here, the piece is zero unless $\chi|_{Z(M(\chi))^\circ}$ is integral. Moreover, $\Eis_{P(\chi)}^{\mc V}$ is fully faithful on $\M'(\Lie M(\chi)/M(\chi))^{-\chi|_{Z(M(\chi))^\circ}}$ and its essential image is independent of a representative of the $W_G$-orbit. 
\end{proposition}

\begin{proof}
If $G$ is connected, this is a special case of \cite[\S 2.2.2]{BuPadurariuToda} with $m=1$, together with Example \ref{ex: MagicCategoryExample}. The disconnected case is a mild generalization (the argument itself, or \cite[Theorem 3.2.3]{BuPadurariuToda}); fix any $W_G$-invariant positive-definite form on $X_*(T)_\Q$ that restricts to the Killing form on the $(G^\circ)_\der$-part. 
\end{proof}

Let $\mrL'(\T^* BM)\subset \Coh(\T^* BM)$ denote the subcategory corresponding to $\MF^{\gr}(\M'(\Lie M/M), 0)$ under the Koszul equivalence. Its summand where $Z(M)^\circ$-weights are trivial is the limit category. 

\begin{lemma}
$\mrL'(\T^* BM)\subset \Coh(\T^* BM)$ is generated under finite colimits by vector bundles on the zero section associated with representations of $M/([M, M]^\circ)$. 
\end{lemma}

\begin{proof}
By \cite[Lemma 2.23]{PTqBPSsymquiverswithpotential}, $\mc F\in \Coh(\T^*BM)$ lies in $\mrL' (\T^* BM)$ if and only if its $*$-pushforward to $BM$ is contained in the category generated by representations of $M/([M,M]^\circ)$. Since the $*$-pushforward is $t$-exact, it is equivalent to the condition that all the truncations of $\mc F$ are inflated from $M/([M,M]^\circ)$, which is equivalent to $\mc F$ being contained in the category generated by vector bundles on the zero section associated with representations of $M/([M,M]^\circ)$. 
\end{proof}

\begin{definition}
We define $\LBas^{\Coh}(\T^{\ast}BG)$ to consist of sheaves $\mc{F}$ such that for each proper semistandard parabolic $P \subset G$ with Levi $M$, we have $\Hom(\mc F, \Eis^{\T^*}_{P}(\mc{F}') ) \simeq 0$ for $\mc{F}'$ of strictly $P$-antidominant weight $\chi \in X^{\ast}(Z(M)^\circ)$. 
\end{definition}

\begin{corollary}{\label{cor: limcatclassifyingstack}}
We have the semiorthogonal decomposition
\[
\Coh(\T^* BG)=\langle \Eis_{P(\chi)}^{\T^*}\mrL' (\T^* BM(\chi))^{-\chi|_{Z(M(\chi))^\circ}} , \chi\in X^*(T)_\Q /W_G\rangle,
\]
and $\Eis_{P(\chi)}^{\T^*}$ is fully faithful on $\mrL' (\T^* BM(\chi))^{-\chi|_{Z(M(\chi))^\circ}}$. In particular, 
\[
\LBas^\Coh(\T^* BG)=\mrL' (\T^* BG), 
\]
and it is generated under finite colimits by vector bundles on the zero section $BG\hookrightarrow \T^* BG$ associated with representations of $G/([G,G]^\circ)$. Similarly, $\mrL (\T^* BG)$ is generated under finite colimits by vector bundles on the zero section $BG\hookrightarrow \T^* BG$ associated with representations of $\pi_0(G)$. 
\end{corollary}

\begin{proof}
With the previous proposition, this is a special case of \cite[Theorem 3.5.6]{BuPadurariuToda}, together with Koszul equivalence and compatibility of Eisenstein series. It is obtained by applying $\MF^{\gr}$ to the previous proposition. The inclusion $\LBas^\Coh(\T^* BG)\subset \mrL' (\T^* BG)$ follows from central weight decomposition and the choice of order in the proof of Proposition \ref{prop: magicategorySOD}, and the other inclusion is analogous to Corollary \ref{cor: LimsubsetBas}.    
\end{proof}

It follows that $\LBas^\Coh(\T^* BG)=\mrL' (\T^* BG)$ is stable under $\D'_{\GS}$ and we call them basic sheaves and write $\Bas^\Coh(\T^* BG)$. 

\subsection{Semiorthogonal decomposition by basic coherent sheaves and Eisenstein series}{\label{ss: cotangenbybasic}}
We will prove part (1) of Theorem \ref{thm: BasicCotangentMain} in the graded setting, using matrix factorizations and limit categories.\footnote{The proof possibly works in the general case, but we treat the general case differently later.} 

\begin{proposition}{\label{prop: cotanbascohSOD}}
    We have 
    \[
    \LBas^\Coh(\T^*\V)=\RBas^\Coh(\T^*\V)
    \]
    and we will simply write $\Bas^\Coh(\T^*\V)$. It is stable under $\D'_\GS$, and we have the following semiorthogonal decompositions of $\Coh(\T^{\ast}\V)$:
    \begin{equation*}
        \langle  \Eis^{\circ}_-,  \Bas^{\Coh}(\T^{\ast}\V) \rangle,  \qquad \langle \Bas^{\Coh}(\T^{\ast}\V),  \Eis^{\circ}_+\rangle.
    \end{equation*}
Note that by $\D'_{\GS}$, each semiorthogonal decomposition implies the other. 
\end{proposition}
\begin{proof}
The statements easily reduce to the case of semisimple groups by Lemma \ref{lem: reductiontosimplyconnectedcover} and its argument to relate Eisenstein series. We shall explain how \cite{BuPadurariuToda} directly implies the above semiorthogonal decompositions in the semisimple case. The Killing form of $\wh G$ induces a $W_G=W_{\wh G}$-invariant positive-definite quadratic form $q$ on $X_*(\wh T)_\Q=X^*(T)_\Q$:
\[
q(\nu)\coloneq \sum_{\alpha\in \Phi (\wh G, \wh T)} \langle \alpha, \nu \rangle^2=\sum_{\alpha \in \Phi (G, T)} \langle \alpha^\vee, \nu \rangle^2. 
\]
Note that if $\nu$, $\nu'$ are both $B$-dominant and $\nu'\geq \nu$ with respect to the $B$-dominance order, then $q(\nu')\geq q(\nu)$. The determinant $\det(V)$ is a character of $H$ and it can be regarded as an element of $X^*(\wh  T)^{W_H}$. We write (only in this proof) $\delta\in X^*(\wh  T)^{W_H}_\Q$ for $-\det(V)/2$. 

We first apply \cite[Theorem 2.1.6]{BuPadurariuToda} to 
\[
\Coh((\T^*V \times \Lie H)/H)=\Perf((\T^*V \times \Lie H)/H)
\]
with these $q$ and $\delta$, and we obtain a semiorthogonal decomposition
\[
\Perf((\T^*V \times \Lie H)/H)=\langle \star_\nu \mathrm{W}_{(\T^*V \times \Lie H)^{\nu=0}/H^{\nu=0}}(\delta_\nu)\mid \nu \in X_*(\wh T)_\Q /W_H \rangle,
\]
where $\star_\nu$ denotes the Hall induction and $\mathrm{W}_{V_\nu/H_\nu}(\delta_\nu)$ is a window subcategory, and $\star_\nu$ is fully faithful on these window subcategories. The Hall induction is $\Eis^{\mc V}$ in our terminology and their window subcategory $\mathrm{W}_{(\T^*V \times \Lie H)/H}(\delta)$ is exactly the magic category $\M((\T^*V \times \Lie H)/H)$. While they formulated their result using cocharacters $\nu$, the relevant piece can be defined only using the character corresponding to $-\nu$ via the identification induced by $q$. Namely, if $\chi\in X^*(\wh T)_\Q$, then there is a unique pair of semistandard parabolic and Levi subgroups $P\supset M$ such that $\chi$ is strictly $P$-antidominant, and we only need these data. We can think of $\mathrm{W}_{(\T^*V \times \Lie H)^{\nu=0}/H^{\nu=0}}(\delta_\nu)$ as a twisted magic category in the sense that the relevant polytope is shifted by a rational character of $Z(\wh M)^\circ$. (In particular, it is nontrivial only if $\chi|_{Z(\wh M)^\circ}$ is integral.) Since any character of $\wh T$ is $W_{\wh G}$-conjugate to a unique $B$-dominant character, we can group the pieces with the same $B$-dominant conjugate and order them by $q$-norm refining the $B$-dominant partial order to obtain a semiorthogonal decomposition. In fact, if two characters are $W_{\wh G}$-conjugate and not $W_H$-conjugate, the pieces are orthogonal to each other; see \cite[\S 2.3.2]{BuPadurariuToda}.

Next, we apply \cite[Theorem 3.5.6]{BuPadurariuToda} (cf.~\cite[\S 4.4.1(i)]{BuPadurariuToda}). It simply means we can obtain a semiorthogonal decomposition by taking the graded matrix factorization $\MF^{\gr}(-,f)$; since $\Eis^{\mc V}$ is fully faithful above, its matrix factorization is the Eisenstein series of a matrix factorization. 
\begin{multline*}
\MF^{\gr}((\T^*V \times \Lie H)/H, f) \\
=\langle \star_\nu \MF^{\gr}(\mathrm{W}_{(\T^*V \times \Lie H)^{\nu=0}/H^{\nu=0}}(\delta_\nu), f)\mid \nu \in X_*(\wh T)_\Q / W_H\rangle    
\end{multline*}
Using the Koszul equivalence and compatibility between $\Eis^{\mc V}$ and $\Eis^{\T^*}$ (Theorem \ref{thm: MFKoszulequiv}, Lemma \ref{lem: MFcotaneiscompatibility}, Proposition \ref{prop: magickoszuleiscompatibility}), we can rewrite the above semiorthogonal decomposition as
\[
\Coh(\T^*\V)=\langle \Eis^{\T^*}_P \mrL'(\T^*\V_M)^{\chi} \rangle, 
\]
where $\mrL'(\T^*\V_M)^{\chi}$ is defined as a full subcategory satisfying
\[
\mrL'(\T^*\V_M)^{\chi} \xrightarrow[\simeq]{\omega\otimes\Psi}\MF^{\gr}(\mathrm{W}_{(\T^*V \times \Lie H)^{\nu=0}/H^{\nu=0}}(\delta_\nu), f),
\]
where $\omega$ is the twist in Lemma \ref{lem: MFcotaneiscompatibility}, and $\chi\in X^*(\wh T)_\Q$ is strictly $P$-antidominant as above and the pieces are taken up to $W_H$-conjugacy of $\chi$; it follows from \cite[(2.1.6.2)]{BuPadurariuToda} that $\mrL'(\T^*\V_M)^{\chi}$ has $\nu$-weight $-q(\nu)$. 
Recalling the inclusion $\mrL(\T^*\V)\subset \LBas^\Coh(\T^*\V)$, we deduce that $\mrL(\T^*\V)=\LBas^\Coh(\T^*\V)$. Since $\mrL(\T^*\V)$ is stable under $\D'_\GS$ (Lemma \ref{lem: limitselfdual}) and $\D'_\GS$ exchanges left-basic and right-basic coherent sheaves (Lemma \ref{lem: Cotangentleftrightcuspidality}), it follows that left- and right- basic coherent sheaves agree. Therefore, we have proven the desired semiorthogonal decomposition.  
\end{proof}

We record the following result proven in Proposition \ref{prop: cotanbascohSOD}.
\begin{proposition}{\label{prop: lim=bas}}
    For $G$ semisimple, we have an equality
    \begin{equation*}
        \mrL(\T^{\ast}\V) = \Bas^{\Coh}(\T^{\ast}\V).
    \end{equation*}
\end{proposition}

\begin{corollary}{\label{cor: LBas=RBasTV}}
We have $\LBas(\T^{\ast}\V) = \RBas(\T^{\ast}\V)$ and we denote this category by $\Bas(\T^{\ast}\V)$. Then $\Bas(\T^{\ast}\V)$ is compactly generated by $\Bas^{\Coh}(\T^{\ast}\V)$.
\end{corollary}

\begin{proof}
These follow from ind-completing the semiorthogonal decompositions of Proposition \ref{prop: cotanbascohSOD}.
\end{proof}

We may now simply refer to basic sheaves on $\T^{\ast} \V$.

\begin{remark}
The above proof via graded matrix factorization also shows a finer  semiorthogonality relation, which we will also prove using the spectral geometric lemma. 
\end{remark}

\subsection{Basic sheaves and window categories}{\label{ss: basicvswindow}}
Consider the diagram \eqref{eqn: naturalitydiagram} with $\U = \U_C$ for some $C \in \mc{C}$.

\begin{lemma}{\label{lem: restriction lemma}}
Let $\mc{F} \in \Coh(\T^{\ast}\V)$ be basic. Then $j^{\ast} \mc{F} \in \mc{G}^{[w_C, w_C+1)}$. Suppose further that $\mc{F}$ is left-cuspidal, and that $C \neq C^{\circ}$ or $\nu_{c^{\circ}} \neq 0$. Then $j^{\ast} \mc{F} \in \mc{W}_C$.
\end{lemma}
\begin{proof}
 First, we consider basic sheaves. By duality (and using Proposition \ref{prop: cotanbascohSOD}), we just need to show one weight inequality. For the first claim, by Lemma \ref{lem: HLEquivalenceWeightConditions}, it suffices to show that $\pi_*i^{!}j^{\ast}\mc{F} \in \QCoh(\T^{\ast}\U_C^{\nu_c = 0})^{< w_C +1}$, where $\pi$ and $i$ are as in \S\ref{ss: windowcat}. By Lemma \ref{lem: eisfiberprod} and base change,
 \[
 \pi_*i^{!}j^{\ast}\mc{F} \simeq j^{=0*}q_* p^!\mc F.
 \]
 By right-basicness, the $\nu_c$-weights of $q_* p^!\mc F$ are $\leq w_C$ as desired. 

For the second claim, note that $\D'_\GS \mc F$ is right-cuspidal by Lemma \ref{lem: Cotangentleftrightcuspidality}. 
The above argument using the right-cuspidality instead implies that $j^* \D'_\GS \mc F\in \mc G^{< w_C}$, and hence $j^* \mc F\in \mc G^{\geq w_C +1}$. 
\end{proof}

\begin{corollary}{\label{cor: basadjunction}}
    For any $\mc{F} \in \Bas^{\Coh}(\T^{\ast}\V)$ and $\mc{F}' \in \Coh(\T^{\ast}\U_{C^{\circ}})$ we have
    \begin{equation*}
       \Hom(j^\circ_{\ush} \mc{F}' , \mc{F}) \simeq \Hom(\mc{F}', j^{\circ\ast} \mc{F}).
    \end{equation*}
\end{corollary}

We can now prove Theorem \ref{thm: ClassificationCuspidalCotangent} (2) for left-cuspidal coherent sheaves:
\begin{proposition}{\label{prop: cuspextopen}}
    Any left-cuspidal coherent sheaf $\mc{F} \in \Coh(\T^{\ast}\V)$ satisfies $\mc{F} \simeq j^{\circ}_{\ush} \mc{V}$ for some $\mc{V} \in \Coh(\T^{\ast}\U_{C^{\circ}})$. In fact, we can take $\mc{V} = j^{\circ \ast} \mc{F}$. Moreover, the restriction functor $j^{\circ \ast}$ is fully faithful on left-cuspidal coherent sheaves.
\end{proposition}
\begin{proof}
    
 We prove this by iterating on the total order $C^{\circ}=C_0, C_1, ..., C_k=\{0\}$. We adopt the notation that $j_i\coloneq j_{C_{i+1}}: \T^{\ast}\U_{C_i} \rightarrow \T^{\ast}\U_{C_{i+1}}$ and $j^i = j_{k-1}\circ \dots \circ j_i$. Firstly, $\U_{C_k}=\V$ and so Lemma \ref{lem: restriction lemma} implies that $\mc{F} \in \mc{W}_{C_k}$. Then Theorem \ref{thm: HL3.3.1} and the construction of $j_{\ush}$ imply that $j^{k-1}_{\ush}j^{k-1, \ast}\mc{F} \simeq \mc{F}$. But now $j^{k-1, \ast} \mc{F} \in \mc{W}_{C_{k-1}}$, so, 
\begin{equation*}
    j_{k-2, \ush} j^{k-2, \ast} \mc{F}\simeq j_{k-2, \ush}j^{\ast}_{k-2} j^{k-1, \ast} \mc{F}  \simeq j^{k-1, \ast} \mc{F}.
\end{equation*}
Applying $j_{k-1, \ush}=j^{k-1}_{\ush}$ then gives
\begin{equation*}
    j^{k-2}_{\ush}j^{k-2, \ast}\mc{F} \simeq j^{k-1}_{\ush}j^{k-1, \ast}\mc{F} \simeq \mc{F}.
\end{equation*}
Iterating this argument proves the claim.

Finally, the full faithfulness follows from the first claim and Corollary \ref{cor: basadjunction}.
\end{proof}
\begin{remark}{\label{rem: cohcleanness}}
    Working instead with the category $\mc{G}^{[w_C, w_C)}$, one can analogously construct a functor 
    \begin{equation*}
        j^{\circ}_{|\ast|}: \Coh(\T^{\ast}\U_{C^{\circ}}) \rightarrow \Coh(\T^{\ast}\V),
    \end{equation*}
which satisfies a right-adjunction property relative to $j^{\circ \ast}$ analogous to the left-adjunction property in Corollary \ref{cor: basadjunction}. One can similarly prove an analogue of Proposition \ref{prop: cuspextopen} showing that 
\begin{equation*}
    j^{\circ}_{|\ast|}j^{\circ \ast}\mc{F} \cong \mc{F},
\end{equation*}
for any right-cuspidal $\mc{F} \in \Coh(\T^{\ast}\V)$. We deduce a coherent analogue of Proposition \ref{prop: constrcleanness}. Namely that for any both left- and right- cuspidal $\mc{F} \in \Coh(\T^{\ast}\V)$ (such as cuspidal perfect complexes), we have that $\mc{F}$ is a clean extension from $\T^{\ast} \U_{C^{\circ}}$ in the sense that there is a sheaf $\mc{V} \in \Coh(\T^{\ast}\U_{C^{\circ}})$ such that
\begin{equation*}
    j^{\circ}_{\ush}\mc{V} \simeq \mc{F} \simeq j^{\circ}_{|\ast|} \mc{V}.
\end{equation*}
We will show that left-cuspidality and right-cuspidality agree. 
\end{remark}

\subsection{Non-existence of cuspidal coherent sheaves}{\label{ss: reductiondiscClassificationCuspidalCotangent}}
We reduce Theorem \ref{thm: ClassificationCuspidalCotangent} (3) in the graded setting to the essentially discrete case.  

\subsubsection{Reduction to reductive stabilizer}{\label{sss: cuspreductiontoreductive}}

In general, it may not be the case that $H_{c^{\circ}}$ is reductive for $c^{\circ} \in C^{\circ}$. If $H_{c^\circ}$ is not reductive, $\nu_{c^{\circ}}$ gives a nontrivial $\Theta$-stratum. We have 
\begin{equation*}
    \begin{tikzcd}
\T^{\ast}\U_{C^{\circ}}^{\nu_{c^{\circ}} = 0} \arrow[d, "j^{=0}"] & \T^{\ast}\U_{C^{\circ}}^{\nu_{c^{\circ}} \geq 0} \arrow[r, "p_c"] \arrow[d, "j^{\geq 0}"] \arrow[l, swap, "q_c"]&  \T^{\ast}\U_{C^{\circ}} \arrow[d, "j"] \\
 \T^{\ast}\V^{\nu_{c^\circ} = 0} & \T^{\ast}\V^{\nu_{c^\circ} \geq 0} \arrow[r, "p"] \arrow[l, swap, "q"] & \T^{\ast}\V .
    \end{tikzcd}
\end{equation*}
where $\T^{\ast}\U_{C^{\circ}}^{\nu_{c^{\circ}} = 0} \simeq \T^{\ast}BH_{M, c^{\circ}}$ and $H_{M, c^{\circ}}$ is the reductive quotient of $H_{c^{\circ}}$ (we have $H_{c^{\circ}} \subset P_{c^{\circ}}$ by Proposition \ref{prop: nuthetastrat}). In this case, $\mc{G}^{[w_{C^{\circ}}, w_{C^{\circ}}+1)} = \Eis^{\T^*}_{P_{c^{\circ}}}\Coh(Z)^0$ and $\mathcal{W}_{C^{\circ}}$ is trivial. Thus, by Lemma \ref{lem: restriction lemma}, any left-cuspidal $\mc{F} \in \Coh(\T^{\ast}\V)$ satisfies 
$j^{\circ\ast} \mc{F} \simeq 0$, and Proposition \ref{prop: cuspextopen} implies that
\[
\mc F\simeq j^\circ_{\ush} j^{\circ*} \mc F \simeq 0. 
\]

\subsubsection{Reduction to finite stabilizer}

In light of \S \ref{sss: cuspreductiontoreductive}, we assume in this subsection that $H_{c^{\circ}}$ is reductive (though possibly disconnected).

Recall that by Corollary \ref{cor: limcatclassifyingstack}, we have that $\Bas^{\Coh}(\T^{\ast}\U_{C^{\circ}})$ agrees with $\mrL'(\T^{\ast}\U_{C^{\circ}})$ and is generated under finite colimits by vector bundles on the $0$-section associated with representations of $H_{c^{\circ}}/[H_{c^{\circ}}, H_{c^{\circ}}]^\circ$.

\begin{proposition}{\label{prop: basrestrictisbas}}
   Suppose $\mc{F} \in \Bas^{\Coh}(\T^{\ast}\V)$. Then $j^{\circ \ast} \mc{F}$ is basic and the $Z(H^{\circ}_{c^{\circ}} \cap \wh{G}_{\der})^{\circ}$-weights of $j^{\circ \ast} \mc{F}$ are trivial.
\end{proposition}
\begin{proof}
We first prove that $j^{\circ\ast} \mc{F}$ is (left-)basic. Let $\nu$ be a cocharacter of $(H_{c^{\circ}})\cap \wh G_{\der}$ defining a parabolic subgroup $P_{\nu, c^{\circ}}$ of $H_{c^{\circ}}$ with unipotent radical $U_{\nu, c^{\circ}}$. Then $\nu$ determines a parabolic subgroup $P_{\nu}$ of $\wh{G}$.   Consider the diagram
\begin{equation}{\label{eqn: eispullbackdiagram}}
    \begin{tikzcd}
\T^{\ast}BH_{c^{\circ}}^{\nu = 0} \arrow[d, "j^{=0}"] \arrow[r, "\sigma_c", bend left = 15]& \T^{\ast}BH_{c^{\circ}}^{\nu \geq 0} \arrow[r, "p_c"] \arrow[d, "j^{\geq 0}"] \arrow[l, swap, "q_c"]&  \T^{\ast}BH_{c^{\circ}} \arrow[d, "j"] \\
 \T^{\ast}\V^{\nu = 0} \arrow[r, "\sigma", bend right = 15] & \T^{\ast}\V^{\nu \geq 0} \arrow[r, "p"] \arrow[l, swap, "q"] & \T^{\ast}\V .
    \end{tikzcd}
\end{equation}
Let $\mc{F}' \in \Coh(\T^{\ast}BH^{\nu = 0}_{c^{\circ}})$ with central weight $\chi$ such that $\langle \nu, \chi \rangle < 0$. Then
\begin{equation*}
    \Hom(j^{\ast}\mc{F}, \Eis_{P_{\nu, c^{\circ}}} \mc{F}') \simeq \Hom(\mc{F}, j_{\ast}\Eis_{P_{\nu, c^{\circ}}} \mc{F}') \simeq \Hom(\mc{F}, \Eis^{\T^*}_{P_{\nu}} j^{=0}_{\ast}\mc{F}'),
\end{equation*}
by Lemma \ref{lem: eisfiberprod}. We have $c^{\circ} \in \Lie M_{\nu}$ since $\nu$ factors through $H_{c^{\circ}}$. Thus $Z(M_{\nu}) \subset Z(H^{\nu = 0}_{c^{\circ}})$ and the central weight of $j^{=0}_{\ast} \mc{F}'$ is $\chi|_{Z(M_{\nu})}$. Hence the Hom spaces in the above equation vanish by basicness of $\mc{F}$.

We now prove the central weights are trivial. Consider a cocharacter $\nu \in X_{\ast}(Z(H^{\circ}_{c^{\circ}} \cap \wh{G}_{\der}))$. Then the centralizer in $\wh{G}$ of $\nu$ is a Levi subgroup $M_{\nu}$. We consider again the diagram \eqref{eqn: eispullbackdiagram}. In this case, we have $H^{\circ}_{c^{\circ}} \subset M_{\nu}$ so the stacks $\T^{\ast}BH_{c^{\circ}}^{\nu = 0}$ and $\T^{\ast}BH_{c^{\circ}}$ are quotients of the same scheme, differing only by the component group of $H_{c^{\circ}} \cap M_{\nu}$ compared to $H_{c^{\circ}}$. Moreover, the kernel of the map $P_{\nu} \rightarrow M_{\nu}$ is unipotent so $q_c\colon \T^{\ast}BH_{c^{\circ}}^{\nu \geq 0} \simeq \T^{\ast}BH_{c^{\circ}}^{\nu = 0}$.

Consider the pullback $q^{-1 \ast}_cp^{\ast}_cj^{\ast} \mc{F} \in \Coh(\T^{\ast}BH_{c^{\circ}}^{\nu = 0})$. Suppose this sheaf has a summand of central weight $\chi$ with $\langle \nu,\chi\rangle\neq 0$ and call it $\mc{F}_{\chi}$. Then up to possibly replacing $\nu$ with $- \nu$, we may assume the pairing $\langle \nu, \chi \rangle < 0$. Note that this does not change $q_c$ which is essentially the identity map. Then we have
\begin{equation*}
    0 \not\simeq \Hom(p^{\ast}_c j^{\ast} \mc{F}, q_c^*\mc{F}_{\chi}),
\end{equation*}    
but on the other hand this $\Hom$ must vanish by Lemma \ref{lem: eisfiberprod} and Lemma \ref{lem: cotangentdomswitch}. This is a contradiction.
\end{proof}

The previous propositions imply that the pullback of a basic sheaf to the open orbit is in the category generated by pushforwards along the $0$-section of reps of $H_{c^{\circ}}/[H_{c^{\circ}}, H_{c^{\circ}}]^\circ$.

\begin{proposition}{\label{prop: cuspreductiontofinite}}
    If $\Coh(\T^{\ast}\V)$ contains a nonzero left-cuspidal coherent sheaf, then $H_{c^{\circ}}$ is finite modulo $Z(\wh{G})$.
\end{proposition}
\begin{proof}
Suppose that $H_{c^{\circ}}$ is not finite modulo $Z(\wh G)$. We first assume that $H_{c^{\circ}}^\circ$ is not a torus. After replacing $c^\circ$ by an $H$-conjugate if necessary, let $B_H \subset H^{\circ}_{c^{\circ}}$ be a Borel subgroup containing $(H^{\circ}_{c^{\circ}}\cap \wh T)^\circ$ and with defining cocharacter $\nu$. Let $\wh{P} \subset \wh{G}$ be the parabolic induced by $\nu$. Then if nonzero $\mc{F} \in \LCusp^\Coh(\T^{\ast}\V)$, we have by Proposition \ref{prop: basrestrictisbas} that $j^{\circ \ast}\mc{F}$ is basic with trivial (up to $\wh{G}_{\ab}$) central weights and hence by Corollary \ref{cor: limcatclassifyingstack} is generated by pushforwards of vector bundles along the $0$-section. Lemma \ref{lem: zerosectiondirectsummand} below shows that such sheaves are generated by Eisenstein series of sheaves with $H_{c^{\circ}}$-basic central weight and hence $j^{\circ \ast}\mc{F}$ admits a nonzero map to such an Eisenstein series. Since in fact the $Z(H^{\circ}_{c^{\circ}} \cap \wh{G}_{\der})^{\circ}$-weights of $j^{\circ \ast}\mc{F}$ are trivial, we cannot have that the central weights $\chi$ of this Eisenstein series satisfy $\langle \nu, \chi \rangle > 0$ for some defining $\nu$. Hence by Lemma \ref{lem: eisfiberprod}, Lemma \ref{lem: cotangentdomswitch}, and adjunction, this contradicts left-cuspidality.

On the other hand, if $H^{\circ}_{c^{\circ}}$ is a non-central torus, which we denote $T_H$, then $Z_{\wh{G}}(T_H) = M_H$ is a Levi subgroup of $\wh{G}$ containing $H^{\circ}_{c^{\circ}}$. Hence we get an Eisenstein series diagram that is essentially the identity as in the second part of the proof of Proposition \ref{prop: basrestrictisbas}. Using a similar argument, we again conclude that $\Coh(\T^{\ast}\V)$ cannot contain a nonzero left-cuspidal sheaf.
\end{proof}

The following lemma was used in the proof of Proposition \ref{prop: cuspreductiontofinite}; we use different notation. 

\begin{lemma}{\label{lem: zerosectiondirectsummand}}
Let $H$ be a possibly disconnected reductive group with maximal torus $T$ and cocharacter $\nu \in X_{\ast}(T)$ such that the associated parabolic and Levi subgroups $M \subset P \subset H$ satisfy $M^{\circ}=T$ and $P^{\circ} \subset H^{\circ}$ is a Borel subgroup. Let $z_M: BM \rightarrow \T^{\ast}BM$ and $z: BH \rightarrow \T^{\ast}BH$ be the zero sections. Then if $\rho$ is a representation of $H$ thought of as a sheaf in $\Coh(BH)$, and on which $H^{\circ}$ acts by a character $\theta$, we have that $z_{\ast} \rho$ is a direct summand of $\Eis^{\T^{\ast}}_P z_{M \ast} \sigma$, where $\sigma = \Res^{H}_{M}\rho$.
\end{lemma}

\begin{proof}
We aim to show that as a complex, $\Eis^{\T^{\ast}}_P z_{M \ast} \sigma$ is a direct sum of its cohomology and the cohomology in degree $0$ is the induction $z_{\ast}((\Ind^{H}_{H^\circ P}\mathbf{1})\otimes \rho)$. This implies the lemma.

Consider the commutative diagram
\begin{equation*}
    \begin{tikzcd}
        (\T^{\ast}BH)^{\nu = 0} &  (\T^{\ast}BH)^{\nu \geq 0} \arrow[l, swap, "q"] \arrow[r, "p"] & \T^{\ast}BH \arrow[d, bend left=30, "\pi"] \\
        BM \arrow[u, "z_M"] & X \arrow[u, "\tilde{z}_M"] \arrow[l, swap, "\tilde{q}"] & BH \arrow[u, "z"],
    \end{tikzcd}
\end{equation*}
where the left square is cartesian. We have by base change
\begin{equation*}
    \Eis^{\T^{\ast}}_P z_{M \ast} \sigma = p_{\ast}q^{\ast} (z_{M \ast} \sigma \otimes \omega_{\V_U}) \simeq p_{\ast}(\tilde{z}_{M, \ast}\tilde{q}^{\ast} \sigma \otimes \omega_{\V_U}),
\end{equation*}
where $U$ is the unipotent radical of $P$ and $\omega_{\V_U}\coloneq\omega_{BP}\otimes \omega_{BM}^{-1}$. We have $(\T^{\ast}BH)^{\nu = 0} = (\Lie M)^{\vee}[-1] / M$ and as in Lemma \ref{lem: eiscotanfiberprodcomparison},
$(\T^{\ast}BH)^{\nu \geq 0} \simeq (\Lie P^{\op})^{\vee}[-1] / P$. Then 
\begin{align*}
  X \simeq (\Lie U^{\op})^\vee[-1] /P  & \simeq  (H \times^P (\Lie U^{\op})^\vee[-1]) /H\\
  & \simeq \T^{\ast}(H/P)[-1] / H \simeq \Spec_{H/P}\Sym(T_{H/P}[1]) / H. 
\end{align*}
Now, $\Eis^{\T^{\ast}}_P z_{M \ast} \sigma$ is computing the pushforward to $\T^{\ast}BH$ of a twist by $\omega_{\V_U}$ of a vector bundle $\mc{L}_{\sigma}$ over $X$. If we push forward along $\pi$, we get the underlying complex attached to $\Eis^{\T^{\ast}}_P z_{M \ast} \sigma$ which is
\begin{equation}{\label{eqn: underlyingrepcohomologyformula}}
    \Gamma(H/P, \Sym(T_{H/P}[1]) \otimes \mc{L}_{\sigma} \otimes \omega_{\V_U}) \simeq \bigoplus\limits_{i=0}^d \Gamma(H/P, \Omega^{d-i}_{H/P} \otimes \mc{L}_{\sigma})[i-d],
\end{equation}
where $d = \dim U$ and the last isomorphism comes from $\omega_{\V_U} \simeq \Omega^d_{H/P}[-d]$. Then $H/P$ is a (disconnected) flag variety and the $H$-equivariant vector bundle $\mc{L}_{\sigma}$ is isomorphic to $\mc{O}_{H/P} \otimes \rho$. It follows from the projection formula that
\begin{equation*}
    H^a(H/P, \Omega^b \otimes \mc{L}_{\sigma}) \cong \Ind^{H}_{H^{\circ}P} (H^a(H^{\circ}/P^{\circ}, \Omega^b)) \otimes \rho.
\end{equation*}
Then $H^a\Gamma(H^{\circ}/P^{\circ}, \Omega^b) \cong 0$ whenever $a \neq b$ since $H^{\circ}/P^{\circ}$ is a connected flag variety whose cohomology is generated by Schubert classes. Hence, the same is true of $H^a\Gamma(H/P, \Omega^b \otimes \mc{L}_{\sigma})$. Since $\Gamma(H^{\circ}/P^{\circ}, \Omega^a)[-a]$ is concentrated in degree $2a$, this means \eqref{eqn: underlyingrepcohomologyformula} is concentrated in even degrees. Finally, when $a=0$, we have 
\begin{equation*}
    H^0\Gamma(H/P, \Omega^0 \otimes \mc{L}_{\sigma}) \cong (\Ind^{H}_{H^\circ P} \mathbf{1}) \otimes \rho.
\end{equation*}

To conclude, we need to show the desired splitting occurs over $\T^{\ast} BH$. For any complex in $\Coh(\T^{\ast}BH)$, its cohomology sheaf lies in the image of $z$, and for $V, W \in \Coh(BH)$ in degree 0, we have
\begin{equation*}
    \Hom_{\T^{\ast} BH}(z_{\ast} V , z_{\ast} W) \simeq \Hom_{BH}(V , W\otimes \Sym((\Lie H)^\vee[-2])),
\end{equation*}
hence lies in even degrees. Any extension class between two cohomology sheaves of $\Eis^{\T^{\ast}}_P z_{M \ast} \sigma$ lies in odd degree and hence is trivial. This proves the desired result.
\end{proof}

\section{Construction of cuspidal vector bundles}{\label{s: CVBs}}
We prove Theorem \ref{thm: ClassificationCuspidalCotangent} (4) in the graded setting by constructing cuspidal vector bundles on $\T^* \V_\lambda$ corresponding to cuspidal representations $\rho$ of $S_\phi$. We work in the graded setting of \S \ref{s: cotangentgeometry} and hence use the notation of that section including $\iota, V,H,\V, H_{c^{\circ}}$. In the general setting, these correspond to $\lambda, V_{\lambda}, H_{\lambda}, \V_{\lambda}, S_{\phi}$.

\subsection{Fundamental cuspidal inequalities}{\label{ss: fundcuspineq}}

We will eventually prove that $\Cusp^{\Coh}(\T^{\ast}\V)$ is generated by cuspidal vector bundles and that these are determined by explicit inequalities on the weights of their underlying representations. These inequalities allow for an explicit classification and also imply strong structural properties of the various categories of cuspidal sheaves. In this subsection, we derive these weight inequalities determining cuspidal vector bundles.

To begin, we discuss vector bundles on $\T^{\ast}\V$.
\begin{lemma}{\label{lem: equivnakayama}}
Consider the maps
\begin{equation*}
   BH \rightarrow \V \xrightarrow{z}  \T^{\ast}\V \rightarrow BH,
\end{equation*}
where the first map is determined by $0\in V$ and the last is the projection. Then the composition is the identity map and any vector bundle on $\V$ or $\T^{\ast}\V$ is given by pullback of a vector bundle on $BH$.
\end{lemma}
\begin{proof}
The first claim is clear, and the second is essentially a corollary of the equivariant Nakayama lemma of Bass and Haboush \cite[Corollary (6.4)]{BassHaboush} as the vector bundle is determined by its restriction to the underlying classical stack. Note that to apply this corollary, we use that $V$ and $\mu^{-1}(0)$ have unique closed orbits (see proof of Proposition \ref{prop: TVthetastrat} (1) to reduce the closed orbit claim for $\mu^{-1}(0)$ to $V$ where it is clear). 
\end{proof}

We now derive weight inequalities satisfied by cuspidal vector bundles. Note $\wh{T}$ is a maximal torus of $H$. For each nonzero $\nu \in X_{\ast}(\wh{T}_{\der})$ (the intersection of $\wh{T}$ with $\wh G_{\der}$), we have an associated parabolic subgroup $P_{\nu} = M_{\nu}U_{\nu}\subset \wh G$ and a diagram of cotangent stacks
\begin{equation*}
\begin{tikzcd}
      \T^{\ast} \V^{\nu = 0} \arrow[r, "\sigma", bend left = 15] & \T^{\ast} \V^{\nu \geq 0} \arrow[r, "p"] \arrow[l, "q"] & \T^{\ast} \V.
\end{tikzcd}
\end{equation*}
By Proposition \ref{prop: basicweightcharacterization}, to determine if a vector bundle $\mc{V} \in \Coh(\T^{\ast}\V)$ is cuspidal, it suffices to check that $p^{\ast} \mc{V} \in \QCoh(\T^{\ast} \V^{\nu \geq 0})^{>\langle \nu, \wt(\omega_{\V_{U_\nu}})\rangle}$. By part (1) of Lemma \ref{lem: HLlemmas}, it suffices to compute the $\nu$-weights of $\sigma^{\ast} p^{\ast} \mc{V}$. We have a diagram
\begin{equation*}
    \begin{tikzcd}
       \T^{\ast}\V^{\nu = 0} \arrow[r, "p \circ \sigma"] \arrow[d, "\pi_{\nu}"] & \T^{\ast} \V \arrow[d, "\pi"] \\
       BH_{M_{\nu}} \arrow[r, "i"] & B H.
    \end{tikzcd}
\end{equation*}
Now, $\mc{V} \cong \pi^{\ast} \gamma$ for some representation $\gamma$ of $H$. Hence,
\begin{equation*}
   (\pi \circ p \circ \sigma)^{\ast} \gamma \simeq (i \circ \pi_{\nu})^{\ast} \gamma \simeq \pi^{\ast}_{\nu}\res^{H}_{H_{M_{\nu}}} \gamma.
\end{equation*}
Then the $\nu$-weights of $\pi^{\ast}_{\nu} \res^{H}_{H_{M_{\nu}}} \gamma$ are given by the pairings $\langle \nu, \chi \rangle$, for $\chi$ a $\wh{T}$-weight of $\gamma$ (\S \ref{sss: derivedtheta strat}). Thus $\mc{V}$ is cuspidal precisely when each weight $\chi$ of $\gamma$ satisfies
\begin{equation*}
    \langle \nu , \chi \rangle > \langle \nu, \wt(\omega_{\V_{U_\nu}}) \rangle.
\end{equation*}
Expanding, we get 
\begin{equation*}
     \langle \nu , \chi \rangle > -\sum\limits_{v \in \wt(V \cap \Lie U_{\nu})} \langle \nu, v \rangle + \sum\limits_{v \in \wt(\Lie H \cap \Lie U_{\nu})} \langle \nu, v \rangle.
\end{equation*}
The weights $\wt(V \cap \Lie U_{\nu})$ are the $\nu$-positive weights of $V$ (roots of $\wh G$). We denote the set of $\nu$-positive roots by $\Phi^{>0}_{\nu}$. By duality, we also get an upper bound for $\langle \nu, \chi \rangle$. Hence, the inequalities are
\begin{equation}\label{eq: CuspidalInequalitiesAll}
 -\sum\limits_{v \in \wt(V) \cap \Phi^{>0}_{\nu}} \langle \nu, v \rangle + \sum\limits_{v \in \wt(\Lie H) \cap \Phi^{>0}_{\nu}} \langle \nu, v \rangle < \langle \nu, \chi \rangle <  -\sum\limits_{v \in \wt(V) \cap \Phi^{<0}_{\nu}} \langle \nu, v \rangle + \sum\limits_{v \in \wt(\Lie H) \cap \Phi^{<0}_{\nu}} \langle \nu, v \rangle.
\end{equation}
We observe that the distance between these bounds is precisely 
\begin{equation}{\label{eqn: bmkleqn}}
    \sum\limits_{v \in \wt(V)} | \langle \nu, v \rangle | - \sum\limits_{v \in \wt(\Lie H)} |\langle \nu, v \rangle|,
\end{equation}
and that the midpoint is given by $\chi$ such that
\begin{equation*}
    \langle \nu, \chi \rangle + [\sum\limits_{v \in \wt(V) \cap \Phi^{>0}_{\nu}} \langle \nu, v \rangle - \sum\limits_{v \in \wt(\Lie H) \cap \Phi^{>0}_{\nu}} \langle \nu, v \rangle] =  [-\sum\limits_{v \in \wt(V) \cap \Phi^{<0}_{\nu}} \langle \nu, v \rangle + \sum\limits_{v \in \wt(\Lie H) \cap \Phi^{<0}_{\nu}} \langle \nu, v \rangle] - \langle \nu, \chi \rangle.
\end{equation*}
Solving, we get
\begin{equation*}
    \langle \nu, \chi \rangle = \frac{1}{2}[- \sum\limits_{v \in \wt(V)} \langle \nu, v \rangle + \sum\limits_{v \in \wt(\Lie H)} \langle \nu, v \rangle].
\end{equation*}
We conclude that $\chi$ is exactly half the total weight $\delta$ of the dualizing complex $\omega_{\V}$.

Finally, as in the proof of \cite[Proposition 2.5]{BMKL}, it suffices to consider $\nu = w^{-1}(\delta_{\alpha})$ for $\delta_{\alpha}$ the fundamental coweight of a simple root $\alpha$ and $w$ in the Weyl group $W(\wh{G}, \wh{T})$. This reduces to a finite number of inequalities of the form
\begin{equation}{\label{eqn: weylcuspineq}}
  -[\sum\limits_{v \in w(\wt(V)) \cap \Phi^{>0}_{\delta_{\alpha}}} v - \sum\limits_{v \in w(\wt(\Lie H)) \cap \Phi^{>0}_{\delta_{\alpha}}} v]_{\alpha}  < [w(\chi)]_{\alpha} < -[\sum\limits_{v \in w(\wt(V)) \cap \Phi^{<0}_{\delta_{\alpha}}}v - \sum\limits_{v \in w(\wt(\Lie H)) \cap \Phi^{<0}_{\delta_{\alpha}}} v]_{\alpha},
\end{equation}
where $[ - ]_{\alpha}\coloneq \langle \delta_\alpha, -\rangle$, $\alpha\in \Delta$ denotes the coefficient vector with respect to the set of simple roots $\Delta$. 

It is shown in \cite{BMKL} that the non-strict version of these inequalities defines a convex polytope in $X^{\ast}(\wh{T}_{\der})_{\R}$ denoted $\mc{P}_{\cusp}$ and that 
  \begin{equation*}
        \mc{P}_{\cusp} = \left(\bigoplus_{\alpha \in \wt(V)} [-\alpha, 0]\right) \bigominus \left(\bigoplus_{\beta \in \Phi(H)} [-\beta, 0]\right),
    \end{equation*}
where $\oplus$ and $\ominus$ denote Minkowski sum and difference and $[-\alpha,0], [-\beta, 0] \subset X^{\ast}(\wh{T}_{\der})_{\R}$ are closed intervals. Then we can restate the condition of cuspidality of $\mc{V}_{\gamma}$ as (the restriction of) all weights of $\gamma$ lie in the interior $\mc{P}^{\circ}_{\cusp}$ of $\mc{P}_{\cusp}$ (where the interior is taken in the topology of $X^{\ast}(\wh{T}_{\der})_{\R}$, so if $\mc{P}_{\cusp}$ has positive codimension, the interior is empty).

\begin{lemma}
If $\mc O_{\T^*\V}\otimes B$ is a vector bundle on $\T^*\V$ such that the weights of $B$ are in $\mc{P}^{\circ}_{\cusp}$ and trivial on $Z(\wh G)^\circ$, then 
\[
\Psi(\mc O_{\T^*\V}\otimes B)\in \MF^{\gr}(\M(\mc V^\vee),f )\simeq \mrL (\T^*\V).
\]
\end{lemma}

\begin{proof}
By the definition of $\Psi$ using the Koszul factorization, it suffices to show that all the weights $\chi$ in the terms of the Koszul resolution of $\mc O_{\T^*\V}\otimes B$ in $\T^*V/H$ satisfy \eqref{eq: MagicAllWeights} for all nonzero $\nu \in X_* (\wh T \cap \wh G_{\der})$, as the weights are assumed to be trivial on $Z(\wh G)^\circ$. Then, it suffices to show
\[
\frac{1}{2}\langle \nu, -\det (V^{\nu >0}) +\det (V^{\nu<0})-\det(V)\rangle
< \langle \nu, \chi \rangle 
< \frac{1}{2}\langle \nu, \det (V^{\nu >0}) -\det (V^{\nu<0})-\det(V)\rangle,
\]
equivalently
\[
\langle \nu, -\det (V^{\nu >0})\rangle
< \langle \nu, \chi \rangle 
<  \langle \nu, -\det (V^{\nu<0})\rangle. 
\]
It follows from \eqref{eq: CuspidalInequalitiesAll} as $\chi= \chi' + \chi_H$ for some multiplicity-free sum $\chi_H$ of weights in $\Lie H$. 
\end{proof}

\begin{remark}\label{rem: CuspidalSubsetLimit}
If $\V$ is essentially discrete, then the same computation implies that it also holds for the vector bundle attached to $B$ on the underlying classical stack of $\T^*\V$ since the derived structure comes only from $\Lie Z(\wh G)$ (together with \cite[Lemma 2.23]{PTqBPSsymquiverswithpotential}). In the notation of the classification theorems, it in particular means $\mc G_\rho$ defined in \S \ref{ss: filtered lifts} for a cuspidal $\rho$ trivial on $Z(\wh G)^\circ$ is in the limit category.  
\end{remark}

\subsection{The simply connected, split, unipotent case}{\label{ss: classificationunipotentcase}}

We now assume $\wh{G}$ is connected, simply connected, split, simple and use \cite{BMKL} to work out the classification of irreducible cuspidal vector bundles. We reduced to this case in \S \ref{ss: reductionurclassificationcuspidalcotangent}. In this case, \cite[Theorem 1.2]{BMKL} and \eqref{eqn: bmkleqn} imply that any cuspidal vector bundle must be supported on $\V$ with distinguished open orbit, which, by Lemma \ref{lem: discequalsdistforunip}, is the same as corresponding to an essentially discrete $L$-parameter. As in Example \ref{ex: Ilambda}, we may choose $\iota$ such that $V = \Lie \wh{G}(2)$. 

We will prove in Proposition \ref{prop: CVBimpliesCLS} that any cuspidal vector bundle on an essentially discrete Vogan stack gives a nontrivial cuspidal local system; the proof is inductive and uses the results in this section for the cases related to proper Levi subgroups. In the present case, these have been classified by Lusztig \cite{LusztigIntersection} and this significantly reduces the number of cases we need to consider. See also \cite{BalaCarterI} (also \cite{BalaCarterII}, \cite[\S8.4]{CollingWoodMcGovern}) where the various distinguished unipotent orbits are classified. We first claim the results of \S\ref{s: cotangentgeometry} imply that each cuspidal local system has at most one lift to a left-cuspidal coherent sheaf. Indeed, 
\begin{itemize}
    \item a local system on the open orbit of $\V$ corresponds to a coherent sheaf $\mc{F}$ on the open $\Theta$-stratum of $\T^{\ast}\V$  (note both the open orbit of $\V$ and open $\Theta$-stratum of $\T^{\ast}\V$ identify with $BH_{c^{\circ}}$ and we are identifying coherent and constructible sheaves on a finite classifying stack),
    \item by the theory of window categories (Theorem \ref{thm: HL3.3.1}, Lemma \ref{lem: restriction lemma}, Proposition \ref{prop: cuspextopen}) there is at most one left-cuspidal coherent sheaf on $\T^{\ast} \V$ whose restriction to the open stratum is $\mc{F}$, 
    \item and it is given by $j_{\ush}\mc{F}$, where $j_{\ush}$ is the exotic left-adjoint of the pullback functor $j^{\ast}$ to the open $\Theta$-stratum.
\end{itemize}
Taking into account that the restriction to the open stratum is fully faithful (Proposition \ref{prop: cuspextopen}), it suffices to produce, for each $\V$ containing a cuspidal local system, the same number of irreducible cuspidal vector bundles.

We now examine each case. In what follows, we normalize $\iota$ to be antidominant.
\subsubsection{A-type}
In the $A$-type case, the unique distinguished orbit is the regular one. In this case, $\wh{G} = \SL_{n+1}$. Let $\wh T \subset \SL_{n+1}$ be the diagonal torus. Then $X^{\ast}(\wh T)$ is given by $\Z^{n+1} / \Z(1, \dots, 1)$ by restricting the standard basis vectors of the diagonal torus of $\GL_{n+1}$. We have $\V \cong \A^n/\mathbb{G}_m^n$ and the weights of $V$ are the negative simple roots. We have $\delta=(1,0, \dots, 0, -1) = \sum\limits^{n}_{i=1} \sigma_i$, where $\sigma_i = e_i - e_{i+1}$.

We consider the inequalities \eqref{eqn: weylcuspineq} for $w=1$. These imply that any cuspidal weight $\mu$ satisfies $0<\mu< \sum\limits^n_{i=1} \sigma_i$ and hence that $\mu = \sum\limits_{i=1}^n r_i\sigma_i$ for $0<r_i<1$. This equality most naturally lives in $X^{\ast}(\wh T)_{\Q} \cong \Q^{n+1} / \Q$ but $\mu$ must be a weight and hence lie in the image of $X^{\ast}(\wh T)$. So written in coordinates, we have some $\lambda \in \Q$ such that
\begin{equation*}
    \mu = (r_1 + \lambda, r_2 -r_1 + \lambda, \dots, r_n - r_{n-1} + \lambda, -r_n + \lambda) = (a_1, \dots, a_{n+1}) \in \Z^{n+1}.
\end{equation*}
The element $\lambda$ is well-defined up to an element of $\Z$ so we choose $0<\lambda<1$ which implies that $a_1=1$.
Adding up all the coefficients gives $a_1+\dots+a_{n+1} = (n+1)\lambda$, and the left-hand side is an integer so $\lambda = \frac{m}{n+1}$. For each $k \in \{ 1, \dots, n\}$ we get an equation $a_1+\dots+a_k = r_k + k\lambda$. Since the left-hand side is an integer, it must equal $[\frac{km}{n+1}]+1$ and we must have that $k\lambda$ is not an integer (otherwise $r_k$ would be an integer). Thus $(m,n+1)=1$. This gives $\phi(n+1)$ (totient function) many choices for $m$. If $[\frac{km}{n+1}]=[\frac{(k-1)m}{n+1}]$ then $a_k=0$ ($2\leq k\leq n$), otherwise $a_k=1$, and $a_{n+1}=0$.

It remains to show that these $\phi(n+1)$ weights are cuspidal. We can use the equation for $\zeta(w)_{\sigma_i}$ in \cite[\S2]{BMKL}. Indeed, the inequalities \eqref{eqn: weylcuspineq} for a fixed $w$, specialized to $\sigma_i$ give a range of possible values for the coefficient of $\sigma_i$ in $w(\mu)$ in order for $\mu$ to be cuspidal. The midpoint of this range is given by $[w(\delta/2)]_{\sigma_i}$. Hence, it suffices to show that each candidate cuspidal weight $\mu$ has the $\sigma_i$ coefficient of $w(\mu)$ close enough to $w(\delta/2)$ to lie within the necessary range. For the weights $\mu$ considered in the previous paragraph, we have $|[\mu - (\delta/2)]_{\sigma_j}|<\frac{1}{2}$ for each $j$. We now fix some $\sigma_i$ and $w$. Then $w(\sigma_j)$ has a nonzero $\sigma_i$ coefficient if and only if $w$ takes one of  $\{j, j+1\}$ to some number $\leq i$ and the other to some number $> i$. This is the same condition that $\sigma_j$ contributes to $\zeta(w)_{\sigma_i}$. The contribution from $\sigma_j$ as a constituent of $\mu-\delta/2$ is then less than $\frac{1}{2}$ in absolute value, while the contribution to $\zeta(w)_{\sigma_i}$ is $1$, and hence the contribution to the ``radius'' around $[w(\delta/2)]_{\sigma_i}$ is exactly $1/2$. This implies the weights are cuspidal.

The group $H$ is a torus and so each cuspidal weight yields an irreducible character of $H$ and hence a cuspidal vector bundle on $\T^{\ast}\V$.
 
\subsubsection{B-type}
The simply connected groups are $\Spin_{2n+1}$ and there are three cases we need to consider.

Firstly consider the partitions $1+3+5+\cdots$ of $2n+1$. This corresponds to distinguished diagrams with weights $2020020\cdots 00$, where the last node corresponding to the short simple root is $0$: the number of $0$s in each block (of the form $2$ followed by some number of $0$s) increases by $1$ to the second to last block, with $n_{k-1}-1$ zeroes and $n_{k-1}$ will be even. Then the last block has $
\frac{n_{k-1}}{2}$ many $0$s. The corresponding unipotent essentially discrete parameter is given by the sum of the $1,3,5, \dots, 2k+1$ dimensional representations of $\SL_2$. Conjugating such that the image of $\begin{pmatrix} t & 0 \\ 0 & t^{-1} \end{pmatrix}$ is antidominant,  we get $k+1$ copies of the $0$-weight in the ``middle'' followed by the other weights arranged symmetrically. Since $k+1$ must be odd, this works. It is easy to check that the dualizing weight $\delta$ is given by $2e_1 + \cdots +2e_{n-\frac{k}{2}}$ and hence $\frac{\delta}{2}$ is indeed integral and therefore cuspidal. We claim that the representation of $H$ with highest weight $\frac{\delta}{2}$ is a character, and hence gives a cuspidal vector bundle. Indeed, this is implied by the fact that $\delta$ is the highest weight of a line bundle (for instance, see the proof of \cite[Lemma 4.1]{BMO}). 

In the second case, we can consider $1+5+9+\cdots=2n+1$, a sum of $k$ numbers. Then following the formula in Bala--Carter, this will give weighted Dynkin diagram equal to $22020200200\cdots 0$ where the penultimate two strings of the form $20\cdots 0$ are of the same length $n_{k-1}$ and the final one (including the final node) is of length $\frac{n_{k-1}}{2} + 1$. So for example $1+5+9$ gives $2202020$. Then $\delta$ will be $e_1+e_2+\cdots$ for each $i$ except the ones contributing to the $0$'s in the final string. The weights of $\Spin_{2n+1}$ have $e_i$ coefficients that are either integral for all $i$ or half-integral for all $i$. Hence $\frac{\delta}{2}$ will not be integral but by adding $\pm \frac{1}{2}e_i$ for each $i > n - \frac{k-1}{2}$, we can make an integral weight and these give a single orbit under the Weyl group of $H$ as can be seen from the weighted Dynkin diagram. 

These weights are cuspidal. To prove this, we use \cite[(3)]{BMKL}. The weights in question differ from $\frac{\delta}{2}$ by $\epsilon = \sum\limits_{i > n - \frac{k-1}{2}} \pm \frac{1}{2}e_i$. Hence, we just need to show that for each $w$ and $\gamma \in \Delta$, we have $\frac{1}{2}\zeta_{B_n}(w)_{\gamma} > \lvert[w(\epsilon)]_{\gamma}\rvert$. The term $\epsilon$ is a half sum of a number $j$ of short roots $e_i$. Each half-short root can contribute at most $\frac{1}{2}$ to $w(\epsilon)_{\gamma}$ and hence the total contribution is at most $\frac{j}{2}$. Moreover, the number of $w(e_i)_{\gamma}$ that are nonzero is exactly $a_m$. Hence when $a_m \neq 0$, $\frac{1}{2}\zeta_{B_n}(w)_{\gamma}\geq a_m > \frac{a_m}{2} \geq \lvert[w(\epsilon)]_{\gamma}\rvert$. When $a_m=0$, the desired bound follows from $\zeta_{B_n}(w)_{\gamma} > 0$.

The convex hull of this Weyl orbit of cuspidal weights does not contain any other weights and the Weyl orbit consists of the weights of a (minuscule) highest weight representation of $H$ corresponding to a cuspidal vector bundle. In fact these half integral weights are related to those of the spin representation of an odd spin group.

In the third case, we consider $3+7+11+\cdots=2n+1$ a sum of $k$ numbers. Then the weighted diagram is of the form $22020200200\cdots 0$ where the penultimate two strings are of different lengths $n_{k-1} - 1$ and $n_{k-1}$ and the final string is length $\frac{n_{k-1}-1}{2} +1$. So for example $3+7+11$ gives $2202020020$. Again, $\delta$ is a sum of the $e_i$ except for the ones contributing to the $0$s of the last string, and again, the cuspidal weights are equal to $\frac{\delta}{2}$ plus $\pm \frac{1}{2}e_i$ for the ``missing'' $e_i$ and the argument that these give a cuspidal vector bundle is essentially the same as before.

\subsubsection{C-type}

The simply connected groups are $\Sp_{2n}$. The relevant partitions of $2n$ are of the form $2+4+6+\cdots$. These correspond to the weighted Dynkin diagram $2020020\cdots 02$ (number of $0$s in each block increases by $1$ and the very last entry is a $2$) where the final node corresponds to the long simple root. The parameter is given by the sum of the $2, 4, 6,\dots,2k$ dimensional representations of $\SL_2$. The weight of a determinant bundle is the sum of the weights of the original bundle. To compute this, we conjugate the torus of $\Sp_{2n}$ such that the weights of the diagonal torus of the Deligne $\SL_2$ appear in increasing order. Then it is easy to check that $\omega_{\V} = \det(V)^{-1}$ is of weight $2e_1+ \cdots + 2e_n$ and its square root has weight $e_1+\cdots+e_n$ and hence is cuspidal and gives a rank one cuspidal vector bundle as in the first $B$-type case. The restriction of this weight to the center gives the nontrivial character precisely when $n$ is odd, matching Lusztig's classification. 

\subsubsection{D-type}
We now consider the situation for simply connected $D_n$, i.e., groups of the form $\Spin_{2n}$ where there are three relevant cases.

In the first case, we take $1+3+5+\cdots=2n$. Then we get a weighted Dynkin diagram $202002000\cdots$ with both tail nodes equal to $0$. Then $\delta$ is $2e_1 + 2e_2 + \cdots$ for each $e_i$ appearing to the left of the $0$s in the final string. Then $\frac{\delta}{2}$ is integral hence cuspidal and gives a rank one cuspidal vector bundle.

In the second case, we consider $1+5+9+\cdots=2n$ (hence $2n$ ranges over half the even triangular numbers) which gives diagram $22020200200\cdots$ where the penultimate strings are the same length and contain $n_{k-1} - 1$ zeroes, and the final string containing the tail nodes (which have weight $0$) has $\frac{1}{2}(n_{k-1})$ or $\frac{1}{2}(n_{k-1} + 1)$ zeroes depending on parity. In this case, we have $\delta = e_1+e_2+\cdots$ for each $i$ except the ones contributing to the $0$s of the final string. Then $\delta/2$ is not integral. We consider the weights given by adding $\pm \frac{1}{2}e_i$ for each $i$ appearing in the $0$ weights of the final string. These weights are integral and we claim they are also cuspidal. Let $\epsilon = \sum\limits_i \pm \frac{1}{2}e_i$, where $i$ is an index belonging to the last string. We need to show that $\frac{1}{2} \zeta_{D_n}(w)_{\gamma} > \lvert[w(\epsilon)]_{\gamma}\rvert$ for each $w$ and $\gamma$. By computing $w(e_i)$, we see the value of $\lvert [w(\epsilon)]_{\gamma}\rvert$ is at most $\frac{1}{2}a_m$ when $\gamma$ is not a tail root and $\frac{1}{4}a_m$ otherwise. On the other hand, $\frac{1}{2}\zeta_{D_n}(w)_{\gamma}$ is at least $a_m$ when $\gamma$ is not a tail root and $\frac{1}{2}a_m$ otherwise. When $a_m=0$, the desired bound follows from $\zeta_{D_n}(w)_{\gamma} > 0$. This proves cuspidality. The group $H$ is a Levi subgroup of $\wh{G}$ and has a factor equal to an even spin group of rank $m$ and these cuspidal weights correspond to two different irreducible representations of dimension $2^{m-1}$. Hence there are two cuspidal vector bundles of rank $2^{m-1}$ in this case. 

The third case is a partition $3+7+11+\cdots=2n$ (hence $2n$ ranges over the other half of the even triangular numbers) where the penultimate strings before the final one now have different lengths. This case is almost identical to the previous and the cuspidal weights are the weights of the two half-spin representations and hence there are two cuspidal vector bundles. However, we remark that the $D_5$ example corresponding to partition $3+7$ has weighted Dynkin diagram $22022$, in particular the tail nodes are $2$s. However, one can check that cuspidal weights in this example agree with the other cases. Namely, they are $\delta/2 \pm \frac{1}{2}e_5$ and correspond to two distinct cuspidal vector bundles of rank $1$.

\subsubsection{Exceptional groups}{\label{sss: exceptionalgroups}}

In \cite{BMKL} it is proven that the distinguished unipotent orbits have nonempty cuspidal polytopes and this immediately implies that the weight $\delta/2$ is cuspidal if it is integral. For $G_2, F_4, E_7, E_8$ and the unique relevant distinguished orbit in each case ($G_2(a_1), F_4(a_3), E_7(a_5), E_8(a_7)$ in \cite{BalaCarterII}), one can compute that these are indeed integral and hence give lifts of the unique cuspidal local system in each case. The remaining case is $E_6(a_3)$ and the cuspidal weights are $\omega_4 - \omega_1$ and $\omega_4 - \omega_6$ regarded as characters of $H$, where $\omega_i$ are the fundamental weights with Bourbaki's numbering. Computer code (written manually) due to Jonathan Leake showing these lie in the cuspidal polytope is given at his webpage \href{https://www.jleake.com/res/distinguished-rho-positivity/#cuspidal-weights}{here}.

\subsubsection{Consequences}

Thus, we have shown that the cuspidal vector bundles are in natural bijection with the cuspidal local systems of \cite{LusztigIntersection}. In each case, one can check that the restriction of the cuspidal vector bundle to the open orbit of $\V$ recovers a cuspidal representation of $H_{c^{\circ}}$ matching Lusztig's classification. 

\section{Generalized coherent and constructible Springer sheaves}{\label{s: Springer}}

We have proven weaker versions of (1) and (2) of Theorem \ref{thm: ClassificationCuspidalCotangent} in the graded setting and also constructed cuspidal vector bundles. It follows by computation or an argument in this section that the $\ast$-pushforward of a cuspidal representation $\rho$ of $S_{\phi}$ along the zero section of $\T^{\ast} BS_{\phi}$ is in the image of the restriction
\begin{equation*}
    j^*\colon \Cusp^\Coh(\T^{\ast} \V_{\lambda}) \rightarrow \Coh(\T^{\ast} BS_{\phi})
\end{equation*}
for the open immersion $j\colon \T^* BS_\phi\hookrightarrow \T^*\V_\lambda$, showing (4). It remains to show that the image of $j^*$ restricted to cuspidal sheaves is generated by such pushforwards. 

To prove (3) in the graded case, we may assume $\Lambda=\C$, and we relate generalized coherent and constructible Springer sheaves. We will directly work with coherent Springer sheaves on cotangent stacks, but they are related to coherent Springer sheaves on $\Par_G$, or rather their $!$-restrictions via Koszul duality. 

The argument is inductive and, along the way, we will also prove Theorem \ref{thm: BasicCotangentMain} (3) in the graded setting, and the full versions of (1) and (2) of Theorem \ref{thm: ClassificationCuspidalCotangent}. 

\subsection{Filtered lifts of cuspidal sheaves on cotangent stacks}{\label{ss: filtered lifts}}
We introduce some notation first, assuming Theorem \ref{thm: ClassificationCuspidalCotangent} holds. Let $\lambda\colon W_F\to \wh G$ be the semisimplification of an essentially discrete $L$-parameter $\phi$, let $\rho$ be a cuspidal representation of $S_\phi$, and let $\mc G_{\rho}$ denote the cuspidal sheaf on $\T^* \V_\lambda$ that restricts to $z_*\rho$ for $z\colon BS_\phi\to \T^* BS_\phi$. It is a vector bundle on the underlying classical stack of $\T^* \V_\lambda$, which identifies with the fiber product
\[
B H_\lambda\times_{B(H_{\lambda}/Z(\wh G)^{\Gamma, \circ})}\T^* (V_\lambda/(H_{\lambda}/Z(\wh G)^{\Gamma, \circ}))\hookrightarrow
\T^* \V_\lambda.  
\]
In the graded setting, $\mc G_\rho$ is in the limit category by Remark \ref{rem: CuspidalSubsetLimit} if $\rho$ is trivial on $Z(\wh G)^\circ$.   

We shall now work in the graded setting and assume that $\rho$ is trivial on $Z(\wh G)^{\circ}$. We write $\overline{\rho}$ for the resulting representation of $S_\phi/Z(\wh G)^\circ$. By (the proof of) Lemma \ref{lem: equivnakayama}, $\mc G_\rho$ is pulled back from a representation of $H_\lambda$. Its restriction to $S_\phi$ is $\rho$, so it is trivial on $Z(\wh G)^{\circ}$; let $\mc V_{\overline\rho}$ be the corresponding vector bundle on $V_\lambda/(H_\lambda/Z(\wh G)^{\circ})$, with the trivial grading. Then $\mc G_\rho$ is obtained from $\mc V_{\overline\rho}|_{\T^* ( V_\lambda/(H_\lambda/Z(\wh G)^{\circ}))}$ on $\T^* ( V_\lambda/(H_\lambda/Z(\wh G)^{\circ}))$ by $*$-pulling back to the classical stack displayed above and $*$-pushing forward to $\T^*\V_\lambda$. By Example \ref{ex: derivedVB}, $\mc V_{\overline\rho}|_{\T^* ( V_\lambda/(H_\lambda/Z(\wh G)^{\circ}))}$ is the Koszul dual of $z_*\mc V_{\overline\rho}$. Let $\widetilde{\mc V}_{\overline\rho}$ be $\mc V_{\overline\rho}$ with the trivial filtration. Then $\ind\widetilde{\mc V}_{\overline\rho}\in F\mc D_c(V_\lambda/(H_\lambda/Z(\widehat G)^{\circ}))$ has the associated graded 
\[
\gr(\ind\widetilde{\mc V}_{\overline\rho})\simeq \mc V_{\overline\rho}|_{\T^* (V_\lambda/(H_\lambda/Z(\widehat G)^{\circ}))},
\]
and its $t$-exact pullback (defined to be $!$-pullback with shift $[\dim Z(\wh G)]$) along $\V_\lambda\to V_\lambda/(H_\lambda/Z(\wh G)^{\circ})$ is a filtered $\mc D$-module $\widetilde{\mc G}_\rho\in F\mc D_c(\V_\lambda)$ lifting $\mc G_\rho$.

\begin{proposition}{\label{prop: CVBimpliesCLS}}
We work in the graded setting. Let $\mc V\neq 0$ be a vector bundle on $V_\lambda / (H_\lambda/ Z(\wh G)^{\circ})$ whose pull-push as above to $\T^*\V_\lambda$ is both left-cuspidal and right-cuspidal. Then, by Proposition \ref{prop: cuspreductiontofinite}, $\V_{\lambda}$ contains an essentially discrete $L$-parameter $\phi=(\lambda, N)$. Write $\widetilde{\mc V}$ for the filtered lift with the trivial filtration. Then, the $t$-exact pullback of the underlying $\mc D$-module of $\ind \widetilde{\mc V}$ is cuspidal on $\V_\lambda$. In particular, $\mc V$ restricts to a cuspidal representation $\rho$ of $S_\phi$ trivial on $Z(\wh G)^{\circ}$, and the $t$-exact pullback of the underlying $\mc D$-module of $\ind \widetilde{\mc V}$ is $\RH^{-1}(\mc{IC}(\mc L_\rho))$.
\end{proposition}

\begin{proof}[Proof of Proposition \ref{prop: CVBimpliesCLS} and Theorem \ref{thm: ClassificationCuspidalCotangent} (4)]
We prove both statements by induction on the semisimple rank of $G$, and for the fixed rank we prove Proposition \ref{prop: CVBimpliesCLS} first, then Theorem \ref{thm: ClassificationCuspidalCotangent} (4) follows by \S \ref{s: CVBs}. 
Note also that the pullback of $\mc V$ to $\T^* (V_\lambda / (H_\lambda/ Z(\wh G)^{\circ}))$ is cuspidal, and hence we may assume that $G$ is semisimple by replacing $\wh G$ by $\wh G/ Z(\wh G)^\circ$ to prove Proposition \ref{prop: CVBimpliesCLS}. 

Take a proper semistandard Levi subgroup $\wh M\subset \wh G$ with a lift $\lambda_M$, the semisimplification of an essentially discrete $\phi_M$. First observe that the inductive hypothesis for $M$ and a cuspidal representation $\rho_{M}$ trivial on $Z(\wh M)^\circ$ implies that the underlying $\mc D$-module of $\widetilde{\mc G}_{\rho_M}$ is cuspidal, and its restriction to $BS_{\phi_M}$ is $\rho_M$. Therefore, its underlying $\mc D$-module is $\RH^{-1}(\mc{IC}(\mc L_{\rho_M}))$ by Theorem \ref{thm: ClassificationCuspidalPerverseSheaves}. 

By Lemma \ref{lem: filteredcusplem} and Theorem \ref{thm: ClassificationCuspidalPerverseSheaves}, it suffices to prove
\[
\Hom (\ind \widetilde{\mc V}, \Eis_{P,\lambda_P}^{F\mc D, !} \widetilde{\mc G}_{\rho_M})\simeq 0
\]
as an object of $\QCoh (\bb A^1 /\Gm)$ for any semistandard proper Levi/parabolic $\wh M\subset \wh P$, a lift $\lambda_P$ and cuspidal $\rho_M$ that is trivial on $Z(\wh M)^{\circ}$.

By Lemma \ref{lem: Dmodulehomequality}, it suffices to show that the derived fiber at $B\Gm$ vanishes. By passing to the associated graded, this is equivalent to
\[
\Hom (\mc V|_{\T^* \V_\lambda}, \Eis_{P,\lambda_P}^{\T^{\ast}} \mc G_{\rho_M}) \simeq 0. 
\]
Since $\rho_M$ is trivial on $Z(\wh M)^{\circ}$, the (left-)cuspidality of $\mc V$ implies the vanishing. 

The final claim follows since the restriction to $\U_{C^{\circ}}\simeq BS_\phi$ is compatible with all the operations in the following diagram (Remark \ref{rem: Dmodulefunctoriality})
\begin{equation*}
    \begin{tikzcd}
        \mc{D}_c(\U_{C^{\circ}})  & F\mc{D}_c(\U_{C^{\circ}}) \arrow[l, swap, "{\un}"] \arrow[r, "{\gr}"]   & \Coh(\T^{\ast}\U_{C^{\circ}}/\Gm)  \\ 
           \mc{D}_c(\V_\lambda) \arrow[u, "j^{\ast}"] & F\mc{D}_c(\V_\lambda) \arrow[l, swap, "{\un}"] \arrow[r, "{\gr}"] \arrow[u, "j^{\ast}"] & \Coh(\T^{\ast}\V_\lambda/\Gm), \arrow[u, "j^{\ast}"]
    \end{tikzcd}
\end{equation*}
and the composition
\[
\un \colon F \mc D(BS_\phi)\to \mc D(BS_\phi)\simeq \QCoh(BS_\phi)
\]
sends the sheaf $\ind \widetilde{\mc V}|_{BS_\phi}$ to $\mc V|_{BS_\phi}$, and the Riemann--Hilbert equivalence is essentially the identity functor on the classifying stack of a finite group. 
\end{proof}

\subsubsection{Comparison of Springer sheaves}
By using Proposition \ref{prop: CVBimpliesCLS}, we can compare Springer sheaves as was already used in its inductive proof. 

Now let $\wh M\subset \wh P\subset  \wh G$ be semistandard Levi and parabolic subgroups, let $\lambda_M$ be the semisimplification of an essentially discrete $L$-parameter $\phi_M$ of $M$, and let $\rho_M$ be a cuspidal representation of $S_{\phi_M}$ trivial on $Z(\wh M)^{\circ}$. 

\begin{lemma}\label{lem: EisIsLimit}
Assume we are in the graded setting. We have $\Eis^{\T^{\ast}}_P \mc{G}_{\rho_M} \in \mrL(\T^*\V_\lambda)$. In particular, it is left-basic and right-basic by Corollary \ref{cor: LimsubsetBas} and the twisted Grothendieck--Serre duality.\footnote{This part can be proven also using Proposition \ref{prop: EisBasisBas}.}
\end{lemma}

\begin{proof}
This follows from Remark \ref{rem: CuspidalSubsetLimit} and Proposition \ref{prop: Limcatfunctorialities}. 
\end{proof}

We again work in the graded setting and, as in the previous subsection, we have the filtered lift $\widetilde{\mc G}_{\rho_M}$. 

\begin{proposition}{\label{prop: filtered springer}}
The filtered $\mc D$-module $\Eis_{P}^{F\mc D, !} \widetilde{\mc G}_{\rho_M}$ is a lift of the graded generalized coherent Springer sheaf $\Eis_{P}^{\T^*} \mc G_{\rho_M}$ whose underlying $\mc D$-module is sent to the generalized constructible Springer sheaf $\Eis_{P}^{\mathrm{cons}, !} \mc{IC}(\mc L_{\rho_M})$ under the Riemann--Hilbert correspondence.     
\end{proposition}

\begin{proof}
This follows from the preceding discussion together with the compatibility of various Eisenstein series. 
\end{proof}

Let us record the following consequence. We will further study generalized coherent Springer sheaves in a sequel. 

\begin{corollary}{\label{cor: coherentspringersummand}}
Suppose $S_\phi$ is finite, i.e., $G$ is semisimple and $\phi$ is essentially discrete. 
If an irreducible representation $\rho$ of $S_{\phi}$ is a summand, up to a shift, of the restriction of some $\Eis_{P}^{\mathrm{cons}, !} \mc{IC}(\mc L_{\rho_M})$ to $BS_{\phi} \simeq \T^{\ast} BS_{\phi}$, then $\rho$ is also a summand, up to a shift, of the restriction of $\Eis_{P}^{\T^*} \mc G_{\rho_M}$. Moreover, every irreducible $\rho$ is a summand up to a shift of the restriction of some $\Eis_{P}^{\T^*} \mc G_{\rho_M}$ including $P=G$. 
\end{corollary}

\begin{proof}
By Proposition \ref{prop: filtered springer}, the restriction of $\Eis_{P}^{F\mc D, !} \widetilde{\mc G}_{\rho_M}$ to $BS_\phi$ is a complex of representations of $S_\phi$ with a finite filtration. Its underlying complex is the restriction of $\Eis_{P}^{\mathrm{cons}, !} \mc{IC}(\mc L_{\rho_M})$, and its associated graded is the restriction of $\Eis^{\T^*}_P \mc G_{\rho_M}$. The first part follows since the category of representations of $S_\phi$ is semisimple. For the second, if $\rho$ is cuspidal take $P=G$. Otherwise, by the proof of Proposition \ref{prop: OrthogonalCuspidalConstructibleSheaves}, $\mc{IC}(\mc L_\rho)$ is a summand up to shift of some $\Eis_{P}^{\mathrm{cons}, !} \mc{IC}(\mc L_{\rho_M})$ with $P$ proper, and we restrict to $BS_\phi$ and apply the first part. 
\end{proof}

\subsection{Generation of basic sheaves}{\label{ss: basgen}}
In the graded setting, Theorem \ref{thm: BasicCotangentMain} (1) and (2) have been proven, and we shall prove (3). In fact, we prove the following stronger form and its proof is largely inspired by \cite[\S6.3]{PTDolbeault}. It partially realizes their expectation \cite[footnote 8, Conjecture 6.4]{PTDolbeault} as a variant of \cite[Theorem 6.13]{PTDolbeault}. 

\begin{theorem}\label{thm: cotanbasgenbycuspeis}
The $\infty$-category $\Bas^\Coh(\T^*\V)$ is generated under finite colimits and retracts by cuspidal sheaves $\mc G_\rho$ and 
$\Eis_P^{\T^*} \Bas^\Coh(\T^*\V_M)^{\kappa(b_M)}$ with proper $P$ and $\kappa (b_M)$ being trivial on $(Z(\wh M) \cap \wh{G}_\der)^\circ$ (in particular, the image of $b_M$ in $B(G)$ is basic). 

For any $c \in C\in \mc C$, with $j^C\colon \T^*\U_C\to \T^* \V$,
\[
j^{C}_{\ush}\Eis_{P_c}^{\T^*}\mrL(Z) \quad (\text{resp.~}j^{C}_{|\ast|} \Eis_{P_c}^{\T^*}\mrL(Z))
\]
lies in $\mrL(\T^*\V)$, and these generate $\Bas^\Coh(\T^*\V)$ in each case if $G$ is semisimple.
\end{theorem}

\begin{remark}{\label{rem: limbasnonss}}
\begin{enumerate}
    \item It also follows by arguing inductively that $\Bas^\Coh(\T^*\V)$ is generated by Eisenstein series of cuspidal sheaves $\mc G_\rho$ with trivial $(Z(\wh M)\cap \wh{G}_\der)^\circ$-weights. By Lemma \ref{lem: EisIsLimit} such sheaves are in $\mrL (\T^*\V)$ if their $Z(\wh G)^\circ$-weights are trivial, and it follows that the subcategory of $\Bas^\Coh(\T^*\V)$ with trivial $Z(\wh G)^\circ$-weights agrees with $\mrL (\T^*\V)$. 
    
    It is also related to \cite[Conjecture 1.12, Remark 1.13]{PTDolbeault} since these cuspidal sheaves $\mc G_\rho$ lift to filtered $\mc D$-modules with cuspidal underlying $\mc D$-modules, and their Eisenstein series generate $\mc D$-modules. 
    \item If one is interested only in the first part of the theorem, whose proof uses the second part, the second part can be replaced by the same claim but with $\mrL(\T^*\V)$ replaced by $\Bas^\Coh(\T^*\V)$; see also the footnote below.
\end{enumerate}
\end{remark}

\begin{proof}
We proceed by induction on the semisimple rank. The case of tori holds trivially. In fact, the statements reduce to the case of semisimple groups, and we do assume $G$ is semisimple from now on. 

Next we will use the $\Theta$-stratification as in \S \ref{sss:ThetaStratificationCotangent}. We shall also use the notation in \S \ref{ss: windowcat}, \ref{ss: basicvswindow}, and \ref{ss: limcat}. 
We first prove

\begin{lemma}{\label{lem: limcatlem}}
Suppose $\mc F\in \mrL(\T^*\V)\subset \Bas^\Coh(\T^*\V)$. 
We have 
\[
j_{C\ush}(\mc F|_{\T^* \U_C\setminus S}), \quad j_{C|*|}(\mc F|_{\T^* \U_C\setminus S})\in \mrL (\T^*\U_C)
\]
for every $C\in \mc C$.\footnote{For any $\mc F\in \Bas^\Coh(\T^*\V)$, the argument can be modified to show that the extensions belong to the subcategory generated under finite colimits by the restriction of $\Bas^\Coh(\T^*\V)$.} 
\end{lemma}

\begin{proof}
It suffices to prove the statement for $j_{\ush}$ by the duality. The claim is vacuous when $C= C^\circ$, so we assume $C\neq C^\circ$. Since $\mc F|_{\T^* \U_C}$ is in the limit category by Proposition \ref{prop: Limcatfunctorialities} (the map is an open immersion), it suffices to show that the cone of the counit
\[
j_{C\ush}(\mc F|_{\T^* \U_C\setminus S})\to
\mc F|_{\T^* \U_C}
\] 
is in the limit category. It has the form $\Eis_{P_c}^{\T^*} \mc F_{M_c}$ with $\mc F_{M_c}\in \Coh(Z)^0$ (Theorem \ref{thm: HL3.3.1} parts (1), (3)). For any $\mc G\in \Coh(Z)^0$, 
\begin{equation}{\label{eqn: HomEiscomparison}}
    \Hom(\mc F|_{\T^* \U_C}, \Eis_{P_c}^{\T^*}\mc G)\simeq 
\Hom(\Eis_{P_c}^{\T^*} \mc F_{M_c}, \Eis_{P_c}^{\T^*}\mc G) \simeq
\Hom( \mc F_{M_c}, \mc G), 
\end{equation}
where the first and second equivalences follow from Theorem \ref{thm: HL3.3.1} (3), (2).
We claim that $\mc F_{M_c}\in \Bas^\Coh (Z)^0$. To verify this, we consider $\mc{G}$ of the form $\Eis^{\T^{\ast}}_{P'} \mc{G}'$ for $P' \subset M_c$ a parabolic subgroup with Levi quotient $M'$ and $\mc{G}'$ with strictly $P'$-antidominant central weight. Applying \eqref{eqn: HomEiscomparison}, we have for a parabolic $P'_c$ of $\wh G$,
\begin{align*}
    \Hom( \mc F_{M_c}, \Eis^{\T^{\ast}}_{P'} \mc{G}') &\simeq \Hom(\mc F|_{\T^* \U_C}, \Eis_{P_c}^{\T^*}\Eis^{\T^{\ast}}_{P'} \mc{G}')\simeq \Hom(\mc F|_{\T^* \U_C}, \Eis_{P'_c}^{\T^*} \mc{G}')\\ &\simeq \Hom(\mc F, j^{C }_{\ast}\Eis_{P'_c}^{\T^*} \mc{G}') \simeq \Hom(\mc F, \Eis_{P'_c}^{\T^*} j^{=0}_{\ast}\mc{G}') \simeq 0, 
\end{align*}
where the vanishing follows since $\mc{F}$ is basic.

By Proposition \ref{prop: Limcatfunctorialities}, it remains to show that $\mc F_{M_c}\in \mrL(Z)$. Since $Z$ is the cotangent stack of a classifying stack, it follows by Corollary \ref{cor: limcatclassifyingstack} that it remains to show $\mc F_{M_c}$ has trivial $Z(H^{\circ}_{M, c})^{\circ}$-weights. We first show that the $Z(M_c)^\circ$-weights of $\mc F_{M_c}$ are trivial, and we show \eqref{eqn: HomEiscomparison} vanishes for every $\mc G\in \Bas^\Coh(Z)^0$ with a single nontrivial $Z(M_c)^\circ$-weight $\chi$. We can find an intermediate Levi/parabolic $\wh M\supset M_c$, $\wh P\supset P_c$ such that $\chi$ restricted to $Z(\wh M)^\circ$ is strictly $P$-antidominant. Indeed, we can write $\nu_c$ as a non-negative linear combination of fundamental coweights and the condition $\langle \nu_c, \chi \rangle =0$ implies either that $\chi$ pairs trivially with each fundamental coweight outside of $M_c$, in which case $\chi$ is trivial contrary to assumption (since we are assuming $G$ is semisimple), or $\chi$ pairs negatively with at least one fundamental coweight, in which case we take $
\wh{P}$ to be the maximal parabolic associated to that coweight. Then, by the left-basicness of $\mc F$ and Lemma \ref{lem: cotangentdomswitch}, the Hom vanishes. 

Next, we want to show the $Z(H_{M,c}^\circ)^\circ$-weights of $\mc F_{M_c}$ are trivial. Write $M'_c\coloneq Z_{\wh G}(Z(H_{M,c}^\circ)^\circ)$. Note that for any semistandard parabolic $M_c \supset \wh P'\supset M'_c$ with Levi subgroup $\wh{M}'\supset M'_c$, the subgroup $H_{M,c}\cap \wh{M}'\subset H_{M,c}$ is of finite index, and the sheaf $\mc F_{M_c}$ appears as a summand of $\Eis_{P'}^{\T^*}\mc F_{\wh{M}'}$ for some $\mc F_{\wh{M}'}$. For a nontrivial $Z(H_{M,c}^\circ)^\circ$-weight $\chi$ trivial on $Z(M_c)^\circ$, we can find $M_c\supset \wh{P}'\supset \wh{M}'\supset M'_c$ such that $\chi$ restricted to $Z(\wh{M}')^\circ$ is strictly $\wh{P}'$-antidominant. Then, using the left-basicness again, we conclude that the $\chi$-summand of $\mc F_{M_c}$ vanishes. This completes the proof. 
\end{proof}

As in the above proof, we have the triangle
\begin{equation}{\label{eqn: limcatfibertriangle}}
    j_{C\ush}(\mc F|_{\T^* \U_C\setminus S})\to
\mc F|_{\T^* \U_C} \to
\Eis_{P_c}^{\T^*} \mc F_{M_c},
\end{equation}
and $\mc F_{M_c}\in \mrL (Z)$ (i.e., by Corollary \ref{cor: limcatclassifyingstack} it is generated by representations of $\pi_0 (H_{M,c})$ on the zero section $BH_{M,c}\hookrightarrow Z$). By the inductive hypothesis (noting that $c$ gives the open orbit in the $\Theta$-stratification of $\V_{M_c}$ and $H_{M,c}$ is reductive hence the Eisenstein series from the center of the $\Theta$-stratum is the identity) $j^\circ_{\ush}\mc F_{M_c}$ is in $\mrL(\T^*\V^{\nu_c=0})$. Then by Proposition \ref{prop: Limcatfunctorialities}, 
\[
\Eis_{P_c}^{\T^*} j_{\ush}^\circ\mc F_{M_c} \in \mrL(\T^*\V)\subset \Bas^\Coh(\T^*\V),
\]
and it restricts, by Lemma \ref{lem: eisfiberprod}, to $\Eis_{P_c}^{\T^*} \mc F_{M_c}$ on $\T^{\ast}\U_C$. 

Now, if $C'>C$ is the element next to $C$, then we apply $j_{C'\ush}$ to the triangle \eqref{eqn: limcatfibertriangle} and note the middle and right terms are in $\mrL(\T^{\ast}\U_{C'})$ by Lemma \ref{lem: limcatlem} (for the right term we use the previous paragraph to show it satisfies the conditions). So we conclude that
\[
j_{C'\ush}j_{C\ush}(\mc F|_{\T^* \U_C\setminus S})\in \mrL(\T^{\ast}\U_{C'}), 
\]
and this holds for every $\mc F$ and $C'>C$. Then arguing similarly (again the right term is the extension by $j_{C'' \ush}j_{C' \ush}$ of the restriction of an element of $\mrL(\T^{\ast}\V)$ hence in the limit category by what we've already proven) we can deduce that
\[
j_{C''\ush}j_{C'\ush}j_{C\ush}(\mc F|_{\T^* \U_C\setminus S})\in \mrL(\T^{\ast}\U_{C''}), 
\]
for $C''$ next to $C'$. Repeating the argument, we obtain
\[
j^{C}_{\ush}j_{C\ush}(\mc F|_{\T^* \U_C\setminus S})\in \mrL(\T^*\V)\subset\Bas^\Coh(\T^*\V).  
\]
\begin{lemma}{\label{lem: extensionlimcat}}
Under the inductive hypothesis as above, 
\[
j^{C}_{\ush}\Eis_{P_c}^{\T^* }\mrL(Z) \subset \mrL(\T^*\V). 
\]
\end{lemma}

\begin{proof}
Suppose $C\neq C^\circ$. Using the inductive hypothesis on Levi subgroups, we have, for any $\mc G\in \mrL(Z)$
\[
\Eis_{P_c}^{\T^* }j^\circ_{\ush}\mc G \in \mrL(\T^*\V),
\]
that restricts to $\Eis_{P_c}^{\T^*}\mc G$ by Lemma \ref{lem: eisfiberprod}, and hence we deduce using the paragraph before this lemma that
\[
j^{C}_{\ush}\Eis_{P_c}^{\T^* }\mrL(Z) \subset \mrL(\T^*\V). 
\]

If $C=C^\circ$, a similar argument using \S \ref{sss: cuspreductiontoreductive} works if $H_{c^\circ}$ is not reductive. If $H_{c^\circ}$ is reductive but not finite, take $\nu\neq 0$ defining a subgroup $H_{c^\circ} \cap P_{\nu} \subset H_{c^{\circ}}$ such that $H_{c^\circ}^\circ \cap {P_{\nu}} \subset H_{c^\circ}^\circ$ is a Borel subgroup. By Lemma \ref{lem: zerosectiondirectsummand}, any representation of $\pi_0(H_{c^\circ})$ as a sheaf on the zero section of $\T^* \U_{C^\circ}$ appears as a summand of an Eisenstein series $\Eis^{\T^{\ast}}_{H_{c^\circ} \cap P_{\nu}}z_{H_{c^\circ} \cap M_\nu *}\rho$ for some representation $\rho$ of $\pi_0(H_{c^\circ} \cap M_\nu)$. Therefore, as in the discussion before this lemma, it appears as a summand of the restriction of $\Eis^{\T^{\ast}}_{P_{\nu}} j^\circ_{\ush} z_{H_{M_\nu}*}\rho$, and we can argue analogously. 
Finally, if $H_{c^\circ}$ is finite, every irreducible representation appears in the restriction of a cuspidal coherent sheaf or a basic generalized coherent Springer sheaf by Corollary \ref{cor: coherentspringersummand}. They are in $\mrL(\T^*\V)$, and the same argument works. 
\end{proof}

Now we proceed by descending induction on the fixed ordering of $\mc C$. For $C\in \mc C$ and a subset $A \subset \lvert \T^* \V\rvert$, we say $A> C$ (resp.~$\geq C$) if every point in the projection of $A$ to $\V$ lies in an orbit $> C$ (resp.~$\geq C$). 

We first prove $\Bas^\Coh (\T^*\V)$ is generated under finite colimits by $j^{C}_{\ush}\Eis_{P_c}^{\T^* }\mrL(Z)$. Suppose the support of $\mc F$ is $\geq C$ for $C\in \mc C$. If $C\neq C^\circ$, we can write
\[
j^{C*}\mc F =\mc F |_{\T^*\U_C} = \Eis_{P_c}^{\T^*} \mc F_{M_c}, \qquad \mc F_{M_c}\in  \Coh(Z)^0, 
\]
and as in the proof of Lemma \ref{lem: limcatlem}, $\mc F_{M_c}$ is in $\mrL (Z)$.
Therefore, by the counit
\[
j^C_{\ush}j^{C*}\mc F \to \mc F
\]
and the inductive hypothesis on $C$, it follows that $\mc F$ is in the desired subcategory. If $C=C^\circ$, then we similarly have
\[
j^{\circ*}\mc F \in \Eis_{P_{c^\circ}}^{\T^*} \mrL(Z), 
\]
so the same argument works. 

Now, to show that $\Bas^\Coh(\T^*\V)$ is generated under finite colimits and retracts by $\Cusp^\Coh(\T^*\V)$ (or cuspidal vector bundles) and 
$\Eis_P^{\T^*} \mrL(\T^*\V_M)$, it suffices to show that $j^C_{\ush}\Eis_{P_c}^{\T^* }\mrL(Z)$ is contained in this category (again by descending induction on $C$). 
Let $\Eis_{P_c}^{\T^*}\mc G \in \Eis_{P_c}^{\T^* }\mrL(Z)$. Then, $\Eis_{P_c}^{\T^*}j_{\ush}^\circ \mc G\in \mrL(\T^*\V)$ restricts to $\Eis_{P_c}^{\T^*}\mc G$, and the cone of the counit
\[
j^C_{\ush}\Eis_{P_c}^{\T^*}\mc G\to 
\Eis_{P_c}^{\T^*}j_{\ush}^\circ \mc G 
\]
is supported on $> C$. Therefore, if $P_c$ is proper, the inductive hypothesis on $C$ implies the claim. If $P_c=\wh G$, then $c=c^\circ$ and $H_{c^\circ}$ is reductive. As in the proof of Lemma \ref{lem: extensionlimcat}, $\Eis_{P_c}^{\T^*}\mc G=\mc G$ appears as a summand of the restriction of a cuspidal sheaf or of an object of $\Eis_P^{\T^{\ast}} \mrL(\T^{\ast}\V_M)$ with $P$ proper, and we are done by the same argument. 
\end{proof}

Part (3) of Theorem \ref{thm: BasicCotangentMain} is proved in the graded setting.

\begin{proof}[Proof of Theorem \ref{thm: ClassificationCuspidalCotangent} except (5)]
(1) reduces to the graded setting and the case $G$ semisimple. But we have proven Theorem \ref{thm: BasicCotangentMain} and the proof also shows we only need cuspidal vector bundles in Part (3). Therefore, inside $\Bas^\Coh(\T^*\V)$, left- and right-cuspidal sheaves are left- and right-orthogonal to $\Bas^\Coh(\T^*\V) \cap \Eis_-$, $\Bas^\Coh(\T^*\V) \cap \Eis_+$ respectively, and they are both generated by cuspidal vector bundles. The first half of Theorem \ref{thm: ClassificationCuspidalCotangent} (1) is now clear, and the second half follows by passing to the ind-completion. 

(1)' follows from (1): recall that
\[
\IndCoh(\T^*\V_\lambda)=\langle \RCusp (\T^*\V_\lambda), \Ind \Eis_+^{\T^*} \rangle,
\]
where $\Eis_+^{\T^*}$ is the category generated by dominant Eisenstein series. By passing to compact objects using (1), we obtain
\[
\Coh(\T^*\V_\lambda)=\langle \Cusp^\Coh (\T^*\V_\lambda), \Eis_+^{\T^*} \rangle,
\]
and the antidominant case follows by Grothendieck--Serre duality together with (1). 

(2) follows as Proposition \ref{prop: cuspextopen} is now applicable to all cuspidal coherent sheaves. 

(3) reduces to the graded setting, and then we have proven that cuspidal sheaves are generated under colimits by $\mc G_\rho$, so (3) is clear. 
\end{proof}

\section{Variants of spectral geometric lemma}{\label{s: spectralgeometriclemma}}

We shall prove orthogonality relations between Eisenstein series on cotangent stacks and some related results by analyzing a certain filtration on the Hom space, reminiscent of the geometric lemma for representations of $p$-adic groups. The idea of such a ``spectral'' geometric lemma appears in \cite{GLCIII}, and plays a key role in \cite{HansenMann} as well with many more details. While we work with cotangent stacks, the argument below using baric structures itself should work directly for $\Par_G$ as well.

Let $\lambda\colon W_F\to \wh G$ be a semisimple $L$-parameter such that $\wh T_H\coloneq (\wh T\cap H_\lambda)^\circ$ is a maximal torus of $H_\lambda$ contained in $\wh T^\Gamma$, chosen as in Lemma \ref{lem: TJcontrol}. We will simply write $\V=\V_\lambda$. Let $\nu_1, \nu_2\colon \Gm\to \wh T_H$ be cocharacters. We can define the diagrams
\[
\T^* \V \xleftarrow{p_1} \T^* \V^{\nu_1\geq 0} \xrightarrow{q_1} \T^*\V^{\nu_1=0}, \quad
\T^* \V \xleftarrow{p_2} \T^* \V^{\nu_2\geq 0} \xrightarrow{q_2} \T^*\V^{\nu_2=0},
\]
as discussed in Example \ref{ex: TVnugeq0} in the graded setting. Note that the constructions there carry through even in this more general setting.  We write $\wh P_1 =\wh M_1 \wh U_1$, $\wh P_2 = \wh M_2 \wh U_2$ for the semistandard parabolic and Levi subgroups and their unipotent radicals associated with $\nu_1$, $\nu_2$ respectively. 

Write
\begin{gather*}
W_{H_\lambda}= N_{H_\lambda}(\wh T_H)/ \wh T_H, \quad
W_{\wh M_1 \cap H_\lambda}= N_{\wh M_1 \cap H_\lambda}(\wh T_H)/ \wh T_H, \\
W_{\wh M_2 \cap H_\lambda}= N_{\wh M_2 \cap H_\lambda}(\wh T_H)/ \wh T_H.
\end{gather*}

The results in \S \ref{sss: derivedtheta strat} and Proposition \ref{prop: basicweightcharacterization} continue to hold in this setting with $W_H$ replaced by $W_{H_\lambda}$ by identical arguments, and we shall use them without further notice.   

\subsection{Orthogonality relations on the cotangent stack}

\begin{lemma}
\[
(\wh P_1 \cap H_\lambda)\backslash H_\lambda / (\wh P_2 \cap H_\lambda) \cong 
W_{\wh M_1 \cap H_\lambda} \backslash W_{H_\lambda} /W_{\wh M_2 \cap H_\lambda}. 
\]
\end{lemma}
\begin{proof}
First note that $H_{\lambda} \cap \wh{P}_i$ is a parabolic subgroup of $H_{\lambda}$ in the sense that it can be defined via the dynamic method. 

When $H_{\lambda}$ is connected reductive and for standard parabolics, the statement of the lemma is \cite[IV.5 Remark 2]{Bou1}. More generally, we have that $H_{\lambda}$ acts on the left of $H_{\lambda} / (\wh{P}_2 \cap H_{\lambda})$ and we claim the $\wh{T}_H$-fixed points are $W_{H_{\lambda}} / W_{\wh{M}_2 \cap H_{\lambda}}$. Indeed, $h(\wh{P}_2 \cap H_{\lambda})$ is a $\wh{T}_H$-fixed point if and only if $h^{-1}$ conjugates $\wh{T}_H$ into $\wh{P}_2 \cap H_{\lambda}$. But then, up to a different choice of $h$, this is equivalent to $h \in N_{H_{\lambda}}(\wh{T}_H)$. Moreover, any two such $h$ in the same coset differ by an element of $N_{\wh{P}_2 \cap H_{\lambda}}(\wh{T}_H)$. 
    
To finish the proof, we need to show that each $\wh{P}_1 \cap H_{\lambda}$ orbit in $H_{\lambda} / (\wh{P}_2 \cap H_{\lambda})$ has a unique $\wh{T}_H$-fixed point (up to action by $N_{\wh{P}_1 \cap H_{\lambda}}(\wh{T}_H)$). The existence follows from the Bruhat decomposition, once we note that the connected components of $H_{\lambda} / (\wh{P}_2 \cap H_{\lambda})$ are of the form $H^{\circ}_{\lambda} / x(\wh{P}_2 \cap H^{\circ}_{\lambda})x^{-1}$ for some $x \in H_{\lambda}$, and that $x(\wh{P}_2 \cap H^{\circ}_{\lambda})x^{-1}$ is a parabolic subgroup of $H^{\circ}_{\lambda}$. For the uniqueness, two $\wh{T}_H$-fixed points $h(\wh{P}_2 \cap H_{\lambda}), h'(\wh{P}_2 \cap H_{\lambda})$ satisfying $h' = ph$ for some $p \in \wh{P}_1 \cap H_{\lambda}$ implies that $\wh{T}_H, p^{-1}
\wh{T}_Hp$ both fix $h(\wh{P}_2 \cap H_{\lambda})$. So both tori lie in the stabilizer of $h(\wh{P}_2 \cap H_{\lambda})$ in $\wh{P}_1 \cap H_{\lambda}$. Hence, modifying $p$ by some element of this stabilizer, we can assume $p \in N_{\wh{P}_1 \cap H_{\lambda}}(\wh{T}_{H})$ as desired.
\end{proof}

We shall use the above identification freely, and choose a representative in $N_{H_\lambda}(\wh T_H)$ implicitly. 

\begin{proposition}{\label{prop: generaleisbasisbas}}
Let $\mc F_1 \in \Coh(\T^{\ast}\V^{\nu_1=0})^{\kappa_{M_1}(b_1)}$ for $b_1 \in B(M_1)_{\bas}$ and $\mc F_2 \in \Coh(\T^* \V^{\nu_2=0})^{\kappa_{M_2}(b_{2})}$ for $b_{2} \in B(M_2)_{\bas}$. Assume either that $\kappa_G(b_1)\neq \kappa_G(b_2)$, or that $\kappa_G(b_1)=\kappa_G(b_2)$, $\mc F_2$ is right-basic, and
\[
\langle \nu_1,\kappa_{M_1}(b_{1})\rangle=\langle \nu_1,\nu_{M_1}(b_1)\rangle
> \langle \nu_1,w(\nu_{M_2}(b_2))\rangle
\]
for every $w\in W_{H_\lambda}$. In both cases, 
\[
\Hom (\Eis^{\T^{\ast}}_{P_1} \mc F_1, \Eis^{\T^{\ast}}_{P_2} \mc F_2)\simeq 0. 
\]
\end{proposition}
\begin{proof}
In the case where $\kappa_G(b_1) \neq \kappa_G(b_2)$, we have $\Hom (\Eis^{\T^{\ast}}_{P_1} \mc F_1, \Eis^{\T^{\ast}}_{P_2} \mc F_2)\simeq 0$ since the Eisenstein series have different central $Z(\wh{G})^{\Gamma}$-weights. 

Suppose $\kappa_G(b_1) = \kappa_G(b_2)$. For $w\in (H_\lambda\cap \wh P_1)\backslash H_\lambda / (H_\lambda \cap \wh P_2)$, represented in $N_{H_\lambda}(\wh T_H)$, define the stack
\begin{equation*}
    \T^{\ast}\V^{\nu_1 \geq 0, w(\nu_2) \geq 0} \coloneq [\{0\} \times_{ (\Lie H_{\lambda})^{\ast, \nu_1 \geq 0, w(\nu_2) \geq 0}} \T^{\ast}V^{\nu_1 \geq 0, w(\nu_2) \geq 0}] / (H_{\lambda} \cap \wh P_1 \cap w(\wh P_2)).
\end{equation*}
Then consider the following commutative diagram where the $(p_1, p_2)$ and $(q_2, w^{-1} \circ q')$ squares are cartesian 
\[
\begin{tikzcd}[column sep=1em]
        && \T^* \V^{\nu_1 \geq 0, w(\nu_2)\geq 0} \arrow[d, "i_w"] \arrow[dddll,bend right,"g_w",swap] \arrow[ddrrr,bend left,"f_w"] \arrow[drr, "h_w"] \\
        && \T^\ast \V^{\nu_1\geq 0}\times_{\T^* \V}\T^\ast \V^{\nu_2\geq 0} \arrow[dl, "p_2"] \arrow[dr, "p_1"]  && Y \arrow[dl] \arrow[dr, "q_2"]\\
        & \T^\ast \V^{\nu_1\geq 0} \arrow[dl,"q_1"] \arrow[dr,"p_1"] && \T^* \V^{\nu_2 \geq 0}  \arrow[dl,"p_2"] \arrow[dr,"q_2"] && (\T^* \V^{w(\nu_2)=0})^{\nu_1 \geq 0} \arrow[dl, "w^{-1}\circ q'"]\\
        \T^* \V^{\nu_1=0} && \T^* \V && \T^* \V^{\nu_2=0}
\end{tikzcd}
\]
As in \cite[Proposition 3.6.1]{HansenMann}, we have a topological stratification of the form
\begin{equation*}
    \lvert\T^\ast \V^{\nu_1\geq 0}\times_{\T^* \V}\T^\ast \V^{\nu_2\geq 0}\rvert = \coprod\limits_{w \in (H_{\lambda} \cap \wh P_1) \backslash H_{\lambda} / (H_{\lambda} \cap \wh P_2)} \lvert\T^* \V^{\nu_1 \geq 0, w(\nu_2)\geq 0}\rvert.
\end{equation*}
We need to show 
\begin{align*}
&\Hom_{\Coh}(p_{1*}q_1^* (\mc F_1\otimes \omega_{\V_{U_1}}), p_{2*}q_2^* (\mc F_2\otimes \omega_{\V_{U_2}})) \\
&\simeq
\Hom_{\IndCoh}(q_1^* (\mc F_1\otimes \omega_{\V_{U_1}}), p_1^! p_{2*}q_2^* (\mc F_2\otimes \omega_{\V_{U_2}})) \\
&\simeq
\Hom_{\IndCoh}(q_1^* (\mc F_1\otimes \omega_{\V_{U_1}}), p_{2*} p_1^! q_2^* (\mc F_2\otimes \omega_{\V_{U_2}})) \\
&\simeq
\Hom_{\IndCoh}(\mc F_1\otimes \omega_{\V_{U_1}},(q_1 \circ p_2)_* p_1^! q_2^* (\mc F_2\otimes \omega_{\V_{U_2}})) 
\end{align*}
vanishes. The functor $p_1^! q_2^*$ admits a filtration indexed by $w\in (H_{\lambda} \cap \wh P_1) \backslash H_{\lambda} / (H_{\lambda} \cap \wh P_2)$ by excision, cf.~the proof of \cite[Theorem 3.6.3]{HansenMann}. We denote the graded pieces of this filtration by $F_w$; they are given by
\begin{equation*}
    F_w \coloneq \wh{i}_{w \ast} \wh{i}^!_w p^!_1 q^{\ast}_2,
\end{equation*}
where $\wh{i}_w$ is the embedding from the formal completion. The functor $F_w$ then admits a further filtration indexed by $n\geq 0$ and corresponding to the excision sequence for the $n$-th infinitesimal thickening of $\T^* \V^{\nu_1 \geq 0, w(\nu_2)\geq 0}$ in $\T^\ast \V^{\nu_1\geq 0}\times_{\T^* \V}\T^\ast \V^{\nu_2\geq 0}$. The graded pieces of this filtration are of the form
\begin{equation*}
    F_{w,n} \simeq \cofib(i_{w, n \ast} i^!_{w,n} p^!_1 q^{\ast}_2 \rightarrow i_{w, n+1 \ast} i^!_{w, n+1} p^!_1 q^{\ast}_2) \simeq i_{w \ast}(\Sym^n \Nm_w \otimes  i^!_w p^!_1 q^{\ast}_2),
\end{equation*}
where $i_{w,n}$ is the inclusion of the $n$-th infinitesimal neighborhood and $\Nm_w$ is the normal bundle of $i_w$. 

We first show 
\begin{align*}
&\Hom_{\IndCoh}( \mc F_1\otimes \omega_{\V_{U_1}}, (q_1 \circ p_2)_* i_{w*} i_w^! p_1^! q_2^* (\mc F_2\otimes \omega_{\V_{U_2}})) \\
&\simeq
\Hom_{\IndCoh}( \mc F_1\otimes \omega_{\V_{U_1}}, (g_w)_* i_w^! p_1^! q_2^* (\mc F_2\otimes \omega_{\V_{U_2}})) \\
&\simeq
\Hom_{\QCoh}( \mc F_1\otimes \omega_{\V_{U_1}}, (g_w)_* i_w^! p_1^! q_2^* (\mc F_2\otimes \omega_{\V_{U_2}})) \\
&\simeq
\Hom_{\QCoh}(g_w^* (\mc F_1\otimes \omega_{\V_{U_1}}), i_w^! p_1^! q_2^* (\mc F_2\otimes \omega_{\V_{U_2}}))
\end{align*}
vanishes, where we are using $\IndCoh^+\simeq \QCoh^+$ on bounded below objects and the compatibility of $\Psi$ and $g_{w*}$, and also abusing notation. On $\T^* \V^{\nu_1 \geq 0, w(\nu_2)\geq 0}$, we can consider the baric structure induced by $\nu_1$, and it is compatible with the one for $\T^* \V^{\nu_1=0}$ via $g_w$, and the one for $(\T^* \V^{w(\nu_2)=0})^{\nu_1\geq 0}$ via $f_w$. It follows that 
\[
g_w^* (\mc F_1\otimes \omega_{\V_{U_1}})
\in \Aperf^*(\T^* \V^{\nu_1 \geq 0, w(\nu_2)\geq 0})^{\geq \langle \nu_1, \wt(\V_{U_1}) + \kappa_{M_1} (b_1) \rangle}.
\]
We will show that, under $\IndCoh^+\simeq \QCoh^+$, 
\[
i_w^! p_1^! q_2^* (\mc F_2\otimes \omega_{\V_{U_2}}) 
\in \QCoh(\T^* \V^{\nu_1 \geq 0, w(\nu_2)\geq 0})^{\leq \langle \nu_1, \wt(\V_{U_1}) + w(\nu_{M_2} (b_2)) \rangle}. 
\]
Since $\mc F_2$ is right-basic, we have by Proposition \ref{prop: basicweightcharacterization},
\[
(w^{-1}\circ q')^! \mc F_2 \in \QCoh((\T^* \V^{w(\nu_2)=0})^{\nu_1 \geq 0})^{\leq \langle \nu_1, \wt (\V_{U_1\cap w(M_2)}) + w(\nu_{M_2} (b_2))\rangle }.
\]
Then 
\begin{equation*}
    i_w^! p_1^! q_2^* \simeq h^!_wq^{\ast}_2(w^{-1} \circ q')^!,
\end{equation*}
so it remains to understand the relative dualizing complex of $h_w$, because both $q_2$ and $f_w=q_2\circ h_w$ are quasi-smooth. Let $\omega$ denote the line bundle whose correct shift is the relative dualizing complex of $h_w$; we can talk about its weight. It suffices to show
\[
\langle \nu_1, \wt (\omega) \rangle=\langle \nu_1, \wt (\V_{U_1})-\wt (\V_{w(U_2)}) - \wt (\V_{U_1\cap w(M_2)})  \rangle. 
\]
The contribution of $V$ to the right-hand side is
\[
-\langle \nu_1, \wt (V^{\nu_1 >0, w(\nu_2)<0}) - \wt(V^{\nu_1 <0, w(\nu_2)>0})\rangle,
\]
and the contribution of $H_{\lambda}$ to the right-hand side is
\[
\langle \nu_1,  \wt (\Lie H^{\nu_1 > 0, w(\nu_2)<0}_{\lambda})- \wt (\Lie H^{\nu_1 < 0, w(\nu_2) > 0}_{\lambda})\rangle
=2 \langle \nu_1, \wt (\Lie H^{\nu_1 > 0, w(\nu_2)<0}_{\lambda})\rangle.
\]

For the left-hand side, the contribution from $\T^* V$ can be computed as
\begin{align*}
\bigl(-\wt(\T^* V^{\nu_1 \geq 0, w(\nu_2)\geq 0}) & +  \wt (\T^* V^{w(\nu_2)\geq 0})\bigr) -  \bigl(-\wt(\T^* V^{w(\nu_2)=0, \nu_1\geq 0}) + \wt(\T^* V^{w(\nu_2)=0})\bigr) \\
&=    
\wt(\T^* V^{\nu_1 < 0, w(\nu_2)\geq 0})  -  \wt(\T^* V^{w(\nu_2)=0, \nu_1 < 0}) \\
&=\wt(\T^* V^{\nu_1 < 0, w(\nu_2) >0}) \\
&=\wt(V^{\nu_1 <0, w(\nu_2) > 0}) + \wt((V^\vee)^{\nu_1 <0, w(\nu_2) > 0}) \\
&=\wt(V^{\nu_1 <0, w(\nu_2) > 0}) - \wt(V^{\nu_1 > 0, w(\nu_2) < 0}),
\end{align*}
and the desired equality holds. To compute the contribution relevant to $H_{\lambda}$, we may assume $V=\{0\}$, and then the relevant spaces look like, for instance, 
\[
((\Lie H_{\lambda})^\vee)^{\nu_1 \geq 0, w(\nu_2)\geq 0}[-1] / H^{\nu_1 \geq 0, w(\nu_2)\geq 0}_{\lambda}, 
\]
and its dualizing weight is
\[
\wt (((\Lie H_{\lambda})^\vee)^{\nu_1 \geq 0, w(\nu_2)\geq 0}) 
+ \wt (\Lie H^{\nu_1\geq 0, w(\nu_2)\geq 0}_{\lambda})=2\wt (\Lie H^{\nu_1\geq 0, w(\nu_2)\geq 0}_{\lambda}). 
\]
Therefore, the contribution to $\omega$ is twice
\begin{align*}
\bigl(\wt (\Lie H^{\nu_1\geq 0, w(\nu_2)\geq 0}_{\lambda}) &- \wt (\Lie H^{\nu_1\geq 0, w(\nu_2)=0}_{\lambda})\bigr) - \bigl(\wt (\Lie H^{w(\nu_2)\geq 0}_{\lambda})-\wt(\Lie H^{w(\nu_2)=0}_{\lambda})\bigr) \\
&=
\wt (\Lie H^{\nu_1\geq 0, w(\nu_2)> 0}_{\lambda}) - \wt (\Lie H^{w(\nu_2)> 0}_{\lambda}) \\
&=-\wt (\Lie H^{\nu_1 <0, w(\nu_2)> 0}_{\lambda})=\wt (\Lie H^{\nu_1 > 0, w(\nu_2)< 0}_{\lambda})
\end{align*}
as desired. 

To treat other pieces of the filtration, it suffices to show that any $\nu_1$-weight of $\Sym^n \Nm_w$ is negative for $n\geq 1$. We shall study the normal bundle via the description using moment maps, and then it suffices to treat $(\Lie H_{\lambda})^\vee$, $*/H_{\lambda}$, $\T^*V$ separately. Their contributions to $\mathrm{Nm}_w$ are, up to shifts,
\[
((\Lie H_{\lambda})^\vee)^{\nu_1 < 0, w(\nu_2)<0}, \quad (\Lie H_{\lambda})^{\nu_1 < 0, w(\nu_2)<0}, \quad 
(\T^* V)^{\nu_1 < 0, w(\nu_2)<0},
\]
which indeed have negative $\nu_1$-weights, and the same is the case for symmetric powers. 
\end{proof}

\begin{proposition}{\label{prop: EisBasisBas}}
Assume $\nu_1$ factors through $\wh G_\der$. 
Let $\mc F_1 \in \Coh(\T^{\ast}\V^{\nu_1=0})^{\kappa_{M_1}(b_1)}$ for $b_1 \in B(M_1)_{\bas}$ and $\mc F_2 \in \Coh(\T^* \V^{\nu_2=0})^{\kappa_{M_2}(b_2)}$ for $b_2 \in B(M_2)_{\bas}$. Assume further that the image of $b_2$ in $B(G)$ is basic, that $\langle \nu_1,\kappa_{M_1}(b_1)\rangle > 0$, and that $\mc F_2$ is right-basic. Then we have 
\[
\Hom (\Eis^{\T^{\ast}}_{P_1} \mc F_1, \Eis^{\T^{\ast}}_{P_2} \mc F_2)=0. 
\]
Consequently, as in the proof of Proposition \ref{prop: basicweightcharacterization}, this implies that $\Eis^{\T^{\ast}}_{P_2} \mc F_2$ is right-basic, and the same holds for left-basic sheaves by the twisted Grothendieck--Serre duality. 
\end{proposition}

\begin{proof}
Since $\nu_1$ factors through $\wh G_\der$ and the image of $b_2$ is basic in $B(G)$, $\langle \nu_1,w(\nu_{M_2}(b_2))\rangle=0$ for any $w\in W_{H_\lambda}$. Hence, the result follows from Proposition \ref{prop: generaleisbasisbas}.
\end{proof}

\begin{corollary}{\label{cor: eisHomnonvanishingconditions}}
Suppose $b_1$ is $P_1$-dominant, $\mc F_2$ is right-basic, and
\[
\Hom (\Eis^{\T^{\ast}}_{P_1} \mc F_1, \Eis^{\T^{\ast}}_{P_2} \mc F_2)\not\simeq 0.
\]
Then $b_1 \leq b_2$ in $B(G)$. Dually, if $b_{2}$ is $P_2$-antidominant, $\mc F_1$ is left-basic, and 
\[
\Hom (\Eis^{\T^{\ast}}_{P_1} \mc F_1, \Eis^{\T^{\ast}}_{P_2} \mc F_2)\not\simeq 0, 
\]
then $b_1\geq b_2$ in $B(G)$. 
\end{corollary}

\begin{proof}
By twisted duality $\D'_{\GS}$, it suffices to show the first statement. 

Recall that $b_1 \leq b_2$ in $B(G)$ precisely when $\kappa_G(b_1) = \kappa_G(b_2)$ and $\nu_G(b_2) - \nu_G(b_1)$ is a $\Q_{\geq 0}$-linear combination of positive coroots.

We fix a Borel $B_1$ such that ${P_1}$ is standard and let $\Delta$ be the associated set of simple roots. The semistandard parabolic $P_2$ may not be standard, but we may choose some $w \in W(G,T)$ such that $w(\nu_{M_2}(b_2))=\nu_{w(M_2)}(w(b_2))$ is $B_1$-dominant. 

We may choose $\nu_1$ so that it factors through $\wh G_\der$. By assumption and Proposition \ref{prop: generaleisbasisbas}, we have $\kappa_G(b_1) = \kappa_G(b_2)$ and 
\[
\langle \nu_1,  w'(\nu_{M_2}(b_2)) -\nu_{M_1}(b_1) \rangle \geq 0
\]
for some $w'\in W_{H_\lambda}$ (note that $W_{H_{\lambda}}$ (non-canonically) maps to $W(G,T) = W_{\wh{G}}$ since an element of $N_{H_{\lambda}}(\wh{T}_H)$ normalizes $\wh{T}$ up to an element of $Z_{\wh{G}}(\wh{T}_H)$ and pairings with cocharacters of $\wh{T}_H$ do not depend on this choice). 
By writing $\Eis^{\T^{\ast}}_{P_1}$ as a composition of Eisenstein series, we may replace $P_1$ by a $\Delta$-standard maximal parabolic containing it. Such a parabolic is determined by a fundamental weight
$\delta_{\sigma}$ for $\sigma \in \Delta$ (the sum over the $\Gamma$-orbit of $\sigma$ if $G$ is not split). Then we have, for some $w'_\delta\in W_{H_\lambda}$, 
\begin{multline*}
    0 \leq \langle  \delta_{\sigma} , w'_\delta(\nu_{M_2}(b_2)) -\nu_{M_1}(b_1)  \rangle \\ 
    \leq \langle  \delta_{\sigma} , w(\nu_{M_2}(b_2)) - \nu_{M_1}(b_1) \rangle = \langle \delta_{\sigma}, \nu_G(b_2) - \nu_G(b_1)  \rangle,
\end{multline*}
where the middle inequality holds because $w(\nu_{M_2}(b_2)) \geq  w'_{\delta}(\nu_{M_2}(b_2))$ and the right equality holds because $w(\nu_{M_2}(b_2))$ and $\nu_{M_1}(b_1)$ are both dominant by assumption.

Now, $\nu_G(b_1) -\kappa_G(b_1)$ and $\nu_G(b_2)  - \kappa_G(b_2)$ lie in $X_{\ast}(T_{\der})_{\Q}$, where we are implicitly identifying $X_{\ast}(A_G)_{\Q}$ and $X^{\ast}(Z(\wh{G})^{\Gamma})_{\Q}$ and $A_G$ is the maximal split torus in the center of $G$. Thus,
\begin{equation*}
    \nu_G(b_2) - \nu_G(b_1) = (\nu_G(b_2) - \kappa_G(b_2)) - (\nu_G(b_1) - \kappa_G(b_1))
\end{equation*}
is a $\Q$-combination of simple coroots and by the computation of the previous paragraph, the coefficients are non-negative for the simple coroots not in $M_1$. For the simple coroots in $M_1$, we have non-negative coefficients as well: write $\nu_G(b_2) - \nu_G(b_1)=\sum_\beta c_\beta \beta^\vee$ (if $G$ is not split, $c_\beta$ is constant on each $\Gamma$-orbit). As $\langle \nu_G(b_1),\alpha\rangle=0$ for $\alpha\in \Delta_{M_1}$ and $\nu_G(b_2)$ is dominant, we have
\[
\sum_{\beta}c_\beta \langle \beta^\vee, \alpha\rangle=
\langle \nu_{G}(b_2), \alpha \rangle \geq 0. 
\]
Since $c_{\beta}\geq 0$ and $\langle \beta^\vee, \alpha\rangle \leq 0$ for $\beta \not\in \Delta_{M_1}$, we obtain $\sum_{\beta\in \Delta_{M_1}}c_\beta \langle \beta^\vee, \alpha\rangle\geq 0$. It follows that $c_\beta \geq 0$ for $\beta\in \Delta_{M_1}$ as well, as the inverse of the Cartan matrix of $M_1$ has non-negative entries. This is the desired result. 
\end{proof}

\begin{corollary}
Suppose $b_1$ is $P_1$-dominant, $b_2$ is $P_2$-antidominant, and  $\mc F_1 \in \Coh(\T^{\ast}\V^{\nu_1 = 0})^{\kappa_{M_1}(b_1)}$ is left-basic and $\mc{F}_2 \in \Coh(\T^{\ast}\V^{\nu_2 = 0})^{\kappa_{M_2}(b_2)}$ is right-basic. If 
\[
\Hom (\Eis^{\T^{\ast}}_{P_1} \mc F_1, \Eis^{\T^{\ast}}_{P_2} \mc F_2)\not\simeq 0, 
\]
then $b_1=b_2$ in $B(G)$.
\end{corollary}

\begin{proof}
This follows by combining the two inequalities above.  
\end{proof}

\subsection{Some variants}
We will use similar arguments to prove some variants. Fix a semisimple $L$-parameter $\lambda\colon W_F\to \wh G$. Let $P_1\supset M_1$, $P_2\supset M_2$ be semistandard parabolic and Levi subgroups through which  $\lambda$ factors; we write $\lambda_1\colon W_F\to \wh M_1$, $\lambda_2\colon W_F\to \wh M_2$ for the lifts to distinguish them. For $b_i\in X^*(Z(\wh M_i)^\Gamma)\cong B(M_i)_{\bas}$, $i \in \{1,2\}$, we let
\[
\mc F_1\in \Coh(\T^*\V_{\lambda_1})^{\kappa_{M_1}(b_1)}, \quad \mc F_2\in \RBas^\Coh(\T^*\V_{\lambda_2})^{\kappa_{M_2}(b_2)}.
\]
If, moreover, $b_1$ is $P_1$-dominant and
\[
\Hom(\Eis^{\T^*}_{P_1}\mc F_1, \Eis^{\T^*}_{P_2}\mc F_2)\not\simeq 0, 
\]
then $\kappa_G(b_1)=\kappa_G(b_2)$ and $b_1\leq b_2$ in $B(G)$ by Corollary \ref{cor: eisHomnonvanishingconditions}, choosing a cocharacter $\nu_1\in X_*(\wh T_H)$ defining $\wh P_1$; the same argument works in the current setup.

\begin{proposition}{\label{prop: EisHomvanishing}}
Suppose $b_i$ is strictly $P_i$-dominant and there exists $w\in W_{\wh G}^\Gamma$ such that $w(\wh M_2)=\wh M_1$, $w(\wh P_2)=\wh P_1$, $w(b_2)=b_1$; such $w$ is unique as a coset $\wh P_1 w \wh P_2=w \wh P_2$. If $\wh P_1 w\wh P_2\not\subset \wh P_1 H_\lambda \wh P_2$, equivalently $w\not\in H_\lambda \wh P_2$,  then
\[
\Hom(\Eis^{\T^*}_{P_1}\mc F_1, \Eis^{\T^*}_{P_2}\mc F_2)\simeq 0.
\] 
Otherwise, $\wh P_1 w\wh P_2$ has a unique lift to $(H_\lambda \cap \wh P_1)\backslash H_\lambda / (H_\lambda\cap \wh P_2)$, and 
\[
\Hom(\Eis^{\T^*}_{P_1}\mc F_1, \Eis^{\T^*}_{P_2}\mc F_2)\simeq \Hom(\mc F_1, w_*\mc F_2). 
\]

The same statements hold if the $b_i$ are strictly $P_i$-antidominant and we assume instead \[
\mc F_1\in \LBas^\Coh(\T^*\V_{\lambda_1})^{\kappa_{M_1}(b_1)}, \quad \mc F_2\in \Coh(\T^*\V_{\lambda_2})^{\kappa_{M_2}(b_2)}. 
\]
\end{proposition}

\begin{proof}
By twisted Grothendieck--Serre duality, it suffices to show the dominant case. We use the notation in the proof of Proposition \ref{prop: generaleisbasisbas} adapted to this setting. In the argument of Proposition \ref{prop: generaleisbasisbas}, if $w'\in (H_\lambda \cap \wh P_1)\backslash H_\lambda /(H_\lambda\cap \wh P_2)$ is not a lift of $w$, then
\[
\langle \nu_1, \kappa_{M_1} (b_1)\rangle=\langle \nu_1, \kappa_{M_1} (w(b_2))\rangle > \langle \nu_1, w'(\nu_{M_2}(b_2))\rangle
\]
and hence the contribution of $w'$ vanishes. In particular, if $\wh P_1w \wh P_2\not \subset \wh P_1 H_\lambda \wh P_2 $, all the contributions vanish. 

Assume $\wh P_1w \wh P_2 \subset \wh P_1 H_\lambda \wh P_2 $. Any element $h$ of $H_\lambda \cap \wh P_1 w\wh P_2$ satisfies $h\wh P_2h^{-1}=\wh P_1$, and hence the double coset $(H_\lambda\cap \wh P_1) h(H_\lambda\cap \wh P_2)$ is independent of $h$. Let us call this double coset $w$ as well by abuse of notation. To compute the contribution of $w$, first note that only the 1st infinitesimal neighborhood contributes.  As $g_w^* (\mc F_1\otimes \omega_{\V_{U_1}})$, $i_w^! p_1^! q_2^* (\mc F_2\otimes \omega_{\V_{U_2}})$ are both pulled back from $\T^*\V^{\nu_1=0}=\T^*\V^{w(\nu_2)=0}$ and have single $\nu_1$-weights, together with the computation of the dualizing complex of $h_w$ as in the proof of Proposition \ref{prop: generaleisbasisbas}, we can rewrite the contribution of $w$ as
\[
\Hom(\mc F_1, (w^{-1})^* \mc F_2)=\Hom(w^*\mc F_1, \mc F_2)=\Hom(\mc F_1, w_*\mc F_2),
\]
as desired. 
\end{proof}

\begin{corollary}{\label{cor: fullyfaithfullbas}}
Suppose $P\supset M$ is semistandard and $b_M\in B(M)_{\bas}$ is strictly $P$-dominant (resp.~strictly $P$-antidominant). Then $\Eis^{\T^*}_P$ is fully faithful on $\RBas^\Coh(\T^*\V_{\lambda_M})^{\kappa_M(b_M)}$ (resp.~$\LBas^\Coh(\T^*\V_{\lambda_M})^{\kappa_M(b_M)}$). 
\end{corollary}

\begin{proof}
In the previous notation, we now have $P_1=P_2=P$ and $w=e$. By the construction of the filtration, the $w=e$ piece induces a map
\[
\Hom(\mc F_1, \mc F_2)\to \Hom(\Eis^{\T^*}_{P_1}\mc F_1, \Eis^{\T^*}_{P_2}\mc F_2)
\]
that agrees with the one from the functoriality, and we have proven it is an equivalence. 
\end{proof}

Therefore, it is important to analyze the case where the images of $b_i$ in $B(G)$ are basic, in which case Eisenstein series preserves basicness. We will show that basic sheaves are generated by Eisenstein series of cuspidal sheaves, and we restrict ourselves to such sheaves from now on. In the following proposition, we work with $\Gm$-equivariant sheaves as this is needed for the applications. A non-equivariant version exists via the same argument but dropping $\Gm$-actions.

\begin{proposition}{\label{prop: Gmequivcusphomcomputation}}
Assume that 
\[
\mc F_1\in \LCusp^\Coh(\T^*\V_{\lambda_1} / \Gm)^{\kappa_{M_1}(b_1)}, \qquad \mc F_2\in \RCusp^\Coh(\T^*\V_{\lambda_2}/\Gm)^{\kappa_{M_2}(b_2)}
\]
and the images of $b_1$ and $b_2$ in $B(G)$ are the same basic element. Then, there is a filtration on (enriched) $\Hom(\Eis^{\T^*}_{P_1}\mc F_1, \Eis^{\T^*}_{P_2}\mc F_2)$ indexed by $w\in (H_{\lambda}\cap \wh P_1)\backslash H_\lambda /(H_{\lambda}\cap \wh P_2)$
such that $\wh M_1=w(\wh M_2)$ for some (equivalently any) representative in $W_{H_\lambda}$, whose graded piece at $w$ is given by
\[
\Hom(\mc F_1, (w^{-1})^* \mc F_2)\langle d_w \rangle[d_w]=\Hom(w^*\mc F_1, \mc F_2)\langle d_w \rangle[d_w]=\Hom(\mc F_1, w_*\mc F_2)\langle d_w \rangle[d_w], 
\]
where $\langle x\rangle[y]$ denotes a grading twist by $t\mapsto t^x$ and a degree shift of $y$ and
\begin{equation*}
    d_w = 2( \dim (\Lie H_{\lambda})^{\nu_1 > 0, w(\nu_2)<0}   - \dim (V_{\lambda})^{\nu_1 > 0, w(\nu_2)<0}).
\end{equation*}

In particular, it vanishes if there does not exist $w\in W_{H_\lambda}$ such that $\wh{M}_1=w(\wh{M}_2)$ and $\Hom (\mc F_1, w_*\mc F_2)\not\simeq 0$. 
\end{proposition}

\begin{proof}
Again we use the notation in the proof of Proposition \ref{prop: generaleisbasisbas} adapted to this setting and note that the construction of the filtration and the description of graded pieces are compatible with $\Gm$-gradings. If $w(\wh M_2)\cap \wh P_1\neq w(\wh M_2)$, then $(w^{-1}\circ q')^! \mc F_2$ has $\nu_1$-weights
\[
< \langle \nu_1, \wt(\V_{U_1\cap w(M_2)}) + w (\nu_{M_2} (b_2)) \rangle 
\]
by Proposition \ref{prop: basicweightcharacterization}. As $\langle \nu_1, w (\nu_{M_2} (b_2))\rangle=\langle \nu_1, \kappa_{M_1} (b_1)\rangle$ under the current assumption on $b_1, b_2$, the argument there shows that the contribution from $w$ vanishes. Note also that only the 1st infinitesimal neighborhood contributes, even if $w(\wh M_2)\cap \wh P_1= w(\wh M_2)$. We can similarly construct another filtration coming from $\mc F_1$ and the twisted Grothendieck--Serre duality, and we see from the left-cuspidality of $\mc F_1$ that the contribution from $w$ vanishes unless $\wh M_1\cap w(\wh P_2)=\wh M_1$. In fact, the two filtrations essentially agree: using that $p_2^* q_1^* (\mc F_1\otimes \omega_{\V_{U_1}})$ is almost perfect and $p_1^{!,\QCoh}q_2^* (\mc F_2 \otimes \omega_{\V_{U_2}})$ is bounded below, we see that the first filtration is obtained from the excision filtration on the internal Hom
\[
\underline{\Hom}_{\QCoh(Z)}(p_2^* q_1^* (\mc F_1\otimes \omega_{\V_{U_1}}), p_1^{!,\QCoh}q_2^* (\mc F_2 \otimes \omega_{\V_{U_2}})),
\]
where $Z\coloneq \T^\ast \V^{\nu_1\geq 0}\times_{\T^* \V}\T^\ast \V^{\nu_2\geq 0}$. Similarly, the second filtration is obtained from, up to a shift,
\[
\underline{\Hom}_{\QCoh(Z)}(p_1^* q_2^* (\D'_{\GS}\mc F_2\otimes \omega_{\V_{U_2}}), p_2^{!,\QCoh}q_1^* (\D'_{\GS}\mc F_1 \otimes \omega_{\V_{U_1}}))
\]
which, by the proof of Proposition \ref{prop: twistedGSEiscompatibility}, is equivalent to
\[
\underline{\Hom}_{\QCoh(Z)}(p_1^* \D_\GS (q_2^* (\mc F_2\otimes \omega_{\V_{U_2}})), p_2^{!,\QCoh}\D_{\GS} (q_1^* (\mc F_1 \otimes \omega_{\V_{U_1}})))
\]
up to a shift. Since
\[
p_1^{!,\QCoh}(-)\simeq \underline{\Hom}_{\QCoh(Z)}(p_1^*\D_{\GS}-,\omega_Z), \quad
p_2^{!,\QCoh}(-)\simeq \underline{\Hom}_{\QCoh(Z)}(p_2^*\D_{\GS}-,\omega_Z)
\]
on coherent sheaves, both internal Homs are equivalent to
\[
\underline{\Hom}_{\QCoh(Z)}(p_2^* q_1^* (\mc F_1\otimes \omega_{\V_{U_1}}) \otimes p_1^* (\D_\GS q_2^* (\mc F_2\otimes \omega_{\V_{U_2}})), \omega_Z), 
\]
and both filtrations are induced from the same excision filtration. 

Assume now that $w(\wh M_2)\cap \wh P_1= w(\wh M_2)$ and $\wh M_1\cap w(\wh P_2)=\wh M_1$, equivalently $\wh M_1=w(\wh M_2)$. In this situation, the contribution from $w$ can be written as
\[
\Hom(g_w^* (\mc F_1\otimes \omega_{\V_{U_1}}), i_w^! p_1^! q_2^* (\mc F_2\otimes \omega_{\V_{U_2}})). 
\]
As $g_w^* (\mc F_1\otimes \omega_{\V_{U_1}})$, $i_w^! p_1^! q_2^* (\mc F_2\otimes \omega_{\V_{U_2}})$ are both pulled back from $\T^*\V^{\nu_1=0}=\T^*\V^{w(\nu_2)=0}$ and have a single $\nu_1$-weight, together with the computation of the dualizing complex of $h_w$ as in the proof of Proposition \ref{prop: generaleisbasisbas}, we can rewrite the contribution of $w$ as
\[
\Hom(\mc F_1, (w^{-1})^* \mc F_2)\langle d_w \rangle[d_w]=\Hom(w^*\mc F_1, \mc F_2)\langle d_w \rangle[d_w]=\Hom(\mc F_1, w_*\mc F_2)\langle d_w \rangle[d_w]. 
\]
Note the shift and twist come from the nontrivial relative dualizing complex of $h_w$ whose computation follows similarly to that of its weight pairing with $\nu_1$ in the proof of Proposition \ref{prop: generaleisbasisbas}, and the shift also comes from $\omega^{-1}_{\V_{U_1}} \otimes w_{\ast}\omega_{\V_{U_2}}$ as well. 
\end{proof}

Let us note the following lemma:

\begin{lemma}{\label{lem: dblcosetlem}}
The natural map 
\[
(H_{\lambda}\cap \wh P_1)\backslash H_\lambda /(H_{\lambda}\cap \wh P_2) \to
\wh P_1 \backslash \wh G /\wh P_2
\]
is injective on the subset consisting of
$w\in (H_{\lambda}\cap \wh P_1)\backslash H_\lambda /(H_{\lambda}\cap \wh P_2)$
such that $w(\wh M_2)=\wh M_1$ for some (equivalently any) representative in $W_{H_\lambda}$.

In particular, if $\wh P\coloneq\wh P_1=\wh P_2$, $\wh M\coloneq \wh M_1=\wh M_2$, then this subset identifies with the subgroup $N_{H_\lambda}(\wh M)/(\wh M \cap H_\lambda)$ of $(N_{\wh G}(\wh M)/\wh M)^\Gamma$. 
\end{lemma}

\begin{proof}
If $w\in W_{H_\lambda}$ satisfies $w (\wh M_2)=\wh M_1$, then 
\[
W_{\wh M_1 \cap H_\lambda}w W_{\wh M_2\cap H_\lambda}=W_{\wh M_1 \cap H_\lambda}w=w W_{\wh M_2 \cap H_\lambda}.
\]
We also have
\[
W_{\wh M_1}w W_{\wh M_2}= W_{\wh M_1}w=w W_{\wh M_2}.
\]

Let $\mathrm{Tr}_{1,2}\coloneq \{g\in \wh G \mid g(\wh M_2)=\wh M_1\}$. The maps 
\[
N_{H_{\lambda}}(\wh T_H) \cap \mathrm{Tr}_{1,2} \to W_{\wh M_1}\backslash W_{\wh G}, \qquad
N_{H_{\lambda}}(\wh T_H) \cap \mathrm{Tr}_{1,2} \to W_{\wh G}/W_{\wh M_2}
\]
induce injective maps
\begin{align*}
N_{H_\lambda \cap \wh M_1}(\wh T_H) \backslash (N_{H_{\lambda}}(\wh T_H) \cap \mathrm{Tr}_{1,2}) &\to W_{\wh M_1}\backslash W_{\wh G}, \\
(N_{H_{\lambda}}(\wh T_H) \cap \mathrm{Tr}_{1,2})/N_{H_\lambda \cap \wh M_2}(\wh T_H) &\to W_{\wh G}/W_{\wh M_2} 
\end{align*}
and both agree with the map in the statement. 
\end{proof}

\subsection{Applications to generalized coherent Springer sheaves}

Keep the notation of the previous subsection. We have applications to generalized coherent Springer sheaves generalizing some results of \cite{BCHNspringer} and \cite{ZhuTame}. 

\begin{proposition}{\label{prop: cusphomfilt}}
Suppose $b_1\in B(M_1)_{\bas}$, $b_2\in B(M_2)_{\bas}$ have the same image in $B(G)$, which is basic, and let $\rho_1$, $\rho_2$ be irreducible cuspidal representations of $S_{\phi_i}$ with central characters $\kappa_{M_1}(b_1)$, $\kappa_{M_2}(b_2)$ respectively. Then, there is a filtration on
$\Hom(\Eis^{\spec, \wedge}_{P_1}\wh{\mc V}_{\rho_1}, \Eis^{\spec, \wedge}_{P_2}\wh{\mc V}_{\rho_2})$
indexed by $w\in (H_\lambda \cap \wh P_1)\backslash H_\lambda / (H_\lambda \cap \wh P_2)$ such that $w(\wh{M}_2)=\wh{M}_1$, whose graded piece at $w$ is given by $\Hom (\wh{\mc V}_{\rho_1}, w_* \wh{\mc V}_{\rho_2})$. In particular, $\Hom(\Eis^{\spec, \wedge}_{P_1}\wh{\mc V}_{\rho_1}, \Eis^{\spec, \wedge}_{P_2}\wh{\mc V}_{\rho_2})$ vanishes if there does not exist $w$ such that $w(\wh{M}_2)=\wh{M}_1$ and $\rho_1\cong {}^w\rho_2$, and it is concentrated in degree $0$ in general. 
\end{proposition}

\begin{proof}
The first part follows from the linear Koszul duality and Lemma \ref{lem: HatCohForgettingGm} by using the natural (or trivial) grading on $\wh{\mc V}_\rho$, and Proposition \ref{prop: Gmequivcusphomcomputation}. Then, the second part follows from the properties of $\wh{\mc V}_\rho$. Note that under sheared Koszul duality, the twist in Proposition \ref{prop: Gmequivcusphomcomputation} cancels the degree shift.
\end{proof}

\begin{corollary}{\label{cor: cohspringerhoms}}
For any cuspidal irreducible representations $\rho_i$ whose central characters are $\kappa_{M_i}(b_i)$ for $b_i\in B(M_i)_{\bas}$ with the same image in $B(G)$, which is basic, 
\[
\Hom(\Eis^\spec_{P_1}\widetilde{\mc V}_{\rho_1}, \Eis^\spec_{P_2}\widetilde{\mc V}_{\rho_2})
\]
is concentrated in degree 0. It is nonzero if and only if $(\lambda_1, \rho_1)$ is $\wh G$-inertially equivalent to $(\lambda_2, \rho_2)$ in the sense that there exist $g\in \wh G$ and $z\in Z(\wh M_2)^{\Gamma, \circ}$ such that $g\wh M_2 g^{-1}=\wh M_1$, ${}^L \lambda_1=g(z\cdot {}^L \lambda_2)g^{-1}$ (at the level of representatives), and $\mc V_{\rho_1}$ and $\mc V_{\rho_2}$ are isomorphic under the identifications 
\[
\V_{\lambda_2}\simeq \V_{z\cdot \lambda_2} \xrightarrow[\simeq]{\Ad(g)}\V_{\lambda_1}.
\]
\end{corollary}

\begin{proof}
We claim it suffices to prove a similar statement after taking $\wh{i}^*_\lambda$ for all $\lambda$. 

Indeed, let us write $f$ for the map to the coarse moduli space
\[
f\colon \Par_G \to X_G. 
\]
It suffices to prove that $f_*\underline{\Hom} (\Eis^\spec_{P_1}\widetilde{\mc V}_{\rho_1},\Eis^\spec_{P_2}\widetilde{\mc V}_{\rho_2})\in \QCoh(X_G)$ is concentrated in degree 0, and it is nonzero if and only if the inertial equivalence holds.\footnote{It is known to be a \emph{bounded} complex of \emph{coherent} sheaves by \cite[Theorem 3.7.3]{HansenMann}, but we do not need it here.} 
It further suffices to prove these claims after $*$-pullback to $\Spec \mathcal{O}_{X_G,x}$, or $\Spec \wh{\mathcal{O}}_{X_G,x}$, for all closed points $x$ of $X_G$. Using coherent completeness as in the proof of Proposition \ref{prop: unipvsformal}, we need only prove these claims for
\[
\Hom(\wh{i}_\lambda^* \Eis^\spec_{P_1}\widetilde{\mc V}_{\rho_1}, \wh{i}_\lambda^* \Eis^\spec_{P_2} \widetilde{\mc V}_{\rho_2})
\]
in $\wh{\Coh}{}^*(\Par_{G,\lambda}^\wedge)$, or dually, 
\[
\Hom(\wh{i}_\lambda^! \Eis^\spec_{P_1}\widetilde{\mc V}_{\rho_1}, \wh{i}_\lambda^! \Eis^\spec_{P_2} \widetilde{\mc V}_{\rho_2})
\]
in $\wh{\Coh}{}^!(\Par_{G,\lambda}^\wedge)$, for any semisimple $L$-parameter $\lambda\colon W_F\to \wh G$. 

It suffices to consider only translates of $\lambda_1$, and we first take $\lambda=\lambda_1$. We only need to consider the case where $w(\wh M_2)=\wh M_1$ for some $w\in W_{H_\lambda}$ by Proposition \ref{prop: cusphomfilt}. It further follows from Proposition \ref{prop: cusphomfilt} together with the description of $\wh{i}^!_{\lambda}\widetilde{\mc V}_{\rho_i}$ that, if this Hom does not vanish, then $(\lambda_1, \rho_1)$ is $\wh M_1$-inertially equivalent to $({}^w \lambda_2, {}^w\rho_2)$, i.e., there exists $z\in Z(\wh M_1)^{\Gamma,\circ}$ such that $\lambda_1=z\cdot {}^w\lambda_2$ and $\rho_1 \cong {}^w \rho_2$ under the identification $S_{\phi_1}\cong S_{{}^w \phi_2}$ induced by the translation by $z$. These conditions together are equivalent to the $\wh G$-inertial equivalence in the statement. Moreover, the nonvanishing and concentration in degree 0 in the $\wh G$-inertial case follow from Proposition \ref{prop: cusphomfilt}. 

In fact, the above argument does not change by replacing $\lambda$ by its translation by some $z'\in Z(\wh M_1)^{\Gamma,\circ}$, and then it follows that the non-vanishing at one $\lambda$ is equivalent to the non-vanishing at all translates of $\lambda$. 
\end{proof}

Fix $P\supset M$ and $\rho$. Recall that by the paragraph after Corollary \ref{cor: CVBproperty}, we have $\End (\widehat{\mc V}_\rho)\simeq \mc O (\Lie Z(\wh M)^\Gamma)^\wedge_0$ naturally. 

\begin{lemma}
The filtration on $\End(\Eis^{\spec,\wedge}_{P}\wh{\mc V}_{\rho})$ is a filtration by bi-$\End (\wh{\mc V}_{\rho})$-modules, and the graded piece at $w$, if nonzero, is isomorphic to the bi-$\End (\wh{\mc V}_{\rho})$-module $\mc O (\Lie Z(\wh M)^\Gamma)^\wedge_0$ whose right action is the natural one and whose left action is twisted by $w$ on $Z(\wh M)^\Gamma$. 
\end{lemma}

\begin{proof}
Only the last claim is not clear from the construction of the filtration. 
After choosing an isomorphism $w_*\widehat{\mc V}_\rho\cong \wh{\mc V}_\rho$, which exists by assumption, we can describe the left action as
\[
\End (\wh{\mc V}_{\rho})\xrightarrow{w_*} \End (w_*\wh{\mc V}_{\rho}) \cong \End (\widehat{\mc V}_\rho). 
\]
Once we identify $\End (\wh{\mc V}_{\rho})$ with $\mc O (\Lie Z(\wh M)^\Gamma)^\wedge_0$, this is simply induced by $w$ on $Z(\wh M)^\Gamma$.  
\end{proof}

\begin{corollary}{\label{cor: EndEisfinitefree}}

\begin{enumerate}
    \item The algebra $\End(\Eis^{\spec, \wedge}_P \wh{\mc V}_\rho)$ is finite free both as left and right $\End(\wh{\mc V}_\rho)$-modules. 
    \item Write $W_{\rho}\subset N_{H_\lambda}(\wh M)/(\wh M \cap H_\lambda) $ for the subgroup (see Lemma \ref{lem: dblcosetlem}) consisting of $w\in (\wh P \cap H_\lambda)\backslash H_\lambda / (\wh P \cap H_\lambda)$ such that $w(\wh M)=\wh M$ and $\rho\cong {}^w \rho$. 
    The center of $\End(\Eis^{\spec, \wedge}_P \wh{\mc V}_\rho)$ is given by $(\mc O(\Lie Z(\wh M)^\Gamma)^{W_{\rho}})^\wedge_0$. 
\end{enumerate}
\end{corollary}

\begin{proof}
(1) This is clear from the previous lemma. Since $W_\rho$ acts faithfully on $\Lie Z(\wh M)^{\Gamma}$, (2) follows from \cite[Lemma 36]{BernsteinRumelhart} (or \cite[Lemma in p.205]{Renard}). 
\end{proof}

The above structure of the endomorphism algebra and its center is reminiscent of those in \cite{BernsteinRumelhart} in the context of Bernstein's compact projective generators of Bernstein blocks. 

\section{Semiorthogonal decompositions}{\label{s: SODs}}
After discussing semiorthogonal decompositions for cotangent stacks beyond the graded setting, we will complete the proofs of Theorems \ref{thm: CVBclassification} and \ref{thm: Generators} and hence Corollary \ref{cor: LCusp=RCusp}. We then deduce the semiorthogonal decompositions of $\Coh(\Par_G)$ and discuss additional consequences. 

\subsection{Semiorthogonal decomposition for cotangent stacks}{\label{ss: SODcotangent}}
We will prove Theorem \ref{thm: ClassificationCuspidalCotangent} (5) and Theorem \ref{thm: BasicCotangentMain} in the general case by induction and we assume they hold for all proper semistandard Levi subgroups. We first introduce some notation. 

Fix $b\in B(G)$. There exists a unique standard Levi subgroup $M_b$ such that $b$ comes from a unique strictly $P_b$-dominant element $b_M \in B(M_b)^{\circ +}_{\bas}$ for $P_b \subset G$ the standard parabolic subgroup with standard Levi $M_b$ (\cite[\S5.1]{KottwitzisoII}). We then attach to $b$ the pair $(P_b, \chi_b)$ where $\chi_b = \kappa_{M_b}(b_M) \in X^{\ast}(Z(\wh{M}_b)^{\Gamma})$. 

Fix a semisimple $\lambda\colon W_F\to \wh G$. By Lemma \ref{lem: TJcontrol}, after replacing $\lambda$ by a $\wh G$-conjugate, we may assume that $(\wh T \cap H_\lambda)^\circ$ is a maximal torus of $H_\lambda$. Recall that $\lambda$ factors through a semistandard $\wh M$ if and only if $Z(\wh M)^{\Gamma,\circ}\subset (\wh T \cap H_\lambda)^\circ$. 

Suppose $b$ is not basic. For any pair $(P, \chi)$ of a proper semistandard parabolic subgroup $P$ with semistandard Levi $M$ and a character $\chi\in X^*(Z(\wh M)^\Gamma)$ such that $(P, \chi)$ is $W_{\wh G}^\Gamma$-conjugate to $(P_b,\chi_b)$ and $\lambda$ factors through $\wh M$ (we then write $\lambda_M$ for $\lambda$), consider
\[
\Eis_{P}^{\T^*} \Bas^\Coh(\T^*\V_{\lambda_M})^{\chi} \subset \Eis^{\circ}_+;
\]
recall that $\Eis_{P}^{\T^*}$ is fully faithful on $\Bas^\Coh(\T^*\V_{\lambda_M})^{\chi}$ by Corollary \ref{cor: fullyfaithfullbas}, so its essential image is closed under finite colimits and retracts. 
Write $\Coh(\T^*\V_\lambda)^{b\sharp}$ for the subcategory generated under finite colimits by these subcategories for all such 
$(P,\chi)$.\footnote{Proposition \ref{prop: EisHomvanishing} implies that different $(P,\chi)$'s produce orthogonal subcategories if they are not $W_{H_\lambda}$-conjugate. Note also that $W_{H_\lambda}$-conjugate pairs give the same subcategory.} Similarly, write $\Coh(\T^*\V_\lambda)^{b!}$ for the subcategory generated under finite colimits by
\[
\Eis^{\T^{\ast}}_{P^\op} \Bas^\Coh(\T^*\V_{\lambda_M})^{\chi}\subset \Eis^{\circ}_-,
\]
for all such $(P,\chi)$. 

If $b$ is basic, we define
\begin{equation*}
    \Coh(\T^{\ast}\V_\lambda)^{b\sharp}\coloneq \RBas^\Coh(\T^*\V_\lambda)^{\kappa(b)}, \quad  \Coh(\T^{\ast}\V_\lambda)^{b!}\coloneq \LBas^\Coh(\T^*\V_\lambda)^{\kappa(b)}. 
\end{equation*}

\begin{proposition}{\label{prop: refinedcotanSOD}}
Assume that Theorem \ref{thm: BasicCotangentMain} holds for all proper semistandard Levi subgroups $M$ of $G$ and all $\lambda_M$ as above.  
\begin{enumerate}
    \item The collection of subcategories 
    \[\{\Coh(\T^*\V_\lambda)^{b\sharp}\}_{b\in B(G)} \qquad  (\text{resp}.~\{\Coh(\T^*\V_\lambda)^{b!}\}_{b\in B(G)})
    \]
    forms a semiorthogonal decomposition of $\Coh(\T^*\V_\lambda)$ with the ascending (resp.~descending) partial order on $B(G)$.
    \item The right-orthogonal in $\Eis_+^{\T^*}$ of the subcategory generated under finite colimits by $\Coh(\T^*\V_\lambda)^{b\sharp}$ for non-basic $b$ is generated under finite colimits and retracts by
\[
\Eis^{\T^{\ast}}_P \Bas^\Coh(\T^*\V_{\lambda_M})^{\chi} \subset \RBas^\Coh(\T^*\V_\lambda)
\]
with $\chi$ corresponding to $b_M\in B(M)_\bas$ whose image in $B(G)$ is basic.
\item The left-orthogonal in $\Eis_-^{\T^*}$ of the subcategory generated under finite colimits by $\Coh(\T^*\V_\lambda)^{b!}$ for non-basic $b$ is given by the same subcategory of $\LBas^\Coh(\T^*\V_\lambda)$.

\item The subcategories $\Coh(\T^* \V_\lambda)^{b\sharp}$ (resp.~$\Coh(\T^* \V_\lambda)^{b!}$) are generated under finite colimits and retracts by 
\[
\Eis^{\T^{\ast}}_P \Cusp^\Coh(\T^*\V_{\lambda_M})^\chi \qquad (\text{resp}.~\Eis^{\T^{\ast}}_{P^\op} \Cusp^\Coh(\T^*\V_{\lambda_M})^\chi),
\]
where $\chi=\kappa(b_M)$ is a $P$-dominant basic element for $b_M\in B(M)_\bas$ whose image in $B(G)$ is $b$. In particular, $\LBas^\Coh(\T^{\ast}\V_\lambda)$, $\RBas^\Coh(\T^{\ast}\V_\lambda)$ are both generated under finite colimits and retracts by $\Eis_P^{\T^*} \Bas^\Coh(\T^*\V_{\lambda_M})^{\chi}$ with $\chi$ basic in $G$ as above and $\Cusp^\Coh(\T^*\V_\lambda)$; the latter is both left and right orthogonal to the former. 
\end{enumerate}
\end{proposition}

\begin{proof}
We begin with the proof of (1) and discuss the semiorthogonal decomposition $\{\Coh(\T^*\V_\lambda)^{b\sharp}\}_{b\in B(G)}$, noting that the other one is analogous, or essentially the twisted Grothendieck--Serre dual. 

We first check the semiorthogonality relations. Suppose that we have $b, b' \in B(G)$ such that $b \nleq b'$. We must show that if $\mc{F} \in \Coh(\T^{\ast}\V_\lambda)^{b \sharp}$ and $\mc{F}' \in \Coh(\T^{\ast}\V_\lambda)^{b' \sharp}$, then $\Hom(\mc{F}, \mc{F}') \simeq 0$. This follows from Corollary \ref{cor: eisHomnonvanishingconditions}. 

We now prove the rest of (1), and we assume we have proven the proposition for all proper Levi subgroups of $G$ and proceed by induction. For (1), we need to show that the categories $\{\Coh(\T^{\ast}\V_{\lambda})^{b\sharp}\}_{b \in B(G)}$ generate $\Coh(\T^{\ast}\V_{\lambda})$. By Theorem \ref{thm: ClassificationCuspidalCotangent} (1)', the subcategories $\Eis^{\T^{\ast}}_+$ and $\Cusp^{\Coh}(\T^{\ast}\V_{\lambda})$ generate $\Coh(\T^{\ast}\V_{\lambda})$, and hence it suffices to show that the subcategory generated by these categories contains $\Eis_+^{\T^*}$. 

Observe that by the inductive assumption, $\Eis_+^{\T^*}$ is generated by $\Eis^{\T^{\ast}}_P \Bas^\Coh(\T^{\ast}\V_{\lambda_M})^{\chi}$ for central weights $\chi$ such that there exists some simple root $\sigma$ omitted by $P$, with respect to some Borel subgroup contained in $P$, and $\langle \wh{\chi}, \delta_{\sigma} \rangle \geq 0$ for $\delta_{\sigma}$, the fundamental weight of $\sigma$. Note that the inductive assumption does not allow us to immediately conclude that $\chi$ is $P$-dominant. We remark that such generators are certainly in $\Eis_+^{\T^*}$. Indeed, if $P_{\sigma} \supset P$ is the semistandard maximal parabolic subgroup omitting $\sigma$, then one can write any sheaf in the above category as $\Eis^{\T^{\ast}}_{P_{\sigma}}$ applied to a sheaf with central weight $\chi_\sigma\coloneq\chi|_{Z(\wh{M}_{\sigma})^{\Gamma}}$, and we have
\begin{equation*}
    0 \leq \langle \wh{\chi}, \delta_{\sigma} \rangle = \langle \wh{\chi_\sigma}, \delta_{\sigma} \rangle,
\end{equation*}
which is equivalent to dominance for the semistandard maximal parabolic $P_{\sigma}$. We denote this condition on $\chi$, property $(\bm{\ast})$; note that it is independent of the choice of Borel subgroup contained in $P$. 

Hence, it remains to show that the categories $\Eis^{\T^{\ast}}_P \Bas^\Coh(\T^{\ast}\V_{\lambda_M})^{\chi}$ for $\chi$ satisfying $(\bm{\ast})$ are contained in the subcategory generated by $\{\Coh(\T^{\ast}\V_{\lambda})^{b\sharp}\}$. Consider $\Eis^{\T^{\ast}}_P\mc{F} \in \Eis^{\T^{\ast}}_P \Bas^\Coh(\T^{\ast}\V_{\lambda_M})^{\chi}$. If $\chi$ is $P$-dominant, we are done (after potentially passing to a larger Levi where $\chi$ is strictly dominant).

Otherwise, choose a simple root $\alpha$ such that $\langle \wh{\chi}, \alpha \rangle < 0$. Let $L \subset Q$ be the Levi and parabolic subgroups containing $M \subset P$ as well as the root subgroup of $\alpha$. By applying this proposition to $L$ (which is allowed by our inductive assumption and $L \neq G$ since otherwise $M$ would be maximal and $(\bm{\ast})$ would imply dominance), we have that $\Eis^{\T^{\ast}}_{P_L} \mc{F}$ (where $P_L  = P \cap L$) has a filtration with graded pieces in $\{\Coh(\T^{\ast}\V_{\lambda_L})^{b\sharp}\}_{b \in B(L)}$. Then the first part of Corollary \ref{cor: eisHomnonvanishingconditions} shows that
\begin{equation*}
    \Hom(\Eis^{\T^{\ast}}_{P'_L} \mc{F}', \Eis^{\T^{\ast}}_{P_L} \mc{F}) \neq 0
\end{equation*}
implies that $b' \leq b_{\chi}$ in $B(L)$, where $b_{\chi} \in B(M)_{\bas}$ is the element corresponding to $\chi$ and $\mc{F}'$ is assumed to have a single central weight $\kappa(b') \in X^{\ast}(Z(\wh{M'_L})^{\Gamma})$. Then only finitely many $b'$ can occur in the filtration. Since $M \subset L$ is the maximal Levi in $L$ such that $b_\chi$ is basic (and it is not basic in $L$), it follows that the top term is $\Eis^{\T^{\ast}}_{P''} \mc{F}''$ for $\mc{F}''$ of central weight $\chi''$ such that $(P'', \chi'')$ is $W_{H_\lambda}$-conjugate to $(P^{\op}_L, \chi)$. Now, the lower terms (induced to $G$) correspond to smaller elements in $B(G)$, which, if they satisfy $(\bm{\ast})$, can inductively be assumed to be expressible in terms of $\{\Coh(\T^{\ast}\V_{\lambda})^{b\sharp}\}$. The top degree term, if it satisfies $(\bm{\ast})$, now pairs negatively with fewer roots of the unipotent radical of the new parabolic than $U_P$, and so by iterating, we can eventually make it dominant. 

Now we are done if we can show that the inductions to $G$ of the graded pieces in the filtration of $\Eis^{\T^{\ast}}_{P_L} \mc{F}$ satisfy $(\bm{\ast})$. Note that $\kappa_L(b') = \kappa_L(b_{\chi})=\chi|_{Z(\wh{L})^{\Gamma}}$. If $\alpha  \neq \sigma$, then each graded piece satisfies $(\bm{\ast})$ via $\sigma$. If $\alpha = \sigma$, then we have 
\begin{equation*}
    \langle \wh{\chi} , \sigma \rangle < 0 \leq \langle \wh{\chi}, \delta_{\sigma} \rangle.
\end{equation*}
Let $\sigma^M$ be the projection of $\sigma$ away from $\Delta_M$. Namely, $\sigma - \sigma^M$ is a combination of simple roots of $M$ and $\langle \gamma^\vee, \sigma^M \rangle = 0$ for each $\gamma \in \Delta_M$. Note that since $\wh{\chi}$ lands in $A_M$, we have $\langle \wh{\chi}, \sigma^M \rangle = \langle \wh{\chi}, \sigma\rangle$. We have $\sigma^M = \sum\limits_{\beta \in \Delta_G \setminus \Delta_M} \langle \beta^{\vee}, \sigma^M \rangle \delta_{\beta}$. The term $\langle \sigma^{\vee}, \sigma^M \rangle \delta_{\sigma}$ pairs non-negatively with $\wh{\chi}$ since $\langle \sigma^{\vee}, \sigma^M \rangle > 0$ (indeed, if we had $\langle \sigma^{\vee}, \sigma^M\rangle  \leq 0$, then $\sigma^M$ would be $L$-antidominant and hence a non-positive weight, but we know that the coefficient of $\sigma$ in $\sigma^M$ is $1$). 
Thus, some other term pairs negatively. But we claim the pairings $\langle \beta^{\vee}, \sigma^M \rangle$ are non-positive. This will imply that $\langle \wh{\chi}, \delta_{\beta}\rangle > 0$ for some $\beta \not\in (\Delta_M \cup \{\sigma\})$ as desired. 

To settle the claim, note that
\begin{equation*}
    \langle \gamma^{\vee}, \sigma - \sigma^M \rangle = \langle \gamma^{\vee}, \sigma \rangle \leq 0,
\end{equation*}
for each $\gamma \in \Delta_M$. Hence $\sigma^M - \sigma$ is $M$-dominant. It follows that we can write 
\begin{equation*}
    \sigma^M = \sigma + \sum\limits_{\gamma \in \Delta_M} c_{\gamma}\gamma,
\end{equation*}
where $c_{\gamma} \geq 0$. The claim that  $\langle \beta^{\vee}, \sigma^M \rangle \leq 0$ for any $\beta \in \Delta_G \setminus (\Delta_M \cup \{\sigma\})$ follows because we have both
\begin{equation*}
    \langle \beta^{\vee}, \sigma\rangle \leq 0, \qquad \langle \beta^{\vee}, \sum\limits_{\gamma \in \Delta_M} c_{\gamma}\gamma \rangle \leq 0.
\end{equation*}

By the semiorthogonal decomposition in (1), we need to compute the intersection of $\Eis^{\T^{\ast}}_+$ and $\{\Coh(\T^*\V_\lambda)^{b\sharp}\}$ for basic $b$. This is provided by Theorem~\ref{thm: BasicCotangentMain} (2) which gives the desired result. Then (3) is the twisted Grothendieck--Serre dual to (2). By combining with the inductive hypothesis, (4) also follows for $G$.  
\end{proof}

This completes the proof of Theorem~\ref{thm: ClassificationCuspidalCotangent} (5) and parts (1) and (3) of Theorem~\ref{thm: BasicCotangentMain} in the general case by passing to the ind-completion.

\subsection{Semiorthogonal decompositions in other settings}

We continue with the notation of \S\ref{ss: SODcotangent}.

\begin{corollary}\label{cor: SODhatCoh}
We have semiorthogonal decompositions $\{\wh{\Coh}{}^!(\wh{\T}\V_{\lambda}[-1])^{b\sharp}\}$ and $\{\wh{\Coh}{}^!(\wh{\T}\V_{\lambda}[-1])^{b!}\}$ of $\wh{\Coh}{}^!(\wh{\T}\V_{\lambda}[-1])$, with analogously defined subcategories, and each piece admits a further decomposition according to ``cuspidal supports''. The same statement holds for $\Par_{G,\lambda}^{\wedge}$. In particular, the image of 
\[
\wh{i}_\lambda^!\colon \Coh(\Par_G) \to \wh{\Coh}{}^!(\Par_{G,\lambda}^\wedge) \simeq \wh{\Coh}{}^!(\wh{\T}\V_{\lambda}[-1])
\]
generates the target under finite colimits and retracts, and the same holds for $\chi$-summands for any $\chi \in X^*(Z(\wh G)^{\Gamma})$. 
\end{corollary}

\begin{proof}
As there are natural generators with graded lifts and Eisenstein series are compatible with linear Koszul duality, the first claim follows from the case of cotangent stacks and Lemma \ref{lem: HatCohForgettingGm}. The second claim then follows from coherent completeness as in the proof of Proposition \ref{prop: unipvsformal}. The last claim holds since, by the first claim refined by cuspidal supports, the target is generated under finite colimits and retracts by $\Eis^{\spec, \wedge}_{P,\lambda_P}\wh{\mc V}_{\rho_M}$; each $\wh{\mc V}_{\rho_M}$ is a summand of $\wh{i}^!_{\lambda_M}\widetilde{\mc V}_{\rho_M}$, so by Corollary \ref{cor: EisComparisonCompletedFiber} $\Eis^{\spec, \wedge}_{P,\lambda_P}\wh{\mc V}_{\rho_M}$ is a summand of $\wh{i}^!_\lambda \Eis^\spec_P \widetilde{\mc V}_{\rho_M}$. 
\end{proof}

\begin{remark}\label{rem: ihatpreservescuspIndCoh}
The last part of Corollary \ref{cor: SODhatCoh} implies that
\[
\wh{i}^!_\lambda\colon \IndCoh(\Par_G)\to \Ind\wh{\Coh}{}^!(\Par^\wedge_{G,\lambda})
\]
preserves right-cuspidality. In Lemma \ref{lem: ihatpreservescusp}, this is shown only for coherent sheaves. 
\end{remark}

We need a slight variant:

\begin{corollary}\label{cor: CohCompletionSOD}
We have semiorthogonal decompositions $\{\Coh(\Par_{G,\lambda}^{\wedge})^{b\sharp}\}$ and $\{\Coh(\Par_{G,\lambda}^{\wedge})^{b!}\}$ of $\Coh(\Par_{G,\lambda}^{\wedge})$, with analogously defined subcategories, and each admits a decomposition according to ``cuspidal supports'' in terms of categories generated by generalized coherent Springer sheaves. In particular, the full subcategory of right-cuspidal sheaves in $\IndCoh(\Par_{G,\lambda}^{\wedge})$ is generated under colimits by cuspidal vector bundles.  
\end{corollary}

\begin{proof}
The second claim follows from the first: any compact object in the right-cuspidal sheaves is a cuspidal coherent sheaf $\mc F\in \Coh(\Par_{G,\lambda}^{\wedge})$. It is therefore cuspidal in $\wh{\Coh}{}^!(\Par_{G,\lambda}^{\wedge})$ and hence is determined by its restriction to $\wh{\T}(BS_\phi)[-1]$ by Theorem \ref{thm: ClassificationCuspidalCotangent} and Koszul duality. Since the restriction of $\mc F$ is in $\Coh(\wh{\T}(BS_\phi)[-1])$, the restriction is in the category generated by irreducible cuspidal representations of $S_\phi$ on the zero section under finite colimits. In turn, this means $\mc F$ is in the category generated by cuspidal vector bundles on $\V_\lambda$ under finite colimits. 

Next, we prove the first claim. We need only prove the refined version with ``right-cuspidal supports''. The relevant categories satisfy the semiorthogonality relations since $\Coh$ is a subcategory of $\wh{\Coh}{}^!$ and the inclusion preserves right-cuspidality. Therefore, by arguing inductively as before, it suffices to show that the subcategory of right-cuspidal sheaves in $\IndCoh(\Par_{G,\lambda}^\wedge)$ is compactly generated by coherent sheaves. This claim reduces to the case that $G$ is semisimple by arguing analogously to the reduction of Theorem \ref{thm: CVBclassification} (3) in \S \ref{ss: reductiondiscsuppcvbclassgraded}--\ref{ss: reductiontosimplyconnecteddiscsuppandCVBclass}. 

Assume $G$ is semisimple. It suffices to show that $(\mc F)^{\rc,\IndCoh}\in \IndCoh(\Par_{G,\lambda}^\wedge)$ for $\mc F\in \Coh(\Par_{G,\lambda}^\wedge)$ is still coherent, where the superscript indicates the ambient category in which we are considering the adjoint functor. Consider the unit
\[
\mc{F} \rightarrow (\mc F)^{\rc, \Ind\wh{\Coh}{}^!}, \qquad (\mc F)^{\rc, \Ind\wh{\Coh}{}^!}\in \wh{\Coh}{}^!(\Par^\wedge_{G,\lambda})
\] 
whose fiber is in the subcategory generated by $\Eis_{P,\lambda_P}^{\spec, \wedge} \mc G$ with $P$-dominant $\mc G \in \wh{\Coh}{}^!(\Par_{M,\lambda_M}^\wedge)$. 
The image of $(\mc F)^{\rc, \Ind\wh{\Coh}{}^!}$ under
\[
\wh{\Coh}{}^!(\Par_{G,\lambda}^\wedge)\subset \Ind\wh{\Coh}{}^!(\Par_{G,\lambda}^\wedge) \to \IndCoh (\Par_{G,\lambda}^\wedge)
\]
is isomorphic to $(\mc F)^{\rc, \IndCoh}$ as $\wh{\Coh}{}^!\subset \IndCoh$ and right-cuspidality and Eisenstein series are compatible with the inclusion. 

On the other hand, by the classification of cuspidal sheaves in the semisimple case (Theorem \ref{thm: ClassificationCuspidalCotangent} and Koszul duality), we know that $(\mc F)^{\rc, \Ind\wh{\Coh}{}^!} \in \wh{\Coh}{}^! (\Par_{G,\lambda}^\wedge)$ is in fact coherent. In particular, its image under the functor 
\[
\Ind\wh{\Coh}{}^!(\Par_{G,\lambda}^\wedge) \to \IndCoh (\Par_{G,\lambda}^\wedge)
\]
is coherent. Therefore, $(\mc F)^{\rc, \IndCoh}$ is coherent. 
\end{proof}

This completes the proofs of Theorems \ref{thm: CVBclassification} and \ref{thm: Generators}, and hence Corollary \ref{cor: LCusp=RCusp}. 

\begin{theorem}{\label{thm: cuspSOD}}
For each $b\in B(G)$, consider the subcategories generated under finite colimits and retracts by 
\[
\Eis^\spec_P \Cusp^\Coh (\Par_M)^{\kappa(b_M)}
\]
where $P$ is standard, $b_M \in B(M)$ is basic, $\kappa(b_M)$ is $P$-dominant, and the image of $b_M$ in $B(G)$ equals $b$. Then, they form a semiorthogonal decomposition of $\Coh(\Par_G)$, and the same holds for ind-coherent sheaves. An analogous statement holds with the $P$-antidominant condition instead. 
\end{theorem}

\begin{remark}
This implies \cite[Theorem 3.7.1(ii)]{HansenMann}. 
\end{remark}

\begin{proof}
Since we have shown (Corollary \ref{cor: LCusp=RCusp}) that $\Cusp(\Par_G)$ is compactly generated by $\Cusp^\Coh(\Par_G)$ for any $G$ and $\Cusp^\Coh(\Par_G)$ is preserved under Grothendieck--Serre duality, it suffices to prove the $P$-dominant version for ind-coherent sheaves. 

There are semiorthogonal decompositions of $\IndCoh(\Par_G)$
\[
\langle \Cusp(\Par_G), \Eis_+\rangle, \qquad \langle \Eis_-, \Cusp(\Par_G)\rangle;
\]
let us emphasize here that $\Eis_+$, $\Eis_-$ are generated by Eisenstein series for standard parabolic subgroups (with suitable weight conditions). 
As in the case of cotangent stacks, this, together with Lemma \ref{lem: cuspSODeisvanishing} below, implies by induction that there are semiorthogonal decompositions given by the subcategories generated by $\Eis^\spec_{P}\Cusp(\Par_{M})^{\kappa(b_M)}$. 
\end{proof}

\begin{lemma}{\label{lem: cuspSODeisvanishing}}
Let $P_1\supset M_1, P_2\supset M_2$ be standard parabolic and Levi subgroups. Let $b_1\in B(M_1)$, $b_2\in B(M_2)$ be basic elements such that $b_1\not\leq b_2$ in $B(G)$ and $b_1$ is $P_1$-dominant. 
Then for any $\mc F_1\in \Cusp(\Par_{M_1})^{\kappa(b_1)}$, $\mc F_2\in \Cusp(\Par_{M_2})^{\kappa(b_2)}$, 
\[
\Hom(\Eis^\spec_{P_1}\mc F_1, \Eis^\spec_{P_2}\mc F_2)\simeq 0,
\]
and the same holds if $b_1\not\geq b_2$, $b_2$ is $P_2$-dominant, and $\Eis^\spec_{P_1^{\op}}$, $\Eis^\spec_{P_2^{\op}}$ instead. 
\end{lemma}

\begin{proof}
We may assume that $\mc F_1, \mc F_2$ are coherent. By Grothendieck--Serre duality, it suffices to prove the first case, as $\Cusp^\Coh(\Par_{M_i})$ is preserved by duality. 
We can check the statement at the completed Vogan stack of each $\lambda$ by an argument as in the proof of Corollary~\ref{cor: cohspringerhoms} and on $\T^{\ast}\V_{\lambda}$ by Koszul duality and Lemma \ref{lem: HatCohForgettingGm}. Then this is Corollary~\ref{cor: eisHomnonvanishingconditions} in the more general setting. 
\end{proof}

\begin{corollary}{\label{cor: BasParGgenbyCusp}}
We have $\LBas(\Par_G)=\RBas(\Par_G)$, and it will be denoted by $\Bas(\Par_G)$. It is compactly generated by
\[
\Bas^\Coh(\Par_G)\coloneq \Bas(\Par_G)\cap \Coh(\Par_G), 
\]
and $\Bas^\Coh(\Par_G)$ is generated under finite colimits and retracts by $\Cusp^\Coh(\Par_G)$ and $\Eis_P^\spec \Cusp^\Coh(\Par_M)^{\kappa(b_M)}$ for proper $P$ and $b_M \in B(M)_{\bas}$ with basic image in $B(G)$; the former is both left and right orthogonal to the latter. 
\end{corollary}
\begin{proof}
We first prove that $\RBas^{\Coh}(\Par_G)$ is generated under finite colimits and retracts by $\Cusp^\Coh(\Par_G)$ and $\Eis^{\spec}_P \Cusp^\Coh(\Par_M)^{\kappa(b_M)}$ as in the statement. This implies that $\RBas^{\Coh}(\Par_G)$ is stable under $\D_{\GS}$ and hence that 
\begin{equation*}
    \RBas^{\Coh}(\Par_G) = \LBas^{\Coh}(\Par_G).
\end{equation*}
The orthogonality statement follows because $\kappa(b_M)$ for $b_M \in B(M)_{\bas}$ with basic image in $B(G)$ is both $P$-dominant and $P$-antidominant.

Now by Theorem \ref{thm: cuspSOD}, we have that $\RBas^{\Coh}(\Par_G)$ is generated by $\Eis^{\spec}_P \Cusp^{\Coh}(\Par_M)^{\kappa(b_M)}$ for (not necessarily proper) standard parabolic subgroups $P$, basic $P$-dominant $b_M \in B(M)$ whose image in $B(G)$ is basic. First, any such generator is indeed basic. As in the proof of Lemma \ref{lem: cuspSODeisvanishing}, it suffices to check that $\Eis^{\T^{\ast}}_P \Cusp^{\Coh}(\T^{\ast}\V_{\lambda_M})^{\kappa(b_M)}$ are basic, and this is Proposition \ref{prop: EisBasisBas}. Conversely, any subcategory of the form $\Eis^{\spec}_P \Cusp^{\Coh}(\Par_M)^{\kappa(b_M)}$  for $P$-dominant $b_M$ whose image in $B(G)$ is not basic is contained in $\Eis^{\circ, \Coh}_+$ by the argument of Lemma \ref{lem: ParGdomswitch}. Then the right-orthogonal complement of $\Eis^{\circ, \Coh}_+$ is simultaneously $\RBas^{\Coh}(\Par_G)$ and generated by the desired terms. 
\end{proof}

\subsection{Some consequences}

\begin{proposition}{\label{prop: formalfiberEisFF}}
Suppose $P$ is semistandard and $b_M\in B(M)_\bas$ is strictly $P$-dominant or strictly $P$-antidominant.  
The functor $\Eis^{\spec,\wedge}_{P,\lambda_P}$ is fully faithful on $\Bas(\Par^{\wedge}_{M, \lambda_M})^{\kappa(b_M)}$.
\end{proposition}

\begin{proof}
It suffices to check the claim on generators, and we have generators with graded lifts and can pass to the linear Koszul dual side by Lemma \ref{lem: HatCohForgettingGm}, where the result is Corollary~\ref{cor: fullyfaithfullbas}.
\end{proof}

\begin{proposition}{\label{prop: formalfiberEisvanishing}}
Suppose $P$ is semistandard and $b_M\in B(M)_\bas$ is strictly $P$-dominant or strictly $P$-antidominant. Let $\lambda_P\neq \lambda'_P$ be two lifts of $\lambda$. 
Then
\[
\Eis^{\spec,\wedge}_{P,\lambda_P}\Bas(\Par^{\wedge}_{M,\lambda_M})^{\kappa(b_M)}, \quad
\Eis^{\spec,\wedge}_{P,\lambda'_P}\Bas(\Par^{\wedge}_{M,\lambda'_M})^{\kappa(b_M)}
\]
are orthogonal to each other. 
\end{proposition}

\begin{proof}
We may replace $\Eis^{\spec,\wedge}_{P,\lambda'_P}\Bas(\Par^{\wedge}_{M, \lambda'_M})^{\kappa(b_M)}$ by $\Eis^{\spec,\wedge}_{w(P),\lambda_P}\Bas(\Par^{\wedge}_{w(M),\lambda_M})^{\kappa(w(b_M))}$ for some $w\in W_{\wh G}^\Gamma$ such that $w\not\in H_\lambda \wh P$. As above, we work with generators and pass to the linear Koszul dual side and apply Proposition \ref{prop: EisHomvanishing}. 
\end{proof}

\begin{corollary}{\label{cor: EisFFonBas}}
If $b_M$ is strictly $P$-dominant or strictly $P$-antidominant, then $\Eis^\spec_P$ is fully faithful on $\Bas (\Par_M)^{\kappa (b_M)}$. 
\end{corollary}

\begin{proof}
This can be checked on each $\lambda$ as in the proof of Corollary \ref{cor: cohspringerhoms}. So, this follows from Propositions \ref{prop: formalfiberEisFF} and \ref{prop: formalfiberEisvanishing}. 
\end{proof}

\begin{theorem}{\label{thm: basSODparG}}
The $\infty$-subcategories 
\[
 \Eis^\spec_{P_b} \Bas (\Par_{M_b})^{\kappa(b_M)} \quad (\textnormal{resp.~$\Eis^\spec_{P_b^\op} \Bas (\Par_{M_b})^{\kappa(b_M)}$})
\]
give rise to semiorthogonal decompositions of $\IndCoh(\Par_G)$ indexed by $B(G)$ (with different orderings), and the same holds for coherent sheaves. 

Moreover, each piece is generated under colimits by 
\[
\Eis^\spec_P \Cusp^\Coh(\Par_M)^{\kappa(b_M)}\quad (\textnormal{resp.~$\Eis^\spec_{P^\op} \Cusp^\Coh(\Par_M)^{\kappa(b_M)}$})
\]
with $P$-dominant $b_M$ whose image in $B(G)$ equals $b$.
\end{theorem}

\begin{proof}
For $\Coh(\Par_G)$, the semiorthogonal decomposition in terms of $\Eis^\spec_{P_b} \Bas^\Coh(\Par_{M_b})^{\chi_b}$ follows from the cuspidal semiorthogonal decomposition (Theorem \ref{thm: cuspSOD}) and the description of generators of $\Bas^{\Coh}(\Par_G)$ (Corollary \ref{cor: BasParGgenbyCusp}); the main new content here is the description of pieces as the essential images of Eisenstein series by Corollary \ref{cor: EisFFonBas}. For $\IndCoh(\Par_G)$ we ind-complete these semiorthogonal decompositions. 
\end{proof}

Here is a refined semiorthogonal decomposition using ``cuspidal supports'':

\begin{theorem}{\label{thm: SODcuspsupp}}
The following objects indexed by pairs $(M, (\phi_M, \rho_M))$ of a standard Levi subgroup $M(\subset P)$ and an irreducible cuspidal representation $\rho_{M}$ of $S_{\phi_{M}}$ whose central character is $P$-dominant, up to $\wh G$-inertial equivalence,
\[
\{ \Eis^\spec_P \widetilde{\mc V}_{\rho_M} \}_{[M,\rho_M]}, \qquad \{ \Eis^\spec_{P^\op} \widetilde{\mc V}_{\rho_M} \}_{[M,\rho_M]}
\]
generate subcategories (under finite colimits and retracts, resp.~colimits) giving rise to semiorthogonal decompositions indexed by the same set. More precisely, they refine the semiorthogonal decompositions in Theorem \ref{thm: basSODparG} under the map taking the central character of $\rho_M$. For a fixed $b \in B(G)$, the resulting subcategories are orthogonal to each other.  
\end{theorem}

\begin{proof}
Given Theorem \ref{thm: basSODparG}, this follows from Corollary \ref{cor: cohspringerhoms}.  
\end{proof}

\begin{proposition}{\label{prop: formalfibereisconservative}}
The functor $\Eis^{\spec,\wedge}_{P,\lambda_P}$ is conservative on $\Cusp(\Par_{M,\lambda_M}^\wedge)^{\chi}$ with $P$-dominant $\chi$. 
\end{proposition}
    
\begin{proof}
We reduce first to the case where $\chi$ is $G$-basic. Choose $P' \supset P$ with Levi subgroup $M'$ such that $M'$ is the largest Levi containing $M$ for which $\chi$ is basic. Then 
\begin{equation*}
    \Eis^{\spec,\wedge}_P = \Eis^{\spec,\wedge}_{P'} \circ \Eis^{\spec,\wedge}_{P \cap M'}.
\end{equation*}
Since $\chi|_{Z(\wh{M'})^{\Gamma}}$ is strictly $P'$-dominant, the conservativity of $\Eis_{P'}^{\spec,\wedge}$ follows from Proposition \ref{prop: formalfiberEisFF}. This completes the reduction.

Now, we may assume that our sheaf is in the subcategory of $\Ind\wh{\Coh}{}^!(\Par_{M,\lambda_M}^\wedge)$ generated under colimits by a single $\wh{\mc V}_{\rho_M}$. Since $\Eis^{\spec,\wedge}_{P,\lambda_P}$ preserves compact objects and colimits, it suffices to show that the base change functor from right $\End(\wh{\mc V}_{\rho_M})$-modules to right $\End(\Eis^{\spec,\wedge}_{P,\lambda_P}\wh{\mc V}_{\rho_M})$-modules is conservative. This follows from Corollary \ref{cor: EndEisfinitefree}.     
\end{proof}

\begin{corollary}{\label{cor: ParGEisconservativity}}
The functor $\Eis^\spec _P$ is conservative on $\Cusp(\Par_M)^{\kappa (b_M)}$ for a $P$-dominant $b_M\in B(M)_\bas$. 
\end{corollary}

\begin{proof}
To show the conservativity, it suffices to work with $\wh{\Coh}{}^!(\Par_{M,\lambda_M}^\wedge)$, and apply Proposition \ref{prop: formalfibereisconservative}. 
\end{proof}

\printbibliography

\end{document}